\documentclass[reqno]{amsart}

\numberwithin{equation}{section}
\newtheorem{theorem}{\sc Theorem}[section]
\newtheorem{lemma}{\sc Lemma}[section]
\newtheorem{proposition}{\sc Proposition}[section]
\newtheorem{remark}{\sc Remark}
\newtheorem{definition}{\sc Definition}[section]

\title[Navier-Stokes equations with surface tension]
      {Unique solvability of the free-boundary Navier-Stokes equations
with surface tension}

\author[D. Coutand and S. Shkoller]{}

\email{coutand@math.ucdavis.edu}
\email{shkoller@math.ucdavis.edu}

\begin{document}

\maketitle

\centerline{\scshape   Daniel Coutand and Steve Shkoller}
 \medskip

  {\footnotesize \centerline{ Department of Mathematics }
  \centerline{ University of California at Davis } \centerline{Davis, CA
   95616 } }

\begin{abstract}
We prove the existence and uniqueness of solutions to the time-dependent 
incompressible Navier-Stokes equations with a free-boundary governed by surface
tension.  The solution is found using a topological fixed-point theorem that 
requires the analysis of a model linear problem consisting of the time-dependent
Stokes equation with linearized mean-curvature forcing on the boundary.
We use energy methods to establish new types of spacetime estimates
which allow us to find a unique weak solution to this linear problem.  
We then prove regularity of the weak solution, and establish the a priori 
estimates required by the topological fixed-point theorem.
\end{abstract}

\tableofcontents

\section{Introduction}
We are concerned with establishing the existence and uniqueness of the 
time-dependent incompressible Navier-Stokes equations with a free-surface
governed by surface tension.
Let $\Omega_0 \subset{\mathbb R}^3$ denote an open bounded domain with boundary
$\Gamma_0:= \partial \Omega_0$.  For each $t\in (0,T]$, we wish to find the
domain $\Omega(t)$, a divergence-free velocity field $u(t, \cdot)$, a pressure
function $p(t, \cdot)$ on $\Omega(t)$, and a volume-preserving transformation
$\eta(t, \cdot):\Omega_0 \rightarrow {\mathbb R}^3$ such that
\begin{subequations}
  \label{ns}
\begin{align}
\Omega(t)&= \eta(t,\Omega_0), 
         \label{ns.a}\\
\eta_t(t,x) &= u(t, \eta(t,x)) \,,
         \label{ns.b}\\
u_t - \nu \triangle u + (u\cdot \nabla)u  &= -\nabla p + f \,,
         \label{ns.c}\\
   \operatorname{div} u &= 0     \,,
         \label{ns.d}\\
\nu \operatorname{Def}u \ n - pn &=  \sigma H n \ \ \text{ on } \ \ \Gamma(t):=
\eta(t,\Gamma_0)  \,,
         \label{ns.e}\\
   u(0,x) &= u_0(x)     \,,
         \label{ns.f}\\
   \eta(0,x) &= x     \,,
         \label{ns.g}
\end{align}
\end{subequations}
where $\nu$ is the kinematic viscosity, $\sigma>0$ denotes the surface tension, 
$n(t,\cdot)$ is the outward pointing unit normal to $\Gamma(t)$, 
$H(t,\cdot):=\frac{1}{2}\operatorname{div}n(t,\cdot)$ 
denotes the mean curvature of 
$\Gamma(t)$, and $\operatorname{Def}u$ is twice the rate of deformation tensor 
of $u$, given in coordinates by $u^i,_j + u^j,_i$.  All Latin indices run 
through $1,2,3$, the Einstein  summation convention is employed, and indices 
after commas denote partial derivatives.

Solonnikov studied the solvability of (\ref{ns})
in \cite{Sol1991,Sol1992}.  His proof did not rely on energy estimates, but
rather on Fourier-Laplace transform techniques, which required
the use of exponentially weighted anisotropic Sobolev-Slobodeskii spaces with
only fractional-order spatial derivatives for the analysis, as well as for the
statement of the main result on solvability.  In particular, he required
the initial data $u_0 \in H^{s} (\Omega_0;{\mathbb R}^3)$ for $s\in (2,2.5)$.
In a more recent article Tani in \cite{Tani1996} used Solonnikov's 
functional framework, but 
applied a new nonlinear iteration procedure in order to prove the solvability
result; unfortunately, the linear problem which he poses in equation (4.3) 
on page 319 of \cite{Tani1996} is not solvable.

In the case that $\Omega_0$ is an infinite horizontal layer of fluid with a 
rigid bottom and a free surface, Beale \cite{Beale1983} established the
stability of the equilibrium state; namely, he showed that small 
perturbation of the flat surface continue for all time.  The analysis made
clever use of the fact that the normal vector is constant on ${\mathbb R}^2$ so
that the boundary terms arising from surface tension involve surface
Laplacians of a scalar field, whereas for a {\it general domain} 
$\Omega_0$ (wherein
the normal does not commute with the surface Laplacian), the operators act
on a vector-valued field.  This is a subtle issue which significantly 
simplifies the stability analysis.  Analytically, Beale used the Laplace 
transform for the time variable; this required the use of fractional-order
derivatives for his analysis as well.

In this paper, we develop energy methods, based on a new type of energy 
inequality for the weak formulation of the fundamental linearized problem,  
arising from Solonnikov's  method of successive approximations . 
We use this inequality to obtain weak solutions of this linear problem,
proceed to develop the necessary regularity theory, and then provide
a topological fixed-point theory for the mapping associated with the method of
successive approximations; this, in turn, yields our unique solution to
(\ref{ns}).  As we shall specify below, we choose the initial velocity
field $u_0$ to be in the space $H_{\operatorname{div}}^{2} 
(\Omega_0;{\mathbb R}^3)$, and
prove the unique solvability of (\ref{ns}) in the natural energy space
$L^2(0,T; H^3_{\operatorname{div}}(\Omega(t);{\mathbb R}^3))$ associated
to this initial condition.  Our proof of regularity illuminates the
subtleties between the smoothing effects of surface tension and the derivative
loss which surface tension appears to induces upon the successive 
approximation scheme.   

In addition to the energy methods that we use for the
basic linear problem, we employ the Tychonoff fixed-point theorem to 
prove existence of solutions to the Navier-Stokes equations; our approach
requires less regularity on the forcing function $f$ than is required by the 
Banach fixed-point theorem used by Solonnikov. 

\section{The basic energy law} To understand
the intrinsic dynamics of interface flow, let us set the external forcing
$f$ to zero.
Let $K(t)=\frac{1}{2}\int_{\Omega(t)}|u(t,x)|^2 dx$ denote the kinetic 
energy of the fluid, and  let
$A(t)$ denote the surface area of $\Gamma(t)$. If there are smooth solutions 
of (\ref{ns}), they must satisfy the following basic energy law:
\begin{equation}\nonumber
\frac{d}{dt} \left[ K(t) + \sigma A(t)\right] = - \nu\int_{\Omega(t)} 
|\operatorname{Def} u|^2 dx \,.
\end{equation}
When $\sigma=0$ but $\nu >0$, the kinetic energy 
decays, whereas when $\nu=0$ and $\sigma >0$, there
is a delicate balance between the kinetic energy  of the fluid and the 
surface area of the moving boundary, the sum achieving a critical point of
the total energy.  In the former case, the problem behaves as
if it was parabolic, whereas in the latter case, it behaves almost as though
it was hyperbolic; 
Beale \cite{Beale1983} viewed the system as being {\it of mixed character}.
We view the system as behaving more as though it was parabolic; the estimates
for the linearized system show that viscosity is necessary  to control the 
higher-order derivatives of the surface tension term on the boundary.  On
the other hand, this coupling is seen only when one linearizes about an
interface for which the normal vector is not constant.  For example, in
the problem that Beale considered in \cite{Beale1983}, smoothing on the
interior due to viscosity, and smoothing on the boundary due to surface tension
decouple at the level of the linearized equations.

\section{Lagrangian formulation of the problem}
Let 
\begin{equation}\label{a}
a(x) = [\nabla\eta(x)]^{-1}, 
\end{equation}
where $(\nabla \eta(x))^i_j = \partial \eta^i/\partial x^j(x)$
denotes the matrix of partial derivatives of $\eta$.  Let $v=u\circ \eta$ denote
the Lagrangian or material velocity field, $q=p \circ \eta$ is the Lagrangian pressure
function, and $F= f \circ \eta$ is the forcing function.  
Then (\ref{ns}) can be written as 
\begin{subequations}
  \label{nsl}
\begin{alignat}{2}
\eta_t &=v&\ &\text{in} \ \ (0,T)\times \Omega_0 \,, 
         \label{nsl.a}\\
v^i_t - \nu (a^j_l a^k_l v^i,_k),_j + a^k_i q,_k &= F^i 
&&\text{in} \ \ (0,T)\times \Omega_0 \,, 
         \label{nsl.c}\\
   a^k_i v^i,_k &= 0     &&\text{in} \ \ (0,T)\times \Omega_0 \,, 
         \label{nsl.d}\\
\nu(v^i,_k a^k_l + v^l,_k a^k_i) a^j_l N_j - q a^j_i N_j &= \sigma 
\triangle_g(\eta)^i
&&\text{on} \ \ (0,T)\times \Gamma_0 \,, 
         \label{nsl.e}\\
   v &= u_0  
 &&\text{on} \ \ \Omega_0\times \{ t=0\} \,, 
         \label{nsl.f}\\
   \eta &= \text{Id}     
 &&\text{on} \ \ \Omega_0\times \{ t=0\} \,, 
         \label{nsl.g}
\end{alignat}
\end{subequations}
where $N$ denotes the outward-pointing unit normal to $\Gamma_0$, and 
$$ \triangle_g(\eta) = (Hn) \circ \eta$$ 
is defined as follows:
denote the Laplacian on $\Gamma_0$ with respect to the induced metric 
$g$ by
$ \triangle_g$, so that in 
local coordinates $\{y^\alpha, \ \alpha=1,2\}$ on $\Gamma_0$, 
\begin{equation}\label{laplacian}
 \triangle_{g(t,\cdot)}  = g^{\alpha \beta}(t,\cdot)\left(\frac{\partial^2 }
{\partial y^\alpha \partial y^\beta} - 
\Gamma^\gamma_{\alpha\beta}(t,\cdot)\frac{\partial }{\partial y^\gamma}\right),
\end{equation}
where $g^{\alpha\beta} = g_{\alpha\beta}^{-1}$,
\begin{equation}\label{metric}
g_{\alpha\beta}= \frac{\partial \eta^i}{\partial y^\alpha} \delta_{ij}
 \frac{\partial \eta^j}{\partial y^\beta}, \ \ \text{and} \ \
\Gamma^\gamma_{\alpha\beta} = \frac{1}{2}g^{\gamma\delta}
\left( \frac{\partial g_{\beta\delta}}{\partial y^\alpha} + 
\frac{\partial g_{\alpha\delta}}{\partial y^\beta} -
\frac{\partial g_{\alpha\beta}}{\partial y^\delta} \right)\,,
\end{equation}
where $\delta_{ij}$ denotes the usual (identity) metric on $ \mathbb{R}  ^3$.
In the computation of the induced metric in (\ref{metric}), the flow
map $\eta$ is restricted to the boundary.
Throughout the paper, all Greek indices run through $1,2$.

\section{Notation and preliminaries} 
We begin by specifying our notation for certain vector and matrix operations.
\begin{itemize}
\item[] We write the Euclidean inner-product between two vectors $x$ and $y$ 
as $x\cdot y$, so that $x\cdot y=x^i\ y^i$.
\item[] The transpose of a matrix $A$ will be denoted by $A^T$, {\it i.e.}, 
$(A^T)^i_j=A^j_i$.
\item[] We write the product of a matrix $A$ and a vector $b$ as $A\ b$, 
{\it i.e}, $(A\ b)^i=A^i_j b^j$.
\item[] The product of two matrices $A$ and $S$ will be denoted by 
$A\cdot S$, {\it i.e.}, $(A\cdot S)^i_j=A^i_k\ S^k_j$.
\item[] The trace of the product of two matrices $A$ and $S$ will be denoted by 
$A: S$, {\it i.e.}, $ A:S=\operatorname{Trace}(A\cdot S)=A^i_j\ S^j_i$.
\end{itemize}

For $s\ge 0$ and a Hilbert space $(X,\|\cdot\|_X)$, $H^s(\Omega; X)$ denotes the
Sobolev space of $X$-valued functions with $s$ distributional derivatives in 
$L^2(\Omega;X)$, the equivalence class of functions 
which are measurable and have finite
$\|\cdot\|_{L^2}$-norm, where 
$\|f\|^2_{L^2(\Omega;X)}= \int_\Omega \|f(x)\|^2_X dx$.

For $T>0$ and integers $m\ge 1$, we set
\begin{align*}
V^1(T)&=\{ w \in L^2(0,T;  H^1(\Omega_0;{\mathbb R}^m)) \ | \ 
w_t \in L^2(0,T;  H^1(\Omega_0;{\mathbb R}^m)') \}, \\
V^2(T)&=\{ w \in L^2(0,T;  H^2(\Omega_0;{\mathbb R}^m)) \ | \ 
w_t \in L^2(0,T;  L^2(\Omega_0;{\mathbb R}^m)) \},  \\
V^3(T)&=\{ w \in L^2(0,T;  H^3(\Omega_0;{\mathbb R}^m)) \ | \ 
w_t \in L^2(0,T;  H^1(\Omega_0;{\mathbb R}^m)) \},  \\
\end{align*}
where for any Hilbert space $X$, we use $X'$ to denote the dual space.
Letting $\Gamma:=\partial \Omega$, we use $H^{-s}(\Gamma; {\mathbb R}^m)$ to
denote the dual space of $H^{s}(\Gamma; {\mathbb R}^m)$.
We shall also need the spaces
\begin{align*}
V^k_{\operatorname{div}}(T)&=\{ w \in V^k(T) \ | \ \operatorname{div}w=0\}, 
\ \ \ k=1,2,3\,,
\end{align*}
as well as the space (of weak solutions)
\begin{align*}
{\mathcal V}(T) = &
\{ w \in L^2(0,T;  H^1(\Omega_0;{\mathbb R}^3)) \ | \ \operatorname{div}w=0,\\
& \qquad \int_0^t N \cdot \nabla_0 
w(r)|_{\Gamma_0}dr \in  L^\infty(0,T;  L^2(\Gamma_0;{\mathbb R}^3)) \},
\end{align*}
where in components $[N\cdot \nabla_0 v]^ \alpha  =  N_i g_0^{ \alpha \beta }
\partial_\beta v^i$.
When there is no time dependence, we shall use ${\mathcal V}$ to denote
the space
$\{ \psi \in H^1(\Omega_0, {\mathbb R}^3) \ | \ \operatorname{div}\psi=0, \ 
N \cdot \nabla_0\psi \in L^2(\Gamma_0; {\mathbb R}^3) \}$.

\begin{lemma}\label{lemma1}
The space ${\mathcal V}$ is separable. 
\end{lemma}
\begin{proof}
The mapping ${ I}: v\rightarrow (\partial_1 v,  \partial_2 v,
\partial_3 v, v,  \nabla_0 v\cdot N)$ 
is an isometry from ${\mathcal V}$ into the
separable space $[L^2(\Omega_0;{\mathbb R}^3)]^4\times  
L^2(\Gamma_0;{\mathbb R}^3)$. As a subset of
a separable space, $I({\mathcal V})$ is separable  
(see, for instance, \cite{Bre1983}), and consequently ${\mathcal V}$ is also a 
separable space.
\end{proof}

\section{The main theorem} 

\begin{theorem}\label{main}
Let $\Omega_0\subset {\mathbb R}^3$ be a smooth, open and bounded subset,
and suppose that $\nu>0$, $\sigma> 0$, and 
$u_0 \in H^2_{\operatorname{div}}(\Omega;{\mathbb R}^3)$
satisfies the compatibility condition 
\begin{equation}\label{compatibility}
[\operatorname{Def}u_0 \ N]_{\operatorname{tan}}=0\ \text{ on } \Gamma_0,
\end{equation}
and that
\begin{equation} \label{f_regularity}
f \in L^2(0,\bar T; H^1({\mathbb R}^3;{\mathbb R}^3)), \ \ \
f_t \in L^2(0,\bar T; H^1({\mathbb R}^3;{\mathbb R}^3)')\,.
\end{equation}

{\it Existence.}  Then there 
exists a $\bar T >0$ depending on $u_0$, $f$, and $\Omega_0$, such that 
there exists a solution $v \in V^3(\bar T)$ and 
$q\in V^2(\bar T)$ 
of the problem (\ref{nsl}).  Furthermore, 
$\eta \in  C^0([0,\bar T]; H^3(\Omega_0; {\mathbb R}^3))$ and
the surface tension term
$\sigma \triangle_g(\eta) \in L^2(0, \bar T; H^{\frac{3}{2}}(\Gamma_0; {\mathbb R}^3))$.

{\it Uniqueness.} Moreover, if there exists $K>0$ such that
\begin{equation}\label{Lip}
\begin{array}{l}
\forall t\le \bar T,\  \forall (x,y)\in {\mathbb R}^3\times {\mathbb R}^3, \\
|f(t,x)-f(t,y)|+ |\nabla f(t,x)-\nabla f(t,y)|+|f_t(t,x)-f_t(t,y)|
\le K\ |x-y|\,,
\end{array}
\end{equation}
i.e, 
$f$, $\nabla f$, and $f_t$ are uniformly Lipschitz continuous in the spatial
variable, then the solution is unique.
\end{theorem}

\begin{remark}
The regularity of our unique solution $v \in V^3(\bar T)$ implies that for
each $t\in [0,\bar T]$ there
is a unique domain $\Omega(t)$ of regularity class $H^3$,  
a unique divergence-free velocity field on $\Omega(t)$,
$u \in L^2(0,T; H^3(\Omega(t);{\mathbb R}^3))$, and a unique pressure function
$p \in L^2(0,T; H^2(\Omega(t);{\mathbb R}))$ 
solving the Eulerian problem (\ref{ns}).  Also, although 
we have stated our results for three-dimensional fluid motion, all
of our results hold for two-dimensional fluid motion as well.
\end{remark}

\begin{remark}
Theorem \ref{main} differs from the existence and uniqueness assertions of
Solonnikov \cite{Sol1992}, which require the initial velocity
field $u_0$ to be taken in $H^{s}_{\operatorname{div}}(\Omega_0;{\mathbb R}^3)$ for $s \in (2,
2.5)$, and does not permit the integer value $s=2$. 
\end{remark}

\begin{remark}
Our theorem also differs from that of Solonnikov's \cite{Sol1992} in that we 
require only the minimal regularity assumptions 
(\ref{f_regularity}) on the forcing function $f$ in order to 
establish the existence of solutions, whereas the additional Lipschitz 
assumption (\ref{Lip}) is needed only for uniqueness.  This is due to our method
of proof which employs the Tychonoff fixed-point theorem instead of the Banach
fixed-point theorem used in \cite{Sol1992}.
\end{remark}

\section{The basic linear problem}
We denote the surface Laplacian on $\Gamma_0$ by
$$ \triangle_0 := \triangle_{g(0,\cdot)}$$
where $\triangle_{g(0,\cdot)}$ is defined in (\ref{laplacian}) and 
(\ref{metric}), and $g(0,\cdot)$ is the induced metric at $t=0$.

We are concerned with the time-dependent linear problem
\begin{subequations}
  \label{linear}
\begin{align}
\bar w_t - \nu \triangle \bar w  &= -\nabla p + \bar f 
& &\text{in} \ \ (0,T)\times \Omega_0 \,, 
         \label{linear.a}\\
   \operatorname{div} \bar w &= \bar a     
& &\text{in} \ \ (0,T)\times \Omega_0 \,, 
         \label{linear.b}\\
\nu \operatorname{Def}\bar w \ N - pN &=  
 \sigma (N\cdot \triangle_0 \int_0^t \bar w(s) ds + \bar B)N +  \bar g 
& &\text{on} \ \ (0,T)\times \Gamma_0 \,, 
         \label{linear.c}\\
   \bar w &= \bar w_0
& &\text{on} \ \ \Omega_0  \times \{t=0\}\,. 
         \label{linear.f}
\end{align}
\end{subequations}

\begin{theorem}\label{thm1}
Given $\nu > 0$, $\sigma > 0$, and
$\bar g \in L^2(0,T; H^{\frac{3}{2}}(\Gamma_0; {\mathbb R}^3))$ with
$\bar g_t \in L^2(0,T; H^{-\frac{1}{2}}(\Gamma_0; {\mathbb R}^3))$,
and $\bar g(0, \cdot)_{\operatorname{tan}}=0$, and the initial data 
$\bar w_0\in H^2_{\operatorname{div}}(\Omega_0;{\mathbb R}^3)$
satisfying the compatibility condition
\begin{equation}\label{comp2}
[\operatorname{Def} \bar w_0\ N]_{\operatorname{tan}} 
= 0 \,,
\end{equation}
$ \bar B \in C^0(0,T; H^{\frac{1}{2}}(\Gamma_0; {\mathbb R})), \
\bar B_t \in L^2(0,T; H^{\frac{1}{2}}(\Gamma_0; {\mathbb R}))$,
with $B(0,\cdot) = 0$,
$\bar a \in V^2(T)$ with $a(0,\cdot)=0$, and $\bar f\in V^1(T)$,
there exists a unique solution
$\bar w \in V^3(T)$  and $p\in V^2(T)$
of (\ref{linear}) for any $T>0$.  Furthermore, we have the estimate
\begin{align}
&\|\bar w\|_{L^2(0,T;H^3(\Omega_0;{\mathbb R}^3))} +
\|N \cdot \triangle_0 \int_0^{\cdot}\bar w(s)ds\|_
{C^0([0,T];H^{\frac{1}{2}}(\Gamma_0;{\mathbb R}))} 
+ \| p\|_{L^2(0,T;H^2(\Omega_0;{\mathbb R}))} \nonumber \\
&+\|\bar w_t\|_{L^2(0,T;H^1(\Omega_0;{\mathbb R}^3))} +
\|N \cdot \triangle_0 \bar w\|_
{L^2(0,T;H^{\frac{1}{2}}(\Gamma_0;{\mathbb R}))} 
+ \| p_t\|_{L^2(0,T;L^2(\Omega_0;{\mathbb R}))} \nonumber \\
&+ \| \nabla \bar w_t\|_{L^2(0,T;H^{-\frac{1}{2}}(\Gamma_0;{\mathbb R}^9))}
+ \| p_t\|_{L^2(0,T;H^{-\frac{1}{2}}(\Gamma_0;{\mathbb R}))} \nonumber \\
& \qquad \le C \left( \|w_0\|_{H^2(\Omega_0; {\mathbb R}^3)}
+ \|g(0,\cdot)\cdot N\|_{H^{\frac{1}{2}}(\Gamma_0; {\mathbb R})} 
+  \|f(0)\|_{L^2(\Omega_0; {\mathbb R}^3)}\right.\nonumber\\
& \qquad\qquad +\|\bar f\|_ {L^2(0,T;H^1(\Omega_0;{\mathbb R}^3))}+
\|\bar f_t\|_ {L^2(0,T;H^1(\Omega_0;{\mathbb R}^3)')} \nonumber \\
& \qquad\qquad
+ \|\bar a\|_ {L^2(0,T;H^2(\Omega_0;{\mathbb R}))} 
+ \|\bar a_t\|_ {L^2(0,T;L^2(\Omega_0;{\mathbb R}))} 
+ \|\bar B\|_ {C^0([0,T];H^\frac{1}{2}(\Gamma_0;{\mathbb R}))}\nonumber \\
& \qquad\qquad  \left.
+ \|\bar B_t\|_ {L^2(0,T;H^{\frac{1}{2}}(\Gamma_0;{\mathbb R}))}
+ \|\bar g\|_ {L^2(0,T;H^\frac{3}{2}(\Gamma_0;{\mathbb R}))}
+ \|\bar g_t\|_ {L^2(0,T;H^{-\frac{1}{2}}(\Gamma_0;{\mathbb R}))}
\right)\,, \label{estimate1}
\end{align}
where the $C$ may depend on $T$, $\Omega_0$, $\nu$, $\sigma$, and 
$C(T)$ remains bounded as $T  \rightarrow  0$.
\end{theorem}

\medskip

\section{Proof of Theorem \ref{thm1}}
\subsection{Weak solutions on $\Omega_0$}  
\subsubsection{The divergence-free linear problem}
We first transform (\ref{linear}) into a divergence-free problem for the
velocity field.  To do so, we shall consider the following elliptic problem:
For a.e.  $t\in [0,T]$, we solve 
\begin{subequations}
  \label{elliptic0}
\begin{alignat}{2}
-\triangle r(t, \cdot) &= \bar a(t, \cdot) & \ &\text{ in } \ \ \Omega_0\,, 
         \label{elliptic0.a}\\
 r(t, \cdot) &= r_0(t, \cdot) &&\text{ on } \ \ \Gamma_0\,, 
         \label{elliptic0.b}\\
-\triangle_0 r_0(t, \cdot) &= \bar a(t, \cdot)  &&\text{ on } \ \ \Gamma_0\,, 
         \label{elliptic0.c}
\end{alignat}
\end{subequations}
where $\triangle_0$ denotes the surface Laplacian on $\Gamma_0$, and define
\begin{equation}\label{v} 
v(t,\cdot) = \nabla r(t, \cdot)\,. 
\end{equation}

\begin{lemma}\label{v_regularity}
$v\in V^3(T)$ and $[\operatorname{Def} v_t \ N] \cdot N \in 
L^2(0,T; H^{\frac{1}{2}}(\Gamma_0;{\mathbb R}))$, and 
\begin{equation}\label{v0}
v(0,\cdot)=0\,.
\end{equation}
\end{lemma}
\begin{proof}
For a.e. $0\le t\le T$, $\bar a(t,\cdot) \in H^2(\Omega_0, {\mathbb R})$ and 
$\bar a_t(t,\cdot) \in L^2(\Omega_0, {\mathbb R})$.  It follows by elliptic
regularity of $\triangle_0$ that for a.e. $0\le t\le T$, 
$r_0(t, \cdot) \in H^{3.5}(\Gamma_0,{\mathbb R})$ and 
$\partial_t r_0(t, \cdot) \in H^{1.5}(\Gamma_0,{\mathbb R})$, so that
by elliptic regularity of the Dirichlet problem, we see that
$r(t, \cdot) \in H^{4}(\Omega_0,{\mathbb R})$ and 
$r_t(t, \cdot) \in H^{2}(\Omega_0,{\mathbb R})$, so that $v \in V^3(T)$.
By the usual compactness argument, $v \in C([0,T]; 
H^2(\Omega_0;{\mathbb R}^3))$, and since $a(0,\cdot)=0$, we see that
$v(0,\cdot) = 0$. 

It remains to show that $[\operatorname{Def} v_t \ N] \cdot N \in 
L^2(0,T; H^{\frac{1}{2}}(\Gamma_0;{\mathbb R}))$ (which is better than
the trace provides).  We may characterize the surface Laplacian as the
trace of the operator $\nabla_0  \nabla_0$, where $\nabla_0$ is the
surface gradient, given by
$$ 
\nabla_0 = (\text{Id} - N \otimes N) \ \nabla \,. 
$$ 
A simple computation shows that on $\Gamma_0$
$$
\triangle = \triangle_0 +  \nabla_0 \cdot [(N \otimes N) \ \nabla]
+ (N \otimes N) : \nabla \nabla \,. 
$$
Because of (\ref{elliptic0.a}) and (\ref{elliptic0.c}), 
\begin{equation}\label{rr}
(N\otimes N) : \nabla \nabla r_t = - \nabla_0 \cdot [(N \otimes N) \ \nabla]
r_t \,, 
\end{equation}
and since
$$
[\operatorname{Def} v_t\ N]\cdot N = 2 (N \otimes N): \nabla v_t =
2( N \otimes N) : \nabla \nabla r_t
$$
(where we have used (\ref{v}) for the second equality), (\ref{rr}) proves
the lemma, since $\nabla r_t \in 
L^2(0,T; H^{\frac{1}{2}}(\Gamma_0; {\mathbb R}^3))$.
\end{proof}

We define the new divergence-free velocity field 
\begin{equation}\label{w}
w= \bar w - v \,, 
\end{equation}
and correspondingly, the modified data:
\begin{align}
f&= \bar f - v_t + \nu\triangle v  
\label{f}\,, \\
 g&= \bar g 
-[\nu \operatorname{Def} v\ N]_{\operatorname{tan}}\,, \nonumber \\
B &= \bar B + N\cdot \triangle_0 \int_0^t v(s) ds - \nu( \operatorname{Def}
v \ N)\cdot N  \,, \nonumber\\
\intertext{and}
w_0&:= w(0) = \bar w_0 \label{w0}\,,
\end{align}
where we have used (\ref{v0}) for (\ref{w0}), and $[\cdot]_{\operatorname{tan}}$
denotes the tangential component.

\begin{lemma}\label{data}
The modified data has the following regularity:
\begin{alignat}{2}
f&\in V^1(T) & \ &\,, \nonumber \\
g &\in L^2(0,T; H^{\frac{3}{2}}(\Gamma_0; {\mathbb R}^3)), 
&g_t &\in L^2(0,T; H^{-\frac{1}{2}}(\Gamma_0; {\mathbb R^3})) \,, \nonumber \\
B &\in C^0([0,T]; H^{\frac{1}{2}}(\Gamma_0; {\mathbb R})), \ \
&B_t &\in L^2([0,T]; H^{\frac{1}{2}}(\Gamma_0; {\mathbb R}))\,, \nonumber \\
\intertext{and}
B(0, \cdot) &= 0 \,,\label{B0} \\
N \cdot g(0, \cdot) &\in H^{\frac{1}{2}}(\Gamma_0; {\mathbb R}) \,.\label{Ng0}
\end{alignat}
\end{lemma}
\begin{proof}
Equations (\ref{B0}) and (\ref{Ng0})
follow from (\ref{v0}) and the fact that $\bar B (0,\cdot)=0$. 
Lemma \ref{v_regularity} also establishes the desired regularity for
$f$, $g$, $g_t$, $B$, and $B_t$.  It remains to show that
$f_t \in L^2(0,T; H^1(\Omega_0;{\mathbb R}^3)')$.  From (\ref{f}) we see
that $f_t = \bar f_t - v_{tt} + \nu\triangle v_t$
has the desired regularity if
$v_{tt} \in L^2(0,T; H^1(\Omega_0;{\mathbb R}^3)')$.  This is true if
$\bar a_{tt} \in L^2(0,T; H^2(\Omega_0;{\mathbb R})')$, which  is the case 
if $\bar w_{tt} \in L^2(0,T; H^1(\Omega_0;{\mathbb R}^3)')$.  This
follows from (\ref{linear.a}) and the a priori estimate (\ref{estimate1}) which
shows that $\nabla p_t  \in L^2(0,T; H^1(\Omega_0;{\mathbb R}^3)')$.  
\end{proof}

\begin{remark}
To be more precise, we may define $\bar a^\epsilon:= J^\epsilon * \bar a$
for a family of Friedrichs mollifiers $J^\epsilon$, and then pass to the
limit as $\epsilon\rightarrow 0$ after obtaining the a priori estimate
(\ref{estimate1}).  The above lemma proves the consistency in allowing
$f \in V^1(T)$ (and hence 
$\nabla p_t  \in L^2(0,T; H^1(\Omega_0;{\mathbb R}^3)')$) with the 
estimate (\ref{estimate1}) which is itself established using that 
$f\in V^1(T)$.
\end{remark}

We now consider the following divergence-free problem:
\begin{subequations}
  \label{linear0}
\begin{alignat}{2}
 w_t - \nu \triangle  w  &= -\nabla p + f 
&&\text{in} \ \ (0,T)\times \Omega_0 \,, 
         \label{linear0.a}\\
   \operatorname{div}  w &= 0
&&\text{in} \ \ (0,T)\times \Omega_0 \,, 
         \label{linear0.b}\\
\nu \operatorname{Def} w \ N - pN &=  
 \sigma (N\cdot \triangle_0 \int_0^t  w(s) ds +  B)N +   g 
&\ \ &\text{on} \ \ (0,T)\times \Gamma_0 \,, 
         \label{linear0.c}\\
    w &=  w_0
&&\text{on} \ \ \Omega_0  \times \{t=0\}\,. 
         \label{linear0.d}
\end{alignat}
\end{subequations}
\subsubsection{Weak solutions of (\ref{linear0})} We shall need
the following 
\begin{lemma}\label{lemma3} For $u\in H^2(\Gamma_0)$ and $v\in H^1(\Gamma_0)$,
$$
\int_{\Gamma_0} (N\cdot \triangle_0 u)(N \cdot v) \, dS   =
-\int_{\Gamma_0} \partial_\alpha u^i g_0^{\alpha\beta}\partial_\beta
(N^jv^jN^i) \, dS \,.
$$
\end{lemma}
\begin{proof}
Letting $\mathfrak{g}= \det{g_0}$, integration by parts using
the formula
%\begin{equation}\label{laplacian}
$$
\mathfrak{g}^{-1}\partial_\alpha (\mathfrak{g} g_0^{ \alpha \beta } 
\partial_\beta u) = g_0^{ \alpha \beta }[ \partial^2_{\alpha\beta} u -
(\Gamma_0)^\gamma _{ \alpha \beta } \partial_ \gamma u]
$$
%\end{equation} 
yields the result for all $u \in C^\infty (\Gamma_0)$, and hence by 
approximation for all $u \in H^2(\Gamma_0)$.
\end{proof}
This allows us to make the following
\begin{definition} \label{def1} Let $w_0\in {L^2_{div}}=\{w\in L^2 (\Omega_0;{\mathbb R}^3)|\ \operatorname {div} w=0 \}$ .
A vector $w \in {\mathcal V}(T)$ with $w_t \in
L^2(0,T;   {\mathcal V}')$ is a weak solution of (\ref{linear0})
provided that 
\begin{align*}
\operatorname{(i)} \ & \langle w_t, v\rangle 
+\frac{ \nu }{2} (\operatorname{Def}w, \operatorname{Def}v)_{L^2(\Omega_0;{\mathbb R}^9)} 
=  \langle f, v\rangle 
+ \langle g , v \rangle_{\Gamma_0}  \\
&\ \  +  \sigma \int_{\Gamma_0} BN\cdot v \, dS 
-\sigma \int_{\Gamma_0} \int_0^t \partial_\alpha w^i(s,\cdot) ds  \,
g_0^{\alpha\beta}\, \partial_\beta (N^jv^jN^i) \, dS \ \ \ \forall v\in 
{\mathcal V}\,,\\
\operatorname{and} \\
\operatorname{(ii)} \ & w(0,\cdot) = w_0, 
\end{align*}
for a.e. $0\le t \le T$, where $\langle \cdot, \cdot \rangle$ denotes the
duality pairing between $H^1(\Omega_0,{\mathbb R}^3)$ and 
$H^1(\Omega_0,{\mathbb R}^3)'$, and 
$\langle \cdot, \cdot \rangle_{\Gamma_0}$ denotes the
duality pairing between $H^{\frac{1}{2}}(\Gamma_0,{\mathbb R}^3)$ and 
$H^{-\frac{1}{2}}(\Gamma_0,{\mathbb R}^3)$.
\end{definition}
%\begin{lemma}\label{c0l2}
%If $w \in {\mathcal V}(T)$ and $w_t \in L^2(0,T; {\mathcal V}')$.  Then
%$u\in C([0,T]; L^2(\Omega_0; {\mathbb R}^3))$, and the mapping
%$t \mapsto \|w(t)\|^2_{L^2(\Omega_0; {\mathbb R}^3)}$ is absolutely 
%continuous, with 
%$$\frac{d}{dt} \|w(t)\|^2_ {L^2(\Omega_0; {\mathbb R}^3)}
%= 2 \langle w_t(t), w(t) \rangle$$
%for a.e. $t\in [0,T]$. 
%\end{lemma}
%\begin{proof}
%Since $w \in{\mathcal V}(T)$, 
%we have $w \in L^2(0,T; {\mathcal V})$,  and since $w_t
%\in L^2(0,T; {\mathcal V}')$, the same argument as in 
%Theorem 3 of Section 5.9.3 of \cite{Evans1998} proves the lemma.
%\end{proof}

\subsubsection{Galerkin approximations}

For 
$$f \in L^2(0,T; H^1(\Omega_0;{\mathbb R}^3)'), \ \
g \in L^2(0,T; H^{- \frac{1}{2}}(\Gamma_0;{\mathbb R}^3)),$$
the Riesz representation theorem asserts the existence of 
$$\tilde f \in L^2(0,T; H^1(\Omega_0;{\mathbb R}^3)), \ \
\tilde g \in L^2(0,T; H^{ \frac{1}{2}}(\Gamma_0;{\mathbb R}^3)),$$
respectively, such that
$$
\langle f, v\rangle = (\tilde f, v)_{H^1(\Omega_0; {\mathbb R}^3)}, \ \
\langle g, v\rangle_{\Gamma_0}
 = (\tilde g, v)_{H^{\frac{1}{2}}(\Gamma_0; {\mathbb R}^3)}.
$$

Lemma \ref{lemma1} guarantees the existence of a basis
$\{ \psi_k\}_{k=1}^\infty$ of $ {\mathcal V}$ which is also an orthonormal basis of $L^2_{div}$ with 
respect to the $L^2(\Omega_0; {\mathbb R}^3)$ inner-product.
Fix $m \in{\mathbb N}$, and let $w_m:[0,T] \rightarrow  {\mathcal V}$ be 
given by

\begin{equation}\label{14} 
w_m(t,x) := \sum_{k=1}^m \lambda^k_m(t) \psi_k(x), 
\end{equation}
where we choose the coefficients $\lambda^k_m(t)$ for $t\in [0,T]$ and
$k=1,...,m$ such that
\begin{equation}\label{15}
\lambda^k_m(0)=(w_0, \psi_k)_{L^2(\Omega_0;{\mathbb R}^3)},
\end{equation}
and, letting $w_m':= \partial_t w_m$,
\begin{align}
& (w_m', \psi_k)_{L^2(\Omega_0;{\mathbb R}^3)} 
+ \frac{ \nu }{2}  (\operatorname{Def}w_m, \operatorname{Def}\psi_k)_{L^2(\Omega_0;{\mathbb R}^3)} 
= (\tilde f, \psi_k)_{H^1(\Omega_0;{\mathbb R}^3)}  
 +(\tilde g, \psi_k)_{H^{\frac{1}{2}}(\Gamma_0;{\mathbb R}^3)}  \nonumber \\
&+  \sigma (B N , \psi_k)_{L^2(\Gamma_0; {\mathbb R}^3)} 
%+\int_{\Gamma_0}B N  \cdot \psi_k \, dS 
-\sigma \int_{\Gamma_0} \partial_\alpha \int_0^t w_m^i(s,\cdot) ds  \,
g_0^{\alpha\beta}\, \partial_\beta (N^j(\psi_k)^jN^i) \, dS \,.
\label{16}
\end{align}

\begin{proposition}\label{prop1}
For each $m\in{\mathbb N}$, there exists a unique function $w_m$ of the form
(\ref{14}) satisfying (\ref{15}) and (\ref{16}). 
\end{proposition}
\begin{proof}
Assuming $w_m$ has the structure of (\ref{14}), we define
$$G_{jk}:=(\operatorname{Def} \psi_j, \operatorname{Def}\psi_k)
_{L^2(\Omega_0; {\mathbb R}^3)}, \ \ 
H_{jk}:= \int_{\Gamma_0} \partial_\alpha (\psi_j)^i \,
g_0^{\alpha\beta}\, \partial_\beta (N^j(\psi_k)^jN^i) \, dS.$$
Letting 
$$F_k(t):= (\tilde f, \psi_k)_{H^1(\Omega_0; {\mathbb R}^3)}
+ \sigma (BN, \psi_k|_{\Gamma_0})_{L^2(\Gamma_0; {\mathbb R}^3)}
+(\tilde g, \psi_k|_{\Gamma_0})_{H^{ \frac{1}{2}}(\Gamma_0; {\mathbb R}^3)}
,$$
we write (\ref{16}) as
\begin{equation}\label{ode1}
{\lambda_m^k}'(t) +  \frac{ \nu }{2}  G_{jk} \lambda^j_m(t) + \sigma H_{jk}
\int_0^t \lambda^j_m(r)dr = F_k(t).
\end{equation}
Let $d^j_m(t) = \int_0^t \lambda^j_m(r)dr$ so that 
\begin{equation}\label{15a}
d^j_m(0)=0.  
\end{equation}

Hence,
${d^j_m}'(t)= \lambda^j_m(t)$, and we may write (\ref{ode1}) as the 
second-order ordinary differential equation
\begin{equation}\label{ode2}
{d_m^k}''(t) +  \frac{ \nu }{2}  G_{kj} {d^j_m}'(t) + \sigma H_{kj}
d^j_m (t) = F_k(t),
\end{equation}
subject to the initial conditions (\ref{15}) and (\ref{15a}).
By the fundamental theorem of ordinary differential equations, there exists
(for each $k=1,...,m$) a unique absolutely continuous function $d^k_m(t)$
satisfying (\ref{ode2}), (\ref{15}) and (\ref{15a}) for a.e. $0\le t\le T$,
and hence $\lambda^k_m(t)$ satisfying (\ref{ode1}).  Thus, $w_m$ defined
by (\ref{14}) solves (\ref{16}) for a.e. $0\le t\le T$.
\end{proof}

\subsubsection{Energy estimates}

\begin{theorem}
There exists a constant $C$, depending only on $\Omega_0$ and $T$, such that
\begin{align}
&\max_{0\le t\le T} \left[\|w_m(t)\|^2_{L^2(\Omega_0; {\mathbb R}^3)}
+ \|N\cdot \nabla_0 \int_0^t w_m(s)ds\|^2_{L^2(\Gamma_0;{\mathbb R})} \right] 
 + \|w_m\|^2_{L^2(0,T; H^1(\Omega_0; {\mathbb R}^3))} 
\nonumber \\
&\qquad \qquad \le
 (\tilde c_3 + C(T))\left[ 
 \|f\|^2_{L^2(0,T;H^1(\Omega_0;{\mathbb R}^3)')} 
    + \|B\|^2_{L^2(0,T;L^2(\Gamma_0;{\mathbb R}))}\right. \nonumber \\
&\qquad \qquad \qquad \qquad\left.
+ \|g\|^2_{L^2(0,T;H^{-\frac{1}{2}}(\Gamma_0;{\mathbb R}^3))} \right]
+ (1 + C(T)) \|w_0\|^2_{L^2(\Omega_0; {\mathbb R}^3)} 
\label{energy1}
\end{align}
for $m=1,2,...$ , where $C(T) = (\tilde c_4 T + \tilde c_5 T^2) 
e^{\tilde c_1 T + \tilde c_2 T^2}$ for some constants
$\tilde c_j>0$, $j=1,2,3,4,5$.
\end{theorem}
\begin{proof}
We multiply (\ref{16}) by $\lambda^k_m(t)$, sum for $k=1,\dots,m$, 
to find that
\begin{align}
& (w_m', w_m)_{L^2(\Omega_0;{\mathbb R}^3)} 
+ \frac{ \nu }{2}  (\operatorname{Def}w_m, \operatorname{Def}w_m)_{L^2(\Omega_0;{\mathbb R}^9)} 
= (\tilde f, w_m)_{H^1(\Omega_0;{\mathbb R}^3)}
+(\tilde g, w_m)_{H^{\frac{1}{2}}(\Gamma_0;{\mathbb R}^3)} \nonumber \\
&\qquad 
+ \sigma \int_{\Gamma_0} B N \cdot w_m \, dS 
-\sigma \int_{\Gamma_0} \partial_\alpha \int_0^t {w_m}^i(s,\cdot) ds  \,
g_0^{\alpha\beta}\, \partial_\beta (N^j{w_m}^jN^i) \, dS
\label{ss}
\end{align}
for a.e. $0\le t\le T$.  We will make use of {\it Korn's  inequality} which
states that for $u\in H^1(\Omega_0; {\mathbb R}^3)$, there exists a 
constant $C>0$ depending only on $\Omega_0$, such that 
$\|\nabla u\|_{L^2(\Omega_0;{\mathbb R}^9)}^2 \le 
C\left[ \| u\|_{L^2(\Omega_0;{\mathbb R}^3)}^2 
+\| \operatorname{Def} u\|_{L^2(\Omega_0;{\mathbb R}^9)}^2  \right]$, so
that the right-hand-side defines a norm equivalent to the $H^1(\Omega_0,
{\mathbb R}^3)$-norm.
For $j\in{\mathbb N}$, we will use $c_j$ to denote a positive constant
which generally only depends on the domain $\Omega_0$ unless otherwise
specified, and $\epsilon >0$ will represent a small number.
 
Using Young's inequality, for a.e. $t\in[0,T]$,
\begin{align*}
(\tilde f, w_m)_ {H^1(\Omega_0; {\mathbb R}^3)}  &=
(\tilde f, w_m)_ {L^2(\Omega_0; {\mathbb R}^3)} +
( \nabla_j\tilde f,  \nabla_jw_m)_ {L^2(\Omega_0; {\mathbb R}^9)} \\
& \le 
c_2(\epsilon) \|\tilde f\|^2_{H^1(\Omega_0; {\mathbb R}^3)} + \epsilon
\|w_m\|^2_{H^1(\Omega_0; {\mathbb R}^3)} \\
& =
c_2(\epsilon) \| f\|^2_{H^1(\Omega_0; {\mathbb R}^3)'} + \epsilon
\|w_m\|^2_{H^1(\Omega_0; {\mathbb R}^3)} 
\end{align*}
Similarly for a.e. $t\in[0,T]$,
\begin{align*}
(\tilde g, w_m)_ {H^{\frac{1}{2}}(\Gamma_0; {\mathbb R}^3)}  &=
(\tilde g, w_m)_ {L^2(\Gamma_0; {\mathbb R}^3)} +
((-\triangle_0)^{\frac{1}{4}}\tilde g,  
(-\triangle_0)^{\frac{1}{4}} w_m)_ {L^2(\Gamma_0; {\mathbb R}^3)} \\
& \le 
c_3(\epsilon) \|\tilde g\|^2_{H^{\frac{1}{2}}(\Omega_0; {\mathbb R}^3)} 
+ \epsilon \|w_m\|^2_{H^{\frac{1}{2}}(\Gamma_0; {\mathbb R}^3)} \\
& \le
c_3(\epsilon) \| g\|^2_{H^{-\frac{1}{2}}(\Gamma_0; {\mathbb R}^3)} 
+ \epsilon \|w_m\|^2_{H^1(\Omega_0; {\mathbb R}^3)} 
\end{align*}
where we have used Young's inequality for the second inequality, 
that $\|g\|_{H^{-\frac{1}{2}}(\Gamma_0; {\mathbb R}^3)} =
\|\tilde g\|_{H^{\frac{1}{2}}(\Gamma_0; {\mathbb R}^3)}$, and
the trace theorem for the last equality.  Also,
$(B N, w_m )_{L^2(\Gamma_0; {\mathbb R}^3)} \le 
c_3(\epsilon) \|B\|^2_{L^2(\Gamma_0; {\mathbb R})} + \epsilon
\|w_m\|^2_{H^1(\Omega_0; {\mathbb R}^3)} $ for a.e. $0\le t\le T$. 
Since $(w'_m, w_m)_{L^2(\Omega_0; {\mathbb R}^3)} = \frac{1}{2}\frac{d}{dt} 
\|w_m\|_{L^2(\Omega_0; {\mathbb R}^3)}^2$,  adding $\nu\| w_m\|^2_
{L^2(\Omega_0; {\mathbb R}^3)}$ to both sides of (\ref{ss}), 
yields the inequality
\begin{align}
&\frac{d}{dt} \|w_m\|_{L^2(\Omega_0; {\mathbb R}^3)}^2
+ c_1 \|w_m\|_{H^1(\Omega_0; {\mathbb R}^3)}^2 \le
c_2(\epsilon) \|f\|^2_{H^1(\Omega_0; {\mathbb R}^3)'}
+c_3(\epsilon) \|g\|^2_{H^{-\frac{1}{2}}(\Gamma_0; {\mathbb R}^3)} \nonumber\\
& \qquad \qquad
+ c_3(\epsilon) \|B\|^2_{L^2(\Gamma_0; {\mathbb R})} 
+ c_4\|w_m\|^2_ {L^2(\Omega_0; {\mathbb R}^3)} 
+ 2 \epsilon\|w_m\|^2_{H^1(\Omega_0; {\mathbb R}^3)} \nonumber\\
& \qquad \qquad
-\sigma \int_{\Gamma_0} \partial_\alpha \int_0^t {w_m}^i(s,\cdot) ds  \,
g_0^{\alpha\beta}\, \partial_\beta (N^j{w_m}^jN^i) \, dS \,. \label{00}
\end{align}

We now study the last term in this inequality.  Expanding the $\partial_\beta$
term, and integrating by parts yields
\begin{align}
&-\int_{\Gamma_0} \int_0^t {w_m}^i,_\alpha(s,\cdot) ds  \,
g_0^{\alpha\beta}\, (N^j{w_m}^jN^i),_\beta \, dS =
-\frac{1}{2} \frac{d}{dt}  \int_{\Gamma_0} \left| \int_0^t N\cdot \nabla_0
w_m dr \right|^2 dS\nonumber  \\
& - \left( N\cdot \nabla_0 \int_0^t w_m\, dr ,\nabla_0 N \cdot w_m
\right)_ {L^2(\Gamma_0; {\mathbb R}^2)} 
+ \left( \nabla_0 N\cdot \int_0^t w_m\, dr , N \cdot \nabla_0 w_m
\right)_ {L^2(\Gamma_0; {\mathbb R}^2)}\nonumber  \\
& + \left( \nabla_0 N\cdot \int_0^t w_m\, dr ,\nabla_0 N \cdot w_m
\right)_ {L^2(\Gamma_0; {\mathbb R}^2)} 
+\int_{\Gamma_0} N^i,_{\alpha\beta} \int_0^t {w_m}^i \, dr g_0^{\alpha\beta}
N^j{w_m}^j \, dS \nonumber \\
& 
+\int_{\Gamma_0} N^i,_{\beta} \int_0^t {w_m}^i \, dr {g_0}^{\alpha\beta},_{
\alpha}
N^j{w_m}^j \, dS\,, \label{01}
\end{align}
where 
$(\nabla_0 f,\nabla_0 h)_{L^2(\Gamma_0; {\mathbb R}^2)}:= 
\int_{\Gamma_0} f,_{\alpha} g^{\alpha\beta}_0 \, h,_{\beta} \,dS$.

Choosing $\epsilon$ sufficiently small and $0 \le t \le \tilde T \le T$, 
integrating (\ref{00}) from
$0$ to $\tilde T$, and using (\ref{01}), we obtain the inequality
\begin{align*}
&\|w_m(\tilde T)\|_{L^2(\Omega_0; {\mathbb R}^3)}^2 
+  \frac{ \sigma }{2} 
\int_{\Gamma_0}\left| \int_0^{\tilde T} N\cdot \nabla_0 w_m(r)dr\right|^2 dS 
+ c_5 \|w_m\|^2_{L^2(0,\tilde T;H^1(\Omega_0; {\mathbb R}^3))} \\
& \ \ \le
\|w_m(0)\|_{L^2(\Omega_0; {\mathbb R}^3)}^2 
 + c_2(\epsilon) \|f\|^2_{L^2(0,\tilde T;H^1(\Omega_0; {\mathbb R}^3)')} \\
& \ \ 
+ c_3(\epsilon) \|g\|^2_{L^2(0,\tilde T;H^{-\frac{1}{2}}(\Gamma_0; {\mathbb R}^3))}  
+ c_3(\epsilon) \|B\|^2_{L^2(0,\tilde T;L^2(\Gamma_0; {\mathbb R}))}  \\
& \ \ + c_4 \|w_m\|^2_{L^2(0,\tilde T;L^2(\Gamma_0; {\mathbb R}^3))}  
-\sigma \int_0^{\tilde T}
\left( N\cdot \nabla_0 \int_0^t w_m\, dr ,\nabla_0 N \cdot w_m
\right)_ {L^2(\Gamma_0; {\mathbb R}^2)}
\, dt
\\
& \ \ +
 \sigma  \int_0^{\tilde T}
\left( \nabla_0 N\cdot \int_0^t w_m\, dr , N \cdot \nabla_0 w_m
\right)_ {L^2(\Gamma_0; {\mathbb R}^2)} \,  dt \\
& \ \ + c_6 \int_0^{\tilde T} \int_{\Gamma_0} \left| \int_0^t w_m(r)dr\right| \,
\left| w_m \right| \, dS \, dt,
\end{align*}
where the last term above bounds the last three terms on the
right-hand-side of equation (\ref{01}).  By Young's inequality and
Korn's inequality,

\begin{align}
&\sigma \int_0^{\tilde T}
\left( N\cdot \nabla_0 \int_0^t w_m\, dr ,\nabla_0 N \cdot w_m
\right)_ {L^2(\Gamma_0; {\mathbb R}^2)}
\,  dt
\nonumber \\
& \ \ \le
c_7(\epsilon)
\int_0^{\tilde T} \int_{\Gamma_0}\left| \int_0^t N\cdot \nabla_0 w_m(r)dr\right|^2 dS 
\, dt +  \epsilon \|w_m\|^2_{L^2(0,{\tilde T};H^1(\Omega_0; {\mathbb R}^3))}
\label{02}
\end{align}
Next, an integration-by-parts in $t$ yields
\begin{align}
& \sigma  \int_0^{\tilde T}
\left( \nabla_0 N\cdot \int_0^t w_m\, dr , N \cdot \nabla_0 w_m
\right)_ {L^2(\Gamma_0; {\mathbb R}^2)} \,  dt \nonumber \\
& \ \ \le -\sigma  \int_0^{\tilde T}
\left( \nabla_0 N\cdot  w_m , N \cdot \nabla_0 \int_0^t w_m\, dr
\right)_ {L^2(\Gamma_0; {\mathbb R}^2)} \,  dt \nonumber \\
& \ \ \ \ \ \ +\sigma  
\left( \nabla_0 N\cdot \int_0^{\tilde T} w_m\, dr , N \cdot \nabla_0 
\int_0^{\tilde T} w_m\, dr \right)_ {L^2(\Gamma_0; {\mathbb R}^2)}. \label{03}
\end{align}
The first term on the right-hand-side of (\ref{03}) is the same as
(\ref{02}).  Young's inequality for the second term gives 
\begin{align}
&\sigma  
\left( \nabla_0 N\cdot \int_0^{\tilde T} w_m\, dr , N \cdot 
\nabla_0 \int_0^{\tilde T} w_m\, dr
\right)_ {L^2(\Gamma_0; {\mathbb R}^2)} \nonumber \\
& \ \  \le \epsilon 
\int_{\Gamma_0} \left| \int_0^{\tilde T}  N \cdot 
\nabla_0 w_m(r)\, dr\right|^2\, dS
+ c_8(\epsilon) \int_0^{\tilde T}\int_{\Gamma_0} |w_m|^2 \, dS\, dt.
\label{04}
\end{align}
By the general trace theorem (see Theorem 5.22 in \cite{Adams1978}), 
there is a constant $C$ such that
$\int_{\Gamma_0}|w_m|^2dS \le C 
\|\nabla w_m\|^2_{L^{1.5}(\Omega_0;{\mathbb R}^3)}$ and 
by the standard interpolation inequality,
$$\|\nabla w_m\|^2_{L^{1.5}(\Omega_0;{\mathbb R}^9)} \le C
\|\nabla w_m\|^{2 r}_{L^2(\Omega_0;{\mathbb R}^9)}
\| w_m\|^{2s}_{L^2(\Omega_0;{\mathbb R}^3)}$$ 
for $r+s=1$.  Thus, with Young's inequality 
\begin{align}
&c_8(\epsilon) \int_0^{\tilde T}\int_{\Gamma_0} |w_m|^2 \, dS\, dt
\le \epsilon \| w_m\|^2_{L^2(0,{\tilde T}; H^1(\Omega_0;{\mathbb R}^3))} 
+ c_9(\epsilon)\int_0^{\tilde T} \|w_m(t)\|^2_{L^2(\Omega_0; {\mathbb R}^3)} \, dt.
\label{05}
\end{align}
Finally, Jensen's inequality together with Young's inequality and
the trace theorem provides the inequality
\begin{align}
& c_6 \int_0^{\tilde T} \int_{\Gamma_0} \left| \int_0^t w_m(r)dr\right| \,
\left| w_m \right| \, dS \, dt \nonumber \\
&\qquad \le
 c_{10}(\epsilon) {\tilde T}\int_0^{\tilde T} \int_{\Omega_0} 
\int_0^t  |\operatorname{Def}w_m(r)|^2 dr \, dS \, dt 
+ \epsilon \| w_m\|^2_{L^2(0,{\tilde T}; H^1(\Omega_0;{\mathbb R}^3))} \nonumber\\
&\qquad \le
 c_{11} {\tilde T}\int_0^{\tilde T} 
\int_0^t  \|\operatorname{Def}w_m(r)\|^2_ {H^1(\Omega_0; {\mathbb R}^9)} 
dr \, dt 
+ \epsilon \| w_m\|^2_{L^2(0,{\tilde T}; H^1(\Omega_0;{\mathbb R}^3))}.
\label{06}
\end{align}
Using (\ref{02}), (\ref{03}), (\ref{04}), (\ref{05}), and (\ref{06}) we 
arrive at the basic inequality
\begin{align}
&\|w_m({\tilde T})\|_{L^2(\Omega_0; {\mathbb R}^3)}^2 
+ c_{12}\int_{\Gamma_0}\left| \int_0^{\tilde T} N\cdot \nabla_0 w_m(r)dr\right|^2 dS 
+ c_{13} \|w_m\|^2_{L^2(0,{\tilde T};H^1(\Omega_0; {\mathbb R}^3))} \nonumber \\
& \ \ \le
\|w_m(0)\|_{L^2(\Omega_0; {\mathbb R}^3)}^2 
 + c_2 \|f\|^2_{L^2(0,{\tilde T};H^1(\Omega_0; {\mathbb R}^3)')} 
+ c_3\|g\|^2_{L^2(0,{\tilde T};H^{-\frac{1}{2}}(\Gamma_0; {\mathbb R}^3))} \nonumber \\
&\ \
+ c_3\|B\|^2_{L^2(0,{\tilde T};L^2(\Gamma_0; {\mathbb R}))}  
+ c_{9} \int_0^{\tilde T} \|w_m(t)\|^2_{L^2(\Omega_0; {\mathbb R}^3)} dt
\nonumber\\
& \ \ + c_{14} \int_0^{\tilde T} \int_{\Gamma_0} \left|
\int_0^t N\cdot \nabla_0 w_m(r)\, dr\right|^2 dS\, dt 
+ c_{11}{\tilde T} \int_0^{\tilde T} \int_0^t \|w_m(r)\|^2_{H^1(\Omega_0; {\mathbb R}^3)} 
dr\, dt \, .
\label{basicenergy}
\end{align}

Letting 
\begin{align}
y_m({\tilde T})&:= \int_0^{\tilde T}\left[ \|w_m(t)\|^2_{L^2(\Omega_0; {\mathbb R}^3)} 
+ \int_{\Gamma_0} \left| \int_0^t N\cdot \nabla_0w_m\, dr\right|^2 \, dS\, dt
\right. \nonumber \\
& \qquad\qquad \left.
+ \int_0^t\|w_m\|^2_{H^1(\Omega_0; {\mathbb R}^3)} dr
\right] dt \label{y} \\
\phi({\tilde T})&:=  
c_2 \|f\|^2_{L^2(0,{\tilde T};H^1(\Omega_0; {\mathbb R}^3)')} 
+ c_3\|g\|^2_{L^2(0,{\tilde T};H^{-\frac{1}{2}}(\Gamma_0; {\mathbb R}^3))} 
+ c_3\|B\|^2_{L^2(0,{\tilde T};L^2(\Gamma_0; {\mathbb R}))}\,, \label{phi}
\end{align}
we may choose new constants $\tilde c_1$, $\tilde c_2$, and $\tilde c_3$,
 and  write (\ref{basicenergy}) as the differential inequality
\begin{align}
y'_m({\tilde T}) &\le [\tilde c_2 {\tilde T} + \tilde c_1]y_m({\tilde T}) 
+ \tilde c_3 \phi(\tilde T)
+ \|w_m(0)\|^2_{L^2(\Omega_0; {\mathbb R}^3)} \nonumber \\
&\le [\tilde c_2 \tilde T + \tilde c_1]y_m(\tilde T) + \tilde c_3 \phi(\tilde T)
+ \|w_0\|^2_{L^2(\Omega_0; {\mathbb R}^3)}.
\label{y'}
\end{align}
It follows that
\begin{align*}
y_m(\tilde T) &\le C e^{\tilde c_1 \tilde T + \tilde c_2 \tilde T^2} 
\int_0^{\tilde T} (\tilde c_3\phi(s)+
\|w_0\|^2_{L^2(\Omega_0; {\mathbb R}^3)}) ds \\
&\le CT e^{\tilde c_1 T + \tilde c_2 T^2} (\tilde c_3 \phi(T)
+\|w_0\|^2_{L^2(\Omega_0; {\mathbb R}^3)}),
\end{align*}
so that 
\begin{align}
y'_m({\tilde T})&\le \left[\tilde c_3 + (\tilde c_4 T + \tilde c_5 T^2) 
e^{\tilde c_1 T + \tilde c_2 T^2} \right]\phi(T)  \nonumber\\
& \qquad + \left[1 + (\tilde c_4 T + \tilde c_5 T^2) 
e^{\tilde c_1 T + \tilde c_2 T^2} \right] 
\|w_0\|^2_{L^2(\Omega_0; {\mathbb R}^3)}\,, \label{y'_gronwall}
\end{align}
for $0\le \tilde T \le T$.
Since $y'_m(\tilde T)$ implies the left-hand-side of inequality (\ref{energy1}),
this proves the theorem.
\end{proof}

\subsubsection{Existence}

We first have the following 
\begin{theorem}
\label{nonuniqueweak}
For $w_0\in L^2_{div}$, there exists a weak solution of (\ref{linear0}) satisfying
\begin{align}
&\sup_{0\le t\le T} \left[\|w(t)\|^2_{L^2(\Omega_0; {\mathbb R}^3)}
+ \|N\cdot \nabla_0 \int_0^t w(s)ds\|^2_{L^2(\Gamma_0;{\mathbb R}^2)} \right]
 + \|w\|^2_{L^2(0,T; H^1(\Omega_0; {\mathbb R}^3))} 
\nonumber\\
&  \le C(T)\left[ \|w_0\|^2_{L^2(\Omega_0;{\mathbb R}^3)} 
+ \|f\|^2_{L^2(0,T;H^1(\Omega_0;{\mathbb R}^3))'}\right.  
\nonumber \\ 
& \qquad\qquad   \left.  
+ \|B\|^2_{L^2(0,T;H^{\frac{1}{2}}(\Gamma_0;{\mathbb R}))}
+ \|g\|^2_{L^2(0,T;H^{-\frac{1}{2}}(\Gamma_0;{\mathbb R}^3))}
 \right]\ .
\label{nonuniqueenergy}
\end{align}

\end{theorem}
\begin{proof}
We first pass to the weak limit as $m\rightarrow \infty$.  The a priori
bounds (\ref{energy1}) show that there exists a subsequence $\{w_{m_l}\}$
such that
\begin{equation}\label{weakconvergence}
w_{m_l} \rightharpoonup w \ \ \text{ in } \ \ {\mathcal V} (T).
\end{equation}
Since equation (\ref{16}) holds for a.e. $0\le t\le T$, we multiply by a
function $\xi \in  {\mathcal D}([0,T])$, 
choose a linear combination of $\{w_k\}_{k=1}^p$, with $p \le m$, which
we denote by $\Psi_p$, and integrate by parts on the time derivative to obtain
\begin{align}
& -\int_0^T(w_{m_l}, \xi'(t)\Psi_p)_{L^2(\Omega_0;{\mathbb R}^3)} dt
+ \frac{ \nu }{2} 
\int_0^T (\operatorname{Def}w_{m_l}, \xi(t)\operatorname{Def}\Psi_p)
_{L^2(\Omega_0;{\mathbb R}^9)} dt \nonumber \\
& \qquad = \int_0^T (\tilde f, \xi(t)\Psi_p)_{H^1(\Omega_0;{\mathbb R}^3)}dt  
+\int_0^T (\tilde g, \xi(t)\Psi_p)_{H^{\frac{1}{2}}(\Gamma_0;{\mathbb R}^3)}dt  
\nonumber \\
& \qquad \qquad
+ \sigma \int_0^T\int_{\Gamma_0} B N  \cdot \xi(t) \Psi_p \, dS dt \nonumber \\
& \qquad \qquad
-\sigma \int_0^T \int_{\Gamma_0} \partial_\alpha \int_0^t w_{m_l}^i(s,\cdot) 
ds  \,
g_0^{\alpha\beta}\, \xi(t) \partial_\beta (N^j(\Psi_p)^jN^i) \, dS dt \,.
\label{integrate16}
\end{align}
From (\ref{weakconvergence}), we see that the limit $w$ satisfies
\begin{align}
& \int_0^T(w', \xi(t)\Psi_p)_{L^2(\Omega_0;{\mathbb R}^3)} dt
+ \frac{ \nu }{2}  
\int_0^T (\operatorname{Def}w, \xi(t)\operatorname{Def}\Psi_p)
_{L^2(\Omega_0;{\mathbb R}^9)} dt \nonumber \\
& \qquad = \int_0^T (\tilde f, \xi(t)\Psi_p)_{H^1(\Omega_0;{\mathbb R}^3)}dt  
+ \int_0^T (\tilde g, \xi(t)\Psi_p)_{H^{\frac{1}{2}}(\Gamma_0;
{\mathbb R}^3)}dt   \nonumber \\
& \qquad \qquad
+ \sigma \int_0^T\int_{\Gamma_0} B \cdot \xi(t) \Psi_p \, dS dt \nonumber \\
& \qquad \qquad
-\sigma \int_0^T \int_{\Gamma_0} \partial_\alpha \int_0^t w^i(s,\cdot) 
ds  \,
g_0^{\alpha\beta}\, \xi(t) \partial_\beta (N^j(\Psi_p)^jN^i) \, dS dt \,.
\label{limitlaw}
\end{align}
By density of the $\Psi_p$ in ${\mathcal V}$, (\ref{limitlaw}) holds
for any $v\in{\mathcal V}$ replacing $\Psi_p$.   By denseness of
${\mathcal D}([0,T]) \otimes {\mathcal V}$ in $L^2(0,T; {\mathcal V})$,
we see that $w_t:= w' \in L^2(0,T; {\mathcal V}')$, and by the inequality
(\ref{energy1}), $w \in {\mathcal V}(T)$, so that condition (i) of 
Definition \ref{def1} is satisfied. 

Let $y$ be given by (\ref{y}) with $w$ replacing $w_m$; we show
that $y'$ satisfies (\ref{y'_gronwall}).  Fix $\tilde T\in [0,T]$, choose
$\delta>0$ small, and integrate (\ref{energy1}) from $\tilde T$ to $\tilde T
+ \delta$ to find
\begin{align}
& \int_{\tilde T}^{\tilde T+\delta}
\left[\|w_m(t)\|^2_{L^2(\Omega_0; {\mathbb R}^3)}
+ \|N\cdot \nabla_0 \int_0^t w_m(s)ds\|^2_{L^2(\Gamma_0;{\mathbb R}^2)} 
+ \|w_m\|^2_{L^2(0,T; H^1(\Omega_0; {\mathbb R}^3))} \right]dt
\nonumber \\
&\qquad \le
 \int_{\tilde T}^{\tilde T+\delta} \left\{
 (\tilde c_3 + C(T))\left[ 
 \|f\|^2_{L^2(0,T;H^1(\Omega_0;{\mathbb R}^3)')} 
    + \|B\|^2_{L^2(0,T;L^2(\Gamma_0;{\mathbb R}))}\right. \right. \nonumber \\
&\qquad \qquad \qquad\left.\left.
+ \|g\|^2_{L^2(0,T;H^{-\frac{1}{2}}(\Gamma_0;{\mathbb R}^3))} \right]
+ (1 + C(T)) \|w_0\|^2_{L^2(\Omega_0; {\mathbb R}^3)} \right\} dt \,.
\label{energy2}
\end{align}
Since $w_m \rightharpoonup w$ in ${\mathcal V}(T)$, by lower semi-continuity
of weak convergence, we have that 
\begin{align}
& \int_{\tilde T}^{\tilde T+\delta}
\left[\|w(t)\|^2_{L^2(\Omega_0; {\mathbb R}^3)}
+ \|N\cdot \nabla_0 \int_0^t w(s)ds\|^2_{L^2(\Gamma_0;{\mathbb R}^2)} 
+ \|w\|^2_{L^2(0,T; H^1(\Omega_0; {\mathbb R}^3))} \right]dt
\nonumber\\
& \qquad \le \liminf_{m \rightarrow \infty}\int_{\tilde T}^{\tilde T+\delta}
\left[\|w_m(t)\|^2_{L^2(\Omega_0; {\mathbb R}^3)}
+ \|N\cdot \nabla_0 \int_0^t w_m(s)ds\|^2_{L^2(\Gamma_0;{\mathbb R}^2)}\right.
\nonumber\\
&\qquad \qquad\qquad\qquad\qquad\left.
+ \|w_m\|^2_{L^2(0,T; H^1(\Omega_0; {\mathbb R}^3))} \right]dt\,. \label{gg}
\end{align}
Since for arbitrary $\delta$ the inequality $\int_{\tilde T}^{\tilde T +\delta}
(f(t)-g(t)) dt \le 0$ implies that $f(\tilde T) - g(\tilde T) \le 0$;
putting together (\ref{gg}) and (\ref{energy2}) shows that
\begin{align}
& \|w(\tilde T)\|^2_{L^2(\Omega_0; {\mathbb R}^3)}
+ \|N\cdot \nabla_0 \int_0^{\tilde T}w(s)ds\|^2_{L^2(\Gamma_0;{\mathbb R}^2)} 
+ \|w\|^2_{L^2(0,T; H^1(\Omega_0; {\mathbb R}^3))} \nonumber\\
&\qquad \le
 (\tilde c_3 + C(T))\left[ 
 \|f\|^2_{L^2(0,T;H^1(\Omega_0;{\mathbb R}^3)')} 
    + \|B\|^2_{L^2(0,T;L^2(\Gamma_0;{\mathbb R}))}\right. \nonumber \\
&\qquad \qquad \qquad\left.
+ \|g\|^2_{L^2(0,T;H^{-\frac{1}{2}}(\Gamma_0;{\mathbb R}^3))} \right]
+ (1 + C(T)) \|w_0\|^2_{L^2(\Omega_0; {\mathbb R}^3)} \,, \label{energyw}
\end{align}
for all $\tilde T \in [0,T]$, which establishes (\ref{nonuniqueenergy}). 

Next, we address the issue of the initial condition.  From (\ref{energyw}),
there is an $M>0$ such that
$$\sup_{t\in[0,T]} \|w(t)\|_{{L^2(\Omega_0; {\mathbb R}^3)}}< M,$$ 
so for each sequence $\{t_n\}_{n=1}^\infty $ with $t_n\rightarrow 0$ as
$n \rightarrow  \infty$, there is a subsequence $\{t_{n_j}\}$ such that
$w(t_{n_j}) \rightharpoonup W$ in $ {L^2(\Omega_0; {\mathbb R}^3)}$ as
$n_j\rightarrow \infty$.  Since $w_t \in L^2(0,T;{\mathcal V}')$, it follows 
that  $w \in C^0([0,T];{\mathcal V}')$, 
so that $\langle w(t_{n_j}), \psi\rangle_{\mathcal V}\rightarrow 
\langle w(0), \psi\rangle_{\mathcal V}$ as $n_j\rightarrow \infty$ for all
$\psi \in {\mathcal V}$. $\langle \cdot, \cdot \rangle_{\mathcal V}$ denotes
the duality pairing between ${\mathcal V}$ and $ {\mathcal V}'$.
Thus $w(0)=W$ in $ {L^2(\Omega_0; {\mathbb R}^3)}$. 

\begin{remark}
The inequality (\ref{y'_gronwall}) shows that 
$$\limsup_{T\rightarrow 0} \|w(T)\|^2_ {L^2(\Omega_0; {\mathbb R}^3)} \le 
\|w_0\|^2_ {L^2(\Omega_0; {\mathbb R}^3)}\,,$$
from which we may infer that $w(t)
\rightarrow w(0)$ strongly in ${L^2(\Omega_0; {\mathbb R}^3)}$ as 
$t\rightarrow 0$. 
\end{remark}

To show that $w(0)= w_0$, we first note from (\ref{limitlaw}), that
\begin{align}
& -\int_0^T(w, v')_{L^2(\Omega_0;{\mathbb R}^3)} dt
+ \frac{ \nu }{2}  \int_0^T (\operatorname{Def}w, \operatorname{Def}v)
_{L^2(\Omega_0;{\mathbb R}^9)} dt = 
\int_0^T (\tilde f, v)_{H^1(\Omega_0;{\mathbb R}^3)}dt  \nonumber  \\
& \qquad
+\int_0^T (\tilde g, v)_{H^{\frac{1}{2}}
(\Gamma_0;{\mathbb R}^3)}dt  
+ \sigma \int_0^T\int_{\Gamma_0} B N \cdot v \, dS dt  \nonumber \\
& \qquad
-\sigma \int_0^T \int_{\Gamma_0} \partial_\alpha \int_0^t w^i(s,\cdot) 
ds  \,
g_0^{\alpha\beta}\, \partial_\beta (N^jv^jN^i) \, dS dt \,. 
+ ( w(0), v(0) )_{L^2(\Omega_0; {\mathbb R}^3)} \label{limit1}
\end{align}
for each $v \in C^1([0,T]; {\mathcal V})$ with $v(T)=0$.	
Similarly, from (\ref{integrate16}) we deduce
\begin{align*}
& -\int_0^T(w_{m_l}, v')_{L^2(\Omega_0;{\mathbb R}^3)} dt
+ \frac{ \nu }{2}  \int_0^T (\operatorname{Def}w_{m_l}, \operatorname{Def}v)
_{L^2(\Omega_0;{\mathbb R}^9)} dt = 
\int_0^T (\tilde f, v)_{L^2(\Omega_0;{\mathbb R}^3)}dt   \\
& \qquad
+\int_0^T (\tilde g, v)_{H^{\frac{1}{2}}(\Omega_0;{\mathbb R}^3)}dt  
+ \sigma \int_0^T\int_{\Gamma_0} B N \cdot v \, dS dt  \\
& \qquad
-\sigma \int_0^T \int_{\Gamma_0} \partial_\alpha \int_0^t w_{m_l}^i(s,\cdot) 
ds  \,
g_0^{\alpha\beta}\, \partial_\beta (N^jv^jN^i) \, dS dt
+ ( w_{m_l}(0), v(0) )_{L^2(\Omega_0; {\mathbb R}^3)} \,.
\end{align*}
Using (\ref{weakconvergence}) and the fact that $w_{m_l}(0)\rightarrow
w_0$ in $ {L^2(\Omega_0; {\mathbb R}^3)}$, we find that 
\begin{align}
& -\int_0^T(w, v')_{L^2(\Omega_0;{\mathbb R}^3)} dt
+ \frac{ \nu }{2}  \int_0^T (\operatorname{Def}w, \operatorname{Def}v)
_{L^2(\Omega_0;{\mathbb R}^3)} dt = 
\int_0^T (\tilde f, v)_{H^1(\Omega_0;{\mathbb R}^3)}dt  \nonumber \\
& \qquad
+\int_0^T (\tilde g, v)_{H^{\frac{1}{2}}(\Gamma_0;{\mathbb R}^3)}dt  
+ \sigma \int_0^T\int_{\Gamma_0} B N \cdot v \, dS dt  \nonumber \\
& \qquad
-\sigma \int_0^T \int_{\Gamma_0} \partial_\alpha \int_0^t w^i(s,\cdot) 
ds  \,
g_0^{\alpha\beta}\, \partial_\beta (N^jv^jN^i) \, dS dt 
+ ( w_0, v(0) )_{L^2(\Omega_0; {\mathbb R}^3)} \label{limit2}\, .
\end{align}
Since $v(0)$ is arbitrary, comparing (\ref{limit1}) with (\ref{limit2}),
we conclude that $w(0)=w_0$, which verifies condition (ii) of Definition
\ref{def1}; thus, $w$ is a weak solution.

%To prove uniqueness, it suffices to check that the only weak solution of
%(\ref{linear0}) with $f\equiv B\equiv g\equiv w_0 \equiv 0$ is $w \equiv 0$.
%We use the basic inequality (\ref{energyw});
%since $\phi=0$  and $w_0=0$, we find that $w=0$.
\subsection{Pressure as a Lagrange multiplier}
Let ${\mathcal W}= \{ \psi \in H^1(\Omega_0, {\mathbb R}^3) \ | \ 
N \cdot \nabla_0\psi \in L^2(\Gamma_0; {\mathbb R}^3) \}$ and
$H = \{ f \in{L^2(\Omega_0; {\mathbb R})} \ | \ \int_{\Omega_0}f(x)dx =0\}$.
\begin{lemma}\label{lemma_lagrange}
For all $p \in H$, there exists a constant
$C>0$ and
$v \in {\mathcal W}$ such that $\operatorname{div} v =p$ and
\begin{equation}\label{v-p} 
\|v\|^2_{H^1(\Omega_0; {\mathbb R}^3)} +  \|N\cdot \nabla_0 v\|^2_
{L^2(\Gamma_0; {\mathbb R}^2)} \le  C\|p\|^2_{L^2(\Omega_0; {\mathbb R})}. 
\end{equation}
\end{lemma}
\begin{proof}
We let $v=\nabla f$ and solve $-\triangle f= -p$ in $\Omega_0$ with
$\partial f/\partial N = 0$ on $\Gamma_0$.  Then $\operatorname{div}v=p$,
and by the standard elliptic estimate $\|v\|^2_{H^1(\Omega_0; {\mathbb R}^3)}
\le C \|p\|^2_{L^2(\Omega_0; {\mathbb R})}$. Since $N\cdot v=0$ on
$\Gamma_0$, then 
$$\|N \cdot \nabla_0 v\|^2_ {L^2(\Gamma_0; {\mathbb R}^2)} 
= \|\nabla_0 N \cdot v \|^2_{L^2(\Gamma_0; {\mathbb R}^2)} 
\le C \|v\|^2_ {H^1(\Omega_0; {\mathbb R}^3)},$$
which establishes (\ref{v-p}).
\end{proof}
We can now follow \cite{SolSca1973}.
Define the linear functional on ${\mathcal W}$
by $(p,\operatorname{div}v)_ {L^2(\Omega_0; {\mathbb R})}$, where $p\in H$.
By the Riesz representation theorem, there is a bounded linear operator
$Q: H\rightarrow  {\mathcal W}$ such that
$$ 
(p,\operatorname{div}v)_ {L^2(\Omega_0; {\mathbb R})}=
(Qp, v)_{\mathcal W}, 
$$ 
where $(\cdot, \cdot)_{\mathcal W}$ denotes the inner-product on 
${\mathcal W}$ (the left-hand-side of (\ref{v-p}) defines the norm on
$ {\mathcal W}$).  
Letting $v=Qp$ shows that 
\begin{equation}\label{Qp0}
\|Qp\|_{\mathcal V} \le C
\|p \|_ {L^2(\Omega_0; {\mathbb R})}
\end{equation}
for some constant $C>0$. 
Using Lemma \ref{lemma_lagrange},  we have the estimate
\begin{equation}\label{Qp}
\|p\|^2_{L^2(\Omega_0; {\mathbb R})} = (p, \operatorname{div} v)
_ {L^2(\Omega_0; {\mathbb R}^3)}\le C\|Qp\|_{\mathcal V}\|v\|_{\mathcal V} 
\le C\|Qp\|_{\mathcal V}\|p\|_{L^2(\Omega_0; {\mathbb R})}. 
\end{equation}
It follows that 
\begin{equation}\label{hodge}
{\mathcal W} = R(Q) \oplus_ {\mathcal W} {\mathcal V}.
\end{equation}
To see this, suppose that $v \in {\mathcal W}\ominus R(Q)$. Then for any
$p \in H$, $(Qp,v)_{\mathcal W}=(p,\operatorname{div}v)_
{L^2(\Omega_0; {\mathbb R}^3)} =0$, so that $\operatorname{div} v$ is a 
constant, and since $v \cdot N=0$ on $\Gamma_0$, the constant must equal
to zero; thus, $v \in {\mathcal V}$.  Finally, by (\ref{Qp0}) and (\ref{Qp}),
$R(Q)$ is closed in ${\mathcal W}$. 
\end{proof}

\begin{lemma} \label{pressure}
A solution $w\in {\mathcal V}(T)$ of (\ref{linear0}) satisfies for a.e. $0\le t \le T$,
\begin{align}
& \langle w_t, v\rangle 
+ \frac{ \nu }{2} 
(\operatorname{Def}w, \operatorname{Def}v)_{L^2(\Omega_0;{\mathbb R}^9)} 
-(p,\operatorname{div} v)_{L^2(\Omega_0; {\mathbb R})}=  \langle f, v\rangle 
+ \langle g , v \rangle_{\Gamma_0} \nonumber  \\
&\ \  + \sigma  \int_{\Gamma_0} BN\cdot v \, dS 
-\sigma \int_{\Gamma_0} \int_0^t \partial_\alpha w^i(s,\cdot) ds  \,
g_0^{\alpha\beta}\, \partial_\beta (N^jv^jN^i) \, dS \ \ \ \forall v\in 
{\mathcal W}\,,\ \label{weakp}
\end{align}
where $p(t,\cdot) \in H$ for a.e. $0\le t\le T$ is termed the pressure 
function, with the energy inequality
\begin{align}
&\sup_{0\le t\le T} \left[\|w(t)\|^2_{L^2(\Omega_0; {\mathbb R}^3)}
+ \|N\cdot \nabla_0 \int_0^t w(s)ds\|^2_{L^2(\Gamma_0;{\mathbb R}^2)}  \right]
\nonumber\\ 
&\ + \|w\|^2_{L^2(0,T; H^1(\Omega_0; {\mathbb R}^3))} 
 + \|p\|^2_{L^2(0,T; L^2(\Omega_0; {\mathbb R}))}  
\nonumber\\
&  \le C(T)\left[ \|w_0\|^2_{L^2(\Omega_0;{\mathbb R}^3)} 
+ \|f\|^2_{L^2(0,T;H^1(\Omega_0;{\mathbb R}^3)')}  
+ \|B\|^2_{L^2([0,T];H^{\frac{1}{2}}(\Gamma_0;{\mathbb R}))}\right.
\nonumber \\ 
& \qquad\qquad   \left.  
+ \|g\|^2_{L^2(0,T;H^{-\frac{1}{2}}(\Gamma_0;{\mathbb R}^3))} \right]\ .
\label{nonuniqueenergyp}
\end{align}
\end{lemma}
\begin{proof}
From (\ref{linear0}) and (\ref{nonuniqueenergy}), we first see that $w_t \in
L^2 (0,T; {\mathcal W}') $, with
\begin{align*}
  \|w_t\|^2_{L^2(0,T; H^1(\Omega_0; {\mathbb R}^3)')} 
 \le & C \left[ \|w\|^2_{L^2(0,T;\mathcal W)}  
+ \|f\|^2_{L^2(0,T;H^1(\Omega_0;{\mathbb R}^3)')}\nonumber\right.\\
&\qquad\qquad\left. 
+ \|B\|^2_{L^2([0,T];H^{\frac{1}{2}}(\Gamma_0;{\mathbb R}))}
+ \|g\|^2_{L^2(0,T;H^{-\frac{1}{2}}(\Gamma_0;{\mathbb R}^3))}
 \right]\ .
\end{align*}
By the decomposition (\ref{hodge}), for $v\in{\mathcal W}$, we let
$v=v_1+v_2$, and $v_1 \in {\mathcal V}$ and $v_2 \in R(Q)$.  We define
$\Lambda\in {\mathcal W}'$ as the difference of the left- and right-hand-sides
of condition (i) in Definition \ref{def1}; 
then $\Lambda\equiv 0$ on ${\mathcal V}$.  It follows that
$$ 
\Lambda(v) = \Lambda(v_2) = ( \psi, v_2)_{\mathcal W} 
= ( \psi, v)_{\mathcal W} \ \text{ for } \ \psi \in R(Q).
$$ 
From Lemma \ref{lemma_lagrange}, for a.e. $t\in[0,T]$, 
$\Lambda(v) = (p,\operatorname{div} v)_{L^2(\Omega_0; {\mathbb R}^3)}$, which establishes (\ref{weakp}), and with the estimate
\begin{align*}
\|p\|_{L^2(\Omega_0; {\mathbb R}^3)}\le & 
 C\left[ \|w_t\|^2_{ H^1(\Omega_0; {\mathbb R}^3))} 
 + \|w\|^2_{L^2(0,T;\mathcal W)}\nonumber\right. \\  
&\qquad\qquad\left.+ \|f\|^2_{H^1(\Omega_0;{\mathbb R}^3)'}+ \|B\|^2_{H^{\frac{1}{2}}(\Gamma_0;{\mathbb R})}
+ \|g\|^2_{H^{-\frac{1}{2}}(\Gamma_0;{\mathbb R}^3)}
 \right]\ ,
\end{align*}
which by integration from $0$ to $T$ gives (\ref{nonuniqueenergyp}).
\end{proof}

While we have proven that $w_t\in L^2(0,T;{\mathcal V}')$, we have not shown
that $w\in L^2(0,T;{\mathcal V})$, and so we cannot employ the 
standard methods to prove uniqueness. 
Using the assumptions on our data for the nonlinear problem, however, we 
have additional regularity for $w_t$, and so uniqueness follows from
 
\begin{theorem}
\label{uniqueweak}
Suppose that
\begin{equation}\label{unique_data}
\begin{array}{c}
f\in L^2 (0,T;  H^1(\Omega_0; {\mathbb R}^3))), \ 
g\in L^2 (0,T;  H^{\frac{3}{2}}(\Gamma_0; {\mathbb R}^3))), \\
f_t\in L^2 (0,T;  H^1(\Omega_0; {\mathbb R}^3))'), \
g_t\in L^2 (0,T; H^{-\frac{1}{2}}(\Gamma_0;{\mathbb R}^3)), \\ 
B_t\in L^2 (0,T; H^{\frac{1}{2}}(\Gamma_0;{\mathbb R}^3))\  B(0)=0, 
\end{array}
\end{equation} 
and that the initial data 
$w_0\in H^2(\Omega_0;{\mathbb R}^3)\cap \mathcal V$ 
satisfies the compatibility condition 
$[\nu \operatorname{Def} (w_0) -g(0)]_{\operatorname{tan}}=0$. 
Then, there exists a unique solution $w$ in $L^2 (0,T;\mathcal V)$ 
of (\ref{linear0}) such that $w_t\in \mathcal V (T)$.
\end{theorem}

\begin{proof} 
Given the regularity assumptions in (\ref{unique_data}), compactness implies 
that
$$f\in C(0,T;L^2 (\Omega_0;{\mathbb R}^3)), \ \
g\in C(0,T;H^{\frac{1}{2}}(\Gamma_0;{\mathbb R}^3)),
$$ 
so that
$f(0)\in L^2(\Omega_0;{\mathbb R}^3)$ and $g(0)\in H^{\frac{1}{2}}(\Gamma_0;{\mathbb R}^3)$.  
Since $w_0\in H^2(\Omega_0;{\mathbb R}^3)$, let $p_0\in H^1(\Omega_0;{\mathbb R}^3)$ be the solution of the Dirichlet problem
\begin{subequations}
  \label{p0}
\begin{gather}
-\triangle p_0 =  -\operatorname{div}( \nu  \triangle  w_0 + f(0)) 
\ \ \text{ in } \ \ \Omega_0,  
         \label{p0.a} \\
p_0 =[\nu\operatorname{Def}w_0 \ N -g(0) ]\cdot N 
\ \ \text{ on } \ \ \Gamma_0\ .
         \label{p0.b}
\end{gather}
\end{subequations}
From our initial compatibility condition $(\nu \operatorname{Def} (w_0) -g(0))_{\operatorname{tan}}=0$, we see that (\ref{p0.b}) implies $$p_0 N=\nu\operatorname{Def}w_0 \ N -g(0)\ ,$$ which will be used later in the proof.
  
Let $\bar{w}_0\in L^2(\Omega_0;{\mathbb R}^3)$ be defined by
$$\bar{w}_0 =\nu \triangle w_0 -\nabla p_0 + f(0)\ .$$
Thanks to (\ref{p0.a}), $\bar w_0 \in L^2_{div}$, and 
standard elliptic estimates for (\ref{p0}) show that there exists a constant
$C>0$ such that
\begin{align}
\| \bar{w}_0 \|^2_{L^2(\Omega_0; {\mathbb R}^3)} & + \| p_0 \|^2_{H^1(\Omega_0; {\mathbb R}^3)} \nonumber\\
& \le C \left(
\| w_0\|^2_{H^2(\Omega_0; {\mathbb R}^3)}
+\| f(0)\|^2_{L^2(\Omega_0; {\mathbb R}^3)} + \| g(0)\|^2_{H^{\frac{1}{2}}(\Gamma_0; {\mathbb R}^3)} \right)\, .
\label{wt(0)}
\end{align}
Now, since $\bar w_0\in L^2_{div}$, from Theorem \ref{nonuniqueweak}, let us define $\bar w\in {\mathcal V} (T)$ to be a solution of
\begin{align*}
&\operatorname{(i)} \  \langle \bar w_t, v\rangle 
 +\frac{ \nu }{2} (\operatorname{Def}\bar w, \operatorname{Def}v)_{L^2(\Omega_0;{\mathbb R}^9)} 
=  \langle f_t, v\rangle 
+ \langle g_t , v \rangle_{\Gamma_0}  \\
&\ \  +  \sigma \int_{\Gamma_0} B_t N\cdot v \, dS 
-\sigma \int_{\Gamma_0} \int_0^t \partial_\alpha \bar w^i(s,\cdot)ds\  g_0^{\alpha\beta}\, \partial_\beta (N^jv^jN^i) dS \\
&\ \ +  \sigma \int_{\Gamma_0}  \partial_\alpha  w^i_0 \ 
g_0^{\alpha\beta}\, \partial_\beta (N^jv^jN^i) \, dS \ \ \ \forall v\in 
{\mathcal V},\ \  \text{a.e. in}\ (0,T),\\
\operatorname{and} \\
&\operatorname{(ii)} \  \bar w(0,\cdot) = \bar w_0, 
\end{align*}
with the energy estimate 
\begin{align}
&\sup_{0\le t\le T} \left[\|\bar w(t)\|^2_{L^2(\Omega_0; {\mathbb R}^3)}
+ \|N\cdot \nabla_0 \int_0^t \bar w(s)ds\|^2_{L^2(\Gamma_0;{\mathbb R}^2)} \right]
 + \|\bar w\|^2_{L^2(0,T; H^1(\Omega_0; {\mathbb R}^3))} 
\nonumber\\
&  \le C(T)\left[ \|\bar w_0\|^2_{L^2(\Omega_0;{\mathbb R}^3)} 
+ \|f_t\|^2_{L^2(0,T;H^1(\Omega_0;{\mathbb R}^3)')}\right.  
\nonumber \\ 
& \qquad\qquad   \left.  
+ \|B_t\|^2_{L^2(0,T;H^{\frac{1}{2}}(\Gamma_0;{\mathbb R}))}
+ \|g_t\|^2_{L^2(0,T;H^{-\frac{1}{2}}(\Gamma_0;{\mathbb R}^3))}
+ \|\nabla_0 w_0\|^2_{H^{\frac{1}{2}}(\Gamma_0;{\mathbb R}^6)}
 \right]\ .
\label{nonuniqueenergybis}
\end{align}

Now, let $\bar p$ be given by Lemma \ref{pressure}, and let us define
$$w(t,x)=w_0(x)+\int_0^t \bar w (s,x) ds\ ,$$
$$p(t,x)=p_0(x)+\int_0^t \bar p (s,x) ds\ .$$
Since $w_t=\bar{w}$, we then have $w_t (0)=\bar{w}(0)=\bar w_0$. We also obviously
have $w(0)=w_0$ and $p(0)=p_0$. Concerning the regularity, we see that $w\in
L^2 (0,T; \mathcal V)$ (since $w_t=\bar w\in \mathcal V (T)$).

By Lemma \ref{pressure}, we know that a.e. in $(0,T)$,
\begin{align*}
& \langle \bar w_t, v\rangle 
+ \frac{ \nu }{2} 
(\operatorname{Def}\bar w, \operatorname{Def}v)_{L^2(\Omega_0;{\mathbb R}^9)} 
-(\bar p,\operatorname{div} v)_{L^2(\Omega_0; {\mathbb R})}=  \langle f_t, v\rangle 
+ \langle g_t , v \rangle_{\Gamma_0} \nonumber  \\
&\ \  + \sigma  \int_{\Gamma_0} BN\cdot v \, dS 
-\sigma \int_{\Gamma_0} \int_0^t \partial_\alpha \bar w^i(s,\cdot) ds  \,
g_0^{\alpha\beta}\, \partial_\beta (N^jv^jN^i) \, dS \nonumber\\
&\ \ + \sigma  \int_{\Gamma_0}  \partial_\alpha w_0^i  \,
g_0^{\alpha\beta}\, \partial_\beta (N^jv^jN^i) \, dS,  \ \ \forall v\in 
{\mathcal W}\ .
\end{align*}
By integrating in time this equality from $0$ to $t\in (0,T)$, we find that
\begin{align*}
& \langle w_t -\bar w_0, v\rangle 
+ \frac{ \nu }{2} 
(\operatorname{Def} (w-w_0), \operatorname{Def}v)_{L^2(\Omega_0;{\mathbb R}^9)} 
-( p -p_0,\operatorname{div} v)_{L^2(\Omega_0; {\mathbb R})}\nonumber\\
& =  \langle f -f(0), v\rangle 
+ \langle g -g(0) , v \rangle_{\Gamma_0} \nonumber  \\
&\ \  + \sigma  \int_{\Gamma_0} BN\cdot v \, dS 
-\sigma \int_{\Gamma_0} \int_0^t \partial_\alpha  (w^i(s,\cdot)-w^i_0) ds  \,
g_0^{\alpha\beta}\, \partial_\beta (N^jv^jN^i) \, dS \nonumber\\
&\ \ + \sigma t \int_{\Gamma_0}  \partial_\alpha  w_0^i  \,
g_0^{\alpha\beta}\, \partial_\beta (N^jv^jN^i) \, dS,  \ \ \forall v\in 
{\mathcal W}\ ,
\end{align*}
and then using the definitions of $\bar w_0$ and $p_0$, we find that
\begin{align*}
& \langle  w_t, v\rangle 
+ \frac{ \nu }{2} 
(\operatorname{Def} w, \operatorname{Def}v)_{L^2(\Omega_0;{\mathbb R}^9)} 
-( p,\operatorname{div} v)_{L^2(\Omega_0; {\mathbb R})}=  \langle f, v\rangle 
+ \langle g , v \rangle_{\Gamma_0} \nonumber  \\
&\ \  + \sigma  \int_{\Gamma_0} BN\cdot v \, dS 
-\sigma \int_{\Gamma_0} \int_0^t \partial_\alpha  w^i(s,\cdot) ds  \,
g_0^{\alpha\beta}\, \partial_\beta (N^jv^jN^i) \, dS, \ \ \forall v\in 
{\mathcal W}\ .
\end{align*}
Consequently, we immediately get that $w\in L^2 (0,T; \mathcal V)$ satisfies 
\begin{align*}
& \langle  w_t, v\rangle 
+ \frac{ \nu }{2} 
(\operatorname{Def} w, \operatorname{Def}v)_{L^2(\Omega_0;{\mathbb R}^9)} 
=  \langle f, v\rangle 
+ \langle g , v \rangle_{\Gamma_0} + \sigma  \int_{\Gamma_0} BN\cdot v \, dS  \nonumber  \\
&\ \  -\sigma \int_{\Gamma_0} \int_0^t \partial_\alpha  w^i(s,\cdot) ds  \,
g_0^{\alpha\beta}\, \partial_\beta (N^jv^jN^i) \, dS, \ \ \forall v\in 
{\mathcal V}\ ,
\end{align*}
which with the condition $w(0)=w_0$ shows that $w$ is a solution of (\ref{linear0}). Furthermore, $w_t\in \mathcal V (T)$. From (\ref{nonuniqueenergybis}), it is
readily seen that
\begin{align}
&\sup_{0\le t\le T} \left[\|w(t)\|^2_{L^2(\Omega_0; {\mathbb R}^3)}
+ \|w_t(t)\|^2_{L^2(\Omega_0; {\mathbb R}^3)}
+ \|N\cdot \nabla_0 \int_0^t w(s)ds\|^2_{L^2(\Gamma_0;{\mathbb R}^2)} 
\right. \nonumber \\
&\
+\left. \|N\cdot \nabla_0 w(t)\|^2_{L^2(\Gamma_0;{\mathbb R}^2)} \right] 
 + \|w\|^2_{L^2(0,T; H^1(\Omega_0; {\mathbb R}^3))} 
 + \|w_t\|^2_{L^2(0,T; H^1(\Omega_0; {\mathbb R}^3))} 
\nonumber\\
&\
 + \|p\|^2_{L^2(0,T; L^2(\Omega_0; {\mathbb R}))} 
 + \|p_t\|^2_{L^2(0,T; L^2(\Omega_0; {\mathbb R}))} 
\nonumber\\
&  \le C(T)\left[ \|w_0\|^2_{H^2(\Omega_0;{\mathbb R}^3)} 
+ \|f(0)\|^2_{L^2(\Omega_0;{\mathbb R}^3)}+  
\|g(0)\|^2_{H^{\frac{1}{2}}(\Gamma_0;{\mathbb R}^3)} \right. \nonumber \\
& \qquad \qquad
+ \|f\|^2_{L^2(0,T;H^1(\Omega_0;{\mathbb R}^3))}  
+ \|f_t\|^2_{L^2(0,T;H^1(\Omega_0;{\mathbb R}^3)')} 
+ \|B\|^2_{C^0([0,T];H^{\frac{1}{2}}(\Gamma_0;{\mathbb R}))}
\nonumber \\ 
& \qquad\qquad   \left.  
+ \|B_t\|^2_{L^2(0,T;H^{\frac{1}{2}}(\Gamma_0;{\mathbb R}))}
+ \|g\|^2_{L^2(0,T;H^{\frac{3}{2}}(\Gamma_0;{\mathbb R}^3))}
+ \|g_t\|^2_{L^2(0,T;H^{-\frac{1}{2}}(\Gamma_0;{\mathbb R}^3))} \right]\ .
\label{H1energy_p}
\end{align}

Now, let us assume that there exists another solution $w'$ to 
(\ref{linear0}), such that $w'\in L^2(0,T;\mathcal V)$ and $w'_t\in 
\mathcal V(T)$. 
By denoting
$\delta w=w-w'$, we see that $\delta w\in L^2(0,T;\mathcal V)$ is solution of 
\begin{align*}
&\operatorname{(i)} \  \langle \delta w_t, v\rangle 
 +\frac{ \nu }{2} (\operatorname{Def}\delta w, \operatorname{Def}v)_{L^2(\Omega_0;{\mathbb R}^9)} 
  \\
&\ \  
= -\sigma \int_{\Gamma_0} \int_0^t \partial_\alpha \delta w^i(s,\cdot)ds\ g_0^{\alpha\beta}\, \partial_\beta (N^jv^jN^i) ds
 \ \ \ \forall v\in 
{\mathcal V},\ \ \text{a.e. in}\ (0,T),\\
\operatorname{and} \\
&\operatorname{(ii)} \  \delta w(0,\cdot) = 0 . 
\end{align*}
Since $\delta w(t\cdot) \in \mathcal V$ a.e. in $(0,T)$, we can use $\delta w$ as
a test function in (i), which gives a.e. in $(0,T)$
\begin{align*}
  \langle \delta w_t, \delta w\rangle 
 &+\frac{ \nu }{2} (\operatorname{Def}\delta w, \operatorname{Def}\delta w)_{L^2(\Omega_0;{\mathbb R}^9)} \\
& =-\sigma \int_{\Gamma_0} \int_0^t \partial_\alpha \delta w^i(s,\cdot) ds\ g_0^{\alpha\beta}\, \partial_\beta (N^j\delta w^jN^i) dS \ .
\end{align*}
Since $\delta w\in L^2 (0,T; \mathcal V)$ and $\delta w_t 
\in L^2 (0,T; {\mathcal V}')$, we have $$\frac{1}{2} \frac{d}{dt}\|\delta w\|^2_{L^2(\Omega_0;{\mathbb R}^3)}=\langle \delta w_t,\delta w\rangle\ . $$
Since $N\cdot \nabla_0 \delta w\in L^2 (0,T; L^2(\Gamma_0;{\mathbb R}^2))$ and
$N\cdot \nabla_0 \int_0^\cdot \delta w\in L^2 (0,T; L^2(\Gamma_0;{\mathbb R}^2))$, we also have  the same relation as (\ref{01}), where $\delta w$ replaces $w_m$.
In a quite similar fashion as we proved (\ref{energyw}) (except that this time there is no limit associated to a Galerkin procedure to consider) we can then infer that
\begin{align}
& \|\delta w(\tilde T)\|^2_{L^2(\Omega_0; {\mathbb R}^3)}
+ \|N\cdot \nabla_0 \int_0^{\tilde T}\delta w(s)ds\|^2_{L^2(\Gamma_0;{\mathbb R}^2)} 
+ \|\delta w\|^2_{L^2(0,T; H^1(\Omega_0; {\mathbb R}^3))} \nonumber\\
&\qquad \le
 (1 + C(T)) \|\delta w (0)\|^2_{L^2(\Omega_0; {\mathbb R}^3)} \,, \label{energydeltaw}
\end{align}
which, with the condition $\delta w(0)=0$, precisely proves that $\delta w=0$,
establishing the uniqueness of such a solution.

\end{proof}

\begin{remark}
We note that the estimate (\ref{H1energy_p}) holds if we replace $\Omega_0$ 
by ${\mathbb R}^3_+$. 
\end{remark}

\subsection{Regularity of weak solutions on ${\mathbb R}^3_+$}
We consider the half-space ${\mathbb R}^3_+ := \{(x^1,x^2,x^3)
\in{\mathbb R}^3 \ | \ x^3>0\}$ with the usual orthonormal 
basis $(e_1,e_2,e_3)$.  The unit normal vector to the horizontal plane
$\{ x^3=0\}$ is $N=(0,0,1)$.  

\begin{definition}\label{differences}
The first-order difference quotient of $w$ of size $h$ at $x$ is given by
$$
D_h w(x) = \frac{w(x+h)-w(x)}{|h|},
$$
where $h$ is any vector orthogonal to $N$.  The second-order difference
quotient of $w$ of size $h$ is defined as $D_{-h}D_h w(x)$, given explicitly by
$$D_{-h}D_hw(x)=\frac{w(x+h)+w(x-h)-2w(x)}{|h|^2}\,.$$
\end{definition}

\begin{lemma}\label{quotients}
Suppose that $w \in{L^2({\mathbb R}^3_+; {\mathbb R}^3)}$.

\noindent
{(i)} If 
$$\|D_h w\|_{L^2( {\mathbb R}^3_+; {\mathbb R}^3)}  \le M$$
for a constant $M$ and for $h=e_1,e_2$, then 
$$\|\nabla_0 w\|_{L^2( {\mathbb R}^3_+; {\mathbb R}^6)}  \le M\,,$$
and $\|D_h w\|_{L^2( {\mathbb R}^3_+; {\mathbb R}^3)}  \le C
\|\nabla_0 w\|_{L^2( {\mathbb R}^3_+; {\mathbb R}^6)}$ for some constant
$C$.

\noindent
{(ii)} If 
$$\|D_{-h}D_h w\|_{L^2( {\mathbb R}^3_+; {\mathbb R}^3)}  \le M$$
for a constant $M$ and for $h=e_1,e_2,(e_1+e_2)/\sqrt{2}$, then 
$$\|\nabla_0\nabla_0 w\|_{L^2( {\mathbb R}^3_+; {\mathbb R}^{12})}  \le M\,,$$
where $\nabla_0 = (\partial_1, \partial_2)$ is the gradient on $ {\mathbb R}^2
\times \{x^3=0\}$. 
\end{lemma}
\begin{proof}
Part (i) is proved on page 277 of \cite{Evans1998}, and part (ii) follows from
page 7 of \cite{CaCa1995}.
\end{proof}

On the half-space ${\mathbb R}^3_+$, the system of equations (\ref{linear0}) 
becomes
\begin{subequations}
  \label{linear7}
\begin{alignat}{2}
 w_t - \nu \triangle  w  &= -\nabla p + f 
& \ &\text{in} \ \ (0,T)\times {\mathbb R}^3_+ \,, 
         \label{linear7.a}\\
   \operatorname{div}  w &= 0
&&\text{in} \ \ (0,T)\times {\mathbb R}^3_+ \,, 
         \label{linear7.b}\\
\nu \operatorname{Def} w \ N - pN &=   \sigma 
\begin{bmatrix}
 0\\
 0\\
 \int_0^t \triangle_0 w^3(r,\cdot) dr + B
\end{bmatrix}
+ 
\begin{bmatrix}
 g^1\\
 g^2\\
 g^3 
\end{bmatrix}
&&\text{on} \ \ (0,T)\times {\mathbb R}^2 \,, 
         \label{linear7.c}\\
    w &=  w_0
&&\text{on} \ \ {\mathbb R}^3_+  \times \{t=0\}\,. 
         \label{linear7.d}
\end{alignat}
\end{subequations}
For the remainder of this section, we replace $\Omega_0$ with ${\mathbb R}^3_+$,
and $\Gamma_0$ with ${\mathbb R}^2 \times \{x^3=0\}$, in the definitions of 
our function spaces,  and we use $ \langle \cdot , \cdot \rangle_{\Gamma_0}$
to denote the duality pairing between 
$H^{\frac{1}{2} }( \mathbb{R}  ^2; \mathbb{R}  ^3)$ and its dual space
$H^{-\frac{1}{2} }( \mathbb{R}  ^2; \mathbb{R}  ^3)$.

\begin{theorem}[Regularity on ${\mathbb R}^3_+$]\label{regularity1}
If $w \in {\mathcal V}(T)$ is a weak solution of (\ref{linear7}), satisfying
the estimate (\ref{H1energy_p}), the weak formulation
\begin{align}
&(w_t, v)_{L^2({\mathbb R}^3_+;{\mathbb R}^3)} 
+ \frac{ \nu }{2}  
(\operatorname{Def}w, \operatorname{Def}v)_{L^2({\mathbb R}^3_+;{\mathbb R}^9)} 
=  \langle f, v\rangle 
+ \langle g , v \rangle_{\Gamma_0} \nonumber  \\
&\ \  +  \sigma (B, v^3)_{L^2( {\mathbb R}^2; {\mathbb R})} 
-\sigma \left( \int_0^t \nabla_0 w^3(r) dr  \ , \
\nabla_0v^3\right)_{L^2( {\mathbb R}^2; {\mathbb R}^2)} \ \ \ \forall v\in 
{\mathcal V}\,, \label{weak0}
\end{align}
as well as the initial condition $w(0,\cdot)=w_0$ and the compatibility
condition 
$$
[\nu \operatorname{Def}w_0 \ N ]_{\operatorname{tan}} = g(0, \cdot )
_{\operatorname{tan}},
$$ 
and if 
$$
\begin{array}{c}
f\in V^1(T), \
g\in L^2(0,T;H^{\frac{3}{2}}({\mathbb R}^2, {\mathbb R}^3)), \
g_t\in L^2(0,T;H^{-\frac{1}{2}}({\mathbb R}^2, {\mathbb R}^3)), \\
B\in C^0([0,T];H^{\frac{1}{2}}({\mathbb R}^2, {\mathbb R})), \
B_t\in L^2(0,T;H^{\frac{1}{2}}({\mathbb R}^2, {\mathbb R}))\,,
\end{array}
$$
and $B(0,\cdot)=0$, 
then $w \in V^3_{\operatorname{div}}(T)$, $p \in V^2(T)$, and we have
the following estimate:
\begin{align}
&\| w\|_{L^2(0,T;H^3({\mathbb R}^3_+;{\mathbb R}^3))} +
\|N \cdot \triangle_0 \int_0^{\cdot} w(s)ds\|_
{C^0([0,T];H^{\frac{1}{2}}({\mathbb R}^2;{\mathbb R}))} 
+ \| p\|_{L^2(0,T;H^2({\mathbb R}^3_+;{\mathbb R}))} \nonumber \\
&+\| w_t\|_{L^2(0,T;H^1({\mathbb R}^3_+;{\mathbb R}^3))} +
\|N \cdot \triangle_0  w\|_
{L^2(0,T;H^{\frac{1}{2}}({\mathbb R}^2;{\mathbb R}))} 
+ \| p_t\|_{L^2(0,T;L^2({\mathbb R}^3_+;{\mathbb R}))} \nonumber \\
&+ \| \nabla  w_t\|_{L^2(0,T;H^{-\frac{1}{2}}({\mathbb R}^2;{\mathbb R}^9))}
+ \| p_t\|_{L^2(0,T;H^{-\frac{1}{2}}({\mathbb R}^2;{\mathbb R}))} \nonumber \\
& \qquad \le C \left( \|w_0\|_{H^2({\mathbb R}^3_+; {\mathbb R}^3)} +  \|f(0)\|_{H^2({\mathbb R}^3_+; {\mathbb R}^3)}+  \|g(0)\|_{H^{\frac{1}{2}}({\mathbb R}^2; {\mathbb R}^3)}
 \right.\nonumber \\
& \qquad\qquad
+ \| f\|_ {L^2(0,T;H^1({\mathbb R}^3_+;{\mathbb R}^3))}+
\| f_t\|_ {L^2(0,T;H^1({\mathbb R}^3_+;{\mathbb R}^3)')}
+ \| B\|_ {C^0([0,T];H^\frac{1}{2}({\mathbb R}^2;{\mathbb R}))}
\nonumber \\
& \qquad\qquad  \left.
+ \| B_t\|_ {L^2(0,T;H^{\frac{1}{2}}({\mathbb R}^2;{\mathbb R}))}
+ \| g\|_ {L^2(0,T;H^\frac{3}{2}({\mathbb R}^2;{\mathbb R}^3))}
+ \| g_t\|_ {L^2(0,T;H^{-\frac{1}{2}}({\mathbb R}^2;{\mathbb R}^3))}
\right)\,. \label{estimate1_flat}
\end{align}

\end{theorem}
\begin{proof}
We will use the superposition principle, first establishing the result
by keeping the data $f$, $g$, and $w_0$, and setting $B=0$ in the first
step, and then keeping $B$ and setting $f=g=w_0=0$ in the second step.
In the third step of the proof, we shall establish the necessary boundary
regularity for the $w_t$ and $p_t$.

1. We first consider the case that $B=0$.  
From the estimate (\ref{H1energy_p}), we see that for a.e. $0\le t\le T$,
$w(t) \in{\mathcal V}$.  We set $v=D_{-h}D_h w$ for $h=e_1,e_2$ in 
(\ref{weak0}).  Since for $f_1,f_2\in L^2({\mathbb R}^3_+;{\mathbb R}^3)$,
$$\int_{ {\mathbb R}^3_+} f_1(x) \cdot D_{-h}D_h f_2(x) dx = 
\int_{ {\mathbb R}^3_+} D_h f_1(x) \cdot D_h f_2(x) dx,$$
we see that for $\epsilon >0$,
\begin{align}
& \frac{1}{2} \frac{d}{dt}\left(
\| D_hw\|^2_{L^2({\mathbb R}^3_+;{\mathbb R}^3)} +
\sigma \left\| D_h\nabla_0 \int_0^tw^3(r) dr
\right\|^2_{L^2( {\mathbb R}^2; {\mathbb R}^2)} \right)
+ \frac{ \nu }{2}  \|D_h\operatorname{Def}w\|^2_{L^2({\mathbb R}^3_+;{\mathbb R}^9)} 
  \nonumber\\
& \qquad =  (f,  D_{-h}D_hw)_{L^2({\mathbb R}^3_+; {\mathbb R}^3)} 
+ ( D_h g , D_h w )_ {L^2({\mathbb R}^2; {\mathbb R}^3)} 
  \nonumber\\
& \qquad \le  C(\epsilon)\|f\|^2_{L^2({\mathbb R}^3_+; {\mathbb R}^3)} +
\epsilon\|D_h \nabla_0 w\|^2_{L^2({\mathbb R}^3_+; {\mathbb R}^6)}
+ C(\epsilon) \|D_h g\|^2_ {H^{-\frac{1}{2}}({\mathbb R}^2; {\mathbb R}^3)} 
\nonumber \\
&\qquad \qquad
+ \epsilon\|D_h  w\|^2_{L^2({\mathbb R}^2; {\mathbb R}^3)}\nonumber\\
& \qquad \le  C\left( \|f\|^2_{L^2({\mathbb R}^3_+; {\mathbb R}^3)} +
\|g\|^2_ {H^{\frac{1}{2}}({\mathbb R}^2; {\mathbb R}^3)}\right)
+ \epsilon\|D_h \nabla_0 w\|^2_{L^2({\mathbb R}^3_+; {\mathbb R}^6)}\nonumber \\
&\qquad \qquad
+ \epsilon\|D_h \operatorname{Def}  w\|^2_{L^2({\mathbb R}^3_+;{\mathbb R}^9)}
+ \epsilon\|D_h  w\|^2_{L^2({\mathbb R}^3_+;{\mathbb R}^3)}\,,
\label{reg1}
\end{align}
where we have used Lemma \ref{quotients} and the trace theorem for the last
inequality.  We choose $\epsilon >0$ sufficiently small, define
$$
y(t):= \|D_h w(t)\|^2_{L^2({\mathbb R}^3_+;{\mathbb R}^3)}+
\sigma \left\| D_h\nabla_0 \int_0^tw^3(r) dr
\right\|^2_{L^2( {\mathbb R}^2; {\mathbb R}^2)}\,,
$$
and
$$
F(t) = C\left(
\|f(t)\|^2_{L^2({\mathbb R}^3_+; {\mathbb R}^3)} 
+\|g(t)\|^2_ {H^{\frac{1}{2}}({\mathbb R}^2; {\mathbb R}^3)}\right)\,,
$$
so that for some constant $c>0$, we may reexpress (\ref{reg1}) as 
$y'(t) - cy(t) \le F(t)$. Gronwall's inequality then shows that for
all $0\le \bar T \le T$, $y(\bar T) \le e^{CT}y(0) + \int_0^T e^{c(T-t)}F(t)dt$.
Since $y(0)=
\|D_h w_0\|^2_{L^2({\mathbb R}^3_+;{\mathbb R}^3)}
\le C\|\nabla_0 w_0\|^2_{L^2({\mathbb R}^3_+;{\mathbb R}^6)}$,
we see that
\begin{align*}
\sup_{0\le t\le T}
\|D_h w(t)\|^2_{L^2({\mathbb R}^3_+;{\mathbb R}^3)}
&\le C(T) \left(
\|\nabla_0 w_0\|^2_{L^2({\mathbb R}^3_+;{\mathbb R}^6)}
+ \|f\|^2_{L^2(0,T; L^2({\mathbb R}^3_+; {\mathbb R}^3))} \right. \\
& \qquad \qquad \qquad \left.
+ \|g\|^2_{L^2(0,T; H^{\frac{1}{2}}({\mathbb R}^2; {\mathbb R}^3))} 
\right)=: \Phi(T) \,,
\end{align*}
where $C(T)\rightarrow +\infty$ as $T\rightarrow+ \infty$. 
Integrating (\ref{reg1}) from $0$ to $T$, and using Korn's inequality, we
find that
$$
\|D_h\nabla w\|^2_{L^2(0,T; L^2({\mathbb R}^3_+; {\mathbb R}^9))}
\le \Phi(T),
$$
possibly readjusting $C(T)$ if necessary.  By Lemma \ref{quotients}, 
\begin{equation}\label{0001}
\|\nabla_0\nabla w\|^2_{L^2(0,T; L^2({\mathbb R}^3_+; {\mathbb R}^{18}))}
\le \Phi(T)\,.
\end{equation}

We next establish a similar bound for $\partial_3\partial_3w^\alpha$.
Since $\operatorname{div}w=0$, $\partial_3\partial_3 w^3 =
-\partial_1\partial_3 w^1-\partial_2\partial_3 w^2$, so that by (\ref{0001}),
\begin{equation}\label{0002}
\| \triangle w^3\|^2_{L^2(0,T; L^2({\mathbb R}^3_+; {\mathbb R}))}
\le \Phi(T)\,.
\end{equation}
Since $\partial_3 p = w^3_t - \nu\triangle w^3 -f^3$, then 
$\partial_3 p \in L^2(0,T; L^2({\mathbb R}^3_+; {\mathbb R}))$ and
satisfies the bound
\begin{equation}\label{0003}
\| \partial_3 p\|^2_{L^2(0,T; L^2({\mathbb R}^3_+; {\mathbb R}))}
\le \Phi(T) + \|w_t\|^2_{L^2(0,T; L^2({\mathbb R}^3_+; {\mathbb R}^3))}\,.
\end{equation}

Using Lemma \ref{pressure}, we see that on ${\mathbb R}^3_+$, $w$ satisfies
\begin{align}
&(w_t, v)_{L^2({\mathbb R}^3_+;{\mathbb R}^3)} 
+ \frac{ \nu }{2}  
(\operatorname{Def}w, \operatorname{Def}v)_{L^2({\mathbb R}^3_+;{\mathbb R}^9)} 
-(p, \operatorname{div}v)_{L^2({\mathbb R}^3_+;{\mathbb R})} \nonumber \\
&\qquad =  \langle f, v\rangle 
+ \langle g , v \rangle_{\Gamma_0} 
-\sigma \left( \int_0^t \nabla_0 w^3(r) dr  \ , \
\nabla_0v^3\right)_{L^2( {\mathbb R}^2; {\mathbb R}^2)} \ \ \ \forall v\in 
{\mathcal W}\,. \label{weak1}
\end{align}
Since the test function $v$ need not satisfy the divergence-free constraint,
the idea is to take $v=(\partial_1 p, \partial_2 p, 0)$ as the test function.
Formally, at least, for this choice of $v$, the boundary term in (\ref{weak1}) 
disappears, and the bound (\ref{0003}) for $\partial_\alpha p$ can be
obtained.  Since we a priori
do not know that $\partial_\alpha p \in{\mathcal W}$, we must use a 
mollification process. 

For each $x_3\ge 0$ and $0\le\epsilon\le \bar \epsilon$, 
let $(x_1,x_2)\mapsto \rho^\epsilon(x_1,x_2,x_3)$ denote
a family of  Friedrichs mollifiers on ${\mathbb R}^2 \times \{ x^3\}$.
It is standard that the $\rho^\epsilon$ satisfy
$\partial_\alpha( \rho^\epsilon * p) =
\partial_\alpha \rho^\epsilon * p =
 \rho^\epsilon * \partial_\alpha p$, and that $\int_{{\mathbb R}^3_+}
(f) (\rho^\epsilon * g) dx = \int_{{\mathbb R}^3_+}
(\rho^\epsilon * f) (g) dx$ for all $f,g \in{L^2({\mathbb R}^3_+; 
{\mathbb R})}$. 

For our test function $v$ in (\ref{weak1}), we
set $v^\alpha= \rho^\epsilon * \partial_\alpha (\rho^\epsilon * p)$ and
$v^3=0$, so that $v\in {\mathcal W}$.  Let $P:=
\rho^\epsilon * p$, and $W:= \rho ^ \epsilon  * w$.  Then (\ref{weak1})
becomes
\begin{align*}
\|\partial_\alpha P\|^2_{L^2({\mathbb R}^3_+; {\mathbb R})} 
&= - \frac{ \nu }{2}  \int_{{\mathbb R}^3_+} \left[ 
2(W^\alpha,_{\beta} + W^\beta,_{\alpha}) P,_{\alpha \beta}
+(W^\alpha,_{3} + W^3,_{\alpha}) P,_{\alpha 3}\right] dx \\
& \qquad 
+ ( \rho ^ \epsilon * f^\alpha -W^\alpha_t, P,_{\alpha})_{L^2({\mathbb R}^3_+; {\mathbb R}^3)}
+ (\rho ^ \epsilon * g^\alpha , P,_{\alpha})_{L^2({\mathbb R}^2; {\mathbb R}^3)} \\
&=  \frac{ \nu }{2}  \int_{{\mathbb R}^3_+} \left[ 
2(W^\alpha,_{\beta\beta} + W^\beta,_{\alpha\beta}) P,_{\alpha }
+(W^\alpha,_{3\alpha} + W^3,_{\alpha\alpha}) P,_{3}\right] dx \\
& \qquad
+ (\rho ^ \epsilon *f^\alpha -W^\alpha_t, P,_{\alpha})_{L^2({\mathbb R}^3_+; {\mathbb R})}
- ( \rho ^ \epsilon * g^\alpha,_\alpha , P)_{L^2({\mathbb R}^2; {\mathbb R})} \,.
\end{align*}
We use Young's inequality (and the trace theorem for the last term), and then
integrate from $0$ to $T$ to find
$$ \|\partial_\alpha P\|^2_{L^2(0,T;L^2({\mathbb R}^3_+; {\mathbb R}))} 
\le C(\varepsilon) \left(
\Phi(T) + \|W_t\|^2_{L^2(0,T; L^2({\mathbb R}^3_+; {\mathbb R}^3))}\right)
+ \varepsilon \|\partial_\alpha P\|^2_{L^2({\mathbb R}^3_+; {\mathbb R})}\,. 
$$
where we have used (\ref{0001}), (\ref{0002}), and (\ref{0003}).
Using the fact that 
$$\|W_t\|^2_{L^2(0,T; L^2({\mathbb R}^3_+; {\mathbb R}^3))}
\le \|w_t\|^2_{L^2(0,T; L^2({\mathbb R}^3_+; {\mathbb R}^3))}$$
and choosing $\varepsilon>0$ sufficiently small, we see that
$$ \|\partial_\alpha P\|^2_{L^2(0,T;L^2({\mathbb R}^3_+; {\mathbb R}))} 
\le C
$$
for some constant $C$.  Thus, $ \partial_\alpha (\rho^\epsilon * p)
 \rightharpoonup \varphi_\alpha$ in 
${L^2(0,T;L^2({\mathbb R}^3_+; {\mathbb R}))}$
as $\epsilon\rightarrow 0$.  Choose $\psi \in {\mathcal D}( {\mathbb R}^3_+,
{\mathbb R})$ with $\text{supp}(\psi) =K$.  Then for each $\alpha =1,2$
$$ 
-\int_0^T \int_{{\mathbb R}^3_+} \partial_\alpha (\rho^\epsilon * \psi) 
\, p \, dx dt
\rightarrow 
\int_0^T \int_{{\mathbb R}^3_+} \varphi_\alpha \, \psi \,  dx dt
$$ 
as $\epsilon\rightarrow 0$.  For each $x^3 \in K$,
$\rho^\epsilon* \partial_\alpha  
\psi(\cdot, x^3)\rightarrow \partial_ \alpha \psi( \cdot, x^3)$ a.e., and
$$
\|\rho^\epsilon * \partial_ \alpha \psi(\cdot, x^3)\|
_{L^ \infty({\mathbb R}^2 \times \{x^3\}; {\mathbb R}) } 
\le \| \rho^\epsilon\|
_{L^1({\mathbb R}^2 \times \{x^3\}; {\mathbb R}) } 
\|\partial_ \alpha \psi(\cdot, x^3)\|
_{L^ \infty({\mathbb R}^2 \times \{x^3\}; {\mathbb R}) } ,
$$
so that
$$
\|\rho^\epsilon * \partial_ \alpha \psi\|
_{L^ \infty({\mathbb R}^3 ; {\mathbb R}) }  \le C
$$
for some constant $C>0$, which follows since the intersection of $K$ with
the $x^3$-axis is compact.
 
Thus, for some limit point $ \xi ^ \alpha \in
L^2({\mathbb R}^3_+; {\mathbb R})$,
$\rho^\epsilon * \partial_\alpha \psi \rightharpoonup \xi_\alpha$ in
$L^2({\mathbb R}^3_+; {\mathbb R})$ and by the a.e. convergence,
$\xi_\alpha= \partial_\alpha \psi$ in $L^2({\mathbb R}^3_+; {\mathbb R})$,
so that
$$
\int_0^T \int_{{\mathbb R}^3_+} \varphi_\alpha \, \psi \,  dx dt
= -
\int_0^T \int_{{\mathbb R}^3_+} p \, \partial_\alpha \psi \,  dx dt.
$$
By lower semi-continuity of weak convergence of the mollification, 
it follows that
\begin{equation}\nonumber
\|\partial_ \alpha  p\|^2_{L^2(0,T;L^2({\mathbb R}^3_+; {\mathbb R}^3))}
\le C\left( \Phi(T) + \|w_t\|^2_{L^2(0,T;L^2({\mathbb R}^3_+; {\mathbb R}^3))}
\right)\,.
\end{equation}
Combining this with (\ref{0003}), we obtain
\begin{equation}\label{0004}
\|\nabla  p\|^2_{L^2(0,T;L^2({\mathbb R}^3_+; {\mathbb R}^3))}
\le C\left( \Phi(T) + \|w_t\|^2_{L^2(0,T;L^2({\mathbb R}^3_+; {\mathbb R}^3))}
\right)\,.
\end{equation}

Since $\nu \partial_3 \partial_3 w^\alpha = -\nu \partial_\beta\partial_\beta
w^\alpha -f^\alpha + w_t^\alpha + \partial_\alpha p$, 
$\partial_3 \partial_3 w^\alpha \in 
L^2(0,T;L^2({\mathbb R}^3_+; {\mathbb R}^3))$ and is bounded by the
right-hand-side of (\ref{0004}).  We obtain the estimate
\begin{equation}\label{H2}
\|w\|_{L^2(0,T;H^2({\mathbb R}^3_+; {\mathbb R}^3))}  
+ \|p\|_{L^2(0,T;H^1({\mathbb R}^3_+; {\mathbb R}))} 
\le C\left( \Phi(T) + \|w_t\|^2_{L^2(0,T;L^2({\mathbb R}^3_+; {\mathbb R}^3))}
\right)\,.
\end{equation}

Next, we set $v=D_{-h}D_h D_{-h}D_h w$ for $h=e_1,e_2, 
\frac{e_1+e_2}{\sqrt{2}}$ 
in (\ref{weak0}).  Just as we obtained the inequality (\ref{reg1}), we
find that
\begin{align}
& \frac{1}{2} \frac{d}{dt}\left(
\| D_{-h}D_hw_t\|^2_{L^2({\mathbb R}^3_+;{\mathbb R}^3)} +
\sigma \left\| D_{-h}D_h\nabla_0 \int_0^tw^3(r) dr
\right\|^2_{L^2( {\mathbb R}^2; {\mathbb R}^2)} \right) \nonumber \\
& 
+ \frac{ \nu }{2}  \|D_{-h}D_h\operatorname{Def}w\|^2_{L^2({\mathbb R}^3_+;{\mathbb R}^9)} 
\le  C\left( \|f\|^2_{H^1({\mathbb R}^3_+; {\mathbb R}^3)} +
\|g\|^2_ {H^{\frac{3}{2}}({\mathbb R}^2; {\mathbb R}^3)}\right)
        \nonumber \\
& + \epsilon\|D_{-h}D_h \nabla_0 w\|^2_{L^2({\mathbb R}^3_+; {\mathbb R}^6)}
+ \epsilon\|D_{-h}D_h \operatorname{Def}  
w\|^2_{L^2({\mathbb R}^3_+;{\mathbb R}^9)}
+ \epsilon\|D_{-h}D_h  w\|^2_{L^2({\mathbb R}^3_+;{\mathbb R}^3)}\,,
\label{reg2}
\end{align}
The identical argument leading to (\ref{0001}) gives
\begin{equation}\label{1001}
\|\nabla_0\nabla_0\nabla w\|^2_{L^2(0,T; L^2({\mathbb R}^3_+; {\mathbb R}^{36}))}
\le \Psi(T)\,.
\end{equation}
where
\begin{equation}\label{psi}
\Psi(T):= C(T) \left(
\|\nabla_0\nabla_0 w_0\|^2_{L^2({\mathbb R}^3_+;{\mathbb R}^{12})}
+ \|f\|^2_{L^2(0,T; H^1({\mathbb R}^3_+; {\mathbb R}^3))}
+ \|g\|^2_{L^2(0,T; H^{\frac{3}{2}}({\mathbb R}^2; {\mathbb R}^3))} 
\right)\,,
\end{equation}
and $C(T)\rightarrow +\infty$ as $T \rightarrow + \infty$.  Since, 
$\|\partial_\alpha \partial_\beta \partial_j w^i\|_
{L^2(0,T; L^2({\mathbb R}^3_+; {\mathbb R}))}\le \Psi(T)$, then
$\|\partial_\alpha \partial_3\partial_3 w^3\|_
{L^2(0,T; L^2({\mathbb R}^3_+; {\mathbb R}))}\le \Psi(T)$
using that $\operatorname{div}w=0$.  It follows that
$$\|\partial_\alpha  \triangle w^3\|_
{L^2(0,T; L^2({\mathbb R}^3_+; {\mathbb R}))}\le \Psi(T),$$ 
and
since $\partial_\alpha \partial_3 p = - \partial_\alpha w^3_t +\partial_\alpha f
+ \partial_\alpha\triangle w^3$, we see that 
\begin{equation}\label{1002}
\| \partial_\alpha\partial_3 p\|^2
_{L^2(0,T; L^2({\mathbb R}^3_+; {\mathbb R}))}
\le \Psi(T) + \|w_t\|^2_{L^2(0,T; H^1({\mathbb R}^3_+; {\mathbb R}^3))}\,.
\end{equation}

Setting $v^\alpha= -\rho^\epsilon * \partial_\alpha \partial_\gamma
\partial_\gamma (\rho^\epsilon * p)$ in (\ref{weak1}) we find that
\begin{align}
&\|\partial_\alpha\partial_\gamma P\|^2_{L^2({\mathbb R}^3_+; {\mathbb R})} 
=  \frac{ \nu }{2}  \int_{{\mathbb R}^3_+} \left[ 
2(W^\alpha,_{\beta\beta\gamma}+W^\beta,_{\alpha\beta\gamma}) P,_{\alpha\gamma }
+(W^\alpha,_{3\gamma\gamma}+W^3,_{\alpha\gamma\gamma}) P,_{3\alpha}\right] dx 
\nonumber \\
& \qquad\qquad\qquad
+ (\partial_\gamma( \rho ^ \epsilon * f^\alpha - W^\alpha_t), 
P,_{\alpha\gamma})_{L^2({\mathbb R}^3_+; {\mathbb R})}
- ( \rho ^ \epsilon * g^\alpha,_{\alpha\gamma} , P,_\gamma)
_{L^2({\mathbb R}^2; {\mathbb R})} \,. \label{pab}
\end{align}
Integrating from $0$ to $T$, using Young's inequality and the trace theorem,
and passing to the limit as $\epsilon\rightarrow 0$, the estimate (\ref{1002})
gives 
\begin{equation}\label{1003}
\|\partial_\alpha\nabla p\|^2_{L^2(0,T;L^2({\mathbb R}^3_+; {\mathbb R}^3))}
\le C\left( \Psi(T) + \|w_t\|^2_{L^2(0,T;H^1({\mathbb R}^3_+; {\mathbb R}^3))}
\right)\,.
\end{equation}
Since $\nu \partial_i\partial_3 \partial_3 w^\beta = 
-\nu \partial_i\partial_\gamma\partial_\gamma w^\beta -
\partial_i f^\beta + \partial_i w_t^\beta + \partial_{i\beta}p$, and
$w^3,_{333} =- w^1_{331} - w^2_{332}$ from the divergence-free condition,
we find that
\begin{equation}\label{H3}
\|w\|_{L^2(0,T;H^3({\mathbb R}^3_+; {\mathbb R}^3))}  
+ \|p\|_{L^2(0,T;H^2({\mathbb R}^3_+; {\mathbb R}))} 
\le C\left( \Psi(T) + \|w_t\|^2_{L^2(0,T;H^1({\mathbb R}^3_+; {\mathbb R}^3))}
\right)\,.
\end{equation}

2. Next, we keep $B$ and set $f=g=w_0=0$. We first assume the following
regularity for $B$:
\begin{equation}\label{Blow}
B \in C^0([0,T]; {L^2({\mathbb R}^3_+; {\mathbb R})})\,, \ \ \ 
B_t\in L^2(0,T; {L^2({\mathbb R}^3_+; {\mathbb R})})\,. 
\end{equation}
During the course of the proof of Theorem \ref{uniqueweak}, we proved that $w_t (0)=\bar w_0$, with $\bar{w}_0$ satisfying (\ref{wt(0)}). Thus, $w_t(0)=0$ with our present choice of initial data and forcings.

We consider the weak
formulation for the time derivative $w_t$ given by 
\begin{align}
&(w_{tt}, v)_{L^2({\mathbb R}^3_+;{\mathbb R}^3)} 
+ \frac{ \nu }{2}  (\operatorname{Def}w_t, \operatorname{Def}v)_{L^2({\mathbb R}^3_+;{\mathbb R}^9)} 
=   \sigma \int_{\Gamma_0}B_t\, v^3 \, dS
-\sigma \int_{\Gamma_0} \partial_\alpha w^3 \partial_\alpha v^3 \, dS
\,, \label{weakt}
\end{align}
for all $v\in L^2(0,T;{\mathcal V})$.  We set $v=D_{-h}D_h w$ in (\ref{weakt})
and integrate from $0$ to $\bar T$ for $0 \le \bar T \le T$ to find that
\begin{align*}
&
 \frac{\nu}{4}
\| D_h \operatorname{Def} w (\bar T)\|^2_{L^2({\mathbb R}^3_+; {\mathbb R}^9)}  
+ 
\sigma \int_0^{\bar T}\int_{\Gamma_0} |D_h\nabla_0 w^3|^2 \, dS\, dt\\
&\qquad  \le
(w_t (\bar T), D_{-h}D_h w(\bar T))_{L^2({\mathbb R}^3_+; {\mathbb R}^3)}
+ C(\epsilon) \int_0^T \|B_t\|^2_{L^2({\mathbb R}^2; {\mathbb R}^3)} dt  \\
&\qquad \qquad
+  \epsilon \int_0^{\bar T}\int_{\Gamma_0} |D_{-h}D_h w^3|^2 \, dS\, dt \, .
\end{align*}
In order to employ Korn's inequality we add $\frac{\nu}{4}\|D_h w(\bar T)\|^2
_{L^2({\mathbb R}^3_+; {\mathbb R}^3)}$ to both sides of the above inequality.
We find that with $ \epsilon >0$ taken sufficiently small,
\begin{align}
&
\| D_h w (\bar T)\|^2_{H^1({\mathbb R}^3_+; {\mathbb R}^3)}  
+ 
\int_0^{\bar T}\int_{\Gamma_0} |D_h\nabla_0 w^3|^2 \, dS\, dt
\le C\left(
\|B_t\|^2_{L^2(0,T; L^2( \mathbb{R}^2  ; \mathbb{R} ))} \right. \nonumber\\
&\qquad \left.
+
\|w_t\|^2_{L^2(0,T; H^1( \mathbb{R}^3_+  ; \mathbb{R}^3  ))}
+
\|w_t(\bar T)\|^2_{L^2( \mathbb{R}^3_+  ; \mathbb{R}^3  )}
+
\|w(\bar T)\|^2_{H^1( \mathbb{R}^3_+  ; \mathbb{R}^3  )}
\right). \label{s1s1}
\end{align}
We need to estimate $\|w(\bar T)\|^2_{H^1( \mathbb{R}^3_+  ; \mathbb{R}^3  )}$.
Letting $v= D_{-h}D_h w$ in (\ref{weak0}), integrating from $0$ to $\bar T$,
and using Young's inequality, we find that
\begin{align*}
&\|D_hw(\bar T)\|^2_{L^2( \mathbb{R}^3_+  ; \mathbb{R}^3  )}
+ 
 \sigma \int_{\Gamma_0} |\int_0^{\bar T}D_h\nabla_0 w^3 \, dt|^2 \, dS \\
& \qquad \le
\frac{1}{2} \|B_t\|^2_{L^2(0,T; L^2( \mathbb{R}^2  ; \mathbb{R}) )} 
+ 
 \frac{\sigma}{2} \int_0^{\bar T} 
 \int_{\Gamma_0} |\int_0^{t}D_h\nabla_0 w^3 \, ds|^2 \, dS\, dt  \\
& \qquad \qquad  +
C( \epsilon ) \|B(\bar T)\|^2_{L^2( \mathbb{R}^2  ; \mathbb{R} )} 
+ \epsilon \int_{\Gamma_0} |\int_0^{\bar T}D_h\nabla_0 w^3 \, dt|^2 \, dS .
\end{align*}
Thus,  for $ \epsilon >0$ sufficiently small,
\begin{align*}
&
\int_{\Gamma_0} |\int_0^{\bar T}D_h\nabla_0 w^3 \, dt|^2 \, dS \\
& \qquad \le 
 C\int_{\Gamma_0} |\int_0^{t}D_h\nabla_0 w^3 \, ds|^2 \, dS\, dt  +
C( \|B_t\|^2_{L^2(0,T; L^2( \mathbb{R}^2  ; \mathbb{R}) )} 
+ 
\|B(\bar T)\|^2_{L^2( \mathbb{R}^2  ; \mathbb{R} )} ).
\end{align*}
Gronwall's inequality  yields the bound
$$
\int_{\Gamma_0} |\int_0^{\bar T}D_h\nabla_0 w^3 \, dt|^2 \, dS 
\le C(1+ \bar T e^{\bar T}) \left(
\|B_t\|^2_{L^2(0,T; L^2( \mathbb{R}^2  ; \mathbb{R}) )} 
+ 
\|B(\bar T)\|^2_{L^2( \mathbb{R}^2  ; \mathbb{R} )}
\right) 
$$
so that
$$
\|D_hw(\bar T)\|^2_{L^2( \mathbb{R}^3_+  ; \mathbb{R}^3  )}
\le C(T) \left(
\|B_t\|^2_{L^2(0,T; L^2( \mathbb{R}^2  ; \mathbb{R}) )} 
+ 
\|B(\bar T)\|^2_{L^2( \mathbb{R}^2  ; \mathbb{R} )}
\right) 
$$
where $C(T) \rightarrow + \infty $ as $T \rightarrow + \infty $. Substituting
this into (\ref{s1s1}) and taking the sup over $0 \le \bar T \le T$, we find
that
\begin{align}
&\sup_{0\le t\le T} \|\nabla_0 w(t) \|^2_{H^1({\mathbb R}^3_+; {\mathbb R}^6)} 
+ \|w^3\|^2_{L^2(0,T;H^2( {\mathbb R}^2, {\mathbb R}))}\nonumber\\
& \qquad \qquad
\le C \left( \|B\|^2_{L^2(0,T;L^2( {\mathbb R}^2, {\mathbb R}))}+
\|B_t\|^2_{L^2(0,T;L^2( {\mathbb R}^2, {\mathbb R}))}
 +\|w_t\|^2_{L^\infty(0,T;L^2( {\mathbb R}^3, {\mathbb R}^3))}\right)\nonumber\\
& \qquad \qquad
\le C \left( \|B\|^2_{C^0([0,T];L^2( {\mathbb R}^2, {\mathbb R}))}
 +\|B_t\|^2_{L^2(0,T;L^2( {\mathbb R}^2, {\mathbb R}))}\right) \,,
\label{BL2}
\end{align}
where we used (\ref{H1energy_p}) for the second inequality, with $\Omega_0$
replaced by ${\mathbb R}^3_+$ and $\Gamma_0$ replaced by ${\mathbb R}^2$.

Having established (\ref{BL2}) under the assumptions (\ref{Blow}), we next
assume that
\begin{equation}\label{Bhigh}
B \in C^0([0,T]; {H^1({\mathbb R}^3_+; {\mathbb R})})\,, \ \ \ 
B_t\in L^2(0,T; {H^1({\mathbb R}^3_+; {\mathbb R})})\,. 
\end{equation}
We take $v=D_{-h}D_h w_t$ in (\ref{weakt}) and obtain 
\begin{align*}
&
\frac{d}{dt} \|D_h w_t\|^2_{L^2({\mathbb R}^3_+; {\mathbb R}^3)} 
+ C\|D_hw_t \|^2_{H^1({\mathbb R}^3_+; {\mathbb R}^3)} 
+\sigma\frac{d}{dt} \|D_h \nabla_0 w^3\|^2_{L^2({\mathbb R}^2;
 {\mathbb R}^2)}  \\
& \qquad \qquad = 2( D_h B_t, D_h w^3_t)_{L^2({\mathbb R}^2; {\mathbb R})} 
+ 2\nu\|D_h w_t\|^2_{L^2({\mathbb R}^3_+; {\mathbb R}^3)} \\
& \qquad \qquad \le C(\epsilon)\| B_t\|^2_{H^1({\mathbb R}^2; {\mathbb R})} 
+ \epsilon\| D_h w_t\|^2_{H^1({\mathbb R}^3_+; {\mathbb R}^3)} 
+  2\nu\|w_t\|^2_{L^2({\mathbb R}^3_+; {\mathbb R}^3)}\,,
\end{align*}
where we have used Young's inequality, the trace theorem and Lemma \ref{quotients}
for the last inequality, as well as Korn's inequality for the first line.
Choosing $\epsilon >0$ sufficiently small, and integrating from $0$ to $\bar T$,
for $0\le \bar T\le T$, we have that
\begin{align*}
&
\|D_h w_t(\bar T)\|^2_{L^2({\mathbb R}^3_+; {\mathbb R}^3)} 
+ 2\|D_hw_t\|^2_{L^2(0,\bar T;H^1({\mathbb R}^3_+;{\mathbb R}^3))} 
+\|D_h \nabla_0 w^3(\bar T)\|^2_{L^2({\mathbb R}^2; {\mathbb R}^2)}  \\
& \qquad \qquad \le 
 C \left( \|B\|^2_{C^0([0,T];L^2( {\mathbb R}^2, {\mathbb R}))}
 +\|B_t\|^2_{L^2(0,T;L^2( {\mathbb R}^2, {\mathbb R}))}\right) \,,
\end{align*}
where we have again used (\ref{H1energy_p}) to bound
$\|w_t\|^2_{L^2(0,T; {L^2({\mathbb R}^3_+; {\mathbb R}^3)})}$. From Lemma
\ref{quotients}, it follows that 
\begin{align}
& \sup_{0\le t\le T} \|\nabla_0  w_t\|^2_{L^2({\mathbb R}^3_+; {\mathbb R}^6)}
+ \|\nabla_0w_t\|^2_{L^2(0,\bar T;H^1({\mathbb R}^3_+;{\mathbb R}^6))}  
\nonumber \\
&\qquad \qquad
\le 
C \left( \|B\|^2_{C^0([0,T];L^2( {\mathbb R}^2, {\mathbb R}))}
 +\|B_t\|^2_{L^2(0,T;L^2( {\mathbb R}^2, {\mathbb R}))}\right) \,.
\label{ss1}
\end{align}

Next, letting $v= D_{-h}D_h D_{-h}D_h w$ in (\ref{weakt}), we get
\begin{align*}
&(D_h w_{tt}, D_h D_{-h}D_h w)_{L^2({\mathbb R}^3_+; {\mathbb R}^3)} 
+\nu \frac{d}{dt} \| D_{-h}D_h \operatorname{Def} w\|^2
_{L^2({\mathbb R}^3_+; {\mathbb R}^9)}  \\
& \qquad
+  \sigma \| D_{-h}D_h \nabla_0 w^3\|^2
_{L^2({\mathbb R}^2; {\mathbb R}^2)} =
(D_h B_t, D_h D_{-h}D_h w^3)_{L^2({\mathbb R}^2; {\mathbb R})} \\
& \qquad \qquad 
\le C(\epsilon)\|B_t\|^2_{H^1({\mathbb R}^2; {\mathbb R})} 
   + \epsilon \|D_{-h}D_h\nabla_0 w^3\|^2_{L^2({\mathbb R}^2; {\mathbb R}^2)} 
\,,
\end{align*}
where we have used Young's inequality and Lemma \ref{quotients} to obtain the
inequality.  Choosing $\epsilon >0$ sufficiently small and integrating from
$0$ to $T$ we get
\begin{align*}
& c_1\|D_{-h}D_h w(T)\|^2_{H^1({\mathbb R}^3_+; {\mathbb R}^3)}
+c_1 \| D_{-h}D_h \nabla_0 w^3\|^2_{L^2(0,T;L^2({\mathbb R}^2; {\mathbb R}^2))}
\\
&
\le \frac{1}{2}   \|D_{-h}D_h w_t\|^2_{L^2(0,T;L^2({\mathbb R}^3_+; {\mathbb R}^3))}
-  (D_h w_t(T), D_hD_{-h}D_h w(T))_{L^2({\mathbb R}^3_+; {\mathbb R}^3)} \\
&
+C \|B_t\|^2_{L^2(0,T;H^1({\mathbb R}^2; {\mathbb R}))}  
+ \nu
\|D_{-h}D_h w(T)\|^2_{L^2({\mathbb R}^3_+; {\mathbb R}^3)}
\\
&
\le \frac{1}{2}   \|D_h w_t\|^2_{L^2(0,T;H^1({\mathbb R}^3_+; {\mathbb R}^3))}
+ C(\epsilon) \|\nabla_0 w_t(T)\|^2_{L^2({\mathbb R}^3_+; {\mathbb R}^6)} \\
&
+C \|B_t\|^2_{L^2(0,T;H^1({\mathbb R}^2; {\mathbb R}))} 
+ \epsilon \|D_{-h}D_h w(T)\|^2_{H^1({\mathbb R}^3_+; {\mathbb R}^3)}
+ \nu \|D_{-h}D_h w(T)\|^2_{L^2({\mathbb R}^3_+; {\mathbb R}^3)} \,.
\end{align*}
Choosing $\epsilon >0$ sufficiently small, using (\ref{ss1}), and Lemma
\ref{quotients}, we find that
\begin{align}
\|w^3\|^2_{L^2(0,T;H^3( {\mathbb R}^2, {\mathbb R}))}
\le C \left( \|B\|^2_{C^0([0,T];L^2( {\mathbb R}^2, {\mathbb R}))}
 +\|B_t\|^2_{L^2(0,T;H^1( {\mathbb R}^2, {\mathbb R}))}\right) \,.
\label{BH1}
\end{align}
Interpolating between (\ref{BL2}) and (\ref{BH1}) we have that
\begin{align}
\|w^3\|^2_{L^2(0,T;H^{2.5}( {\mathbb R}^2, {\mathbb R}))}
\le C \left( \|B\|^2_{C^0([0,T];L^2( {\mathbb R}^2, {\mathbb R}))}
 +\|B_t\|^2_{L^2(0,T;H^{\frac{1}{2}}( {\mathbb R}^2, {\mathbb R}))}\right) \,.
\label{B2.5}
\end{align}

Having established the estimate (\ref{B2.5}), we now return to 
the original assumptions of the theorem  on $B$ so that
$$
B \in C^0([0,T]; {H^{\frac{1}{2}}({\mathbb R}^2; {\mathbb R})})\,, \ \ \ 
B_t\in L^2(0,T; {H^{\frac{1}{2}}( {\mathbb R}^2; {\mathbb R})})\,. 
$$
Equations (\ref{linear7.a}) and (\ref{linear7.c}) take the form
\begin{subequations}
  \label{linear8}
\begin{alignat}{2}
 w_t - \nu \triangle  w  &= -\nabla p
&& \text{in} \ \ (0,T)\times {\mathbb R}^3_+ \,, 
         \label{linear8.a}\\
 \partial_3 w^\alpha &= - \partial_\alpha w^3
&&\text{on} \ \ (0,T)\times {\mathbb R}^2 \,, 
         \label{linear8.b}\\
 2 \nu \partial_3 w^3 -p &= \sigma \left(\int_0^t\triangle_0 w^3(r)\, dr + B \right)
&\ \ &\text{on} \ \ (0,T)\times {\mathbb R}^2 \,. 
         \label{linear8.c}
\end{alignat}
\end{subequations}
Taking the $L^2$ inner-product of (\ref{linear8.a}) with $v \in 
L^2(0,T;{\mathcal V})$ we find that
\begin{align*} 
(w_t, v)_{L^2({\mathbb R}^3_+; {\mathbb R}^3)} +
 \nu ( \nabla w, \nabla v)_{L^2({\mathbb R}^3_+; {\mathbb R}^9)} 
& = \int_{\Gamma_0} \left[ \nu \partial_3 w^i \, v^i - p\, v^3\right] dS \\
& = \int_{\Gamma_0} \left[- \nu \partial_\alpha w^3 \, v^\alpha 
+  ( 2 \nu \partial_3 w^3 - p)\, v^3\right) dS \,, 
\end{align*} 
where we have used (\ref{linear8.b}) for the second equality.
Setting $v=D_{-h}D_hD_{-h}D_h w$,
\begin{align*}
&\frac{1}{2}\frac{d}{dt}\|D_{-h}D_h w\|^2_{L^2({\mathbb R}^3_+; {\mathbb R}^3)} 
+  \nu  \| D_{-h}D_h 
 \nabla w\|^2_{L^2({\mathbb R}^3_+; {\mathbb R}^9)} \\
&\qquad =
-\int_{\Gamma_0} D_{-h}D_h w^3,_{\alpha} \, D_{-h}D_h w^\alpha\, dS
+\int_{\Gamma_0} D_h ( 2 \nu w^3,_3-p) \, D_hD_{-h}D_h w^3\, dS\,.
\end{align*}
Using Korn's inequality and integrating from $0$ to $T$, we have that
\begin{align*}
&
\nu \| D_{-h}D_h w\|^2_{L^2(0,T;H^1({\mathbb R}^3_+; {\mathbb R}^3))} 
\le  -\int_0^T
\int_{\Gamma_0} D_{-h}D_h w^3,_{\alpha} \, D_{-h}D_h w^\alpha\, dS \, dt \\
&\qquad + \int_0^T
\int_{\Gamma_0} D_h (w^3,_3-p) \, D_hD_{-h}D_h w^3\, dS\, dt 
+ \nu \| D_{-h}D_h w\|^2_{L^2(0,T;L^2({\mathbb R}^3_+; {\mathbb R}^3))} \,.
\end{align*}
We look at the right-hand-side term by term.  First,
$\nu \| D_{-h}D_h w\|^2_{L^2(0,T;L^2({\mathbb R}^3_+; {\mathbb R}^3))}$ is
bounded by (\ref{BL2}); second,
\begin{align*}
&\int_0^T\int_{\Gamma_0} D_{-h}D_h w^3,_{\alpha} \, D_{-h}D_h w^\alpha\, dS\,dt
\\
&\qquad 
\le \int_0^T \| D_{-h}D_h w^3,_{\alpha}\|
_{H^{-\frac{1}{2}}({\mathbb R}^2; {\mathbb R})} 
\| D_{-h}D_h w^\alpha\|
_{H^{\frac{1}{2}}({\mathbb R}^2; {\mathbb R})}  \, dt \\
&\qquad 
\le \int_0^T C\| D_{-h}D_h \triangle_0^{-\frac{1}{4}}w^3,_{\alpha}\|
_{L^2({\mathbb R}^2; {\mathbb R})} 
\| D_{-h}D_h w^\alpha\|
_{H^{\frac{1}{2}}({\mathbb R}^2; {\mathbb R})}  \, dt \\
&\qquad 
\le  C(\epsilon) \|w^3\|^2_{L^2(0,T;H^{2.5}({\mathbb R}^2; {\mathbb R}))} 
+ \epsilon \|D_{-h}D_h w^\alpha\|^2_{L^2(0,T;H^{\frac{1}{2}}
({\mathbb R}^2; {\mathbb R}))}  \\
&\qquad 
\le C \left( \|B\|^2_{C^0([0,T];L^2( {\mathbb R}^2, {\mathbb R}))}
 +\|B_t\|^2_{L^2(0,T;H^{\frac{1}{2}}( {\mathbb R}^2, {\mathbb R}))}\right) \\
&\qquad \qquad \qquad
+ \epsilon \|D_{-h}D_h w\|^2_{L^2(0,T;H^1
({\mathbb R}^3_+; {\mathbb R^3}))}  \,,
\end{align*}
where we have used the duality for the first inequality, the fact that
$\|\triangle_0^{\frac{1}{4}} ( \cdot) \|^2_ {L^2({\mathbb R}^2;{\mathbb R})}$
defines an equivalent norm on $H^{\frac{1}{2}}({\mathbb R}^2; {\mathbb R})$ for
the second inequality, Young's inequality for the third inequality, and
the trace theorem together with the estimate (\ref{B2.5}) for the last
inequality.

Finally, using the identical arguments, we find that
\begin{align*}
& \int_0^T \int_{\Gamma_0} D_h ( 2 \nu w^3,_3-p) \, D_hD_{-h}D_h w^3\, dS\, dt \\
&\qquad 
\le \int_0^T \| D_hD_{-h}D_h w^3\|
_{H^{-\frac{1}{2}}({\mathbb R}^2; {\mathbb R})} 
\| D_h ( 2 \nu w^3,_{3} -p)\|
_{H^{\frac{1}{2}}({\mathbb R}^2; {\mathbb R})}  \, dt \\
&\qquad 
\le  C(\epsilon) \| w^3\|^2_{L^2(0,T;H^{2.5}({\mathbb R}^2; {\mathbb R}))} 
+ \epsilon
\|\partial_\alpha 2 \nu \partial_3 w^3-\partial_\alpha p)\|
^2_{L^2(0,T;H^1({\mathbb R}^3_+; {\mathbb R}))}
\,.
\end{align*}

Thus, choosing $\epsilon >0$ sufficiently small, combining the previous
estimates,  and using the fact that for some $C>0$, 
$$
\|\partial_\alpha 2 \nu \partial_3 w^3-\partial_\alpha p\|
^2_{L^2(0,T;H^1({\mathbb R}^3_+; {\mathbb R}))}
\le C\left( \|w\|^2_{L^2(0,T;H^3({\mathbb R}^3_+; {\mathbb R}^3))} 
+ \|p\|^2_{L^2(0,T;H^2({\mathbb R}^3_+; {\mathbb R}))}\right)\,, 
$$
we see that for some (possibly readjusted) $\epsilon >0$, and using Lemma
\ref{quotients}
\begin{align}
\| \nabla_0 \nabla_0 w\|^2_{L^2(0,T;H^1({\mathbb R}^3_+; {\mathbb R}^{12}))} 
&\le C \left( \|B\|^2_{C^0([0,T];L^2( {\mathbb R}^2, {\mathbb R}))}
 +\|B_t\|^2_{L^2(0,T;H^{\frac{1}{2}}( {\mathbb R}^2, {\mathbb R}))}\right) 
\nonumber \\
&\qquad + \epsilon
\left( \|w\|^2_{L^2(0,T;H^3({\mathbb R}^3_+; {\mathbb R}^3))} 
+ \|p\|^2_{L^2(0,T;H^2({\mathbb R}^3_+; {\mathbb R}))}\right) \,.
\label{sss2}
\end{align}
Let $\Xi(T)$ equal the right-hand-side of (\ref{sss2}).
Since $\partial_3\partial_3\partial_\alpha w^3 = - \partial_1\partial_\alpha
\partial_3 w^1 - \partial_2\partial_\alpha \partial_3 w^2$, we see that
$$\| \partial_\alpha \triangle w^3\|^2
_{L^2(0,T;L^2({\mathbb R}^3_+; {\mathbb R}))} \le \Xi(T);$$ 
hence, from (\ref{linear8.a}),
$$ \| \partial_\alpha \partial_3 p\|^2
_{L^2(0,T;L^2({\mathbb R}^3_+; {\mathbb R}))} \le c( \Xi(T)
+ \|w_t\|^2_{L^2(0,T;{H^1({\mathbb R}^3_+; {\mathbb R}^3)} )})
\le C \Xi(T) \,,
$$
where we used (\ref{H1energy_p}) for the second inequality.
From (\ref{linear8.a}), we see that $\nu \partial_3\partial_3\partial_3
 w^\alpha = \partial_3 w^\alpha_t - \nu \partial_3\triangle_0 w^\alpha
+ \partial_3 \partial_\alpha p$, so that
$$ \| \partial_3\partial_3\partial_3 w^\alpha \|^2
_{L^2(0,T;L^2({\mathbb R}^3_+; {\mathbb R}))} \le C \Xi(T)\,.
$$
Next, setting $f=g=0$ in (\ref{pab}) and again using the estimate
(\ref{H1energy_p}) we see that
\begin{equation}\nonumber
\|\partial_\alpha\nabla p\|^2_{L^2(0,T;L^2({\mathbb R}^3_+; {\mathbb R}^3))}
\le C \Xi(T) \,.
\end{equation}
Since $\nu \partial_i\partial_3 \partial_3 w^\beta = 
-\nu \partial_i\partial_\gamma\partial_\gamma w^\beta 
+ \partial_i w_t^\beta + \partial_i\partial_\beta p$, and
$w^3,_{333} =- w^1,_{331} - w^2,_{332}$, we see that
\begin{equation}\nonumber
\|w\|_{L^2(0,T;H^3({\mathbb R}^3_+; {\mathbb R}^3))}  
+ \|p\|_{L^2(0,T;H^2({\mathbb R}^3_+; {\mathbb R}))} 
\le C \Xi(T) \,.
\end{equation}
Choosing $\epsilon >0$ sufficiently small in (\ref{sss2}) we obtain the
estimate
\begin{align}
&\|w\|_{L^2(0,T;H^3({\mathbb R}^3_+; {\mathbb R}^3))}  
+ \|p\|_{L^2(0,T;H^2({\mathbb R}^3_+; {\mathbb R}))} \nonumber \\
&\qquad\qquad \qquad
\le C \left( \|B\|^2_{C^0([0,T];L^2( {\mathbb R}^2, {\mathbb R}))}
 +\|B_t\|^2_{L^2(0,T;H^{\frac{1}{2}}( {\mathbb R}^2, {\mathbb R}))}\right) 
\,.
\label{H3_part2}
\end{align}

Combining (\ref{H1energy_p}),(\ref{psi}), (\ref{H3}), and (\ref{H3_part2}),
it follows that
\begin{align}
& \|w\|^2_{L^2(0,T; H^3({\mathbb R}^3_+; {\mathbb R}^3))} 
 +\|w_t\|^2_{L^2(0,T; H^1({\mathbb R}^3_+; {\mathbb R}^3))} 
 + \|p\|^2_{L^2(0,T; H^2({\mathbb R}^3_+; {\mathbb R}))} 
\nonumber\\
&\
 + \|p_t\|^2_{L^2(0,T; L^2({\mathbb R}^3_+; {\mathbb R}))}\nonumber\\ 
& \le C(T)\left[ \|w_0\|^2_{H^2({\mathbb R}^3_+;{\mathbb R}^3)}+ \|f(0)\|^2_{L^2({\mathbb R}^3_+;{\mathbb R}^3)}+  \|g(0)\|^2_{H^{\frac{1}{2}}({\mathbb R}^2;{\mathbb R}^3)}    \right. \nonumber\\
& \qquad\qquad
+ \|f\|^2_{L^2(0,T;H^1(\Omega_0;{\mathbb R}^3))}  
+ \|f_t\|^2_{L^2(0,T;H^1(\Omega_0;{\mathbb R}^3)')}  
+ \|B\|^2_{C^0([0,T];H^{\frac{1}{2}}(\Gamma_0;{\mathbb R}))}
\nonumber\\
& \qquad\qquad 
   \left. 
+ \|B_t\|^2_{L^2(0,T;H^{\frac{1}{2}}(\Gamma_0;{\mathbb R}))} 
+ \|g\|^2_{L^2(0,T;H^{\frac{3}{2}}(\Gamma_0;{\mathbb R}^3))}
+ \|g_t\|^2_{L^2(0,T;H^{-\frac{1}{2}}(\Gamma_0;{\mathbb R}^3))} \right]\ ,
\nonumber
\end{align}
where $C(T)$ stays bounded as $T\rightarrow 0$.

3. We next want to prove that 
$$\nabla w_t \in L^2(0,T; {H^{-\frac{1}{2}}({\mathbb R}^2; {\mathbb R}^9)}), \ \
p_t \in L^2(0,T; {H^{-\frac{1}{2}}({\mathbb R}^2; {\mathbb R})})$$ 
with the required bounds.   Let $\phi \in H^{\frac{1}{2}}
({\mathbb R}^2; {\mathbb R})$ and for $\alpha=1,2$ consider the duality pairing
\begin{align*}
\int_0^T \langle \partial_\alpha w_t, \phi\rangle_{\Gamma_0}  dt
&= 
-\int_0^T \langle w_t, \partial_\alpha \phi\rangle_{\Gamma_0} dt\\
&\le \|w_t\|^2_{L^2(0,T; {H^{\frac{1}{2}}({\mathbb R}^2; {\mathbb R}^3)})}
 \|\partial_\alpha \phi\|^2_
{L^2(0,T; {H^{-\frac{1}{2}}({\mathbb R}^2; {\mathbb R})})}\\
&\le \|w_t\|^2_{L^2(0,T; {H^1({\mathbb R}^3_+; {\mathbb R}^3)})}
 \|\phi\|^2_
{L^2(0,T; {H^{\frac{1}{2}}({\mathbb R}^2; {\mathbb R})})}\,,
\end{align*}
where we have used the trace theorem for the last inequality. Thus,
\begin{equation}
\label{dualagain}
 \|\partial_\alpha w_t\|^2_
{L^2(0,T; {H^{-\frac{1}{2}}({\mathbb R}^2; {\mathbb R}^3)})}
\le \|w_t\|^2_{L^2(0,T; {H^1({\mathbb R}^3_+; {\mathbb R}^3)})}\,.
\end{equation}
Using the divergence-free condition $\partial_3 w^3_t =
-\partial_1 w^1_t -\partial_2 w^2_t$, we see that
$$
 \|\partial_3 w^3_t\|^2_
{L^2(0,T; {H^{-\frac{1}{2}}({\mathbb R}^2; {\mathbb R})})}
\le \|w_t\|^2_{L^2(0,T; {H^1({\mathbb R}^3_+; {\mathbb R}^3)})}
$$
as well.  It remains to verify this bound for $\partial_3 w^\beta_t$ for
$\beta=1,2$.  This follows from the fact that
$$\partial_3 w^\beta_t = - \partial_\beta w^3_t + g^\beta_t$$
so that
\begin{equation}\label{grad_wt}
 \|\nabla w_t\|^2_
{L^2(0,T; {H^{-\frac{1}{2}}({\mathbb R}^2; {\mathbb R}^9)})}
\le \|w_t\|^2_{L^2(0,T; {H^1({\mathbb R}^3_+; {\mathbb R}^3)})}
+ 
 \|g_t\|^2_
{L^2(0,T; {H^{-\frac{1}{2}}({\mathbb R}^2; {\mathbb R}^3)})}\,.
\end{equation}
Differentiating (\ref{linear8.c}) with respect to $t$, we see that
on the boundary ${\mathbb R}^2 \times \{x^3=0\}$, 
$$
p_t = 2 \nu \partial_3 w^3_t	- \sigma\triangle_0 w^3 -B_t \,, 
$$
so that using the estimate (\ref{B2.5}) we have verified the estimate
(\ref{estimate1_flat}), and have thus concluded the proof of the theorem.
\end{proof}

\subsection{Regularity of weak solutions on $\Omega_0$}
For the remainder of the paper, we shall make use of the
Gagliardo-Nirenberg inequality (\cite{Nirenberg1959}):
{\it
Suppose}
$$
\frac{1}{p} = \frac{i}{n} + a \left( \frac{1}{r} -
\frac{m}{n}
\right) + (1 - a ) \frac{1}{q}
$$
{\it where $ i/m \leq a
\leq 1$ (if $m - i - n/r $ is
an integer $ \geq 1 $, only $
a < 1 $ is allowed).
Then for $ f : \Omega \subset {\mathbb R}^n \rightarrow
\mathbb{R}^k $},
\begin{equation} \label{GN}
\|D^if\|_{L^{p}(\Omega; {\mathbb R}^k)} \leq
C\|D^m f\|^a_{L^{r}(\Omega; {\mathbb R}^k)}
\|f\|^{1-a}_{L^{q}(\Omega; {\mathbb R}^k)} \,.
\end{equation}

For an open subset $U,V$ in ${\mathbb R}^3$, we write
$$
U \subset \subset V
$$
if $U \subset \bar U \subset V$ and $\bar U$ is compact.
%For the purposes of this section, we may set $\nu=1$ and $\sigma=1$.

We shall need the following function spaces for the remainder of the paper.
\begin{definition}\label{XT}
Let
\begin{align*}
X_T= \{& (u,p)\in V^3 (T)\times V^2 (T)\ |\ \\
       &\ \ (\nabla u_t, p_t)\in 
{L^2(0,T; H^{-\frac{1}{2}}({\mathbb R}^2;{\mathbb R}^9))}
\times {L^2(0,T;H^{-\frac{1}{2}}({\mathbb R}^2;{\mathbb R}))}\}\ ,
\end{align*}
 endowed with its natural Hilbert norm
\begin{align*}
\|(u,p)\|^2_{X_T}=&\|  u \|^2_ {L^2(0,T;H^3({\mathbb R}^3_+;{\mathbb R}^3))}+ 
\|u_t\|^2_ {L^2(0,T;H^1({\mathbb R}^3_+;{\mathbb R}^3))}\nonumber \\
& + \|  p \|^2_ {L^2(0,T;H^2({\mathbb R}^3_+;{\mathbb R}))}+ 
\| p_t\|^2_ {L^2(0,T;L^2({\mathbb R}^3_+;{\mathbb R}))}\nonumber \\
& +\| \nabla u_t \|^2_ {L^2(0,T;H^{-\frac{1}{2}}({\mathbb R}^2;{\mathbb R}^9))}+ 
\| p_t\|^2_ {L^2(0,T;H^{-\frac{1}{2}}({\mathbb R}^2;{\mathbb R}))}\ ,
\end{align*}
and  let
\begin{align*}
Y_T=\{& (u,p)\in V^2 (T)\times V^1 (T)| \nonumber\\
&\ \  \ (u_t,p_t)\in L^2(0,T;H^1({\mathbb R}^3_+;{\mathbb R}^3))\times L^2(0,T;L^2({\mathbb R}^3_+;{\mathbb R}^3))  \nonumber\\
& \ \ (\nabla u_t, p_t)\in {L^2(0,T; H^{-\frac{1}{2}}({\mathbb R}^2;{\mathbb R}^9))}\times {L^2(0,T;H^{-\frac{1}{2}}({\mathbb R}^2;{\mathbb R}))}\}\ ,
\end{align*}
 endowed with its natural Hilbert norm
\begin{align*}
\|(u,p)\|^2_{Y_T}=&\|  u \|^2_ {L^2(0,T;H^2({\mathbb R}^3_+;{\mathbb R}^3))}+ 
\|u_t\|^2_ {L^2(0,T;H^1({\mathbb R}^3_+;{\mathbb R}^3))}\nonumber \\
& + \|  p \|^2_ {L^2(0,T;H^1({\mathbb R}^3_+;{\mathbb R}))}+ 
\| p_t\|^2_ {L^2(0,T;L^2({\mathbb R}^3_+;{\mathbb R}))}\nonumber \\
& +\| \nabla u_t \|^2_ {L^2(0,T;H^{-\frac{1}{2}}({\mathbb R}^2;{\mathbb R}^9))}+ 
\| p_t\|^2_ {L^2(0,T;H^{-\frac{1}{2}}({\mathbb R}^2;{\mathbb R}))}\ .
\end{align*}
\end{definition}

\begin{theorem}[Regularity on $\Omega_0$]\label{regularity2}
Let the data satisfy
$$
\begin{array}{c}
f\in V^1(T), \
g\in L^2(0,T;H^{\frac{3}{2}}(\Gamma_0, {\mathbb R}^3)), \
g_t\in L^2(0,T;H^{-\frac{1}{2}}(\Gamma_0, {\mathbb R}^3)), \\
B\in C^0([0,T];H^{\frac{1}{2}}(\Gamma_0, {\mathbb R})), \
B_t\in L^2(0,T;H^{\frac{1}{2}}(\Gamma_0, {\mathbb R}))\,,
\end{array}
$$
and $B(0, \cdot) =0$,
\begin{equation}
\label{trco1}
g(0,\cdot)_{\operatorname{tan}}=0\ ,
\end{equation}
and  let $w_0\in H^2 (\Omega_0, {\mathbb R}^3)\cap \mathcal V$ satisfy the associated compatibility condition
\begin{equation}
\label{trco2}
 (\operatorname{Def}( w_0)\ N)_{\operatorname{tan}}=0\ \ \text{on}\ \Gamma_0\ .
\end{equation}
Then the unique weak solution $w$ to (\ref{linear0}) satisfying the estimate (\ref{H1energy_p}) has the regularity $w \in V^3_{\operatorname{div}}(T)$, $p \in V^2(T)$, with 
the following estimate:
\begin{align}
&\| w\|_{L^2(0,T;H^3(\Omega_0;{\mathbb R}^3))} +
\|N \cdot \triangle_0 \int_0^{\cdot} w(s)ds\|_
{C^0([0,T];H^{\frac{1}{2}}(\Gamma_0;{\mathbb R}))} 
+ \| p\|_{L^2(0,T;H^2(\Omega_0;{\mathbb R}))} \nonumber \\
&+\| w_t\|_{L^2(0,T;H^1(\Omega_0;{\mathbb R}^3))} +
\|N \cdot \triangle_0  w\|_
{L^2(0,T;H^{\frac{1}{2}}(\Gamma_0;{\mathbb R}))} 
+ \| p_t\|_{L^2(0,T;L^2(\Omega_0;{\mathbb R}))} \nonumber \\
&+ \| \nabla  w_t\|_{L^2(0,T;H^{-\frac{1}{2}}(\Gamma_0;{\mathbb R}^9))}
+ \| p_t\|_{L^2(0,T;H^{-\frac{1}{2}}(\Gamma_0;{\mathbb R}))} \nonumber \\
& \qquad \le C(T) \Bigl( \|w_0\|_{H^2(\Omega_0; {\mathbb R}^3)} + 
\| f\|_ {L^2(0,T;H^1(\Omega_0;{\mathbb R}^3))}
+ \| f_t\|_ {L^2(0,T;H^1(\Omega_0;{\mathbb R}^3)')} \nonumber \\
& \qquad\qquad
+ \|f(0,\cdot)\|_{L^2(\Omega_0; {\mathbb R}^3)}
+ \| B\|_ {C^0([0,T];H^\frac{1}{2}(\Gamma_0;{\mathbb R}))}
+ \| B_t\|_ {L^2(0,T;H^{\frac{1}{2}}(\Gamma_0;{\mathbb R}))}
\nonumber \\
& \qquad\qquad  
+ \| g(0,\cdot)\|_ {H^\frac{1}{2}(\Gamma_0;{\mathbb R}^3)}
+\| g\|_ {L^2(0,T;H^\frac{3}{2}(\Gamma_0;{\mathbb R}^3))}
+ \| g_t\|_ {L^2(0,T;H^{-\frac{1}{2}}(\Gamma_0;{\mathbb R}^3))}
\Bigr)\,, \label{estimate1_omega}
\end{align}
where $C(T)$ remains bounded as $T \rightarrow 0$.
\end{theorem}
\begin{proof}

\smallskip
\noindent
{\bf 1. Interior regularity.} 
Choose an open set $U\subset \subset \Omega_0$  and an open set 
$W$ such that $U \subset\subset W \subset \subset \Omega_0$. 
Let $\zeta$ be a smooth cutoff function satisfying
$$
\left\{
\begin{array}{l}
\zeta \equiv 1 \ \text{ on } \ U, \ \ \zeta \equiv 0 \ \text{ on }
{\mathbb R}^3 - W, \\
0 \le \zeta \le 1.
\end{array}
\right.
$$

Let $w$ denote the weak solution of (\ref{linear0}), and set
\begin{equation}\label{Uw}
\chi= \zeta^3 w, \ \ \theta= \zeta^3 p
\end{equation}
Then $\chi$ is a solution of
\begin{subequations}
  \label{linear9}
\begin{align}
 \chi_t - \nu\triangle  \chi +\nabla \theta &= \zeta^3 f 
- 6\nu\zeta^2 \nabla w \cdot\nabla\zeta 
-3 \nu\left[ 2\zeta|\nabla \zeta|^2  +  \zeta^2\triangle \zeta \right]w 
+ 3 \zeta^2 \nabla \zeta  p \nonumber \\
& \qquad \qquad \qquad \qquad \qquad
 \text{in} \ \ (0,T)\times \Omega_0 \,, 
         \label{linear9.a}\\
\operatorname{div} \chi &= 3\zeta^2\nabla \zeta \cdot w 
\ \ \text{on} \ \ (0,T)\times \Omega_0 \,, 
         \label{linear9.b}\\
\chi&= \zeta^3w_0
\ \ \text{on} \ \ \Omega_0 \times \{t=0\} \,.
         \label{linear9.c}
\end{align}
\end{subequations}
It follows that (\ref{linear9}) is an initial-boundary value problem on
$(0,T]\times W$ with boundary value $\chi=0$ on $[0,T] \times \partial W$. This is
the classical Stokes problem for $(\chi,\theta)$ 
with Dirichlet boundary conditions, and
is well-known to have regularity properties.  We bootstrap applying
a new cut-off function whose support is contained strictly inside the
support of $\zeta$. Since initially, 
$w \in L^2(0,T; H^1(\Omega_0;{\mathbb R}^3))$, we find that  
$\chi \in L^2(0,T; H^2(\Omega_0;{\mathbb R}^3))$ so that by (\ref{Uw}), 
$w \in L^2(0,T; H^2(\Omega_0;{\mathbb R}^3))$; this in turn, shows that
$\chi \in L^2(0,T; H^3(\Omega_0;{\mathbb R}^3))$ 
and satisfies the 
estimate (\ref{estimate1_flat}) with all of the boundary terms set
to zero.

\smallskip
\noindent
{\bf 2: Boundary regularity.}  
We will localize our analysis by employing a partition of unity which is
subordinate to a specific collection of coordinate neighborhoods.  We begin
the analysis by constructing one such coordinate patch for a neighborhood of
the boundary $\Gamma _0:= \partial \Omega_0 $.

\begin{definition}[Exponential map]
For $x'\in \mathbb{R}^2  $ and $r>0$, let $D(0,r)$ denote the open disc of 
radius $r$ in $\mathbb{R}^2  $ centered at $x$.  With the induced metric
$g_0$, the pair $(\Gamma_0, g_0)$ is a compact, smooth, Riemannian manifold.
For $x_0 \in \Gamma_0$,  we denote the exponential map
$$\exp_{x_0}:   D(0,r) \subset T_{x_0} \Gamma _0 \rightarrow \Gamma _0$$
by
$ \exp_{x_0}(v) = \gamma(1,x_0,v)$, where $\gamma(1,x_0,v)$ is the time one
geodesic eminating from $x_0 \in \Gamma _0$ with velocity vector $v \in 
\mathbb{R}  ^2$.
\end{definition}

Using this definition, we can cover $\partial \Omega _0$ with finitely 
many {\it normal} or {\it geodesic} balls which are  images of $D(0,r)$
under the exponential map.  In particular, for $x_0 \in \Gamma_0$, we
set
$$
 \mathcal{D} (x_0,r) := \exp_{x_0}( D(0,r)),
$$
where $r>0$ is taken small enough to ensure that $\exp_{x_0}:D(0,r)
\rightarrow  \mathcal{D} (x_0,r)$ is a smooth diffeomorphism, and that
\begin{equation}\label{no_perp}
N(x')\not\perp N(x_0) \ \ \forall x' \in \mathcal{D} (x_0,r).
\end{equation} 

Next, with $B(0,r)$ denoting the open ball in $ \mathbb{R}  ^3$ of radius $r$,
let
$$
B^+(0,r) := \{ x \in  B(0,r) \ | x^3 >0\}
$$
denote the upper hemisphere. 
\begin{definition}[Coordinate map]
With $x'=(x^1,x^2)$ and $x=(x',x^3)$,
let $\Phi: B^+(0,r) \rightarrow \mathbb{R}  ^3$ be given by
$$
\Phi(x', x^3) = \exp_{x_0}(x') + x^3 N(\exp_{x_0}(x')).
$$
\end{definition}
\begin{definition}[Coordinate neighborhood]
Let 
$$U(x_0,r):= \Phi( B^+(0,r))$$ 
denote the local coordinate neighborhood of the boundary point $x_0 \in
\Gamma _0$.  We let $(\tilde x^1,\tilde x^2, \tilde x^3) := 
(\Phi^1, \Phi^2, \Phi^3)$ represent the coordinates on $U(x_0,r)$.
\end{definition}

Since $\exp_{x_0}(0)= x_0$ and $ \nabla \exp_{x_0}(0)= \text{Id}$, we 
can immediately record the following
\begin{lemma}\label{exp_lemma}
$\Phi(0) = x_0$ and $ \nabla \Phi(0,0)= \operatorname{Id}$.
Furthermore, for $x \in U(x_0,r)$, 
\begin{equation}\label{determinant}
\operatorname{Det} \nabla \Phi(x)= 
\operatorname{Det} \nabla \exp_{x_0}(x).
\end{equation} 
\end{lemma}

Let us use $\tilde x^i$ as our Cartesian coordinates on $\Omega_0$; we
now localize the problem to $U(x_0,r)$.  Let $\zeta$ be a smooth cutoff 
function satisying
$$
\left\{
\begin{array}{l}
\zeta \equiv 1 \ \text{ on } \ U(x_0,\frac{r}{4}), 
\ \ \zeta \equiv 0 \ \text{ on }
U(x_0, \frac{r}{2})^c , \\
0 \le \zeta \le 1 .
\end{array}
\right.
$$
With $w(\tilde x)$ denoting the solution of (\ref{linear0}), and once again 
setting
$\chi= \zeta^3 w$ and $\theta= \zeta^3 p$ as in (\ref{Uw}), we see that  
the pair ($\chi(t, \tilde x), \theta(t, \tilde x))$ solves
\begin{subequations}
  \label{linear10}
\begin{align}
 \chi_t - \nu  \triangle\chi +  \nabla \theta &=  \mathfrak{F}(\zeta,w,
 \nabla w)  \ \  \text{ in }  (0,T)\times U(x_0,r) \,, 
         \label{linear10.a}\\
\operatorname{div} \chi &=  \mathfrak{A}(\zeta,w)
\ \ \text{ in }  (0,T)\times U(x_0,r) \,, 
         \label{linear10.b}\\
\nu(\chi^i,_j+\chi^j,_i)N^j - \theta N^i
&-\sigma N^k \triangle_0\int_0^t \chi^k(r)dr\,  N^i \nonumber\\
&= 
 \sigma \mathfrak{B}(\zeta,w, \nabla _0w)N^i + \mathfrak{G}^i 
\ \ \text{ on }  (0,T) \times \partial U(x_0,r)  \,, 
         \label{linear10.c}\\
\chi&= \zeta^3w_0
\ \text{ on }  U(x_0,r) \times \{t=0\} \,,
         \label{linear10.d}
\end{align}
\end{subequations}
where
\begin{subequations}
\begin{align}
\mathfrak{F}^i(\zeta,w, \nabla w) &=
\zeta^3f^i - 3\zeta^2\zeta,_i p 
 + \nu [ 6\zeta \zeta,_j
\zeta,_j w^i + 3 \zeta^2 \zeta,_{jj}w^i + 3 \zeta^2 \zeta,_j w^i,_j]
\label{force1} \\ 
\mathfrak{A}(\zeta,w) &=  3 \zeta^2 \zeta,_k)w^k
\label{force2}\\
\mathfrak{B}(\zeta,w, \nabla _0 w)&=
- N^i g^{\alpha\beta}_0 \int_0^t \left[
6 \zeta \zeta,_ \alpha  \zeta,_{\beta} w^i + 3 \zeta^2 \zeta,_\alpha w^i,_\beta
+3 \zeta ^2 w^i \zeta ,_{ \alpha \beta }  \right. \nonumber \\
& \qquad \left.
- 3 \zeta^2 (\Gamma_0)^\gamma_{\alpha\beta} \zeta,_\gamma w^i \right]dr
+ \zeta^3 B
\label{force3}\\
\mathfrak{G}^i(\zeta,w)&=
\zeta^3g^i
+ 3\nu \zeta^2 [(w\cdot N)\ \zeta,_i + (\nabla\zeta\cdot N)\ w^i] ,
\label{force4}
\end{align}
\end{subequations}
and where we continue to use the notation $g_0$ and $\Gamma_0$ to denote
the metric and connection on $\Gamma_0$, respectively, at time $t=0$ 
when $\eta(0,\cdot)
= \operatorname{Id}$ (see (\ref{metric})).  The problem 
(\ref{linear10}) is localized to $U(x_0,r)$.

We next transfer the problem (\ref{linear10}) from $(U(x_0,r), \{ \tilde x\})$ 
to $(B^+(0,r),\{x^i\})$.
We set
$$
W(t,x) = \chi(t,\tilde x), \ \ Q(t,x)=\theta(t, \tilde x) ,
$$
where $\tilde x = \Phi(x)$.  With
$$
\begin{array}{c}
A(x) = [ \nabla \Phi(x)]^{-1}  \ \ \forall \ x\in B^+(0,r), \\
\mathcal{A} (x')
 = [ \nabla \exp_{x_0}(x')]^{-1}  \ \ \forall \ x'\in D(0,r), \\
\end{array}
$$
the pair $(W(t,x), Q(t,x))$ satisfies, in the weak sense, the  following
system of equations:
\begin{subequations}
  \label{linear11}
\begin{align}
W^i_t - \nu A^s_k(A^r_k  W^i,_r),_s  + A^k,_i Q,_k
&= \mathfrak{F}^i(\Phi) \ \  \text{ in }  (0,T)\times B^+(0,r) \,, 
         \label{linear11.a}\\
A^j_i W^i,_j & = \mathfrak{A}(\Phi)
\ \ \text{ on }  (0,T)\times B^+(0,r) \,, 
         \label{linear11.b}\\
[\nu (W^i,_r A^r_j + W^j,_r A^r_i) N_j(\Phi) - Q N_i(\Phi)]\, | A\  e_3|
 & \nonumber\\
= \sigma N_k(\Phi) 
\mathfrak{g} ^{-1} [\mathfrak{g}g_0^{\alpha\beta} &\int_0^t W^k,_\gamma(r)\, dr
\mathcal{A}^\gamma_\beta],_\delta
\mathcal{A}^\delta_\alpha \,  N_i(\Phi)\nonumber\\
 + \sigma [\mathfrak{B}N^i](\Phi) &+ \mathfrak{G}^i(\Phi)
\ \ \text{ on }  (0,T) \times  \partial B^+(0,r)  \,, 
         \label{linear11.c}\\
W&= \zeta^3 w_0(\Phi)
\ \text{ on }  B^+(0,r) \times \{t=0\} \,,
         \label{linear11.d}
\end{align}
\end{subequations}
where $\mathfrak{g}= \det(g_0)$.
Our notation 
$f(\Phi)$ denotes $f(t,\Phi(x))$.  We have made use of (\ref{determinant}) to
remove the Jacobian factors in the boundary condition (\ref{linear11.c}).
Note that $N(\Phi(x_0))= e_3$.

\begin{remark}
Since $\zeta$ has compact support in $B^+(0,r)$,
and since the forcing functions in (\ref{linear11}) vanish on the complement 
of $B^+(0,r)$, we can extend $W$ and $Q$ to
the entire half-space $ \mathbb{R}^3  _+$, and consider (\ref{linear11}) as a
system of PDE on $[0,T] \times \mathbb{R}^3  _+$.
\end{remark}

The next step is to
show that there exists a unique solution in $ X_T$ to (\ref{linear11}) 
(after first proving that it belongs to $Y_T$). Since it is
readily seen that the existence and uniqueness of a weak solution to
(\ref{linear11}), extended to ${\mathbb R}^3_+$, can be established
exactly in the same fashion as we established Theorem \ref{uniqueweak}, our
weak solution $(W,Q)$ will thus be in $X_T$, which will establish the desired
result. 

We initiate an iteration sequence starting with 
$(v_0,q_0)$, where $q_0=0$ and $v_0$ is the solution of the parabolic problem
\begin{subequations}
  \begin{align}
\partial_tv_0- \nu \triangle {{v_0}}&= 0
\ \ \text{in} \ \ (0,T)\times {\mathbb R}^3 \,, 
         \\
{v_0} &= \overline{(\zeta^3\ w_0)(\Phi)}   
 \ \ \text{on} \ \ {\mathbb R}^3\times \{ t=0\} \,, 
         \label{inicombis}
\end{align}
\end{subequations}
where $\overline{(\zeta^3\ w_0)(\Phi)}$ denotes an $H^2$-extension of ${(\zeta^3\ w_0)(\Phi)}$ to ${\mathbb R^3}$. We set the problem on
${\mathbb R}^3$ so as to avoid any compatibility condition at time $t=0$.

By using the same techniques on this problem as we developed in the 
previous section (which is easier since the problem is set on the full 
space $ \mathbb{R}^3  $), we find that 
\begin{equation}
\label{inicom}
\|(v_0,0)\|_{X_T}\le C \|w_0\|_{H^2(\Omega_0; {\mathbb R}^3)}\ .
\end{equation}

Then, setting
\begin{equation}\label{Gforce}
G  = \Pi_2\Pi_0(\Phi) \left[ \sigma \mathfrak{B} N(\Phi) + \mathfrak{G}(\Phi)\right]
+e_3 \cdot \left[ \sigma \mathfrak{B} N(\Phi) + \mathfrak{G}(\Phi)\right]  e_3 \,,
\end{equation} 
given $(v_n,q_n) \in X_T$, we define $(v_{n+1},q_{n+1})$ 
to be the solution of
\begin{subequations}
  \label{linear12}
\begin{align}
\partial_t v_{n+1}- \nu \triangle {{v_{n+1}}} +
\nabla q_{n+1}&= \mathfrak{F}(\Phi) + \bar F (v_n, q_n)
\ \ \text{in} \ \ (0,T)\times {\mathbb R}^3_+ \,, 
         \label{linear12.a}\\
\operatorname{div} (v_{n+1})&= \bar A (v_n) + \mathfrak{A}(\Phi) 
\ \ \text{in} \ \ (0,T)\times {\mathbb R}^3_+\,, 
         \label{linear12.b}\\
S({v_{n+1}},q_{n+1})\ e_3\ - 
\sigma\ e_3&\cdot 
[ g_0^{\alpha\beta} (x_0)\ \partial^2_{ \alpha \beta} \int_0^tv_{n+1}(r) dr]
\ e_3 \nonumber\\
&=G+ \bar G_1 (v_n)+[\bar G_2 (v_n) 
+  \sigma \bar B (v_n)]\ e_3 \nonumber\\
&\qquad\qquad\qquad \text{on} \ \ (0,T)\times [{\mathbb R}^2\times\{0\}]\,, 
\label{linear12c}\\ 
        \\
{v_{n+1}} &= (\zeta^3\ w_0)(\Phi)   
 \ \ \text{on} \ \ {\mathbb R}^3_+\times \{ t=0\} \,, 
         \label{linear12d}
\end{align}
\end{subequations}
where $S(v,q)=\nu \operatorname{Def} v - q N$, 
$\Pi_2$ denotes the tangential projection onto 
${\mathbb R}^2\times \{0\}$, $\Pi_0$ continues to denote the projection 
onto the tangent plane of $\Gamma_0$, and
\begin{subequations}
  \label{linear13}
\begin{align}
\bar F (v_n, q_n) &=
- \nu (\triangle {v_n} - \operatorname{div}
[ \nabla v_n \cdot A \cdot A^T]) + \nabla
q_{n}-A^T \nabla q_{n}
\ \ \text{in} \ \ (0,T)\times {\mathbb R}^3_+ \,, 
  \label{linear13.a}       \\
\bar A (v_n)     &=
  \operatorname{div} (v_{n})-A^j_i
\frac{\partial {v_n}^i}{\partial x^j } 
\ \ \text{in} \ \ (0,T)\times {\mathbb R}^3_+ \,, 
  \label{linear13.b}       \\
 \bar G_1 (v_n)
&= 
\Pi_2\left[S(v_{n}, q_n)\ e_3\right]
 - \Pi_2\Pi_0(\Phi)\left( \nu \left(\nabla v_n \cdot A
+ [ \nabla v_n \cdot A]^T
\right)\cdot  N(\Phi)\ |A\ e_3|\right)\nonumber\\
 &\qquad+\Pi_2\Pi_0(\Phi)\left(  q_n N(\Phi) \ |A \ e_3|\right)  
 \ \ \text{in} \ \ (0,T)\times {\mathbb R}^2\times\{0\}\,, 
  \label{linear13.c}       \\
\bar G_2 (v_n)
&= 
e_3\cdot \left(S(v_{n}, q_n)\ e_3\right) - e_3\cdot 
\left( \nu \left(\nabla v_n \cdot A + [ \nabla v_n \cdot A]^T
\right)\cdot N(\Phi) \ |A e_3|\right)\nonumber\\
&\qquad\qquad\qquad\qquad +
e_3\cdot\left(  q_n  N(\Phi) \ |A \ e_3| \right)  
\ \ \text{in} \ \ (0,T)\times {\mathbb R}^2\times\{0\}\,, 
  \label{linear13.d}       \\
\bar B (v_n)
&= N^i(\Phi) \mathfrak{g} ^{-1} \partial_ \delta \left[ \mathfrak{g}
g_0^{\alpha\beta} \int_0^t \partial_ \gamma v_n^i(r) dr
\mathcal{A} ^ \gamma _ \beta \right]
\mathcal{A} ^\delta_\alpha 
\ N(\Phi)\cdot e_3\nonumber\\
& -e_3\cdot g_0^{\alpha\beta} (x_0)
\int_0^t\partial^2_{\alpha \beta}v_{n}(r)\ dr 
\ \ \text{on} \ \ (0,T)\times {\mathbb R}^2\times\{0\} \,. 
  \label{linear13.e}  
\end{align}
\end{subequations}

\begin{remark}
Note the presence of the projector $\Pi_0$ in the forcing functions
$G$ defined in (\ref{Gforce}) and  $\bar G_1$ defined
in (\ref{linear13.b}).  This projector is used to remove the less regular
terms, such as $B$ and $\mathfrak{B}$, from the tangential component of the
forcing.
\end{remark}

For each $n\ge 1$, we define
$$
\begin{array}{l}
{\delta v}_n=v_{n+1}-v_{n} \\
{\delta q}_n=q_{n+1}-q_n 
\end{array}
\begin{array}{l}
\text{(difference of velocities)}\\
\text{(difference of pressures)\, . }
\end{array}
$$
Since for each iteration $n\ge 0$, the problem 
(\ref{linear12}) is linear, the pair 
$(\delta v_{n+1}, \delta q_{n+1})$ satisfies
\begin{subequations}
  \label{linear15}
\begin{align}
\partial_t \delta v_{n+1}- \nu \triangle {{\delta v_{n+1}}} +
\nabla \delta q_{n+1}&=  \bar F (\delta v_n,\delta q_n)
\ \ \text{in} \ \ (0,T)\times {\mathbb R}^3_+ \,, 
         \\
\operatorname{div} (\delta v_{n+1})&
= \bar A (\delta v_n)    \ \ \text{in} \ \ (0,T)\times {\mathbb R}^3_+\,, \\
 \ S({\delta v_{n+1}},\delta q_{n+1})\ e_3\ - \sigma\ e_3&
\cdot ( g_0^{\alpha\beta} (x_0)\ 
\int_0^t \partial^2_{\alpha\beta}\delta v_{n+1}(r) dr)\ e_3 \nonumber\\ 
& = \bar G_1 (\delta v_n)+(\bar G_2 (\delta v_n) 
+  \sigma \bar B (\delta v_n))\ e_3  \nonumber \\
&\qquad\qquad\qquad 
\text{on} \ \ (0,T)\times {\mathbb R}^2\times\{0\} \,, 
        \\
{\delta v_{n+1}} &= 0   
 \ \ \text{on} \ \ {\mathbb R}^3_+\times \{ t=0\} \,, 
        \end{align}
\end{subequations}
where the forcing terms 
appearing on the right-hand-side of this system are defined 
by the same relations as (\ref{linear13}) with $(\delta v_n,\delta q_n)$ 
replacing $(v_n,q_n)$.

For each $n\ge 1$, the initial data satisfies the compatibility condition 
$$\Pi_2(\operatorname{Def} \left({\delta v}_{n+1}(0,\cdot)\right)N)
=0\ \text{on}\ {\mathbb R}^2\times \{0\},$$
associated to $\bar{G_1}(\delta v_n)(0,\cdot)=0$, 
$\bar{G_2}(\delta v_n)(0,\cdot)=0$ and $\bar{B}(\delta v_n)(0,\cdot)=0$. 
Thus, (\ref{linear15}) satisfies all the assumptions of the basic linear
problem (\ref{linear7}), with the elliptic operator $v^3,_{\alpha\alpha}$ 
replaced  by the constant coefficient elliptic operator 
$g_0^{\alpha\beta} (x_0)\ v^3,_{\alpha\beta}$ (which is also coercive); thus,
using the $H^2$-type estimates of the previous section, we  find that
\begin{align}
&\| ({\delta v}_{n+1},{\delta q}_{n+1})\|_{Y_T} \nonumber\\
&\le   C \ \left(
\| \bar F ({\delta v}_{n},{\delta q}_{n})\|_ 
{L^2(0,T;L^2({\mathbb R}^3_+;{\mathbb R}^3))}+ 
\| \bar F ({\delta v}_{n},{\delta q}_{n})_t\|_ 
{L^2(0,T;H^1({\mathbb R}^3_+;{\mathbb R}^3)')}
 \right.\nonumber \\
& \qquad\qquad  
+ \|\bar A ({\delta v}_{n})\|_ {L^2(0,T;H^1({\mathbb R}^3_+;{\mathbb R}^9))} 
+ \|\bar A ({\delta v}_{n})_t\|_ {L^2(0,T;H^1({\mathbb R}^3_+;{\mathbb R}^9)')} \nonumber \\
& \qquad\qquad  
+ \|\bar B ({\delta v}_{n})\|_ {L^2(0,T;H^{-\frac{1}{2}}({\mathbb R}^2;{\mathbb R}))}
+ \|\bar B({\delta v}_{n})_t\|_ {L^2(0,T;H^{-\frac{1}{2}}({\mathbb R}^2;{\mathbb R}))}\nonumber \\
&\qquad\qquad +\|\bar G_1({\delta v}_{n},{\delta q}_{n})\|_ {L^2(0,T;H^\frac{1}{2}({\mathbb R}^2;{\mathbb R}^3))}
 +\|\bar G_1({\delta v}_{n},{\delta q}_{n})_t\|_ {L^2(0,T;H^{-\frac{1}{2}}({\mathbb R}^2;{\mathbb R}^3))}\nonumber\\
& \qquad\qquad \left. +\|\bar G_2({\delta v}_{n},{\delta q}_{n})\|_ {L^2(0,T;H^\frac{1}{2}({\mathbb R}^2;{\mathbb R}))}
 +\|\bar G_2({\delta v}_{n},{\delta q}_{n})_t\|_ {L^2(0,T;H^{-\frac{1}{2}}({\mathbb R}^2;{\mathbb R}))}\right)
\ . \nonumber
\end{align}
Consequently, taking the support of $\zeta$ in a small enough neighborhood 
of $x_0$,
$\nabla\Phi$ is sufficiently close to Id in $L^\infty$, so that the 
differences between  quantities appearing in (\ref{linear13}), which are
evaluated at $x$ and at $0$,  are sufficiently small so as to allow us 
to deduce that
\begin{equation}
\label{sequencecontraction}
\| ({\delta v}_{n+1},{\delta q}_{n+1})\|_{Y_T}\le \frac{1}{2}\ \| 
({\delta v}_{n},{\delta q}_{n})\|_{Y_T} \ ;
\end{equation}
by the contraction mapping principle, the sequence 
$(v_n,q_n)_{n\in {\mathbb N}}$ is convergent in $Y_T$ to a limit $(\bar w,r)$ which 
is also the unique solution of the extension of (\ref{linear11}) to ${\mathbb R}^3_+$. To see this (which is not
completely straightforward at first sight), we first notice that $(\bar w,r)$ is the
unique solution in $X_T$ of the problem
 
\begin{subequations}
  \label{linear12bis}
\begin{align}
\partial_t \bar w- \nu \triangle {\bar w} +
\nabla r&= \mathfrak{F}(\Phi) + \bar F (\bar w, r)
\ \ \text{in}\ \ (0,T)\times {\mathbb R}^3_+ \,, 
         \label{linear12bis.a}\\
\operatorname{div} (\bar w)&= \bar A (\bar w) + \mathfrak{A}(\Phi)
\ \ \text{in} \ \ (0,T)\times {\mathbb R}^3_+\,, 
         \label{linear12bis.b}\\
S(\bar w,r)\ e_3\ - \sigma\ e_3&\cdot (\int_0^t {g_0^{\alpha\beta} 
(x_0)\ \bar w,_{\alpha\beta}\ dt'})\ e_3 \nonumber\\
&=G+ \bar G_1 (\bar w)+(\bar G_2 (\bar w) 
+  \sigma \bar B (\bar w))\ e_3 \nonumber\\
&\qquad\qquad\qquad \text{on} \ \ (0,T)\times [{\mathbb R}^2\times\{0\}]\,, 
\label{linear12bis.c}\\ 
        \\
\bar w &= (\zeta^3\ w_0)(\Phi)   
 \ \ \text{on} \ \ {\mathbb R}^3_+\times \{ t=0\} \,, 
         \label{linear12bis.d}
\end{align}
\end{subequations}
From the definitions of  $\bar A$, $\bar F$, $\bar G$, $\bar B$, 
$\mathfrak{A}(\Phi)$, $\mathfrak{F}(\Phi)$, and $G$, we have that
\begin{subequations}
  \label{linear12ter}
\begin{align}
\partial_t\bar w- \nu \operatorname{div}
[ \nabla \bar w\cdot A \cdot A^T] +
A^T \nabla r &= \mathfrak{F}(\Phi) 
\ \ \text{in}\ \ (0,T)\times {\mathbb R}^3_+ \,, 
         \label{linear12ter.a}\\
A^j_i
\frac{\partial {\bar w}^i}{\partial x^j }&= \mathfrak{A}(\Phi)   
\ \ \text{in} \ \ (0,T)\times {\mathbb R}^3_+\,, 
         \label{linear12ter.b}\\
\Pi_2 \Pi_0 (\Phi) (\Sigma) + (\Sigma\cdot e_3)\ e_3 + \lambda (N(\Phi)\cdot e_3)\ e_3 
&=0\  \text{on} \ \ (0,T)\times [{\mathbb R}^2\times\{0\}]\,, 
\label{linear12ter.c}\\ 
\bar w &= (\zeta^3\ w_0)(\Phi)   
 \ \ \text{on} \ \ {\mathbb R}^3_+\times \{ t=0\} \,, 
         \label{linear12ter.d}
\end{align}
\end{subequations}
where
\begin{align*}
\Sigma=& \left(\nu \left( \nabla  \bar w \cdot A
+ [ \nabla \bar w \cdot A]^T
\right)\cdot N(\Phi)- r N(\Phi)\right) |A \ e_3| - \mathfrak{G}(\Phi)-\sigma \mathfrak{B}(\Phi)\,,
\end{align*}
and
$$
\lambda=-\sigma  N^i(\Phi)\mathfrak{g} ^{-1}  
\partial_\delta [ \mathfrak{g} g_0^{ \alpha \beta } \int_0^t \partial_ \gamma 
\bar w ^i(r)dr \mathcal{A} ^\gamma_\beta]
\mathcal{A}^\delta_\alpha \ N(\Phi)\cdot e_3 \,.
$$

By (\ref{no_perp}), $N(\Phi)$ and $e_3=N(x_0)$ are not orthogonal, so that
\begin{equation*}
\Pi_2\Pi_0 (\Phi)(\Sigma)=0\ \ \text{and}\ \  \Sigma\cdot e_3 + \lambda (N(\Phi)\cdot e_3) 
=0\  \text{on} \ \ (0,T)\times [{\mathbb R}^2\times\{0\}]\ ,
\end{equation*}
which in turn implies that
\begin{equation*}
\Pi_0 (\Phi)(\Sigma)=\kappa e_3\ (\kappa\in {\mathbb R})\  \ \text{and}\ \  \Sigma\cdot e_3 + \lambda (N(\Phi)\cdot e_3) 
=0\  \text{on} \ \ (0,T)\times [{\mathbb R}^2\times\{0\}]\ .
\end{equation*}
Again,  since $N(\Phi)$ and $e_3$ are not orthogonal,
\begin{equation*}
\Pi_0 (\Phi)(\Sigma)=0\ \  \ \text{and}\ \  \Sigma\cdot e_3 + \lambda (N(\Phi)\cdot e_3) 
=0\  \text{on} \ \ (0,T)\times [{\mathbb R}^2\times\{0\}]\ ,
\end{equation*}
which provides us with
\begin{equation*}
\Sigma=\kappa' N(\Phi)\ (\kappa'\in {\mathbb R})\ \  \ \text{and}\ \  \Sigma\cdot e_3 + \lambda (N(\Phi)\cdot e_3) 
=0\  \text{on} \ \ (0,T)\times [{\mathbb R}^2\times\{0\}]\ ,
\end{equation*}
and finally,
\begin{equation}
\label{geometry}
\Sigma=-\lambda N(\Phi)\  \text{on} \ \ (0,T)\times [{\mathbb R}^2\times\{0\}]\ ,
\end{equation}
By using the definitions of $\Sigma$ and $\lambda$ in (\ref{geometry}), we find
that the limit of our contractive sequence $(\bar w,r)$ is the unique solution 
of (\ref{linear11}), extended to the half-space ${\mathbb R}^3_+$.

Since $(W,Q)$ is the unique weak solution of (\ref{linear11}), we have shown
that $(W,Q)=(\bar w, r)$. We can then estimate the norms of $(W,Q)$ in $Y_T$ 
and $X_T$ as follows.
By summing the 
inequalities (\ref{sequencecontraction}) from $n=1$ to $\infty$ we see that
\begin{equation*}
\| (W,Q)\|_{Y_T} \le \ 
 \| ({ v}_{1}, {q}_{1})\|_{Y_T} + \|(v_0,0)\|_{Y_T}\ ,
\end{equation*}
which, thanks to (\ref{inicom}), shows that
\begin{equation}
\label{sumbis}
\| (W,Q)\|_{Y_T} \le \ 
 \| ({ v}_{1}, {q}_{1})\|_{Y_T} + C \|w_0\|_{H^2 (\Omega_0;{\mathbb R}^3)}\ .
\end{equation}

To estimate the right-hand-side of (\ref{sumbis}), we  apply our previous
estimate again on the system
\begin{align}
\partial_t{{v_{1}}}- \nu \triangle {{v_{1}}} +
\nabla q_{1}&= \mathfrak{F} (\Phi) +\bar{F}(v_0,0)\ \ \text{in} \ \ (0,T)\times {\mathbb R}^3\,, 
   \nonumber      \\
  \operatorname{div} (v_{1})&=  \mathfrak{A} (\Phi) + \bar{A}(v_0)  \ \ \text{in} \ \ (0,T)\times {\mathbb R}^3\,, 
    \nonumber      \\
\ S({v_{1}},q_{1})\ e_3\  -
\sigma\ \left(e_3\cdot \int_0^t {g^{\alpha\beta}_0 (x_0)\ {v_1},_{\alpha\beta}}\right)\ e_3&= G + \bar{G_1}(v_0)+ (\bar{G_2}(v_0)+\sigma \bar{B}(v_0) )e_3\ \nonumber\\ & \qquad \text{on} \ \ (0,T)\times \ {\mathbb R}^2\times\{0\} \,, 
      \nonumber    \\
  {v_{1}} &= (\zeta^3\ w_0)(\Phi)    
 \ \ \text{on} \ \ {\mathbb R}^3_+ \times \{ t=0\} \nonumber\,, 
\end{align}
whose initial data satisfy the compatibility condition 
 at time $t=0$, given by
 \begin{equation}\label{ssss1}
\nu\ \Pi_2 (\operatorname{Def} \left(v_{1}(0,\cdot)\right)\ e_3)=\Pi_2 \left(G(0,\cdot)+\bar{G_1}(v_0(0,\cdot))\right)\ \ \ \text{on}\ 
{\mathbb R}^2\times \{0\}.
\end{equation} 
To see that (\ref{ssss1}) indeed holds, let us start from the condition
$$\Pi_0 (\operatorname{Def}(w_0))=0\ \text{on}\ \Gamma_0\ .$$
Since $W_0=(\zeta^3 w_0)({ \Phi})$, the previous condition is 
equivalent to
\begin{align}
&\Pi_0 (\Phi) \left(|A\ e_3|\ \nu \left(\nabla W_0 \cdot A + [ \nabla W_0 \cdot A ]^T
\right)\cdot   N (\Phi) 
\right)\nonumber\\
&\qquad\qquad\qquad\qquad\qquad= \Pi_0 (\Phi) \left(3\nu \zeta^2 (\Phi) [(w_0\cdot N)\ \nabla\zeta+ (\nabla\zeta\cdot N)\ w_0](\Phi)\right) \nonumber \\
& \qquad \qquad \qquad\qquad \qquad \qquad
\ \text{ on }  (0,T) \times {\mathbb R}^2\times\{0\}\ .
\label{transformcc}
\end{align}
Since $W_0=v_0$ at $t=0$,
$$\Pi_2 \left(G(0,\cdot)+
\bar{G_1}(v_0(0,\cdot))\right)= \Pi_2 (G(0,\cdot)+\bar{G_1}(W_0))\ .$$
Using the definitions of $G$ and $\bar{G_1}$ together with the condition
(\ref{trco1}), we have that
\begin{align*}
\Pi_2 (G+\bar{G_1}(v_0(0,\cdot))&= \Pi_2\left(S(W_0,0)\ e_3\right)\nonumber\\
&\qquad - \Pi_2\Pi_0(\Phi)\left( |A\ e_3|\ \nu \left(\nabla W_0 \cdot A
+ [ \nabla W_0 \cdot A]^T
\right)\cdot N(\Phi)\right)\nonumber\\
& \qquad + \Pi_2 \Pi_0 (\Phi) \left(3\nu \zeta^2 (\Phi) [(w_0\cdot N)\ \nabla\zeta+ (\nabla\zeta\cdot N)\ w_0](\Phi)\right)
\ ,
\end{align*}
which allows us to infer from (\ref{transformcc}) that 
$$\Pi_2 \left(G(0,\cdot) +\bar{G_1}(v_0(0,\cdot))\right)= \Pi_2\left(S(W_0,0)\ e_3\right)\ ,$$
so that the compatibility condition (\ref{ssss1}) is satisfied.

It follows that
\begin{align}
\label{v1q1bis}
&\| ({v}_{1}, {q}_{1})\|_{Y_T} \nonumber \\
 & \qquad \le C \ \left(\ \|(\zeta^3 w_0)(\Phi)\|
 _{H^1({\mathbb R}^3_+;{\mathbb R}^3)}+ 
\|  \mathfrak{F} (\Phi) +\bar{F} (v_0)\|_ {L^2(0,T;L^2({\mathbb R}^3_+;{\mathbb R}^3))} \right.\nonumber\\
&\qquad\qquad + \|  (\mathfrak{F} (\Phi))_t +\bar{F} (v_0,0)_t\|_ {L^2(0,T;H^1({\mathbb R}^3_+;{\mathbb R}^3)')} + \|  \mathfrak{A} (\Phi) +\bar{A}(v_0)\|_ {L^2(0,T;H^1({\mathbb R}^3_+;{\mathbb R}^3))} \nonumber\\
&\qquad\qquad +
\|  (\mathfrak{F} (\Phi))_t +\bar{A}(v_0)_t\|_ {L^2(0,T;H^1({\mathbb R}^3_+;{\mathbb R}^3)')}\nonumber\\
&\qquad\qquad +\|\tilde G(0,\cdot) +\bar{G_1}(v_0)(0,\cdot)+\bar{G_2}(v_0)(0,\cdot)e_3\|_{H^\frac{1}{2}({\mathbb R}^2;{\mathbb R}^3)}\nonumber\\
&\qquad\qquad + \| \ \tilde G+\bar{G_1}(v_0) + \bar{G_2}(v_0)\ e_3\|_ {L^2(0,T;H^\frac{1}{2}({\mathbb R}^2;{\mathbb R}^3))}\nonumber\\
&\qquad\qquad +
\| \ \tilde G_t + \bar{G_1}(v_0)_t + \bar{G_2}(v_0)_t\ e_3\|_ {L^2(0,T;H^{-\frac{1}{2}}({\mathbb R}^2;{\mathbb R}^3))}
\nonumber\\
&\qquad\qquad\left. + \| \tilde B \|_ {L^2(0,T;H^{-\frac{1}{2}}({\mathbb R}^2;{\mathbb R}))}
+ \|\ \tilde B_t\ \|_ {L^2(0,T;H^{-{\frac{1}{2}}}({\mathbb R}^2;{\mathbb R}))}\right)\ ,  
\end{align}
with
\begin{align*}
\tilde B &= \Pi_2 \Pi_0 (\Phi) (\sigma \mathfrak{B} (\Phi) N (\Phi)) 
+ e_3\cdot \sigma \mathfrak{B} (\Phi) N (\Phi)\  e_3\ ,\\
\tilde G &= \Pi_2 \Pi_0 (\Phi) (\mathfrak{G} (\Phi) )
+ e_3\cdot \mathfrak{G} (\Phi)\  e_3\ ,
\end{align*}
so that $\tilde B (0,\cdot)=0$ and $G=\tilde B+ \tilde G$.
Thanks to this splitting of the tangential forcing, we may infer from 
(\ref{sumbis}) and (\ref{v1q1bis}) that
\begin{align}
\label{wqh2bis}
\|(W,Q)\|_{Y_T} \le & C\ (\ \|w_0\|_{H^2(\Omega_0;{\mathbb R}^3)}
+ \|f\|_{L^2(0,T;L^2(\Omega_0;{\mathbb R}^3))}
+\|f_t\|_{L^2(0,T;H^1(\Omega_0;{\mathbb R}^3)')}\nonumber\\
&  
+\|  g(0,\cdot)\|_ {H^\frac{1}{2}(\Gamma_0;{\mathbb R}^3)} 
+\|  g\|_ {L^2(0,T;H^\frac{1}{2}(\Gamma_0;{\mathbb R}^3))}+
\|  g_t\|_ {L^2(0,T;H^{-\frac{1}{2}}(\Gamma_0;{\mathbb R}^3))}
\nonumber\\
& +\|  f(0,\cdot)\|_ {L^2(\Omega_0;{\mathbb R}^3)} 
+ \|  B \|_ {L^2(0,T;H^{-\frac{1}{2}}(\Gamma_0;{\mathbb R}))}
+ \|\  B_t\ \|_ {L^2(0,T;H^{-{\frac{1}{2}}}(\Gamma_0;{\mathbb R}))}\nonumber\\
& +\|w\|_{L^2(0,T;H^1(\Omega_0;{\mathbb R}^3))}
+\|\nabla w_t \|_{L^2(0,T;L^2(\Omega_0;{\mathbb R}^3))}\nonumber\\
&  + \|\partial_{\alpha} w\|_{L^2(0,T;L^2(\Gamma_0;{\mathbb R}^3))}
+ \|(v_0,0)\|_{Y_T})\ .
\end{align}
Since $w$ and $w_t$ satisfy the estimate
\begin{align*}
&\|w\|_{L^2(0,T;H^1(\Omega_0;{\mathbb R}^3))}
+\|\nabla w_t\|_{L^2(0,T;L^2(\Omega_0;{\mathbb R}^3))}
+ \|\partial_{\alpha} w\|_{L^2(0,T;L^2(\Gamma_0;{\mathbb R}^3))}\ \nonumber\\
&\qquad\le C\ (\ \|w_0\|_{L^2(\Omega_0;{\mathbb R}^3)}+\|f\|_{L^2(0,T;L^2(\Omega_0;{\mathbb R}^3))}+\|f_t\|_{L^2(0,T;H^1(\Omega_0;{\mathbb R}^3)')}\nonumber\\
&\qquad +\|  g(0,\cdot)\|_ {H^\frac{1}{2}(\Gamma_0;{\mathbb R}^3)}
 +\| \  g\|_ {L^2(0,T;H^\frac{1}{2}(\Gamma_0;{\mathbb R}^3))}+
\|  g_t\|_ {L^2(0,T;H^{-\frac{1}{2}}(\Gamma_0;{\mathbb R}^3))}
\nonumber\\
&\qquad+\|  f(0,\cdot)\|_ {L^2(\Omega_0;{\mathbb R}^3)} + \|  B \|_ {L^2(0,T;H^{-\frac{1}{2}}(\Gamma_0;{\mathbb R}))}
+ \|\  B_t\ \|_ {L^2(0,T;H^{-{\frac{1}{2}}}(\Gamma_0;{\mathbb R}))}\ )\ ,\nonumber
 \end{align*}
we then get from (\ref{wqh2bis}) and (\ref{inicom}) that
\begin{align*}
\|(W,Q)\|_{Y_T} \le & C\ (\ \|w_0\|_{H^2(\Omega_0;{\mathbb R}^3)}+\|f\|_{L^2(0,T;L^2(\Omega_0;{\mathbb R}^3))}+\|f_t\|_{L^2(0,T;H^1(\Omega_0;{\mathbb R}^3)')}\nonumber\\
& +\|  g(0,\cdot)\|_ {H^\frac{1}{2}(\Gamma_0;{\mathbb R}^3)}
 +\| \  g\|_ {L^2(0,T;H^\frac{1}{2}(\Gamma_0;{\mathbb R}^3))}+
\| \  g_t\|_ {L^2(0,T;H^{-\frac{1}{2}}(\Gamma_0;{\mathbb R}^3))}
\nonumber\\
&+\|  f(0,\cdot)\|_ {L^2(\Omega_0;{\mathbb R}^3)} + \|  B \|_ {L^2(0,T;H^{-\frac{1}{2}}(\Gamma_0;{\mathbb R}))}
+ \|\  B_t\ \|_ {L^2(0,T;H^{-{\frac{1}{2}}}(\Gamma_0;{\mathbb R}))}\ )\ .
\end{align*}

Thus, 
\begin{align}
\label{wqh2qua}
\|\zeta^3(\Phi) (w,q)(\Phi)\|_{Y_T} \le & C\ (\ \|w_0\|_{H^2(\Omega_0;{\mathbb R}^3)}+\|  g(0,\cdot)\|_ {H^\frac{1}{2}(\Gamma_0;{\mathbb R}^3)}+\|  f(0,\cdot)\|_ {L^2(\Omega_0;{\mathbb R}^3)}\nonumber\\
& +\|f\|_{L^2(0,T;L^2(\Omega_0;{\mathbb R}^3))}+\|f_t\|_{L^2(0,T;H^1(\Omega_0;{\mathbb R}^3)')}\nonumber\\
&  +\| \  g\|_ {L^2(0,T;H^\frac{1}{2}(\Gamma_0;{\mathbb R}^3))}+
\| \  g_t\|_ {L^2(0,T;H^{-\frac{1}{2}}(\Gamma_0;{\mathbb R}^3))}
\nonumber\\
&+ \|  B \|_ {L^2(0,T;H^{-\frac{1}{2}}(\Gamma_0;{\mathbb R}))}
+ \|\  B_t\ \|_ {L^2(0,T;H^{-{\frac{1}{2}}}(\Gamma_0;{\mathbb R}))}\ )\ .
\end{align}

We will now bootstrap, using the gain in regularity of $(w,q)$ in
the system (\ref{linear11}) to deduce
the appropriate regularity gain for $(W,Q)$.
The problem occurs in the region where $\zeta =0$, wherein
(\ref{wqh2qua}) introduces singular terms of the type 
$\frac{\nabla\zeta}{\zeta}$. To avoid this singular behavior, we 
localize the problem to an even smaller neighborhood of $x_0$ 
by introducing another smooth
cutoff function $\zeta'$ satisfying the same properties as $\zeta$ but whose support is strictly included in the support of 
$\zeta$. We then introduce $(W',Q')=\zeta'^3 (w,q)$ which  exactly satisfies
the same problem as $(W,Q)$ in (\ref{linear11}), with $\zeta$ replaced by 
$\zeta'$. Next, we  introduce a sequence $(v_n',q_n')$ defined exactly in the 
same
manner as the sequence $(v_n,q_n)$ with $\zeta'$ replacing $\zeta$. 
Our compatibility conditions are still satisfied , which allows us to use our basic energy estimates, but now
in $X_T$, to get a contractive sequence in $X_T$ (assuming, of course, that
the support of $\zeta'$ is chosen to be in a sufficiently small neighborhood 
of $x_0$), together with the estimate
\begin{equation}
\label{sumter}
\| (W',Q')\|_{X_T} \le \ 
\| ({ v'}_{1}, {q'}_{1})\|_{X_T} \ .
\end{equation}
We proceed to estimate the right-hand-side of (\ref{sumter}) in the same 
way as we did for $(v_1,q_1)$, but now in the space $X_T$ since our forcing 
terms have the additional regularity.  We arrive at the estimate
\begin{align}
\label{qui}
&\| (W',Q')\|_{X_T} \nonumber \\
 & \qquad \le  C\ (\ \|w_0\|_{H^2(\Omega_0;{\mathbb R}^3)}+\|f\|_{L^2(0,T;H^1(\Omega_0;{\mathbb R}^3))}+\|f_t\|_{L^2(0,T;H^1(\Omega_0;{\mathbb R}^3)')}\nonumber\\
& \qquad \qquad
+\|  g(0,\cdot)\|_ {H^\frac{1}{2}(\Gamma_0;{\mathbb R}^3)} +\| \ g\|_ {L^2(0,T;H^\frac{3}{2}(\Gamma_0;{\mathbb R}^3))}+
\| \  g_t\|_ {L^2(0,T;H^{-\frac{1}{2}}(\Gamma_0;{\mathbb R}^3))}
\nonumber\\
&\qquad\qquad +\|  f(0,\cdot)\|_ {L^2(\Omega_0;{\mathbb R}^3)} + \|  B \|_ {L^2(0,T;H^{\frac{1}{2}}(\Gamma_0;{\mathbb R}))}
+ \|\  B_t\ \|_ {L^2(0,T;H^{{\frac{1}{2}}}(\Gamma_0;{\mathbb R}))}\ )\ .
\end{align}

Since $\zeta'=1$ in a neighborhood of $x_0$, we may emply a finite-covering
argument and deduce, using  (\ref{qui}) and the interior regularity result, 
that the basic linear estimate (\ref{estimate1_omega}) is satisfied.
\end{proof}

\subsection{Proof of the a priori estimate for the basic linear problem} 
To finish the proof of Theorem \ref{thm1}, it remains only to include
the dependence on $\bar a$ in the a priori estimate (\ref{estimate1_omega}).

From (\ref{v}) and (\ref{w}), we need only substitute $w= \bar w -\nabla r$ 
into (\ref{estimate1_omega}).  From (\ref{elliptic0}), we have that
$$v=\nabla r = \nabla (- \triangle)^{-1} \bar a.$$ 
Thus, it is clear that 
\begin{align*}
\|v\|^2_{L^2(0,T;H^3(\Omega_0; {\mathbb R}^3)} & =
\|\bar a\|^2_{L^2(0,T;H^2(\Omega_0; {\mathbb R})} \\
\intertext{and that}
\|v_t\|^2_{L^2(0,T;H^1(\Omega_0; {\mathbb R}^3)} & =
\|\bar a_t\|^2_{L^2(0,T;L^2(\Omega_0; {\mathbb R})}.
\end{align*}
Also,
\begin{align*}
\| N\cdot \triangle_0 \int_0^\cdot v(s,\cdot)ds \|^2_
{C^0([0,T];H^{\frac{1}{2}}(\Gamma_0; {\mathbb R})} 
&\le C(T) \| \triangle_0 \nabla (- \triangle)^{-1}\bar  a\|^2_
{L^2(0,T;H^{ \frac{1}{2}}(\Gamma_0; {\mathbb R})} \\
&\le C(T) \| \nabla (- \triangle)^{-1}\bar  a\|^2_
{L^2(0,T;H^3(\Omega_0; {\mathbb R})} \\
&\le C(T) \| \bar a\|^2_
{L^2(0,T;H^2(\Omega_0; {\mathbb R})} \,,
\end{align*}
where we have made use of the trace theorem.  The term
$\| N\cdot \triangle_0 v \|^2_{L^2(0,T;H^{\frac{1}{2}}(\Omega_0; {\mathbb R})}$
has a similar estimate, so we have finished the proof.

\section{The fixed-point scheme for the nonlinear problem}

We will make use of the Tychonoff Fixed-Point Theorem  in our procedure (see, for example, \cite{Deimling}).  Recall that
this states that
for a reflexive separable Banach space $X$, and $C\subset X$ a closed,
convex, bounded subset,  if $F:C \rightarrow C$ is weakly
sequentially continuous into $X$, then $F$ has at least one fixed-point.

Our analysis will be performed using the function space $X_T$, but this time
on the domain $\Omega_0$.
\begin{definition}\label{XT_omega0}
Let
\begin{align*}
X_T= \{& (u,p)\in V^3 (T)\times V^2 (T)\ |\ \\
       &\ \ (\nabla u_t, p_t)\in 
{L^2(0,T; H^{-\frac{1}{2}}(\Gamma_0;{\mathbb R}^9))}
\times {L^2(0,T;H^{-\frac{1}{2}}(\Gamma_0;{\mathbb R}))}\}\ ,
\end{align*}
 endowed with its natural Hilbert norm
\begin{align*}
\|(u,p)\|^2_{X_T}=&\|  u \|^2_ {L^2(0,T;H^3(\Omega_0;{\mathbb R}^3))}+ 
\|u_t\|^2_ {L^2(0,T;H^1(\Omega_0;{\mathbb R}^3))}\nonumber \\
& + \|  p \|^2_ {L^2(0,T;H^2(\Omega_0;{\mathbb R}))}+ 
\| p_t\|^2_ {L^2(0,T;L^2(\Omega_0;{\mathbb R}))}\nonumber \\
& +\| \nabla u_t \|^2_ {L^2(0,T;H^{-\frac{1}{2}}(\Gamma_0;{\mathbb R}^9))}+ 
\| p_t\|^2_ {L^2(0,T;H^{-\frac{1}{2}}(\Gamma_0;{\mathbb R}))}\,.
\end{align*}
\end{definition}

We define a mapping $\Theta_T$ from $X_T$ into itself, which to a given
pair $(v,q) \in  X_T $ associates the pair
$(\tilde{v},\tilde{q})\in  X_T$, a
solution of a linear problem defined as follows. Let
$\eta\in V^3 (T)\cap C([0,T]; {H^3(\Omega_0; {\mathbb R}^3)})$ denote
the Lagrangian flow map of the velocity $v$, defined by 
\begin{equation}\label{eta}
\eta_t =v \ \ \text{in} \ \ (0,T)\times \Omega_0 \ \ \text{ with } \ \
\eta = \text{Id}     \ \ \text{on} \ \ \Omega_0\times \{ t=0\} \,,
\end{equation}
let $a$ be given by (\ref{a}), and 
define $(\tilde{v},\tilde{q})\in\{(u,p)\in X_T|\ u(0)=u_0\}$ to be the 
solution of 
\begin{subequations}
  \label{nlfp}
\begin{align}
{\tilde{v}}^i_t - \nu (a^j_l a^k_l {\tilde{v}}^i,_k),_j + a^k_i \tilde{q},_k&= F^i 
\ \ \text{in} \ \ (0,T)\times \Omega_0 \,, 
         \label{nlfp.c}\\
   a^k_i {\tilde{v}}^i,_k &= 0     \ \ \text{in} \ \ (0,T)\times \Omega_0 \,, 
         \label{nlfp.d}\\
\Pi_0\left(\Pi_{\eta} D_{\eta} (\tilde{v})\right)\cdot a^T\ N &=0
\ \ \text{on} \ \ (0,T)\times \Gamma_0 \,, 
         \label{nlfp.e}\\
N\cdot (\ S_{\eta}(\tilde{v},\tilde{q})\cdot a^T\ N\ ) -
\sigma\ N\cdot \triangle_{g(t)}(\int_0^t {\tilde{v}(t',\cdot)\ dt'})&= 
\sigma\ N\cdot \triangle_{g(t)}(\text{Id})\nonumber\\
& \qquad \ \ \text{on} \ \ (0,T)\times \Gamma_0 \,, 
         \label{nlfp.g}\\
   \tilde{v} &= u_0    
 \ \ \text{on} \ \ \Omega_0\times \{ t=0\} \,, 
         \label{nlfp.h}
\end{align}
\end{subequations}
where $N$ denotes the outward-pointing unit normal to $\Gamma_0$, $\Pi_0$ and
$\Pi_{\eta}$ are, respectively, the projectors onto the tangent planes 
to the surfaces $\Gamma_0$ and $\eta (\Gamma_0)$, 
\begin{subequations}
\begin{align}
D_{\eta}(\tilde{v})^i_l&= 
({\tilde{v}}^i,_k a^k_l + {\tilde{v}}^l,_k a^k_i) 
 \ \ \text{in} \ \ (0,T)\times \Omega_0 \,, 
         \\
S_{\eta}(\tilde{v},\tilde{q})&= \nu D_{\eta}(\tilde{v})- \tilde{q}\ \text{Id}
 \ \ \text{in} \ \ (0,T)\times \Omega_0 \,, 
 \end{align} 
    \end{subequations}
and $\triangle_{g(t)}$ denotes the Laplacian with respect to the induced
metric $g(t)$ defined by the equations (\ref{laplacian}) and (\ref{metric}).
Clearly, $D_{\eta}(\tilde{v})_i^l= D_{\eta}(\tilde{v})_l^i$.
$D_\eta$ is the Lagrangian version of the deformation tensor, and $S_\eta$
is Lagrangian version of the stress tensor.

In the following, we will prove that (\ref{nlfp}) has a unique solution for
a short time $T$ by using an iteration scheme founded upon the basic linear 
problem (\ref{linear}); this will establish that
the mapping $\Theta_T$ is well-defined.  Using the  
Tychonoff fixed point theorem, we will then establish that the map 
$\Theta_T:(v,q) \mapsto (\tilde v,\tilde q)$ has a fixed point,
and that this fixed point is indeed the unique
solution of the system (\ref{nsl}).  This will prove our main theorem. 

In order to do so, let us consider a given element 
$(v,q)\in \{ X_T \ | \ u(0)=u_0\} $.   To apply our fixed point procedure,
we shall define the following closed
convex subset of $X_T$.

\begin{definition}[The convex subset]\label{Cm}
Let 
$$C_T=\{(u,p)\in X_T|\ u(0)=u_0\ \text{and}\ \|(u,p)\|_{X_T}\le M(T)\}$$
with
\begin{align*}
 M(T)= 2\ C_m \Bigl[&|f(0,\cdot)\|_{L^2({\mathbb R}^3;{\mathbb R}^3)}
 + \| f \|_ {L^2(0,T;H^1({\mathbb R}^3;{\mathbb R}^3))}\ (1+c_6 c_3\ \|u_0\|_{H^1(\Omega_0;{\mathbb R}^3)})\\
& + \|f_t\|_ {L^2(0,T;H^1({\mathbb R}^3;{\mathbb R}^3)')}
+\|u_0\|_{H^2(\Omega_0;{\mathbb R}^3)}+\|\triangle_0 \operatorname{Id}
\cdot N\|_{H^\frac{1}{2}(\Gamma_0;{\mathbb R})}\Bigr]\,,
\end{align*}
where we fix $T_m>0$ with $T \le T_m$ and let
$C_m$ denotes the function $C(T)$ in the estimate (\ref{estimate1_omega}) with
$T_m$ replacing $T$. (Note that $C_m$ is an increasing function of time, so that
we can replace $C$ by $C_m$ in the estimate (\ref{estimate1_omega}).) In this
definition, $c_3$ and $c_6$ are the Sobolev constants satisfying for any
$v\in H^1(\Omega_0;{\mathbb R}^3)$ the inequalities
\begin{equation}
\label{C_3}
\|v\|_{L^3(\Omega_0;{\mathbb R}^3)}\le c_3\ \|v\|_{H^1(\Omega_0;{\mathbb R}^3)}\ ,
\end{equation}
\begin{equation}
\label{C_6}
\|v\|_{L^6(\Omega_0;{\mathbb R}^3)}\le c_6\ \|v\|_{H^1(\Omega_0;{\mathbb R}^3)}\ .
\end{equation}

\end{definition}

We note that of course, $C_T$ is well defined since $V^3(T)\times V^2(T)$ is 
continuously embedded into
$C([0,T];H^2(\Omega_0;{\mathbb R}^3)\times H^1(\Omega_0;{\mathbb R}))$ 
(see, for example, \cite{Evans1998}). We will solve (\ref{nlfp}) for a 
short time by an iteration procedure and obtain a unique solution 
$(\tilde v,\tilde q)$ belonging to $C_T$. Note that this short time 
restriction is necessary not only in order to get a solution in $C_T$, but 
also to ensure the very existence of a solution of
(\ref{nlfp}) in $X_T$.  

\begin{remark}
It is somewhat surprising, given that (\ref{nlfp}) is a linear problem, 
that we obtain only a short-time solution; nevertheless, this restriction
arises from the iteration scheme we must employ.  We use
this type of procedure 
since, due to the (non-natural) boundary conditions (\ref{nlfp.e}) and 
(\ref{nlfp.g}), 
we cannot write the variational form of (\ref{nlfp.e}).
\end{remark}

In the following, $C$ denotes a generic constant which may depend on $\Omega_0$ 
and the coefficients of (\ref{nsl}).
The following lemma concerns the $L^\infty$ control of $a$, 
and asserts that $a$ remains near the identity in an appropriate norm, and that 
the
unit normal $N$ to $\Gamma_0$ at $x\in\Gamma_0$ and the unit normal $n$ to
$\Gamma(t)=\eta(t,\Gamma_0)$ at $\eta(t,x)$ are not orthogonal for at
least a short time.

\begin{lemma}
\label{control}
There exists $K>0$, $T_0>0$ such that if $0<T\le T_0$, then, for any $(v,q)\in C_T$,
we have that
\begin{align}
\label{control1}
 \| a^T-\text{Id} \|_ {L^\infty(0,T;C^0(\overline{\Omega_0};{\mathbb R}^9))}
&\le K\ \sqrt{T}\ ,
\\
\label{control2}
\frac{a^T\ N}{|a^T\ N|} \cdot N 
&\ge \frac{1}{2}\ ,\ \ \operatorname{in}\  [0,T]\times\Gamma_0\ ,
\\
\label{control3}
\operatorname{det}\ a 
&\ge \frac{1}{2} ,\ \ \operatorname{in}\  [0,T]\times\Omega_0\ ,
\end{align}
\end{lemma}

\begin{proof}
Using (\ref{eta}), we have that
$$\|\eta (t,\cdot)-\text{Id}\|_{H^3(\Omega_0;{\mathbb R}^3)} =  
\|\int_0^t v(t',\cdot)\ dt'\|_{H^3(\Omega_0;{\mathbb R}^3)}\ ,
$$ 
which shows by Jensen's inequality for the Bochner integral that
$$\|\eta (t, \cdot)-\text{Id}\|_{H^3(\Omega_0;{\mathbb R}^3)} \le  
\int_0^t \|v(t', \cdot)\|_{H^3(\Omega_0;{\mathbb R}^3)}\ dt'\ ,
$$ 
and thus from the Cauchy-Schwarz inequality,
\begin{equation*}
\|\eta-\text{Id} \|_ {L^\infty(0,T;H^3(\Omega_0;{\mathbb R}^3))} \le  
C\sqrt{T} \ \|v\|_{L^2(0,T;H^3(\Omega_0;{\mathbb R}^3))}\le C\sqrt{T}\ \|(v,q)\|_{X_T},
\end{equation*}
and since $(v,q)\in C_T$,
\begin{equation}
\label{etah3bis}
\|\eta-\text{Id} \|_ {L^\infty(0,T;H^3(\Omega_0;{\mathbb R}^3))} \le  
C\sqrt{T} \ \ M(T)\le C\ \sqrt{T}\ M(T_0)\ .
\end{equation}
From the continuity of the embedding of 
$H^3({\Omega_0})$ into $C^1(\overline{\Omega_0})$, we infer that
$$\|\nabla\eta-\text{Id} \|_ {L^\infty(0,T;C^0(\overline{\Omega_0};{\mathbb R}^9))} \le  
 C\ \sqrt{T} \ M(T_0) \ ,
$$
which implies that there exists $T_0>0$ such that for 
$T\le T_0$, $\nabla\eta$ is invertible and is in a neighborhood of the identity,
and thus from the $C^{\infty}$ regularity of the mapping $M\mapsto M^{-1}$ 
from a neighborhood of $\text{Id}$ in
$H^2(\Omega_0;{\mathbb R}^9)$ into $H^2(\Omega_0;{\mathbb R}^9)$, we
have that
$$\|a^T-\text{Id} \|_ {L^\infty(0,T;C^0(\overline{\Omega_0};{\mathbb R}^9))} 
\le  C\ \sqrt{T}\ \ M(T_0) 
$$ 
which proves (\ref{control1}). Since $a$ is close to Id in 
$L^\infty(0,T;C^0(\overline{\Omega_0};{\mathbb R}^9))$, we then easily
deduce (\ref{control2}) and (\ref{control3}) for $T_0 >0$ sufficiently
small; this is due to the fact that when $a=\text{Id}$, we have that
both $a^T \ N/| a^T \ N|\cdot N =1$ and $\operatorname{det}  a=1$.
\end{proof}

Henceforth, we take $T < T_0$, and we allow our generic constant $C$ to
depend on $T_0$.  With this convention, (\ref{etah3bis}) reads 
\begin{equation}
\label{etah3}
\|\eta-\text{Id} \|_ {L^\infty(0,T;H^3(\Omega_0;{\mathbb R}^3))} \le  
C\sqrt{T}  \ \ \forall  (v,q) \in C_T.
\end{equation} 
We next state the following lemma about the $L^\infty$-in-time control 
of elements of $C_T$.
\begin{lemma}
\label{linftybis}
There exists a constant $C$ such that for any $(w,r)\in C_T$ we have
\begin{equation}
\label{wrh2h1}
 \|w\|_{L^\infty(0,T;H^2(\Omega_0;{\mathbb R}^3))} \le C \ .
\end{equation}
\end{lemma}
\begin{proof}
From the standard interpolation inequality (see \cite{Evans1998} for instance)
\begin{align*}
 \|w\|_{L^\infty(0,T;H^2(\Omega_0;{\mathbb R}^3))}\ \le 
 &\|u_0\|_{H^2(\Omega_0;{\mathbb R}^3)}\nonumber\\
&+ C\ \left( \|w_t\|_{L^2 (0,T;H^1(\Omega_0;{\mathbb R}^3))} + \|w\|_{L^2 (0,T;H^3(\Omega_0;{\mathbb R}^3))}\right)\ ,
\end{align*}
we infer that
$$ \|w\|_{L^\infty(0,T;H^2(\Omega_0;{\mathbb R}^3))}\ \le 
C[ \|u_0\|_{H^2(\Omega_0; \mathbb{R} ^3)} + M(T)]\ ,$$
which finishes the proof of the lemma, using the fact that 
$M(T)\le M(T_0)$ (and again allowing $C$ to depend on $T_0$).
\end{proof}

We will also need the following 
\begin{lemma}
\label{linfty}
There exists a constant $C$ such that for any $(v,q)\in C_T$,
\begin{equation}
\label{timeah1}
 \|a_t\|_{L^\infty(0,T;H^1(\Omega_0;{\mathbb R}^3))} \le C \ .
\end{equation}
\end{lemma}
\begin{proof}
We first need a few estimates on $a$ and its time derivative. From (\ref{etah3}), we infer that
$$\|\nabla\eta-\text{Id} \|_ {L^\infty(0,T;H^2(\Omega_0;{\mathbb R}^9))}   
\le C\ \sqrt{T}  \ ,
$$
and by the Sobolev embeddings
$$\|\nabla\eta-\text{Id} \|_ {L^\infty(0,T;C^0(\overline{\Omega_0};{\mathbb R}^9))}   
\le C\ \sqrt{T}  \ .
$$
Thus $\nabla\eta$ stays in a neighborhood of Id in ${\mathbb R}^9$ on $[0,T]\times \overline{\Omega_0}$. From the $C^{\infty}$ regularity of the mapping $M\mapsto M^{-1}$ from a neighborhood of $\text{Id}$ in
$H^2(\Omega_0;{\mathbb R}^9)$ into $H^2(\Omega_0;{\mathbb R}^9)$, we infer from (\ref{etah3}) that
\begin{equation}
\label{ah2}
\|a-\text{Id} \|_ {L^\infty(0,T;H^2({\Omega_0};{\mathbb R}^9))}   
\le C\ \sqrt{T}  \ .
\end{equation}
Now, since $a_t=-a\cdot\nabla\eta_t\cdot a=-a\cdot \nabla v\cdot a$, 
we see from the Cauchy-Schwarz inequality that
\begin{align*}
\|a_t (t,\cdot)\|_{H^1({\Omega_0};{\mathbb R}^9)}\le &\|a(t,\cdot)\|^2_{L^\infty({\Omega_0};{\mathbb R}^9)}\ \|\nabla v(t,\cdot)\|_{H^1({\Omega_0};{\mathbb R}^9)}\\ &+
\|a(t,\cdot)\|^2_{W^{1,4}({\Omega_0};{\mathbb R}^9)}\ \|\nabla v(t,\cdot)\|_{L^4({\Omega_0};{\mathbb R}^9)}\ .
\end{align*}
Thus by the Sobolev embedding theorem,
$$\|a_t (t,\cdot)\|_{H^1({\Omega_0};{\mathbb R}^9)}\le C\ \ \|a(t,\cdot)\|^2_{H^2({\Omega_0};{\mathbb R}^9)}\ \|v(t,\cdot)\|_{H^2({\Omega_0};{\mathbb R}^9)}\ ,$$
and consequently from (\ref{ah2}),
\begin{equation}
\label{y1}
\|a_t(t,\cdot)\|_{H^1({\Omega_0};{\mathbb R}^9)}\le C\ \|v(t,\cdot)\|_{H^2({\Omega_0};{\mathbb R}^9)}\ .
\end{equation}

We then conclude from (\ref{y1}) and (\ref{wrh2h1}) that
\begin{equation}
\label{timeah1bis}
 \|a_t\|_{L^\infty(0,T;H^1(\Omega_0;{\mathbb R}^3))}\ \le C \,.
\end{equation}
\end{proof}
We will also need to estimate the following boundary forcing terms:
\begin{lemma}
We have for any $(v,q)\in C_T$,
\begin{align}
\|   N\cdot g^{\alpha \beta}\frac{\partial^2 \text{Id} }
{\partial y^\alpha \partial y^\beta}\|
_{L^2(0,T;H^\frac{3}{2}(\Gamma_0;{\mathbb R}))}
+& \|\partial_t(N\cdot g^{\alpha \beta}\frac{\partial^2 \text{Id} }
{\partial y^\alpha \partial y^\beta})\ 
\|_ {L^2(0,T;H^{-\frac{1}{2}}(\Gamma_0;{\mathbb R}))}\nonumber\\
& \le \epsilon(T)\  M(T) \,, \label{again}
\end{align}
where $\epsilon (T)\rightarrow 0$ as $T\rightarrow 0$.
\end{lemma}
\begin{proof}
For the first term of this inequality, we have that
\begin{align*}
\|   N\cdot g^{\alpha \beta}\frac{\partial^2 \text{Id} }
{\partial y^\alpha \partial y^\beta}\|_ {L^2(0,T;H^\frac{3}{2}(\Gamma_0;{\mathbb R}))}\le & \|   N\cdot (g^{\alpha \beta}-g^{\alpha,\beta}(0,\cdot))\frac{\partial^2 \text{Id} }
{\partial y^\alpha \partial y^\beta}\|_ {L^2(0,T;H^\frac{3}{2}(\Gamma_0;{\mathbb R}))}\nonumber\\
& +\|   N\cdot g^{\alpha \beta}(0,\cdot)\frac{\partial^2 \text{Id} }
{\partial y^\alpha \partial y^\beta}\|_ {L^2(0,T;H^\frac{3}{2}(\Gamma_0;{\mathbb R}))}\
 .
\end{align*}
Since $\eta(0,\cdot)=\text{Id}$, we infer from (\ref{etah3}) that
\begin{equation}
\label{again1}
\|   N\cdot g^{\alpha \beta}\frac{\partial^2 \text{Id} }
{\partial y^\alpha \partial y^\beta}\|_ {L^2(0,T;H^\frac{3}{2}(\Gamma_0;{\mathbb R}))}\le C\ \sqrt{T}\ .
\end{equation}
For the second term on the left-hand-side of (\ref{again}), 
the smoothness of $\Gamma_0$ implies that
\begin{equation*}
\|\partial_t(N\cdot g^{\alpha \beta}\frac{\partial^2 \text{Id} }
{\partial y^\alpha \partial y^\beta})\ \|_ {L^2(0,T;H^{-\frac{1}{2}}(\Gamma_0;{\mathbb R}))}\le C\ \|\partial_t g^{\alpha \beta}\|_{L^2(0,T;H^{-\frac{1}{2}}(\Gamma_0;{\mathbb R}))}\ ,
\end{equation*}
and from (\ref{etah3}) and Lemma \ref{againlemma} that we prove next,
\begin{equation*}
\|\partial_t(N\cdot g^{\alpha \beta}\frac{\partial^2 \text{Id} }
{\partial y^\alpha \partial y^\beta})\ \|_ {L^2(0,T;H^{-\frac{1}{2}}(\Gamma_0;{\mathbb R}))}\le C\ \|\nabla_0 v\|_{L^2(0,T;H^{-\frac{1}{2}}(\Gamma_0;{\mathbb R}^9))}\ .
\end{equation*}
Since there is only a surface gradient involved in the right-hand-side of 
this inequality, from the same type of duality argument as we used for the 
proof of (\ref{dualagain}), we find that
\begin{equation*}
\|\partial_t(N\cdot g^{\alpha \beta}\frac{\partial^2 \text{Id} }
{\partial y^\alpha \partial y^\beta})\ \|_ {L^2(0,T;H^{-\frac{1}{2}}(\Gamma_0;{\mathbb R}))}\le C\ \|\nabla v\|_{L^2(0,T;L^2(\Omega_0;{\mathbb R}^9))}\ ,
\end{equation*}
and from (\ref{wrh2h1}), we then have in turn that
\begin{equation}
\label{again2}
\|\partial_t(N\cdot g^{\alpha \beta}\frac{\partial^2 \text{Id} }
{\partial y^\alpha \partial y^\beta})\ \|_ {L^2(0,T;H^{-\frac{1}{2}}(\Gamma_0;{\mathbb R}))}\le C\ \sqrt{T}\ .
\end{equation}
Form (\ref{again1}) and (\ref{again2}), we then deduce the desired inequality (\ref{again}).

\end{proof}

We next have the following simple result concerning the 
$H^{-\frac{1}{2}}(\Gamma_0;{\mathbb R})$ norm of a product.
\begin{lemma}
\label{againlemma}
For any $\alpha\in H^{-\frac{1}{2}}(\Gamma_0;{\mathbb R})$ and $\beta$ which is the
trace of an element of $\beta\in W^{1,4}(\Gamma_0;{\mathbb R})$, we have that
\begin{equation}
\label{dual}
\|\alpha\beta\|_{H^{-\frac{1}{2}}(\Gamma_0;{\mathbb R})} \le C\ \|\alpha\|_{H^{-\frac{1}{2}}(\Gamma_0;{\mathbb R})}\ \|\beta\|_{W^{1,4}(\Omega_0;{\mathbb R})}\ .
\end{equation}
\end{lemma}

\begin{proof}
Let $\phi\in {H^{\frac{1}{2}}(\Omega_0;{\mathbb R})}$ and let $R$ be a linear and continuous mapping from ${H^{\frac{1}{2}}(\Omega_0;{\mathbb R})}$ into ${H^1 (\Omega_0;{\mathbb R})}$ such that $\phi$ is the trace of $R(\phi)$. 

 We then have 
$$ \langle \alpha\beta, \phi\ \rangle_{\Gamma_0}=\langle\ \alpha,\beta\phi 
\rangle_{\Gamma_0}\ ,$$
where $\langle \cdot,\cdot \rangle_{\Gamma_0}$ denotes the duality 
product between ${H^{-\frac{1}{2}}(\Gamma_0;{\mathbb R})}$ and 
${H^{\frac{1}{2}}(\Gamma_0;{\mathbb R})}$. By the trace theorem,
$$ \langle \alpha\beta, \phi \rangle_{\Gamma_0}\le C
\|\alpha\|_{H^{-\frac{1}{2}}(\Gamma_0;{\mathbb R})} 
\ \|\beta R(\phi)\|_{H^1 (\Omega_0;{\mathbb R})}\ .$$
Since by the Cauchy-Schwarz inequality 
$$\|\beta R(\phi)\|_{H^1 (\Omega_0;{\mathbb R})}\le \|\beta\|_{L^\infty (\Omega_0;{\mathbb R})}\ \|\phi\|_{H^1 (\Omega_0;{\mathbb R})}+ 
\|\beta\|_{W^{1,4} (\Omega_0;{\mathbb R})}\ \|R(\phi)\|_{L^4 (\Omega_0;{\mathbb R})}\ ,$$
we then have by the Sobolev embedding theorem,
$$\|\beta R(\phi)\|_{H^1 (\Omega_0;{\mathbb R})}\le 
C\ \|\beta\|_{W^{1,4} (\Omega_0;{\mathbb R})}\ \|R(\phi)\|_{H^1 (\Omega_0;{\mathbb R})}\ ,$$
and thus
$$ \langle \alpha\beta, \phi \rangle_{\Gamma_0}\le 
\|\alpha\|_{H^{-\frac{1}{2}}(\Gamma_0;{\mathbb R})} 
\ C\ \|\beta\|_{W^{1,4} (\Omega_0;{\mathbb R})}
\ \|\phi\|_{H^{\frac{1}{2}} (\Omega_0;{\mathbb R})}\ ,$$
which proves the lemma.
\end{proof}
For the uniqueness of the solution to (\ref{nsl}), we will need the following
\begin{lemma}
For any $(v,q)\in C_T$ and $(\tilde v,\tilde q)\in C_T$,
\begin{equation}
\label{deltaah2}
\|a-\tilde{a}\|_{L^\infty (0,T; H^2 (\Omega_0;{\mathbb R}^9)}\le C\ \sqrt{T}\ ,
\end{equation}
where our notation means that $a$ and $\tilde a$ are formed from the
Lagrangian flow maps of $v$ and $\tilde v$, respectively.
\end{lemma}
\begin{proof}
Note that from (\ref{control3}), $a$ and $\tilde{a}$ are invertible, 
so that we may use the decomposition
\begin{equation*}
a(t,\cdot)-\tilde{a}(t,\cdot)=\tilde{a}(t,\cdot)\cdot 
[\tilde{a}^{-1}(t,\cdot)\cdot a(t,\cdot)-\text{Id}]\ ,
\end{equation*}
and thus, since $H^2 (\Omega_0;{\mathbb R}^9)$ is a Banach algebra, 
(\ref{ah2}) enables us to assert that
\begin{equation*}
\|a(t,\cdot)-\tilde{a}(t,\cdot)\|_{H^2 (\Omega_0;{\mathbb R}^9)}
\le C\ \|\tilde{a}^{-1}(t,\cdot)\cdot 
a(t,\cdot)-\text{Id}\|_{H^2 (\Omega_0;{\mathbb R}^9)}\ .
\end{equation*}
Hence,
\begin{align}
\label{delta1}
\|a(t,\cdot)-\tilde{a}(t,\cdot)\|_{H^2 (\Omega_0;{\mathbb R}^9)}
\le & C\ \|\tilde{a}^{-1}(t,\cdot)\|_{H^2 (\Omega_0;{\mathbb R}^9}\ \|
 a(t,\cdot)-\text{Id}\|_{H^2 (\Omega_0;{\mathbb R}^9)}\nonumber \\
& + C\ \|\tilde{a}^{-1}(t,\cdot)-\text{Id}\|_{H^2 (\Omega_0;{\mathbb R}^9)}\ .\end{align}
Since $\tilde{a}^{-1}(t,\cdot)=\nabla\tilde{\eta}$, we infer from (\ref{ah2}) that
\begin{equation*}
\|\tilde{a}^{-1}(t,\cdot)-\text{Id}\|_{H^2 (\Omega_0;{\mathbb R}^9)}\le C\ \sqrt{T}\ ,
\end{equation*}
and consequently, with (\ref{delta1}), we obtain the desired inequality.
\end{proof}

With regards to the forcing $F=f\circ\eta$, we have the following
\begin{lemma}
For any $(v,q)\in C_T$,
\begin{align}
&\|F\|_{L^2 (0,T; H^1 (\Omega_0;{\mathbb R}^3))}+ 
\|F_t\|_{L^2 (0,T; H^1 (\Omega_0;{\mathbb R}^3)')}\nonumber\\
& \le 
 (1+ \epsilon(T)) 
\left(\ (1+c_6 c_3 \|u_0\|_{H^1 (\Omega_0;{\mathbb R}^3)})\ \|f\|_{L^2 (0,T; H^1 ({\mathbb R}^3;{\mathbb R}^3))}+
\|f_t\|_{L^2 (0,T; H^1 ({\mathbb R}^3;{\mathbb R}^3)')}\right)\,,
\label{forcingcmbis}
\end{align}
where $\epsilon(T)\rightarrow 0$ as $T\rightarrow 0$ and $c_3$ and $c_6$ 
are the Sobolev constants defined in (\ref{C_3}) and (\ref{C_6}).
\end{lemma}

\begin{proof}
From the chain rule, 
$$\nabla F =[\nabla f \circ \eta] \cdot \nabla\eta.$$ 
The estimate (\ref{etah3}) provides us with  $L^\infty$ spacetime control 
of both $\nabla\eta$ and $ \nabla \eta ^{-1}$ for
$T\le T_0$, so that we obtain 
\begin{equation}
\label{forcingcm1}
\|F\|_{L^2 (0,T; H^1 (\Omega_0;{\mathbb R}^3))}\le  \ 
\|f\|_{L^2 (0,T; H^1 ({\mathbb R}^3;{\mathbb R}^3))}\ (1+\epsilon(T))\ ,
\end{equation}
where $\epsilon(T)\rightarrow 0$ as $T\rightarrow 0$.

The chain rule also shows that
\begin{equation*}
F_t =f_t \circ \eta +[\nabla f \circ\eta] \ v \ .
\end{equation*}

By duality, 
\begin{equation}\label{SsS}
\|F_t (t,\cdot)\|_{ H^1 (\Omega_0;{\mathbb R}^3)'}\le 
\|[f_t \circ  \eta](t,\cdot) \|_{ H^1 (\Omega_0;{\mathbb R}^3)'}
+c_6\|[\nabla f \circ \eta](t,\cdot)\ 
v(t,\cdot)\|_{ L^{\frac{6}{5}} (\Omega_0;{\mathbb R}^3)}\ .
\end{equation}
From H\"older's inequality and the Sobolev embedding theorem, 
\begin{align*}
\|F_t (t,\cdot)\|_{ H^1 (\Omega_0;{\mathbb R}^3)'}\le & 
\|[f_t \circ  \eta](t,\cdot) \|_{ H^1 (\Omega_0;{\mathbb R}^3)'}\\
& +c_6 c_3\ \|[\nabla f \circ \eta](t,\cdot)\ 
\|_{ L^2 (\Omega_0;{\mathbb R}^3)}\ \|v(t,\cdot)\ 
\|_{ H^1 (\Omega_0;{\mathbb R}^3)}.
\end{align*}
Next, note that
$$
\| v (t,\cdot) \|_{H^1(\Omega_0; \mathbb{R}^3)}
\le \| u_0 \|_{H^1(\Omega_0; \mathbb{R}^3)}+
\int_0^t \| v_t (t',\cdot)\|
_{H^1(\Omega_0; \mathbb{R}^3)} dt'\ ,
$$
which by means of the Cauchy-Schwarz inequality and the fact that $(v,q)\in C_T$ leads us to
$$
\| v (t,\cdot) \|_{L^\infty (0,T; H^1(\Omega_0; \mathbb{R}^3))}
\le \| u_0 \|_{H^1(\Omega_0; \mathbb{R}^3)}+
\sqrt{T}\ M(T) .
$$
Again using (\ref{etah3}), we arrive at the estimate
\begin{align*}
 \|F_t\|_{L^2 (0,T; H^1 (\Omega_0;{\mathbb R}^3)')}
 \le & \|f_t\|_{L^2 (0,T; H^1 ({\mathbb R}^3;{\mathbb R}^3)')}\ (1+\epsilon(T))\\
& + c_6 c_3\ \|f\|_{L^2 (0,T; H^1 ({\mathbb R}^3;{\mathbb R}^3))}\ \|u_0\|_{H^1(\Omega_0; \mathbb{R}^3)}\ (1+\epsilon (T))\ ,
\end{align*}
where $\epsilon(T)\rightarrow 0$ as $T\rightarrow 0$.
Together with (\ref{forcingcm1}), the lemma is proved.
\end{proof}

For the proof of the uniqueness result, which relies on the Lipschitz condition 
(\ref{Lip}), we will need the following
\begin{lemma}
Suppose that $f$ satisfies (\ref{Lip}) with $T$ replacing $\bar T$.
Then, there exists $C>0$ such that for any $(v,q)$ and $(\tilde v,\tilde q)$ 
both in $C_T$, 
\begin{align}
\| f\circ\eta- f\circ\tilde\eta \|_ {L^2(0,T;H^1(\Omega_0;{\mathbb R}^3))}+& 
\|(f\circ\eta- f\circ\tilde\eta )_t \|_{L^2(0,T;H^1(\Omega_0;{\mathbb R}^3)')}
\nonumber\\
&\le C\ T^{\frac{1}{8}} \|(v-\tilde {v},q-\tilde {q})\|_{X_{T}}\ . 
\label{fcirceta}
\end{align}
\end{lemma}
\begin{proof}
We first notice that
$$\|\eta (t,\cdot)-\tilde\eta (t,\cdot)\|_{H^3(\Omega_0;{\mathbb R}^3)} =  
\left\|\int_0^t (v-\tilde v)(t',\cdot)\ dt'
\right\|_{H^3(\Omega_0;{\mathbb R}^3)}\ ,
$$ 
which shows by Jensen's inequality that
$$\|\eta (t, \cdot)-\tilde\eta (t,\cdot)\|_{H^3(\Omega_0;{\mathbb R}^3)} \le  
\int_0^t \|(v-\tilde v)(t', \cdot)\|_{H^3(\Omega_0;{\mathbb R}^3)}\ dt'\ ,
$$ 
and thus from the Cauchy-Schwarz inequality,
\begin{equation}
\label{fcirceta1}
\|\eta-\tilde\eta \|_ {L^\infty(0,T;H^3(\Omega_0;{\mathbb R}^3))}   
\le C\sqrt{T}\ \|(v-\tilde v,q-\tilde q)\|_{X_T}.
\end{equation}
On the other hand, since $v(0,\cdot)=\tilde v(0,\cdot)$, a standard
interpolation argument (see \cite{Evans1998}) yields
$$\|v-\tilde v\|_{L^\infty(0,T;H^2(\Omega_0;{\mathbb R}^3))}\le C\ (\|v-\tilde v\|_{L^2(0,T;H^3(\Omega_0;{\mathbb R}^3))}+\|v_t-\tilde v_t\|_{L^2(0,T;H^1(\Omega_0;{\mathbb R}^3))})\ ,$$
and thus 
\begin{equation}
\label{fcirceta3}
\|v-\tilde v\|_{L^\infty(0,T;H^2(\Omega_0;{\mathbb R}^3))}\le C\ \|(v-\tilde v,q-\tilde q)\|_{X_T}\ .
\end{equation}
By the Cauchy-Schwarz inequality we also have that
$$\|v-\tilde v\|_{L^\infty(0,T;H^1(\Omega_0;{\mathbb R}^3))}
\le C\ \sqrt{T}\ \|v_t-\tilde v_t\|_{L^2(0,T;H^1(\Omega_0;{\mathbb R}^3))}\ ,$$
so that
 \begin{equation}\label{SsS2}
\|v-\tilde v\|_{L^\infty(0,T;H^1(\Omega_0;{\mathbb R}^3))}\le C\ \sqrt{T}\ 
\|(v-\tilde v,q-\tilde q)\|_{X_T}\,.
\end{equation} 
We may interpolate between  (\ref{fcirceta3}) and (\ref{SsS2}) using
the Gagliardo-Nirenberg inequality to find that
\begin{equation}
\label{fcirceta4}
\|v-\tilde v\|_{L^\infty(0,T;L^\infty(\Omega_0;{\mathbb R}^3))}\le 
C\ T^{\frac{1}{8}}\ \|(v-\tilde v,q-\tilde q)\|_{X_T}\ .
\end{equation}
Since
\begin{equation}\label{SsS3}
\nabla f(\eta)\cdot \nabla\eta-\nabla f(\tilde\eta)\cdot \nabla\tilde\eta
=\nabla f(\eta)\cdot ( \nabla\eta-\nabla\tilde\eta) +
(\nabla f(\eta)-\nabla f(\tilde\eta))\cdot \nabla\tilde\eta\ ,
\end{equation} 
using (\ref{fcirceta1}) and (\ref{etah3}) together with (\ref{Lip}) for
the second term on the right-hand-side of (\ref{SsS3}), we obtain
\begin{align*}
&\|\nabla f(\eta)\cdot \nabla\eta-\nabla f(\tilde\eta)\cdot 
\nabla\tilde\eta\|_{L^2(0,T;L^2(\Omega_0;{\mathbb R}^9))}\nonumber\\
&\qquad\qquad\qquad\qquad\qquad\le 
C\ \sqrt{T}\ (K+\|\nabla f\|_{L^2(0,T;L^2({\mathbb R}^3;{\mathbb R}^9))})\ 
\|(v-\tilde v,q-\tilde q)\|_{X_T}\ .
\end{align*}
The Lipschitz continuity of $f$ from (\ref{Lip}) 
together with (\ref{fcirceta1}) provides 
us with the simple estimate 
$$\|f\circ\eta-f\circ\tilde\eta\|_{L^2(0,T;L^2(\Omega_0;{\mathbb R}^3))}
\le
C\ K\ \sqrt{T}\ \|(v-\tilde v,q-\tilde q)\|_{X_T}\ ,$$
from which it follows that
\begin{equation}
\label{fcirceta2}
\| f\circ\eta- f\circ\tilde\eta \|_ {L^2(0,T;H^1(\Omega_0;{\mathbb R}^3))}
\le C\ T^{\frac{1}{2}} \|(v-\tilde {v},q-\tilde {q})\|_{X_{T}}\ . 
\end{equation}
To estimate the time-derivative term, we have that
$$(f\circ\eta-f\circ\tilde\eta)_t=f_t(\eta)-f_t(\tilde\eta)+\nabla f(\eta)\ v-
\nabla f(\tilde\eta)\ \tilde v\ .$$
By using the Lipschitz condition (\ref{Lip}) 
on $f_t$ and $\nabla f$ together with (\ref{etah3})
and (\ref{fcirceta4}) we see that
$$\|(f\circ\eta-f\circ\tilde\eta)_t\|_{L^2(0,T;L^2(\Omega_0;{\mathbb R}^3))}
\le C (K\ \sqrt{T} + T^{\frac{1}{8}}\ \|\nabla f\|_{L^2(0,T;L^2({\mathbb R}^3;{\mathbb R}^9))})\  \|(v-\tilde {v},q-\tilde {q})\|_{X_{T}}\ ,$$
and hence that
\begin{equation}
\label{fcirceta5}
\|(f\circ\eta-f\circ\tilde\eta)_t\|_{L^2(0,T;H^1(\Omega_0;{\mathbb R}^3)')}
\le C\ T^{\frac{1}{8}}\  \|(v-\tilde {v},q-\tilde {q})\|_{X_{T}}\ ,
\end{equation}
which together with (\ref{fcirceta2}) completes the proof.
\end{proof}

We next have the following existence and uniqueness result which
shows that the convex  subset $C_T$ is closed under the map $\Theta_T$.
\begin{lemma}
\label{stability}
There exists $\overline{T}_0\in (0,T_0)$ such that for 
$0<T\le \overline{T}_0$ and for any $(v,q)\in C_T$, there exists a unique
solution $(\tilde v,\tilde q)\in C_T$ to the problem (\ref{nlfp}).
\end{lemma}

\begin{proof}
The linear problem (\ref{nlfp}) will be solved by the following sequence of
problems: We initiate the iteration with $(v_0,q_0)=(0,0)$; then, given
$(v_n,q_n) \in X_T$ we define $(v_{n+1},q_{n+1})$ 
as the solution of
\begin{subequations}
  \label{inlfp}
\begin{align}
(v_{n+1})_t- \nu \triangle {{v_{n+1}}} +
\nabla q_{n+1}&= F + \bar f (v_n, q_n)
\ \ \text{in} \ \ (0,T)\times \Omega_0 \,, 
         \label{inlfp.c}\\
  \operatorname{div} (v_{n+1})&= \bar a (v_n)     \ \ \text{in} \ \ (0,T)\times \Omega_0 \,, 
         \label{inlfp.d}\\
\Pi_0\left(\operatorname{Def} (v_{n+1})N\right) &= \bar g_1 (v_n)
\ \ \text{on} \ \ (0,T)\times \Gamma_0 \,, 
         \label{inlfp.e}\\
N\cdot (\ S({v_{n+1}},q_{n+1})\ N\ ) 
&=
\sigma\ N\cdot \triangle_0(\int_0^t {v_{n+1}}) 
+ \bar g_2 (v_n,q_n)+ \sigma \bar B (v_n)\nonumber\\
& \qquad + \sigma\ N\cdot \triangle_{g(t)}(x)
 \ \ \text{on} \ \ (0,T)\times \Gamma_0 \,, 
         \label{inlfp.g}\\
  {v_{n+1}} &= u_0    
 \ \ \text{on} \ \ \Omega_0\times \{ t=0\} \,, 
         \label{inlfp.h}
\end{align}
\end{subequations}

where

\begin{subequations}
  \label{ilfp}
\begin{align}
\bar f (v_n, q_n) &=
- \nu (\triangle {{v_{n}}} - (a^j_l a^k_l {v_n},_k),_j) + \nabla
q_{n}-a^T \nabla q_{n}
\ \ \text{in} \ \ (0,T)\times \Omega_0 \,, 
         \label{ilfp.c}\\
\bar a (v_n)     &=
  \operatorname{div} (v_{n})-a^k_i {v_n^i}_{,k}
\ \ \text{in} \ \ (0,T)\times \Omega_0 \,, 
         \label{ilfp.d}\\
 \bar g_1 (v_n)
&= 
\Pi_0\left(\operatorname{Def} (v_{n})\ N-\Pi_{\eta(t,\cdot)}
( D_{\eta(t,\cdot)}(v_n)\cdot a^T\ N)\right) 
\ \ \text{on} \ \ (0,T)\times \Gamma_0 \,, 
         \label{ilfp.e}\\
\bar g_2 (v_n,q_n)\,, 
&=
N\cdot (\ S({v_{n}},q_{n})\ N\ ) - 
N\cdot (\ S_{\eta(t,\cdot)}({v_{n}},q_{n})\cdot a^T\ N\ )
\ \ \text{on} \ \ (0,T)\times \Gamma_0 \,, 
         \label{ilfp.g}\\
\bar B (v_n)
&=
 N\cdot (\triangle_{g(t)}-\triangle_0)(\int_0^t {v_n})
\ \ \text{on} \ \ (0,T)\times \Gamma_0 \ . 
         \label{ilfp.h}
\end{align}
\end{subequations}

We note that the system (\ref{inlfp}) is well defined since by construction $\bar{g}_1 (v_n)$ belongs
to the tangent plane to the surface $\Gamma_0$. This is precisely the reason
why the iteration proposed by Tani \cite{Tani1996} fails, since the 
analogue of $g_1(v_n)$ in equation (4.3) on page 319 
does not necessarily belong to this tangent plane as it must.

For each $n\ge 1$ define
$$
\begin{array}{l}
{\delta v}_n=v_{n+1}-v_{n} \\
{\delta q}_n=q_{n+1}-q_n 
\end{array}
\begin{array}{l}
\text{(difference of velocities)}\\
\text{(difference of pressures)\, . }
\end{array}
$$
Since for each iteration $n\ge 0$, the problem 
(\ref{inlfp}) is linear, the pair $(\delta v_{n+1}, \delta q_{n+1})$ satisfies
\begin{subequations}
  \label{dnlfp}
\begin{align}
(\delta v_{n+1})_t- \nu \triangle {{\delta v}_{n+1}} +
\nabla {\delta q}_{n+1}&= \bar f ({\delta v}_{n},{\delta q}_{n})
\ \ \text{in} \ \ (0,T)\times \Omega_0 \,, 
         \label{dnlfp.c}\\
  \operatorname{div} ({\delta v}_{n+1})&=\bar a ({\delta v}_{n})     \ \ \text{in} \ \ (0,T)\times \Omega_0 \,, 
         \label{dnlfp.d}\\
\Pi_0(\operatorname{Def} \left({\delta v}_{n+1}\right)N)&= \bar g_1 ({\delta v}_n)
\ \ \text{on} \ \ (0,T)\times \Gamma_0 \,, 
         \label{dnlfp.e}\\
N\cdot (\ S({{\delta v}_{n+1}},{\delta q}_{n+1})\ N\ ) &=
\sigma\ N\cdot \triangle_0\left(\int_0^t {\delta v}_{n+1}\right) \nonumber \\
&\qquad
+ \bar g_2 ({\delta v}_n,{\delta q}_{n}) +
\sigma \bar B ({\delta v}_n) \ \ \text{on} \ \ (0,T)\times \Gamma_0 \,, \label{dnlfp.g}\\
 {\delta v}_{n+1} &= 0    
 \ \ \text{on} \ \ \{t=0\}\times \Omega_0 \ , 
         \label{dnlfp.h}
\end{align}
\end{subequations}
where the forcing appearing in the right-hand-side of this system are defined 
by the same relations as (\ref{ilfp}) with $(\delta v_n,\delta q_n)$ replacing
$(v_n,q_n)$. We note also that $\delta q_n (0,\cdot)=0$ for any $n\ge 1$.

For each $n\ge 1$, the initial boundary forcings $\bar{g_1}(\delta v_n)(0,\cdot)$, $\bar{g_2}(\delta v_n)(0,\cdot)$ and $\bar{B}(\delta v_n)(0,\cdot)$ satisfy the condition (\ref{trco1}) since they are all zero (from (\ref{dnlfp.h}) and the initial velocity satisfies the compatibility condition 
$$\Pi_0(\operatorname{Def} \left({\delta v}_{n+1}(0,\cdot)\right)N)
=0\ \text{on}\ \Gamma_0.$$
Thus, (\ref{dnlfp}) satisfies all of the assumptions of the basic linear
problem (\ref{linear}), and so we can use
estimate (\ref{estimate1})  to find that
\begin{align}
&\| ({\delta v}_{n+1},{\delta q}_{n+1})\|_{X_T} \nonumber\\
&\le   C \ \left(
\| \bar f ({\delta v}_{n},{\delta q}_{n})\|_ {L^2(0,T;H^1(\Omega_0;{\mathbb R}^3))}+ 
\|\bar f ({\delta v}_{n},{\delta q}_{n})_t\|_{L^2(0,T;H^1(\Omega_0;{\mathbb R}^3)')} \right.\nonumber \\
& \qquad\qquad \qquad 
+ \|\bar a ({\delta v}_{n})\|_ {L^2(0,T;H^2(\Omega_0;{\mathbb R}^9))} 
+ \|\bar a ({\delta v}_{n})_t\|_ {L^2(0,T;L^2(\Omega_0;{\mathbb R}^9))} \nonumber \\
& \qquad\qquad \qquad 
+ \|\bar B ({\delta v}_{n})\|_ {L^2(0,T;H^\frac{1}{2}(\Gamma_0;{\mathbb R}))}
+ \|\bar B({\delta v}_{n})_t\|_ {L^2(0,T;H^{\frac{1}{2}}(\Gamma_0;{\mathbb R}))}\nonumber \\
& \qquad\qquad \qquad +\|\bar g_1({\delta v}_{n},{\delta q}_{n})\|_ {L^2(0,T;H^\frac{3}{2}(\Gamma_0;{\mathbb R}^3))}
+ \|\bar g_1({\delta v}_{n},{\delta q}_{n})_t\|_ {L^2(0,T;H^{-\frac{1}{2}}(\Gamma_0;{\mathbb R}^3))}\nonumber\\
&\qquad\qquad\qquad\left. +\|\bar g_2({\delta v}_{n},{\delta q}_{n})\|_ {L^2(0,T;H^\frac{3}{2}(\Gamma_0;{\mathbb R}))}
+ \|\bar g_2({\delta v}_{n},{\delta q}_{n})_t\|_ {L^2(0,T;H^{-\frac{1}{2}}(\Gamma_0;{\mathbb R}))}\ . \right)\nonumber
\end{align}

Hence, as a consequence of Lemma \ref{estimates} (which directly follows
this proof),
we infer that there exists  $T_1\in (0,T_0)$ such that for any $T\in (0,T_1)$
$$
\|( {\delta v}_{n+1}, {\delta q}_{n+1})\|_{X_T} \nonumber \\
\le C \ T^{\frac{1}{16}}\  \| (v,q)\|_{X_T}\ 
\|( {\delta v}_{n}, {\delta q}_{n+1})\|_{X_T} \ .
$$

Since $v\in C_T$, we see from the previous inequality that there exists ${T}_2\in (0,T_1)$ such that for $T\le {T}_2$, we have 
\begin{equation*}
\| ({\delta v}_{n+1},{\delta q}_{n+1})\|_{X_T}\le  \frac{1}{2} T^{\frac{1}{32}}\ \|( {\delta v}_{n},{\delta q}_{n}\|_{X_T} 
\ .\end{equation*}
which by the contraction mapping principle (since we can assume $T_2\in (0,1)$) shows that the sequence $(v_n,q_n)_{n\in {\mathbb N}}$ is convergent in $X_T$ to a limit $(\tilde{v},\tilde{q})$ which is the unique solution in  $X_T$ of (\ref{nlfp}). By summing the previous inequalities from $n=1$ to $\infty$, and taking into account that the first iterate $(v_0,q_0)=(0,0)$, we see that
\begin{equation}
\label{sum}
\| ({\tilde v},{\tilde q})\|_{X_T} \le \ \frac{\frac{1}{2}\ T^{\frac{1}{32}}}{1-\frac{1}{2}\ T^{\frac{1}{32}}}\ 
\| ({ v}_{1}, {q}_{1})\|_{X_T} + \| ({ v}_{1}, {q}_{1})\|_{X_T}\ .
\end{equation}
In local coordinates $y^\alpha$, $\alpha=1,2$, on $\Gamma_0$, we have that
\begin{equation*}
N\cdot \triangle_g \text{Id}  
= N\cdot g^{\alpha \beta}\left(\frac{\partial^2 \text{Id} }
{\partial y^\alpha \partial y^\beta} - 
\Gamma^\gamma_{\alpha\beta}\frac{\partial \text{Id} }
{\partial y^\gamma}\right)\ ;
\end{equation*}
since $N$ is the unit normal to $\Gamma_0$,
and $\frac{\partial \text{Id}}{\partial y^\alpha}$ is tangent to
$\Gamma_0$,  we have that
\begin{equation}
\label{surface0}
 N\cdot \frac{\partial \text{Id} }{\partial y^\gamma}=
0\ \text{on}\ \Gamma_0 .
\end{equation}

Now, to estimate the right-hand side of (\ref{sum}), we simply apply our basic 
energy estimate (\ref{estimate1_omega}) to the following system, with zero initial tangential boundary forcing and an initial normal boundary forcing in $H^{\frac{3}{2}}(\Gamma_0;{\mathbb R})$, whose initial data satisfies the initial compatibility condition $$\Pi_0 (\operatorname{Def} \left(v_{1}(0,\cdot)\right)N)=0\ \text{on}\ \Gamma_0,$$ (since $v_{1}(0,\cdot)=u_0$):
\begin{subequations}
\begin{align}
\partial_t{{v_{1}}}- \nu \triangle {{v_{1}}} +
\nabla q_{1}&= F\ \ \text{in} \ \ (0,T)\times \Omega_0 \,, 
   \nonumber      \\
  \operatorname{div} (v_{1})&= 0    \ \ \text{in} \ \ (0,T)\times \Omega_0 \,, 
    \nonumber      \\
\Pi_0(\operatorname{Def} (v_{1})N) &=  0
\ \ \text{on} \ \ (0,T)\times \Gamma_0 \,, 
     \nonumber     \\
N\cdot (\ S({v_{1}},q_{1})\ N\ ) -
\sigma\ N\cdot \triangle_0(\int_0^t {v_{1}})&= \sigma\ N\cdot \triangle_{g(t)}(\text{Id})
 \ \ \text{on} \ \ (0,T)\times \Gamma_0 \,, 
      \nonumber    \\
  {v_{1}} &= u_0    
 \ \ \text{on} \ \ \Omega_0\times \{ t=0\} \,, 
         \end{align}
\end{subequations}
to get
\begin{align}
\label{v1q1}
&\| ({v}_{1}, {q}_{1})\|_{X_T} \nonumber \\
 & \le C_m \ \left( \|u_0\|_{H^2(\Omega_0;{\mathbb R}^3)} + 
\|  F \|_ {L^2(0,T;H^1(\Omega_0;{\mathbb R}^3))}+ 
\|F_t\|_ {L^2(0,T;H^1(\Omega_0;{\mathbb R}^3)')} \right.\nonumber \\
& 
\qquad\qquad+ \|\triangle_0\text{Id}\cdot N\|_{H^\frac{1}{2}(\Gamma_0;{\mathbb R})}+\| \  N\cdot g^{\alpha \beta}\frac{\partial^2 \text{Id} }
{\partial y^\alpha \partial y^\beta}\|_ {L^2(0,T;H^\frac{3}{2}(\Gamma_0;{\mathbb R}))}\nonumber\\
&\left.\qquad\qquad+ \|F(0,\cdot)\|_{L^2 (\Omega_0;{\mathbb R}^3)}+ \| \partial_t(N\cdot g^{\alpha \beta}\frac{\partial^2 \text{Id} }
{\partial y^\alpha \partial y^\beta})\ 
\|_ {L^2(0,T;H^{-\frac{1}{2}}(\Gamma_0;{\mathbb R}))}\right)\,,
\end{align}
where $C_m$ is defined in Definition \ref{Cm}.

Then, combining (\ref{again}), (\ref{sum}), (\ref{v1q1}) and 
(\ref{forcingcmbis}), we find that
\begin{equation*}
\| ({\tilde v},{\tilde q})\|_{X_T} \le \ \frac{1}{2}\  \left( 1+ \frac{ T^{\frac{1}{32}}}{1-\frac{1}{2}\ T^{\frac{1}{32}}}\right) \ M(T)\ (1+\epsilon (T))\ ,\end{equation*}
where $\epsilon (T)\rightarrow 0$ as $T\rightarrow 0$.
Thus, choosing 
$\overline{T}_0\in (0,T_2)$, such that $0<T\le \overline{T}_0$, implies that
\begin{equation*}
\| ({\tilde v}, {\tilde q})\|_{X_T} \le M (T) \ . 
\end{equation*}
so that $({\tilde v},{\tilde q})\in C_T$, which proves the stability lemma.
\end{proof}

We now prove the following estimate which was used to prove 
Lemma \ref{stability}.
\begin{lemma}
\label{estimates}
There exists $T_1\in (0,T_0)$ such that for any $(v,q)\in C_T$, 
and $(w,r)\in X_T$ such that $(w(0),r(0))=(0,0)$, we have that
\begin{align*}
&\| \bar f (w,r)\|_ {L^2(0,T;H^1(\Omega_0;{\mathbb R}^3))}+ 
\|\bar f (w,r)_t\|_ {L^2(0,T;H^1(\Omega_0;{\mathbb R}^3)')} \nonumber \\
& \qquad\qquad  
+ \|\bar a (w)\|_ {L^2(0,T;H^2(\Omega_0;{\mathbb R}^3))} 
+ \|\bar a (w)_t\|_ {L^2(0,T;L^2(\Omega_0;{\mathbb R}^3))} \nonumber \\
& \qquad\qquad  
+ \|\bar B (w)\|_ {L^2(0,T;H^\frac{1}{2}(\Gamma_0;{\mathbb R}))}
+ \|\bar B(w)_t\|_ {L^2(0,T;H^{\frac{1}{2}}(\Gamma_0;{\mathbb R}))}\nonumber \\
& \qquad\qquad + \Sigma_{\alpha=1}^2 \left(\|\bar g_\alpha (w,r)\|_ {L^2(0,T;H^\frac{3}{2}(\Gamma_0;{\mathbb R}))}
+ \|\bar g_\alpha(w,r)_t\|_ {L^2(0,T;H^{-\frac{1}{2}}(\Gamma_0;{\mathbb R}))}\right)\nonumber \\
& \qquad\qquad\qquad\qquad \le C \ T^{\frac{1}{16}}\  \| (v,q)\|_{X_T} \  \| (w,r)\|_{X_T}\end{align*}
where
\begin{align}
\bar f (w, r) &=
- \nu (\triangle {w} - (a^j_l a^k_l {w},_k),_j) + \nabla
r-a^T \ \nabla r
\ \ \operatorname{in} \ \ (0,T)\times \Omega_0 \,, 
         \nonumber\\
\bar a (w)   &=  
  \operatorname{div} (w)-a^i_k {w}^k_{,i}
\ \ \operatorname{in} \ \ (0,T)\times \Omega_0 \,, 
         \nonumber\\
\bar g_1 (w) &=
\Pi_0\left(\operatorname{Def} (w)N-\Pi_{\eta(t,\cdot)}\left( D_{\eta(t,\cdot)}
(w)\cdot a^T\ N\right)\ \right) 
\ \ \operatorname{on} \ \ (0,T)\times \Gamma_0 \,, 
         \nonumber\\
\bar g_2 (w,r) &=
N\cdot (\ S({w},r)\ N\ ) - N\cdot(\ S_{\eta(t,\cdot)}(w,r)\cdot a^T\ N\ ) 
\ \ \operatorname{on} \ \ (0,T)\times \Gamma_0 \,,
         \nonumber\\
\bar B (w) &=
 N\cdot (\triangle_{g(t)}-\triangle_0)(\int_0^t {w})
 \ \ \operatorname{on} \ \ (0,T)\times \Gamma_0 \ . 
         \nonumber
\end{align}
\end{lemma}

\begin{proof} 
In order to simplify the notation, we shall 
omit the explicit dependence on time $t\in (0,T)$ when we evaluate certain 
Sobolev norms.  For instance, letting $t$ denote a generic element of $(0,T)$, 
$\|w\|_{H^3(\Omega_0;{\mathbb R}^3)}$ will stand for 
$\|w(t,\cdot)\|_{H^3(\Omega_0;{\mathbb R}^3)}$.
For the sake of clarity, we divide the proof into ten steps.

{\bf 1. Estimate of $\bar{a}$:}\hfill\break
Since $\bar{a}(w)=(\delta^i_k-a^i_k)\ w^k,_{i}$, we have that
\begin{align}
& \|\bar{a}(w)\|_{H^2(\Omega_0;{\mathbb R}^9)}\nonumber\\
&\le C\ \left(
 \|\delta^i_k-a^i_k\|_{L^\infty(\Omega_0;{\mathbb R}^9)}\  \|w^k,_{i}\|_{H^2(\Omega_0;{\mathbb R}^9)}
+  \|\delta^i_k-a^i_k\|_{H^2 (\Omega_0;{\mathbb R}^9)}\  \|w^k,_{i}\|_{L^\infty (\Omega_0;{\mathbb R}^9)}\right)\ ,
\end{align}
which by the continuous embedding of $H^2$ into $L^\infty$ shows that
$$ \|\bar{a}(w)\|_{H^2(\Omega_0;{\mathbb R}^9)}\le C\ \|\delta^i_k-a^i_k\|_{H^2(\Omega_0;{\mathbb R}^9)}\  \|w^k,_{i}\|_{H^2(\Omega_0;{\mathbb R}^9)}
,
$$
and consequently from (\ref{ah2}),
\begin{equation}
\label{e1}
\|\bar{a}(w) \|_ {L^2(0,T;H^2 ({\Omega_0};{\mathbb R}^3))}\le 
C\ \sqrt{T}\  
\| w \|_ {L^2 (0,T;H^3({\Omega_0};{\mathbb R}^3))}\ .
\end{equation}

{\bf 2. Estimate of $\bar{a}_t$:}\hfill\break
Since $\bar{a}(w)=(\delta^i_k-a^i_k)\ w^k,_{i}$, we have by time differentiation,
$$ \bar{a}(w)_t=(\delta^i_k-a^i_k)\ (w^k,_{i})_t+
(a^i_k)_t \ w^k,_{i}\ ,
$$
which together with (\ref{ah2}) and (\ref{timeah1}) shows that
\begin{align}
\|\bar{a}(w)_t \|_ {L^2(0,T;L^2 ({\Omega_0};{\mathbb R}^9))}\le & 
C\ \sqrt{T}\  
\|w_t \|_ {L^2 (0,T;H^1({\Omega_0};{\mathbb R}^3))}  \nonumber\\
&+ \| a_t\cdot \nabla w \|_ {L^2 (0,T;L^2({\Omega_0};{\mathbb R}^9))}
.
\label{l1}
\end{align}
Now, we notice that
\begin{equation}
\label{l2}
 \| a_t\cdot \nabla w \|_ {L^2 (0,T;L^2({\Omega_0};{\mathbb R}^9))}\le
 \|  a_t \|_ {L^\infty (0,T;L^2({\Omega_0};{\mathbb R}^9))}
 \| \nabla w \|_ {L^2 (0,T;L^\infty({\Omega_0};{\mathbb R}^9))}\ .
\end{equation}
From the Sobolev embedding theorem, we have for all $t$,
\begin{equation*}
 \| \nabla w (t,\cdot)\|_ {L^\infty({\Omega_0};{\mathbb R}^9)}\ \le
\| \nabla w (t,\cdot)\|_ {W^{1,4}({\Omega_0};{\mathbb R}^9)}\ ,
\end{equation*}
and thus the Gagliardo-Nirenberg inequality (\ref{GN}) show that
\begin{equation}
\label{l3}
 \| \nabla w (t,\cdot) \|_ {L^\infty({\Omega_0};{\mathbb R}^9))}\  \le
C \| \nabla w (t,\cdot) \|_ {L^2({\Omega_0};{\mathbb R}^9))}^{\frac{1}{8}}\ 
\| \nabla w (t,\cdot)\|_ {H^2({\Omega_0};{\mathbb R}^9))}^{\frac{7}{8}}\ .
\end{equation}
Since $w(0)=0$, we have $\nabla w(t,\cdot)=\int_0^t \nabla w_t (t',\cdot)\ dt'$, and thus from the Cauchy-Schwarz inequality and Fubini's theorem,
\begin{equation*}
 \| \nabla w (t,\cdot)\|_ {L^2 ({\Omega_0};{\mathbb R}^9)}  \le
\sqrt{t}\ 
\| \nabla w_t \|_ {L^2(0,t; L^2({\Omega_0};{\mathbb R}^9)} \ ,
\end{equation*}
and then
\begin{equation}
\label{l4}
 \| \nabla w (t,\cdot)\|_ {L^2 ({\Omega_0};{\mathbb R}^9)}  \le
\sqrt{T}\ 
\| \nabla w_t \|_ {L^2(0,T; L^2({\Omega_0};{\mathbb R}^9))}\
.
\end{equation}
From (\ref{l3}), and (\ref{l4}), we deduce that
\begin{equation*}
\| \nabla w \|_ {L^2 (0,T;L^\infty({\Omega_0};{\mathbb R}^9))}\le C T^{\frac{1}{16}}\ \|  w_t \|_ {L^2(0,T; H^1({\Omega_0};{\mathbb R}^9))}^{\frac{1}{8}}
\left( \int_0^T \| w(t,\cdot) \|_ {H^3({\Omega_0};{\mathbb R}^9))}^{\frac{7}{4}} dt\right)^{\frac{1}{2}}\ ,
\end{equation*}
and from H\"{o}lder's inequality,
\begin{equation}
\label{l5}
\| \nabla w \|_ {L^2 (0,T;L^\infty({\Omega_0};{\mathbb R}^9))}\le C T^{\frac{1}{8}}\  
 \| w \|_ {L^2 (0,T; H^3({\Omega_0};{\mathbb R}^3))}^{\frac{7}{8}}\ 
\|  w_t \|_ {L^2(0,T; H^1({\Omega_0};{\mathbb R}^3))}^{\frac{1}{8}}
\ .
\end{equation}
From (\ref{timeah1}), (\ref{l2}) and (\ref{l5}), and Young's inequality 
we infer that
\begin{align}
 \| a_t\cdot \nabla w &\|_ {L^2 (0,T;L^2({\Omega_0};{\mathbb R}^3))}\nonumber\\
&\le C D(T)\ T^{\frac{1}{8}} 
\left(\ \| w(t,\cdot) \|_ {L^2 (H^3({\Omega_0};{\mathbb R}^3))} + \|  w_t \|_ {L^2(0,T; H^1({\Omega_0};{\mathbb R}^3))} \ \right)\ .\nonumber
\end{align}

From this inequality and (\ref{l1}), we then deduce the existence of 
$T_{\bar{a}_t}\in(0,T_0)$ such that  for any $T \in (0,T_{\bar a_t})$,
\begin{equation}
\label{l6}
\|a(w)_t\|_{X_T} \le C T^{\frac{1}{8}}\ \|(w,r)\|_{X_T}\ .
\end{equation}

{\bf 3. Estimate of $\bar{f}$:}\hfill\break
On $(0,T)$, 
\begin{equation}
\label{f1}
{\|\bar{f}(w,r)\|}_{H^1(\Omega_0;{\mathbb R}^3)}\le \nu\ {\|\operatorname{div}(\nabla w \cdot (\text{Id}-a\cdot a^T))\|}_{H^1(\Omega_0;{\mathbb R}^3)}+ {\|(\text{Id}-a)\ \nabla r \|}_{H^1(\Omega_0;{\mathbb R}^3)}\ .
\end{equation}
From the usual product formula for the derivatives and the
Cauchy-Schwarz inequality we see that
\begin{align}
 {\|(\text{Id}-a)\nabla r\|}_{H^1(\Omega_0;{\mathbb R}^3)}\le C\ &
 {\|\text{Id}-a\|}_{W^{1,4}(\Omega_0;{\mathbb R}^9)}  {\|\nabla r\|}_{L^4(\Omega_0;{\mathbb R}^3)}\nonumber\\
&+ C\ {\|\text{Id}-a\|}_{L^{4}(\Omega_0;{\mathbb R}^9)}  {\|\nabla r\|}_{W^{1,4}(\Omega_0;{\mathbb R}^3)}
\ ,\nonumber
\end{align}
which by the continuity of the embedding of $H^2$ into $W^{1,4}$ shows that
\begin{equation*}
 {\|(\text{Id}-a)\nabla r\|}_{H^1(\Omega_0;{\mathbb R}^9)}\le C\ {\|\text{Id}-a\|}_{H^2 (\Omega_0;{\mathbb R}^9)}  {\|r\|}_{H^2 (\Omega_0;{\mathbb R})}\ ,
\end{equation*}
and consequently from (\ref{ah2}),
\begin{equation}
\label{f2}
\|(\text{Id}-a)\nabla r\|_ {L^2(0,T;H^1 ({\Omega_0};{\mathbb R}^3))}\le 
C\ \sqrt{T}\ \  
\| r \|_ {L^2 (0,T;H^2({\Omega_0};{\mathbb R}^3))})\ .
\end{equation}
For the other term on the right-hand-side of (\ref{f1}) we similarly have
\begin{align}
 {\|\operatorname{div}(\nabla w \cdot (\text{Id}-a\cdot a^T))\|}_{H^1(\Omega_0;{\mathbb R}^9)}\le C\ &
 {\|\text{Id}-a\cdot a^T\|}_{L^\infty (\Omega_0)}  {\| w\|}_{H^3(\Omega_0;{\mathbb R}^3)}\nonumber\\
&+ C\ {\|\text{Id}-a\cdot a^T\|}_{W^{1,4}(\Omega_0;{\mathbb R}^9)}  {\|r\|}_{W^{2,4}(\Omega_0;{\mathbb R})}
\ ,\nonumber
\end{align}
which by the continuity of the embeddings of $H^2$ into $W^{1,4}$, of
$H^3$ into $W^{2,4}$ and of $H^2$ into $L^\infty$ shows that
\begin{equation}
\label{f3}
 {\|\operatorname{div}(\nabla w \cdot (\text{Id}-a\cdot a^T))\|}_{H^1(\Omega_0;{\mathbb R}^9)}\le C\ 
 {\|\text{Id}-a\cdot a^T\|}_{H^2(\Omega_0;{\mathbb R}^9)}  {\| w\|}_{H^3(\Omega_0;{\mathbb R}^3)} \ .
\end{equation}
Now, from (\ref{ah2}) and the fact that $H^2$ is a Banach algebra, we deduce
that
\begin{equation*}
\|a\cdot a^T-\text{Id} \|_ {L^\infty(0,T; H^2(\overline{\Omega_0};{\mathbb R}^9))} \le  
 C\ \sqrt{T}  \ ,
\end{equation*}
which combined with (\ref{f3}) gives us the estimate
\begin{equation}
\label{f4}
 {\|\operatorname{div}(\nabla w \cdot (\text{Id}-a\cdot a^T))\|}_{L^2 (0,T; H^1({\Omega_0};{\mathbb R}^3))}\le C\ \sqrt{T}\ \  
 {\| w\|}_ {L^2(0,T; H^3({\Omega_0};{\mathbb R}))} \ .
\end{equation}
Thus, (\ref{f1}), (\ref{f2}) and (\ref{f4}) show that
\begin{equation}
\label{fbis}
 {\|\bar{f}(w,r)\|}_{L^2 (0,T; H^1({\Omega_0};{\mathbb R}^3))}\le C\ \sqrt{T}\  \left( {\| w\|}_ {L^2(0,T; H^3({\Omega_0};{\mathbb R}^3))}+
 {\| r\|}_ {L^2(0,T; H^2({\Omega_0};{\mathbb R}))}\right) \ .
\end{equation}

{\bf 4. Estimate of $\bar{f}_t$:}\hfill\break
We have 
\begin{align}
\label{ft1}
{\|\bar{f}(w,r)_t\|}_{H^1(\Omega_0;{\mathbb R}^3)'}
&\le \nu\ {\|\operatorname{div}(\nabla w_t \cdot 
(\text{Id}-a\cdot a^T))\|}_{H^1(\Omega_0;{\mathbb R}^3)'}
\nonumber\\
&\qquad 
+ \nu\ {\|\operatorname{div}(\nabla w \cdot 
(\text{Id}-a\cdot a^T)_t)\|}_{H^1(\Omega_0;{\mathbb R}^3)'} \nonumber \\
&\qquad 
+ {\|(\text{Id}-a^T)\nabla r_t\|}_{H^1(\Omega_0;{\mathbb R}^3)'}
+ {\|(\text{Id}-a^T)_t\nabla r\|}_{H^1(\Omega_0;{\mathbb R}^3)'}\, .
\end{align}
Now, we see that for any $\phi\in  H^1({\Omega_0};{\mathbb R}^3)$,
\begin{equation*}
\langle (\text{Id}-a^T)\nabla r_t,\phi \rangle_1 = 
\langle \nabla r_t,(\text{Id}-a) \phi \rangle_1 \
\end{equation*}
where $\langle \cdot ,\cdot \rangle_1$ denotes the duality product 
between $H^1({\Omega_0};{\mathbb R}^3)$ and 
$ H^1({\Omega_0};{\mathbb R}^3)'$. Consequently,
\begin{align}
 |\langle (\text{Id}-a^T)\nabla r_t,\phi \rangle_1|&\le   
 {\|(\text{Id}-a)\phi\|}_{H^1({\Omega_0};{\mathbb R}^3)}  {\|\nabla r_t\|}_{H^1({\Omega_0};{\mathbb R}^3)'}\nonumber\\
& \le C\ {\|\text{Id}-a\|}_{H^2({\Omega_0};{\mathbb R}^3)} 
{\|\nabla r_t\|}_{H^1({\Omega_0};{\mathbb R}^3)'} 
\|\phi\|_{H^1({\Omega_0};{\mathbb R}^3)}\ ,\nonumber
\end{align}
which shows that
\begin{equation*}
\|(\text{Id}-a^T)\nabla r_t\|_{H^1({\Omega_0};{\mathbb R}^3)'}\le C\ 
{\|\text{Id}-a\|}_{H^2({\Omega_0};{\mathbb R}^3)} {\|\nabla r_t\|}_{H^1({\Omega_0};{\mathbb R}^3)'}\ ,
\end{equation*}
and thus from (\ref{ah2}),
\begin{equation}
\label{ft2}
\|(\text{Id}-a^T)\nabla r_t\|_{L^2(0,T;H^1({\Omega_0};{\mathbb R}^3)')}\le C\ \sqrt{T}\ {\|\nabla r_t\|}_{L^2(0,T;H^1({\Omega_0};{\mathbb R}^3)')}\ .
\end{equation}

For the fourth term on the right-hand-side of (\ref{ft1}), 
we notice that for any $\phi\in  H^1({\Omega_0};{\mathbb R}^3)$,
$$ \langle\ (\text{Id}-a^T)_t\nabla r,\phi\ \rangle_1 
= - \int_\Omega a^T_t \nabla r\cdot \phi .$$
Thus,$$ |\langle (\text{Id}-a^T)_t\nabla r,\phi \rangle_1| 
\le  \|a_t\|_{L^2({\Omega_0};{\mathbb R}^9)} 
\|\nabla r\|_{L^4({\Omega_0};{\mathbb R};{\mathbb R}^3)} 
\|\phi\|_{L^4({\Omega_0};{\mathbb R}^3)}\ ,$$
and by Sobolev's embedding theorem and the Gagliardo-Nirenberg inequality,
we find that
$$ |\langle (\text{Id}-a^T)_t\nabla r,\phi \rangle_1| 
\le  C\ \|a_t\|_{L^2({\Omega_0};{\mathbb R}^3)} 
{\| r\|_{H^2({\Omega_0};{\mathbb R})}}^{\frac{7}{8}} 
{\| r\|_{L^2({\Omega_0};{\mathbb R})}}^{\frac{1}{8}}
\|\phi\|_{H^1({\Omega_0};{\mathbb R}^3)}\ ,$$
which shows that
$$ \|(\text{Id}-a^T)_t\nabla r\|_{H^1({\Omega_0};{\mathbb R}^9)'} \le  C\ \|a_t\|_{L^2({\Omega_0};{\mathbb R}^9)} {\| r\|_{H^2({\Omega_0};{\mathbb R})}}^{\frac{7}{8}} 
{\| r\|_{L^2({\Omega_0};{\mathbb R})}}^{\frac{1}{8}}\ .$$
and consequently from (\ref{timeah1}),
\begin{equation}
\label{ft3}
\|(\text{Id}-a^T)_t\nabla r\|_{H^1({\Omega_0};{\mathbb R}^3)'}\le C\  
{\| r\|_{H^2({\Omega_0};{\mathbb R})}}^{\frac{7}{8}} 
{\| r\|_{L^2({\Omega_0};{\mathbb R})}}^{\frac{1}{8}}
\ .
\end{equation}
Since $r(0)=0$, we have $r(t,\cdot)=\int_0^t r_t (t',\cdot)\ dt'$, 
and thus from the Cauchy-Schwarz inequality and Fubini's theorem,
\begin{equation*}
 \| r (t,\cdot)\|_ {L^2 ({\Omega_0};{\mathbb R}))}  \le
\sqrt{t}\ 
\|  r_t \|_ {L^2(0,t; L^2({\Omega_0};{\mathbb R}))} \ ,
\end{equation*}
and thus
\begin{equation*}
 \| r (t,\cdot)\|_ {L^2 ({\Omega_0};{\mathbb R}))}  \le
\sqrt{T}\ 
\|  r_t \|_ {L^2(0,T; L^2({\Omega_0};{\mathbb R}))} \ ,
\end{equation*}
which with (\ref{ft3}) gives
\begin{align*}
\|(\text{Id}-a^T)_t\nabla r&\|_{L^2(0,T;H^1({\Omega_0};{\mathbb R}^3)')}\nonumber\\
\le &C\ T^{\frac{1}{16}}\ \|  r_t \|_ {L^2(0,T; L^2({\Omega_0};{\mathbb R}))}^{\frac{1}{8}}\  
\left( \int_0^T \| r(t,\cdot) \|_ {H^2({\Omega_0};{\mathbb R}))}^{\frac{7}{4}} dt\right)^{\frac{1}{2}}
\ .\nonumber
\end{align*}
Now, from H\"{o}lder's inequality,
\begin{equation*}
\|(\text{Id}-a^T)_t\nabla r\|_{L^2(0,T;H^1({\Omega_0};{\mathbb R}^3)')}\le C\  T^{\frac{1}{8}}\ \|  r_t \|_ {L^2(0,T; L^2({\Omega_0};{\mathbb R}))}^{\frac{1}{8}}\ \| r \|_ {L^2(0,T; H^2({\Omega_0};{\mathbb R}))}^{\frac{7}{8}}\ ,
\end{equation*}
and Young's inequality finally give us the estimate
\begin{equation}
\label{ft4}
\|(\text{Id}-a^T)_t\nabla r\|_{L^2(0,T;H^1({\Omega_0};{\mathbb R}^3)')}\le C\ T^{\frac{1}{8}}\ [\|  r_t \|_ {L^2(0,T; L^2({\Omega_0};{\mathbb R}))}+ \| r \|_ {L^2(0,T; H^2({\Omega_0};{\mathbb R}))} ]\ .
\end{equation}

For the first term on the right-hand-side of (\ref{ft1}), we first notice that
for any $\phi\in  H^1({\Omega_0};{\mathbb R}^3)$ we have
\begin{align}
 \langle \operatorname{div}(\nabla w_t \cdot 
(\text{Id}-a\cdot a^T)),\phi \rangle_1 =&
-\int_{\Omega_0} (\text{Id}-a\cdot a^T)^j_k \partial_k w_t\cdot \partial_j \phi\ dx \nonumber\\
&\qquad
+ \langle \partial_j w_t,
(\text{Id}-a\cdot a^T)^j_k N_j \phi \rangle_{\Gamma_0} ,\nonumber\ .
\end{align}
 Consequently,
\begin{align}
 \langle \operatorname{div}(\nabla w_t \cdot (\text{Id}-a\cdot a^T)),
\phi \rangle_1 \le &
\|\text{Id}-a\cdot a^T\|_{L^\infty ({\Omega_0};{\mathbb R}^9)} 
\|\nabla w_t\|_{L^2 ({\Omega_0};{\mathbb R}^9)} \|\phi\|_{H^1 ({\Omega_0};{\mathbb R}^3)}\nonumber\\
&+  \| \partial_j w_t\|_{H^{-\frac{1}{2}}({\Gamma_0};{\mathbb R}^3)}
\|(\text{Id}-a\cdot a^T)^j_k N_j \phi\|_{H^{\frac{1}{2}}({\Gamma_0};{\mathbb R}^3)}\ .\nonumber
\end{align}
If $\tilde N$ denotes a smooth extension of $N$ to $\Omega_0$, we infer that
\begin{align}
 \langle \operatorname{div}(\nabla w_t \cdot (\text{Id}-a\cdot a^T)),
\phi \rangle_1 &\le 
\|\text{Id}-a\cdot a^T\|_{L^\infty ({\Omega_0};{\mathbb R}^9)} 
\|\nabla w_t\|_{L^2 ({\Omega_0};{\mathbb R}^9)} \|\phi\|_{H^1 ({\Omega_0};{\mathbb R}^3)}\nonumber\\
&+  C\ \| \partial_j w_t\|_{H^{-\frac{1}{2}}({\Gamma_0};{\mathbb R}^3)}
\|(\text{Id}-a\cdot a^T)^j_k \tilde{N}_j \phi\|_{H^1({\Omega_0};{\mathbb R}^3)}\ ,\nonumber
\end{align}
and by the Cauchy-Schwarz inequality,
\begin{align}
 \langle \operatorname{div}(\nabla w_t \cdot (\text{Id}&-a\cdot a^T ))
,\phi \rangle_1\nonumber\\ \le\  
&\|\text{Id}-a\cdot a^T\|_{L^\infty ({\Omega_0};{\mathbb R}^9)} 
\|\nabla w_t\|_{L^2 ({\Omega_0};{\mathbb R}^3)} \|\phi\|_{H^1 ({\Omega_0};{\mathbb R}^3)}\nonumber\\
&+  C\ \| \partial_j w_t\|_{H^{-\frac{1}{2}}({\Gamma_0};{\mathbb R}^3)}
\|\text{Id}-a\cdot a^T\|_{L^\infty ({\Omega_0};{\mathbb R}^9)}\ 
\|\phi\|_{H^1({\Omega_0};{\mathbb R}^3)}\nonumber\\
&+  C\ \| \partial_j w_t\|_{H^{-\frac{1}{2}}({\Gamma_0};{\mathbb R}^3)}
\|\text{Id}-a\cdot a^T\|_{W^{1,4}({\Omega_0};{\mathbb R}^9)}\ 
\|\phi\|_{L^4({\Omega_0};{\mathbb R}^3)}
\ ,\nonumber
\end{align}
followed by Sobolev's embedding theorem,
\begin{align}
 \langle \operatorname{div}(\nabla w_t \cdot &(\text{Id}-a\cdot a^T)),
\phi \rangle_1 \nonumber\\  \le &\ 
\|\text{Id}-a\cdot a^T\|_{H^2({\Omega_0};{\mathbb R}^9)} 
\|\nabla w_t\|_{L^2 ({\Omega_0};{\mathbb R}^3)} \|\phi\|_{H^1 ({\Omega_0};{\mathbb R}^3)}\nonumber\\
&+  C\ \| \partial_j w_t\|_{H^{-\frac{1}{2}}({\Gamma_0};{\mathbb R}^3)}
\|\text{Id}-a\cdot a^T\|_{H^2 ({\Omega_0};{\mathbb R}^9)}\ 
\|\phi\|_{H^1({\Omega_0};{\mathbb R}^3)}
\ .\nonumber
\end{align}
Consequently,
\begin{align}
 \|\operatorname{div}(\nabla w_t \cdot (\text{Id}-a\cdot a^T))\|_{H^1({\Omega_0};{\mathbb R}^3)'} \le & 
\|\text{Id}-a\cdot a^T\|_{H^2({\Omega_0};{\mathbb R}^9)} 
\|\nabla w_t\|_{L^2 ({\Omega_0};{\mathbb R}^3)} \nonumber\\
&+  C\ \| \partial_j w_t\|_{H^{-\frac{1}{2}}({\Gamma_0};{\mathbb R}^3)}
\|\text{Id}-a\cdot a^T\|_{H^2 ({\Omega_0};{\mathbb R}^9)}\ 
\ .\nonumber
\end{align}
From (\ref{ah2}), we then deduce that
\begin{align}
 \|\operatorname{div}(\nabla w_t \cdot (\text{Id}-a\cdot a^T))\|_{H^1({\Omega_0};{\mathbb R}^3)'} \le & 
C\ \sqrt{ T}\  
\|\nabla w_t\|_{L^2 ({\Omega_0};{\mathbb R}^3)}\nonumber\\ 
&+ C\ \sqrt{ T}\ \| \partial_j w_t\|_{H^{-\frac{1}{2}}({\Gamma_0};{\mathbb R}^3)} 
\ ,\nonumber
\end{align}
from which it follows that
\begin{align}
\label{ft5}
 \|\operatorname{div}(\nabla w_t \cdot (\text{Id}-a\cdot a^T))\|_{L^2(0,T;H^1({\Omega_0};{\mathbb R}^3)')} \le & 
C\ \sqrt{ T}\  
\|\nabla w_t\|_{L^2(0,T; L^2 ({\Omega_0};{\mathbb R}^3))}\nonumber\\ 
&+ C\sqrt{ T}\ \| \nabla w_t\|_{L^2(0,T;H^{-\frac{1}{2}}({\Gamma_0};{\mathbb R}^3))} 
\ .\end{align}

Similarly, for the second term on the right-hand-side of (\ref{ft1}), 
we notice that
for any $\phi\in  H^1({\Omega_0};{\mathbb R}^3)$ we have that
\begin{equation*}
 \langle \operatorname{div}(\nabla w \cdot (\text{Id}-a\cdot a^T)_t),
\phi \rangle_1 =-
\int_{\Omega_0} \left((a_t\cdot a^T+a\cdot a^T_t)^j_k  w,_k\right),_j\cdot  \phi \ dx\ .
\end{equation*}
Then, by H\"{o}lder's inequality,
\begin{equation*}
 \langle \operatorname{div}(\nabla w \cdot (\text{Id}-a\cdot a^T)_t),\phi\ \rangle_1 \le 
\left\|\left((-a_t\cdot a^T-a\cdot a^T_t)^j_k  w,_k\right),_j\right\|
_{L^{\frac{6}{5}} ({\Omega_0};{\mathbb R}^3)} 
 \|\phi\|
_{L^{6} ({\Omega_0};{\mathbb R}^3)}
\ .
\end{equation*}
Consequently, by the continuity of the embedding of $H^1 ({\Omega_0};{\mathbb R}^3)$ into $L^{6} ({\Omega_0};{\mathbb R}^3)$, we have that
\begin{equation}
\label{ft6}
 \|\operatorname{div}(\nabla w \cdot (\text{Id}-a\cdot a^T)_t)\|_{H^1 ({\Omega_0};{\mathbb R}^3)'} \le 
\|\left((a_t\cdot a^T+a\cdot a^T_t)^j_k  w,_k\right),_j\|
_{L^{\frac{6}{5}} ({\Omega_0};{\mathbb R}^3)} 
\ .
\end{equation}
On the one hand, we have from H\"{o}lder's inequality that
\begin{equation*}
\|(a_t\cdot a^T+a\cdot a^T_t)^j_k  w,_{kj}\|
_{L^{\frac{6}{5}} ({\Omega_0};{\mathbb R}^3)}\le
 \|(a_t\cdot a^T+a\cdot a^T_t)^j_k\|_{L^2 ({\Omega_0};{\mathbb R}^3)}  
\|w,_{kj}\|_{L^3({\Omega_0};{\mathbb R}^3)}\ ,
\end{equation*}
and by (\ref{ah2}), we deduce that
\begin{equation*}
\|(a_t\cdot a^T+a\cdot a^T_t)^j_k  w,_{kj}\|
_{L^{\frac{6}{5}} ({\Omega_0};{\mathbb R}^3)}\le
C\ \|a_t\|_{L^2 ({\Omega_0};{\mathbb R}^9)}  
\|w,_{kj}\|_{L^3({\Omega_0};{\mathbb R}^3)}\ .
\end{equation*}
With (\ref{timeah1}) and the Gagliardo-Nirenberg inequality, we then obtain
the estimate
\begin{equation*}
\|(a_t\cdot a^T+a\cdot a^T_t)^j_k  w,_{kj}\|
_{L^{\frac{6}{5}} ({\Omega_0};{\mathbb R}^3)}\le
C\   
\|w\|^{\frac{1}{4}}_{H^1({\Omega_0};{\mathbb R}^3)}
\|w\|^{\frac{3}{4}}_{H^3({\Omega_0};{\mathbb R}^3)}\ ,
\end{equation*}
which with (\ref{l4}) provides 
\begin{equation}
\label{ft7}
\|(a_t\cdot a^T+a\cdot a^T_t)^j_k  w,_{kj}\|
_{L^{\frac{6}{5}} ({\Omega_0};{\mathbb R}^3)}\le C\ 
T^{\frac{1}{8}}\  
\|w_t\|^{\frac{1}{4}}_{L^2(0,T;H^1({\Omega_0};{\mathbb R}^3))}
\|w\|^{\frac{3}{4}}_{H^3({\Omega_0};{\mathbb R}^3)}\ .
\end{equation}

On the other hand, from H\"{o}lder's inequality
\begin{equation*}
\|((a_t\cdot a^T+a\cdot a^T_t)^j_k),_j  w,_k\|
_{L^{\frac{6}{5}} ({\Omega_0};{\mathbb R}^3)}\le
 \|((a_t\cdot a^T+a\cdot a^T_t)^j_k),_j\|_{L^2 ({\Omega_0};{\mathbb R}^3)}  
\|w,_k\|_{L^3({\Omega_0};{\mathbb R}^3)}\ ,
\end{equation*}
(\ref{ah2}) allows us to assert that
\begin{align}
\|((a_t\cdot a^T+a\cdot a^T_t)^j_k),_j  w,_{k}\|
_{L^{\frac{6}{5}} ({\Omega_0};{\mathbb R}^3)}\le &
C\ \ \left(\|\nabla a_t\|_{L^2 ({\Omega_0};{\mathbb R}^{27})}  \|w,_k\|_{L^3({\Omega_0};{\mathbb R}^3)}\right.\nonumber\\
&\left. + \|\nabla a\|_{L^4 ({\Omega_0};{\mathbb R}^{27})} \|a_t\|_{L^4 ({\Omega_0};{\mathbb R}^9)} 
\|w,_k\|_{L^3({\Omega_0};{\mathbb R}^3)}\right)\ ,\nonumber
\end{align}
and thus from the Sobolev embedding theorem, and another use of (\ref{ah2}),
we find that
\begin{equation*}
\|((a_t\cdot a^T+a\cdot a^T_t)^j_k),_j  w,_{k}\|
_{L^{\frac{6}{5}} ({\Omega_0};{\mathbb R}^3)}\le 
C\ \|\nabla a_t\|_{L^2 ({\Omega_0};{\mathbb R}^{27})}  
\|w,_k\|_{L^3({\Omega_0};{\mathbb R}^3)}\ .
\end{equation*}
From the Gagliardo-Nirenberg inequalities, and (\ref{l4}), we then get
\begin{align}
\label{ft8}
&\|((a_t\cdot a^T+a\cdot a^T_t)^j_k),_j  w,_k\|
_{L^{\frac{6}{5}} ({\Omega_0};{\mathbb R}^3)}\nonumber\\
&\qquad\qquad\le   
C\|\nabla a_t\|_{L^2 ({\Omega_0};{\mathbb R}^{27})} T^{\frac{3}{8}}\|w_t\|^{\frac{3}{4}}_{L^2(0,T;H^1({\Omega_0};{\mathbb R}^3))} \|w\|^{\frac{1}{4}}_{H^3({\Omega_0};{\mathbb R}^3)}\ .
\end{align}
Consequently, from (\ref{ft6}), (\ref{ft7}) and (\ref{ft8}), we infer that
\begin{align*}
& \|\operatorname{div}(\nabla w \cdot (\text{Id}-a\cdot a^T)_t)\|_{L^2(0,T;H^1({\Omega_0};{\mathbb R}^3)')} \le\nonumber\\
 &C\ T^{\frac{1}{8}}\|w_t\|^{\frac{1}{4}}_{L^2(0,T;H^1({\Omega_0};{\mathbb R}^3))} \left(\int_0^T \|w(t,\cdot)\|^{\frac{3}{2}}_{H^3({\Omega_0};{\mathbb R}^3)}dt\right)^{\frac{1}{2}} \nonumber\\ 
 &+C\ T^{\frac{3}{8}}\|w_t\|^{\frac{3}{4}}_{L^2(0,T;H^1({\Omega_0};{\mathbb R}^3))} \left(\int_0^T \|w(t,\cdot)\|^{\frac{1}{2}}_{H^3({\Omega_0};{\mathbb R}^3)}  \|\nabla a_t(t,\cdot)\|^2_{L^2({\Omega_0};{\mathbb R}^3)} dt\right)^{\frac{1}{2}} \ .
\end{align*}
From H\"{o}lder's inequality and (\ref{timeah1}), we then have that
\begin{align*}
& \|\operatorname{div}(\nabla w \cdot (\text{Id}-a\cdot a^T)_t)\|_{L^2(0,T;H^1({\Omega_0};{\mathbb R}^3)')} \nonumber\\
 &\le C\ T^{\frac{1}{4}}\|w_t\|^{\frac{1}{4}}_{L^2(0,T;H^1({\Omega_0};{\mathbb R}^3))} \|w\|^{\frac{3}{4}}_{L^2(0,T;H^3({\Omega_0};{\mathbb R}^3))} \nonumber\\ 
 &\qquad
+ C\  T^{\frac{3}{4}}\|w_t\|^{\frac{3}{4}}_{L^2(0,T;H^1({\Omega_0};{\mathbb R}^3))} \|w\|^{\frac{1}{4}}_{L^2(0,T;H^3({\Omega_0};{\mathbb R}^3))}\ ,
\end{align*}
which by Young's inequalities gives
\begin{align}
\label{ft9}
& \|\operatorname{div}(\nabla w \cdot (\text{Id}-a\cdot a^T)_t)\|
_{L^2(0,T;H^1({\Omega_0};{\mathbb R}^3)')} \nonumber\\
 &\le C T^{\frac{1}{4}} \left(\|w_t\|_{L^2(0,T;H^1({\Omega_0};{\mathbb R}^3))}+ \|w\|_{L^2(0,T;H^3({\Omega_0};{\mathbb R}^3))}\right)\ .
\end{align}
Finally, from (\ref{ft1}), (\ref{ft2}), (\ref{ft4}), (\ref{ft5}), and 
(\ref{ft9}), we see that there exists $T_{\bar{f_t}}\in (0,T_0)$ such that 
for any $T\in (0,T_{\bar{f_t}})$ we have that
\begin{equation}
\label{ft}
\|f_t\|_{L^2(0,T;H^1({\Omega_0};{\mathbb R}^3)')} \le C\ T^{\frac{1}{8}}\ \|(w,r)\|_{X_T}\ .
\end{equation}
\hfill\break
Let us define $\bar{g}=\bar{g}_1+ \bar{g}_2 N$. We next estimate $\bar{g}$ and 
its time derivative.\hfill\break

{\bf 5. Estimate of $\bar g$:}\hfill\break

By definition,
\begin{align}
\bar{g}(w,r)=&S(w,r)\ N- S_{\eta(t,\cdot)}(w,r)\cdot a^T\ N+\nonumber\\
&\Pi_0 \left((a^T\ N)\cdot( S_{\eta(t,\cdot)}(w,r)\cdot a^T\  N)\ \frac{a^T\ N}{|a^T\ N|^2} \right)\ .
\label{g0}
\end{align}
Letting $\tilde{N}$ denote a smooth extension of $N$ to $\Omega_0$, we have that
$$ \|S(w,r)\ N- S_{\eta(t,\cdot)}(w,r)\cdot a^T\ N\|_{H^{\frac{3}{2}}({\Gamma_0};{\mathbb R}^3)}\le C \|S(w,r)\ \tilde{N}- S_{\eta(t,\cdot)}(w,r)\cdot a^T\ \tilde{N}\|_{H^2({\Omega_0};{\mathbb R}^3)}\ .$$
From
\begin{align}
 S(w,r)\ \tilde{N}&- S_{\eta(t,\cdot)}(w,r)\cdot a^T\ \tilde{N}\nonumber\\
&=
((\nabla w) + (\nabla w)^T)\tilde{ N} -r \tilde{N} - (a^T\cdot \nabla w + (a^T\cdot\nabla w)^T)\cdot a^T \tilde{N} +r a^T\ \tilde{N}\ ,
\end{align}
we infer that
\begin{align*}
&\|S(w,r)\ \tilde{N}- S_{\eta(t,\cdot)}(w,r)\cdot a^T\ \tilde{N}\|_{H^{\frac{3}{2}}({\Gamma_0};{\mathbb R}^3)}\le 
C\  \|r (\text{Id}-a^T)\ \tilde{N}\|_{H^2({\Omega_0};{\mathbb R}^3)}\\
& 
+C\ \|\nabla w \cdot (\text{Id}-a^T)\ \tilde{ N}\|
_{H^2({\Omega_0};{\mathbb R}^3)}
+C\ \|(\text{Id}-a^T) \cdot \nabla w \cdot a^T \ \tilde{ N}\|
_{H^2({\Omega_0};{\mathbb R}^3)}\\
&+C\ \|(\nabla w)^T\cdot (\text{Id}-a )\ \tilde{N}\|
_{H^2({\Omega_0};{\mathbb R}^3)} 
+ C\ \|(\nabla w)^T\cdot a\cdot (\text{Id}-a^T )\ \tilde{N}\|
_{H^2({\Omega_0};{\mathbb R}^3)}\ ,
\end{align*}
and since $H^2({\Omega_0};{\mathbb R}^3)$ is a Banach algebra,
\begin{align*}
\|S(w,r)\ \tilde{N}- S_{\eta(t,\cdot)}(w,r)&\cdot a^T\ \tilde{N}\|_{H^{\frac{3}{2}}({\Gamma_0};{\mathbb R}^3)}\nonumber\\
\le &  C \|\text{Id}-a^T\|_{H^2({\Gamma_0};{\mathbb R}^9)}
\ \|\nabla w \|_{H^2({\Omega_0};{\mathbb R}^9)}\\
&+ C\ \|a\|_{H^2({\Omega_0};{\mathbb R}^9)}\ \|\text{Id}-a^T\|_{H^2({\Gamma_0};{\mathbb R}^9)}
\ \|\nabla w \|_{H^2({\Omega_0};{\mathbb R}^9)}\\
&+ C\  \|r\|_{H^2({\Omega_0};{\mathbb R})}\  \|\text{Id}-a^T\|_{H^2({\Omega_0};{\mathbb R}^9)}\ .
\end{align*}
From (\ref{ah2}) we then get
\begin{align*}
&\|S(w,r)\ \tilde{N}- S_{\eta(t,\cdot)}(w,r)\cdot a^T\ \tilde{N}\|_{H^{\frac{3}{2}}({\Gamma_0};{\mathbb R}^3)} \\
&\qquad\qquad\qquad\qquad \le C\ \sqrt{T} \  
( \| w \|_{H^3({\Omega_0};{\mathbb R}^3)}
+  \|r\|_{H^2({\Omega_0};{\mathbb R})})\ ,
\end{align*}
and by integrating in time
\begin{align}
\label{g1}
&\|S(w,r)\ \tilde{N}- S_{\eta(t,\cdot)}(w,r)\cdot a^T\ \tilde{N}\|_{L^2(0,T; H^{\frac{3}{2}}({\Gamma_0};{\mathbb R}^3))} \nonumber \\
&\qquad\qquad\qquad\le C\ \sqrt{T} \   
\left( \| w \|_{L^2(0,T; H^3({\Omega_0};{\mathbb R}^3))}+ \|r\|_{L^2(0,T; H^2({\Omega_0};{\mathbb R}))}\right)\ .
\end{align}
On the other hand, we also have
\begin{align*}
&\Pi_0 \left((a^T N)\cdot (\ S_{\eta(t,\cdot)}(w,r)\cdot a^T\ N)\ a^T N\right)
\\
&\qquad\qquad\qquad = 
(a^T N)\cdot (\ S_{\eta(t,\cdot)}(w,r)\cdot a^T\ N\ ) a^T\ N\\
&\qquad\qquad\qquad\qquad
-(a^T N)\cdot(\ S_{\eta(t,\cdot)}(w,r)\cdot a^T\ N\ )\ (a^T N)\cdot (N)\ N\ ,
\end{align*}
and thus,
\begin{align}
&\Pi_0 \left((a^T N)\cdot(\ S_{\eta(t,\cdot)}(w,r)\cdot a^T\ N\ )\ a^T N\right)\nonumber\\
&\qquad\qquad\qquad= (a^T N)\cdot (\ S_{\eta(t,\cdot)}(w,r)\cdot a^T\ N\ ) (a^T- ((a^T\ N)\cdot N)\ \text{Id}) N\ .\nonumber
\end{align}
Since 
$$a^T- (a^T\ N\cdot N) \text{Id}=a^T-\text{Id}+((\text{Id}-a^T)\ N\cdot N)\ \text{Id}\ ,$$
we have in turn that
\begin{align}
\label{pio}
&\Pi_0 \left((a^T N)\cdot(\ S_{\eta(t,\cdot)}(w,r)\cdot a^T\ N\ )\ a^T N\right)\nonumber\\
&\qquad\qquad\qquad = (a^T N)\cdot (S_{\eta(t,\cdot)}(w,r)\cdot a^T\ N)\ (a^T-\text{Id}) N \nonumber\\
&\qquad\qquad\qquad +(a^T N)\cdot (\ S_{\eta(t,\cdot)}(w,r)\cdot a^T\ N\ ) ((\text{Id}-a^T)\ N)\cdot (N)\  N\ .
\end{align}
From the trace properties, we then find that
\begin{align}
\label{g2}
&\left\|\Pi_0 ((a^T N)\cdot (\ S_{\eta(t,\cdot)}(w,r)\cdot a^T\ N\ )\ \frac{a^T N}{|a^T N|^2})\right\|_{H^{\frac{3}{2}}({\Gamma_0};{\mathbb R}^3)}\nonumber\\ 
&\qquad\qquad\le \ \| \frac{a^T N}{|a^T N|^2}\cdot (\ S_{\eta(t,\cdot)}(w,r)\cdot  a^T\ N\ )\ (a^T-\text{Id}) N\|_{H^2({\Omega_0};{\mathbb R}^3)}\nonumber\\
&\qquad\qquad\qquad
+\|\frac{a^T N}{|a^T N|^2}\cdot(\ S_{\eta(t,\cdot)}(w,r)\cdot a^T\ N\ )\ ((\text{Id}-a^T)\ N)\cdot (N)\ N\|_{H^2 ({\Gamma_0};{\mathbb R}^3)}\ .
\end{align}
To estimate the first term on the right-hand-side of (\ref{g2}), we notice that
since $H^2({\Omega_0};{\mathbb R}^3)$ is a Banach algebra, we have that
\begin{align*}
&\left\| \frac{a^T N}{|a^T N|^2}\cdot  (\ S_{\eta(t,\cdot)}(w,r)\cdot a^T\ N\ )\ (a^T-\text{Id})\ N\right\|_{H^2({\Omega_0};{\mathbb R}^3)} \\
&\le\ C \left\|\frac{1}{|a^T N|^2}\right\|_{H^2({\Omega_0};{\mathbb R})} \|a\|_{H^2({\Omega_0};{\mathbb R}^9)}
\|S_{\eta(t,\cdot)}(w,r)\cdot a\|_{H^2({\Omega_0};{\mathbb R}^9)}\
\| \text{Id}-a^T\|_{H^2 ({\Gamma_0};{\mathbb R}^9)}\ .
\end{align*}
From (\ref{ah2}) and (\ref{control2}), we infer that
\begin{align}
&\left\| \frac{a^T N}{|a^T N|^2}\cdot (\ S_{\eta(t,\cdot)}(w,r)\cdot a^T\ N\ ) (a^T-\text{Id}) N\right\|_{H^2({\Omega_0};{\mathbb R}^3)}\nonumber\\
&\qquad\qquad\qquad\qquad\qquad\qquad\le 
C \ \sqrt{T}\ \|S_{\eta(t,\cdot)}(w,r)\cdot a^T\|_{H^2({\Omega_0};{\mathbb R}^9)}
\ .
\label{g3}
\end{align}
Once again using that 
$H^2({\Omega_0};{\mathbb R}^3)$ is a Banach algebra we can conclude that
\begin{equation*}
\|S_{\eta(t,\cdot)}(w,r)\cdot a^T\|_{H^2({\Omega_0};{\mathbb R}^9)}\le
C\ \|a^T\cdot\nabla w+ (a^T\cdot\nabla w)^T\|_{H^2({\Omega_0};{\mathbb R}^9)}
\|a\|_{H^2({\Omega_0};{\mathbb R}^9)}
\ ,
\end{equation*}
and that
\begin{equation*}
\|S_{\eta(t,\cdot)}(w,r)\cdot a^T\|_{H^2({\Omega_0};{\mathbb R}^9)}\le
C\ (\|w\|_{H^3({\Omega_0};{\mathbb R}^3)}+\|r\|_{H^2({\Omega_0};{\mathbb R}^3)})\ \|a\|^2_{H^2({\Omega_0};{\mathbb R}^9)}
\ ,
\end{equation*}
which with (\ref{ah2}), provides us with the estimate
\begin{equation}
\label{g4}
\|S_{\eta(t,\cdot)}(w,r)\cdot a^T\|_{H^2({\Omega_0};{\mathbb R}^9)}\le
C\  (\|w\|_{H^3({\Omega_0};{\mathbb R}^3)}+\|r\|_{H^2({\Omega_0};{\mathbb R}^3)})
\ .
\end{equation}
Thus, from (\ref{g3}) and (\ref{g4}), we find that
\begin{align}
\label{g5}
&\left\|\frac{a^T N}{|a^T N|^2}\cdot  (\ S_{\eta(t,\cdot)}(w,r)\cdot a^T\ N\ )\ (a^T-\text{Id}) N\right\|_{H^2({\Omega_0};{\mathbb R}^3)}\nonumber\\
&\qquad\qquad\qquad\qquad\qquad\qquad\le C\  \sqrt{T}\ (\|w\|_{H^3({\Omega_0};{\mathbb R}^3)}+\|r\|_{H^2({\Omega_0};{\mathbb R}^3)})
\ .
\end{align}
In an identical fashion, we also have that
\begin{align}
&\left\|\frac{a^T N}{|a^T N|^2}\cdot(S_{\eta(t,\cdot)}(w,r)\cdot a^T\ N ) ((\text{Id}-a^T)\ N)\cdot (N)\ N \right\|_{H^2({\Omega_0};{\mathbb R}^3)}\nonumber\\
&\qquad\qquad\qquad\qquad\qquad\qquad  \le C\ \sqrt{T}\ (\|w\|_{H^3({\Omega_0};{\mathbb R}^3)}+\|r\|_{H^2({\Omega_0};{\mathbb R}^3)})
\ .
\label{g6}
\end{align}
From (\ref{g0}), (\ref{g5}) and (\ref{g6}), we may finally conclude that
\begin{equation}
\label{g}
\|\bar{g} \|_{L^2(0,T;H^{\frac{3}{2}}({\Gamma_0};{\mathbb R}^3))}\le C\ \sqrt{T}\ (\|w\|_{L^2(0,T; H^3({\Omega_0};{\mathbb R}^3))}+\|r\|_{L^2(0,T;H^2({\Omega_0};{\mathbb R}^3))})
\ .
\end{equation}

\medskip
\medskip

{\bf 6. Estimate of $g_t$:}\hfill\break
We next estimate the time derivatives of the terms in (\ref{g0}). 
First we have that
\begin{align*}
&\|{(S(w,r)\ \tilde{N}- S_{\eta(t,\cdot)}(w,r)\cdot a^T \tilde{N})}_t\|
_{H^{-\frac{1}{2}}({\Gamma_0};{\mathbb R}^3)} \nonumber \\
&\le 
C\  \|r_t (\text{Id}-a^T)\ \tilde{N}\|
_{H^{-\frac{1}{2}}({\Gamma_0};{\mathbb R}^3)}
+ C \|r a^T_t\ \tilde{N}\|
_{H^{-\frac{1}{2}}({\Gamma_0};{\mathbb R}^3)}\nonumber \\
& +C\ \| \nabla w_t \cdot (\text{Id}-a^T)\ \tilde{ N}\|
_{H^{-\frac{1}{2}}({\Gamma_0};{\mathbb R}^3)} 
+C\ \| \nabla w \cdot a^T_t \ \tilde{ N}\|
_{H^{-\frac{1}{2}}({\Gamma_0};{\mathbb R}^3)}\\
&+C\ \|(\text{Id}-a^T) \cdot \nabla w \cdot a_t^T\ \tilde{ N}\|
_{H^{-\frac{1}{2}}({\Gamma_0};{\mathbb R}^3)}
+ C\ \|(\text{Id}-a^T) \nabla w_t \cdot a^T\ \tilde{ N}\|
_{H^{-\frac{1}{2}}({\Gamma_0};{\mathbb R}^3)}\\
& +C\ \|a_t^T\cdot\nabla w\cdot a^T\ \tilde{ N}\|
_{H^{-\frac{1}{2}}({\Gamma_0};{\mathbb R}^3)}
+C\ \|(\nabla w_t)^T\cdot (\text{Id}-a )\ \tilde{N}\|
_{H^{-\frac{1}{2}}({\Gamma_0};{\mathbb R}^3)}\\
&+C\ \|(\nabla w)^T\cdot a_t\ \tilde{N}\|
_{H^{-\frac{1}{2}}({\Gamma_0};{\mathbb R}^3)}
+ C\ \|(\nabla w_t)^T\cdot a\cdot (\text{Id}-a^T )\ \tilde{N}\|
_{H^{-\frac{1}{2}}({\Gamma_0};{\mathbb R}^3)}\\
&+ C\ \|(\nabla w)^T\cdot a_t\cdot (\text{Id}-a^T )\ \tilde{N}\|
_{H^{-\frac{1}{2}}({\Gamma_0};{\mathbb R}^3)} 
+  C\ \|(\nabla w)^T\cdot a\cdot a^T_t\ \tilde{N}\|_
{H^{-\frac{1}{2}}({\Gamma_0};{\mathbb R}^3)}\ .
\end{align*}
From (\ref{dual}), we find that
\begin{align*}
&\|{(S(w,r)\ \tilde{N}- S_{\eta(t,\cdot)}(w,r)\cdot a^T\ \tilde{N})}_t\|_{H^{-\frac{1}{2}}({\Gamma_0};{\mathbb R}^3)} \nonumber \\
&\le C\  \|r_t\|_{H^{-\frac{1}{2}}({\Gamma_0};{\mathbb R})}\ \|\text{Id}-a \|_{W^{1,4}({\Omega_0};{\mathbb R}^9)}+
C\  \| a_t \|_{H^{-\frac{1}{2}}({\Gamma_0};{\mathbb R}^9)}\ 
 \|r \|_{W^{1,4}({\Omega_0};{\mathbb R})}\nonumber \\
& +C\ \|\nabla w_t\|_{H^{-\frac{1}{2}}({\Gamma_0};{\mathbb R}^9)}\ \|\text{Id}-a\|_{W^{1,4}({\Omega_0};{\mathbb R}^9)} +C\ \|a_t\|_{H^{-\frac{1}{2}}({\Gamma_0};{\mathbb R}^9)}\ \|\nabla w\|_{W^{1,4}({\Omega_0};{\mathbb R}^9)}\\
&+C\ \|a_t\|_{H^{-\frac{1}{2}}({\Gamma_0};{\mathbb R}^9)}\ \|\nabla w\|_{W^{1,4}({\Omega_0};{\mathbb R}^9)}\ \|\text{Id}-a \|_{W^{1,4}({\Omega_0};{\mathbb R}^9)}\\
&+C\ \|\nabla w_t\|_{H^{-\frac{1}{2}}({\Gamma_0};{\mathbb R}^9)}\ \|a\|_{W^{1,4}({\Omega_0};{\mathbb R}^9)}\ \|\text{Id}-a \|_{W^{1,4}({\Omega_0};{\mathbb R}^9)}\\
&+C\ \|a_t\|_{H^{-\frac{1}{2}}({\Gamma_0};{\mathbb R}^9)}\ \|a\|_{W^{1,4}({\Omega_0};{\mathbb R}^9)}\ \|\nabla w\|_{W^{1,4}({\Omega_0};{\mathbb R}^9)}
\ ,
\end{align*}
and with the Sobolev embedding theorem and the Gagliardo-Nirenberg inequality,
we obtain
\begin{align*}
\|(S(w,r)\ \tilde{N}&- S_{\eta(t,\cdot)}(w,r)\cdot a^T\ \tilde{N})_t\|_{H^{-\frac{1}{2}}({\Gamma_0};{\mathbb R}^3)} \nonumber \\
\le &C\  \|r_t\|_{H^{-\frac{1}{2}}({\Gamma_0};{\mathbb R})}\ \|\text{Id}-a \|_{H^{2}({\Omega_0};{\mathbb R}^9)}\\
& +C\   \| a_t \|_{H^1({\Omega_0};{\mathbb R}^9)}
 \|r \|^{\frac{1}{8}}_{L^2({\Omega_0};{\mathbb R})}\  \|r \|^{\frac{7}{8}}_{H^2({\Omega_0};{\mathbb R})}\nonumber \\
& +C\ \|\nabla w_t\|_{H^{-\frac{1}{2}}({\Gamma_0};{\mathbb R}^9)}\ \|\text{Id}-a\|_{H^2({\Omega_0};{\mathbb R}^9)}\\
& +C\ \| a_t \|_{H^1({\Omega_0};{\mathbb R}^9)}\ \|\nabla w\|^{\frac{1}{4}}_{H^1({\Omega_0};{\mathbb R}^9)}\ \|\nabla w\|^{\frac{3}{4}}_{H^2({\Omega_0};{\mathbb R}^9)}\\
&+C\   \| a_t \|_{H^1({\Omega_0};{\mathbb R}^9)}\ \| w\|_{H^3({\Omega_0};{\mathbb R}^3)}\ \|\text{Id}-a \|_{W^{1,4}({\Omega_0};{\mathbb R}^9)}\\
&+C\ \|\nabla w_t\|_{H^{-\frac{1}{2}}({\Gamma_0};{\mathbb R}^9)}\ \|a\|_{H^2({\Omega_0};{\mathbb R}^9)}\ \|\text{Id}-a \|_{H^2({\Omega_0};{\mathbb R}^9)}\\
&+C\ \|a_t\|_{H^1({\Omega_0};{\mathbb R}^9)}\ \|a\|_{H^2({\Omega_0};{\mathbb R}^9)}\ \|\nabla w\|^{\frac{1}{4}}_{H^1({\Omega_0};{\mathbb R}^9)}\ \|\nabla w\|^{\frac{3}{4}}_{H^2({\Omega_0};{\mathbb R}^9)}\ .
\end{align*}
From (\ref{ah2}), (\ref{timeah1}), and the fact that $r(0)=0$ and $w(0)=0$, we then infer that
\begin{align*}
\|(S(w,r)\ \tilde{N}&- S_{\eta(t,\cdot)}(w,r)\cdot a^T\ \tilde{N})_t\|_{H^{-\frac{1}{2}}({\Gamma_0};{\mathbb R}^3)} \nonumber \\
& \le C\  \sqrt{T}\  \|r_t\|_{H^{-\frac{1}{2}}({\Gamma_0};{\mathbb R}^3)}\ 
+C\  \sqrt{T}\ \|\nabla w_t\|_{H^{-\frac{1}{2}}({\Gamma_0};{\mathbb R}^9)}\\
&+C  \ T^{\frac{1}{16}}\  
 \|r_t \|^{\frac{1}{8}}_{L^2(0,T; L^2({\Omega_0};{\mathbb R}))}\  \|r \|^{\frac{7}{8}}_{H^2({\Omega_0};{\mathbb R})}\nonumber \\
& +C\  T^{\frac{1}{8}}\|\nabla w_t\|^{\frac{1}{4}}_{L^2(0,T;H^1({\Omega_0};{\mathbb R}^9))}\ \|\nabla w\|^{\frac{3}{4}}_{H^2({\Omega_0};{\mathbb R}^9)}\\
&+C\ \sqrt{T}\ \| w\|_{H^3({\Omega_0};{\mathbb R}^3)} 
+C\ \sqrt{T}\ \|\nabla w_t\|_{H^{-\frac{1}{2}}({\Gamma_0};{\mathbb R}^9)}\\ 
&+C\  T^{\frac{1}{8}}\ \|\nabla w_t\|^{\frac{1}{4}}_{L^2(0,T;H^1({\Omega_0};{\mathbb R}^9))}\ \|\nabla w\|^{\frac{3}{4}}_{H^2({\Omega_0};{\mathbb R}^9)}\ ,
\end{align*}
and from Young's inequality,
\begin{align*}
&\|{(S(w,r)\ \tilde{N}- S_{\eta(t,\cdot)}(w,r)\cdot a^T)\ \tilde{N})}_t\|_{H^{-\frac{1}{2}}({\Gamma_0};{\mathbb R}^3)} \nonumber \\
&\le \ C\  \sqrt{T}\  \|r_t\|_{H^{-\frac{1}{2}}({\Gamma_0};{\mathbb R})}\ 
+C  \sqrt{T}\ \|\nabla w_t\|_{H^{-\frac{1}{2}}({\Gamma_0};{\mathbb R}^9)}\\
&+C \ T^{\frac{1}{16}}\  
 \left( \|r_t \|_{L^2(0,T; L^2({\Omega_0};{\mathbb R}))}\ + \|r \|_{H^2({\Omega_0};{\mathbb R})}\right)\nonumber \\
& +C \ T^{\frac{1}{8}}\left(\|\nabla w_t\|_{L^2(0,T;H^1({\Omega_0};{\mathbb R}^9))}\ + \|\nabla w\|_{H^2({\Omega_0};{\mathbb R}^9)}\right)\\
&+C\  \sqrt{T}\ \| w\|_{H^3({\Omega_0};{\mathbb R}^3)} 
+C\  \sqrt{T}\ \|\nabla w_t\|_{H^{-\frac{1}{2}}({\Gamma_0};{\mathbb R}^9)}\\ 
&+C\  T^{\frac{1}{8}}\ \left(\|\nabla w_t\|_{L^2(0,T;H^1({\Omega_0};{\mathbb R}^9))}\ + \|\nabla w\|_{H^2({\Omega_0};{\mathbb R}^9)}\right)\ ,
\end{align*}
Then, we see that there exists $T'_1\in (0,T_0)$ such that if $0<T<T'_1$, we 
have that
\begin{equation}
\label{gt1}
\|{(S(w,r)\ \tilde{N}- S_{\eta(t,\cdot)}(w,r)\cdot\ a^T\ \tilde{N})}_t\|_{L^2(0,T;H^{-\frac{1}{2}}({\Gamma_0};{\mathbb R}^3))}\le C\  T^{\frac{1}{16}}\ \|(w,r)\|_{X_T}\ .
\end{equation}

To estimate the remaining term in (\ref{g0}), we use (\ref{pio}) to obtain the
estimate
\begin{align*}
&\left\|\ \left(\Pi_0 \left(\ \frac{a^T N}{|a^T N|^2}\cdot (\  S_{\eta(t,\cdot)}(w,r)\cdot a^T\ N\ )\ a^T N\ \right)\right)_t\ \right\|_{H^{-\frac{1}{2}}({\Gamma_0};{\mathbb R}^3)} \\
& \qquad\qquad\le\  \left\|\left(\ \frac{a^T N}{|a^T N|^2}\cdot (\ S_{\eta(t,\cdot)}(w,r)\cdot a^T\ N)\ \right)_t\ (a^T-\text{Id})\ N)\right\|_{H^{-\frac{1}{2}}({\Gamma_0};{\mathbb R}^3)}\\
 & \qquad\qquad + \left\|\frac{a^T N}{|a^T N|^2}\cdot (\ S_{\eta(t,\cdot)}(w,r)\cdot a^T\ N\ )\ a^T_t N\right\|_{H^{-\frac{1}{2}}({\Gamma_0};{\mathbb R}^3)}\\
& \qquad\qquad +\left\|\left(\ \frac{a^T N}{|a^T N|^2}\cdot (\ S_{\eta(t,\cdot)}(w,r)\cdot a^T\ N)\ \right)_t ((\text{Id}-a^T)\ N\cdot N) N\right\|_{H^{-\frac{1}{2}}({\Gamma_0};{\mathbb R}^3)}\\
& \qquad\qquad + \left\|\frac{a^T N}{|a^T N|^2}\cdot(\ S_{\eta(t,\cdot)}(w,r)\cdot a^T\ N)\ (a^T_t\ N)\cdot( N)\ N\right\|_{H^{-\frac{1}{2}}({\Gamma_0};{\mathbb R}^3)}\ .
\end{align*}
From (\ref{dual}), we then find that
\begin{align*}
&\left\|\ \left(\Pi_0 \left(\ \frac{a^T N}{|a^T N|^2}\cdot (\  S_{\eta(t,\cdot)}(w,r)\cdot a^T\ N\ )\ a^T N\ \right)\right)_t\ \right\|_{H^{-\frac{1}{2}}({\Gamma_0};{\mathbb R}^3)} \\
&\qquad\le\ C\ \left\|\left(\ \frac{a^T N}{|a^T N|^2}\cdot (\ S_{\eta(t,\cdot)}(w,r)\cdot a^T\ N)\ \right)_t\right\|_{H^{-\frac{1}{2}}({\Gamma_0};{\mathbb R}^3)}\ \| a-\text{Id} \|_{W^{1,4}({\Omega_0};{\mathbb R}^9)}\\
& \qquad +C\  \left\|\frac{a^T N}{|a^T N|^2}\cdot (\ S_{\eta(t,\cdot)}(w,r)\cdot a^T\ N\ )\right\|_{W^{1,4}({\Omega_0};{\mathbb R})}\  \|a_t\ N\|_{H^{-\frac{1}{2}}({\Gamma_0};{\mathbb R}^3)}\ .
\end{align*}
Using (\ref{dual}) for the first term on the right-hand-side of the 
previous inequality together with the fact that 
$W^{1,4}({\Omega_0};{\mathbb R}^3)$ is a Banach algebra, provides us with
the inequality
\begin{align*}
&\left\|\ \left(\Pi_0 \left(\ \frac{a^T N}{|a^T N|^2}\cdot (\  S_{\eta(t,\cdot)}(w,r)\cdot a^T\ N\ )\ a^T N\ \right)\right)_t\ \right\|_{H^{-\frac{1}{2}}({\Gamma_0};{\mathbb R}^3)} \left\|\frac{1}{|a^T N|^2}\right\|^{-1}_{W^{1,4}({\Omega_0};{\mathbb R})}\\
&\le C\ \|(S_{\eta(t,\cdot)}(w,r))_t\|_{H^{-\frac{1}{2}}({\Gamma_0};{\mathbb R}^9)} \|a \|^2_{W^{1,4}({\Omega_0};{\mathbb R}^9)}\| a-\text{Id} \|_{W^{1,4}({\Omega_0};{\mathbb R}^9)}\\
&+ C\ \|S_{\eta(t,\cdot)}(w,r)\|_{W^{1,4}({\Omega_0};{\mathbb R}^9)}\ \|a \|_{W^{1,4}({\Gamma_0};{\mathbb R}^9)}\ \|a_t\|_{H^{-\frac{1}{2}}({\Gamma_0};{\mathbb R}^9)}\
\| a-\text{Id} \|_{W^{1,4}({\Omega_0};{\mathbb R}^9)}\\
&+ C\ \|(a^T N)\cdot (\ S_{\eta(t,\cdot)}(w,r)\cdot a^T\ N\ )\|_{W^{1,4}({\Omega_0};{\mathbb R})}\  \|a^T_t\ N\|_{H^{-\frac{1}{2}}({\Gamma_0};{\mathbb R}^3)}\ .
\end{align*}
By taking into account (\ref{dual}), we deduce that
\begin{align*}
&\left\|\ \left(\Pi_0 \left(\ \frac{a^T N}{|a^T N|^2}\cdot 
(\  S_{\eta(t,\cdot)}(w,r)\cdot a^T\ N\ )\ a^T N\ \right)\right)_t\ \right\|
_{H^{-\frac{1}{2}}({\Gamma_0};{\mathbb R}^3)} 
\left\|\frac{1}{|a^T N|^2}\right\|^{-1}_{W^{1,4}({\Omega_0};{\mathbb R})}\\
&\qquad
\le 
C\ \|\nabla w_t\|_{H^{-\frac{1}{2}}({\Gamma_0};{\mathbb R}^9)}\ \|a \|^3
_{W^{1,4}({\Omega_0};{\mathbb R}^9)}\ 
\| a-\text{Id} \|_{W^{1,4}({\Omega_0};{\mathbb R}^9)}\\
&\qquad 
+ C\ \|r_t\|_{H^{-\frac{1}{2}}({\Gamma_0};{\mathbb R})}\ \|a \|^2
_{W^{1,4}({\Omega_0};{\mathbb R}^9)}\ 
\| a-\text{Id} \|_{W^{1,4}({\Omega_0};{\mathbb R}^9)}\\
&\qquad
+ C\ \|a_t\|_{H^{-\frac{1}{2}}({\Gamma_0};{\mathbb R}^9)}\ \|\nabla w\|
_{W^{1,4}({\Omega_0};{\mathbb R}^9)}
\ \|a \|^2_{W^{1,4}({\Gamma_0};{\mathbb R}^9)}\
\| a-\text{Id} \|_{W^{1,4}({\Omega_0};{\mathbb R}^3)}\\
&\qquad
+ C\ \|a_t\|_{H^{-\frac{1}{2}}({\Gamma_0};{\mathbb R}^9)}\ \|r\|
_{W^{1,4}({\Omega_0};{\mathbb R})}\ 
\|a \|_{W^{1,4}({\Gamma_0};{\mathbb R}^9)}\
\| a-\text{Id} \|_{W^{1,4}({\Omega_0};{\mathbb R}^9)}\\
&\qquad+ C\ \|\nabla w\|_{W^{1,4}({\Omega_0};{\mathbb R}^9)}\ \|a \|^2_{W^{1,4}({\Omega_0};{\mathbb R}^9)}\ \|a_t\|_{H^{-\frac{1}{2}}({\Gamma_0};{\mathbb R}^9)}\
\| a-\text{Id} \|_{W^{1,4}({\Omega_0};{\mathbb R}^3)}\\
&\qquad+ C\ \|r\|_{W^{1,4}({\Omega_0};{\mathbb R})}\ \|a \|_{W^{1,4}({\Omega_0};{\mathbb R}^9)}\ \|a_t\|_{H^{-\frac{1}{2}}({\Gamma_0};{\mathbb R}^9)}\
\| a-\text{Id} \|_{W^{1,4}({\Omega_0};{\mathbb R}^9)}\\
&\qquad+ C\ \|a\|^3_{W^{1,4}({\Omega_0};{\mathbb R}^9)}\ \|\nabla w\|_{W^{1,4}({\Omega_0};{\mathbb R}^9)}\ \|a_t \|_{H^{-\frac{1}{2}}({\Gamma_0};{\mathbb R}^9)}\\
&\qquad+ C\ \|a\|^2_{W^{1,4}({\Omega_0};{\mathbb R}^9)}\ \|r\|_{W^{1,4}({\Omega_0};{\mathbb R})}\ \|a_t \|_{H^{-\frac{1}{2}}({\Gamma_0};{\mathbb R}^9)}\ ,
\end{align*}
and by the Sobolev embedding theorem and the Gagliardo-Nirenberg inequality,
\begin{align*}
&\left\|\ \left(\Pi_0 \left(\ \frac{a^T N}{|a^T N|^2}\cdot (\  S_{\eta(t,\cdot)}(w,r)\cdot a^T\ N\ )\ a^T N\ \right)\right)_t\ \right\|_{H^{-\frac{1}{2}}({\Gamma_0};{\mathbb R}^3)} \left\|\frac{1}{|a^T N|^2}\right\|^{-1}_{W^{1,4}({\Omega_0};{\mathbb R})}\\
&\qquad \le\ C\ \|\nabla w_t\|_{H^{-\frac{1}{2}}({\Gamma_0};{\mathbb R}^9)}\ \|a \|^3_{H^{2}({\Omega_0};{\mathbb R}^9)}\ \| a-\text{Id} \|_{H^{2}({\Omega_0};{\mathbb R}^9)}\\
&\qquad+ C\ \|r_t\|_{H^{-\frac{1}{2}}({\Gamma_0};{\mathbb R})}\ \|a \|^2_{H^{2}({\Omega_0};{\mathbb R}^9)}\  
\| a-\text{Id} \|_{H^{2}({\Omega_0};{\mathbb R}^9)}\\
&\qquad+ C\ \|a_t\|_{H^1({\Omega_0};{\mathbb R}^9)}\ \|w\|_{H^{3}({\Omega_0};{\mathbb R}^3)}\ \|a \|^2_{H^{2}({\Omega_0};{\mathbb R}^3)}\ 
\| a-\text{Id} \|_{H^{2}({\Omega_0};{\mathbb R}^9)}\\
&\qquad+ C\ \|a_t\|_{H^1({\Omega_0};{\mathbb R}^9)}\ \|r\|_{H^{2}({\Omega_0};{\mathbb R})}\ \|a \|_{H^{2}({\Omega_0};{\mathbb R}^9)}\ 
\| a-\text{Id} \|_{H^{2}({\Omega_0};{\mathbb R}^9)}\\
&\qquad+ C\ \| w\|_{H^{3}({\Omega_0};{\mathbb R}^3))}\ \|a \|^2_{H^2 ({\Omega_0};{\mathbb R}^9)}\ \|a_t\|_{H^1 ({\Omega_0};{\mathbb R}^9)}\
\| a-\text{Id} \|_{H^2({\Omega_0};{\mathbb R}^9)}\\
&\qquad+ C\ \| r\|_{H^{2}({\Omega_0};{\mathbb R}))}\ \|a \|_{H^{2}({\Omega_0};{\mathbb R}^9))}\ \|a_t\|_{H^1 ({\Omega_0};{\mathbb R}^9)}\
\| a-\text{Id} \|_{H^2({\Omega_0};{\mathbb R}^9)}\\
&\qquad+ C\ \|a\|^3_{H^2({\Omega_0};{\mathbb R}^9)}\ \| w\|^{\frac{1}{8}}_{H^1 ({\Omega_0};{\mathbb R}^3)}\  \| w\|^{\frac{7}{8}}_{H^3 ({\Omega_0};{\mathbb R}^3)} \ \|a_t \|_{H^1({\Omega_0};{\mathbb R}^9)}\\
&\qquad+ C\ \|a\|^2_{H^2({\Omega_0};{\mathbb R}^9)}\ \| r\|^{\frac{1}{8}}_{L^2 ({\Omega_0};{\mathbb R})}\  \| r\|^{\frac{7}{8}}_{H^2 ({\Omega_0};{\mathbb R})} \ \|a_t \|_{H^1({\Omega_0};{\mathbb R}^9)}\ .
\end{align*}
Finally, from (\ref{control2}), (\ref{ah2}) and (\ref{timeah1}), we infer that
\begin{align*}
&\left\|\ \left(\Pi_0 \left(\ \frac{a^T N}{|a^T N|^2}\cdot 
(\  S_{\eta(t,\cdot)}(w,r)\cdot a^T\ N\ )\ 
a^T N\ \right)\right)_t\ \right\|
_{H^{-\frac{1}{2}}({\Gamma_0};{\mathbb R}^3)} \\
&\le 
C\  \sqrt{T}\ \left(\|\nabla w_t\|
_{H^{-\frac{1}{2}}({\Gamma_0};{\mathbb R}^9)}+ 
\|r_t\|_{H^{-\frac{1}{2}}({\Gamma_0};{\mathbb R})}+ 
\|w\|_{H^{3}({\Omega_0};{\mathbb R}^3)}  \right)\\
& + 
C\  \sqrt{T}\  
\left( \|r\|_{H^{2}({\Omega_0};{\mathbb R})}+  
\| w\|_{H^{3}({\Omega_0};{\mathbb R}^3)}+\| r\|
_{H^{2}({\Omega_0};{\mathbb R})}\right)\\
& 
+ C\ T^{\frac{1}{8}} \left( \| w_t\|^{\frac{1}{8}}
_{L^2(0,T; H^1 ({\Omega_0};{\mathbb R}^3))}\  
\| w\|^{\frac{7}{8}}_{H^3 ({\Omega_0};{\mathbb R}^3)}+ 
\| r_t\|^{\frac{1}{8}}_{L^2(0,T; L^2 ({\Omega_0};{\mathbb R}^3))}\  
\| r\|^{\frac{7}{8}}_{H^2 ({\Omega_0};{\mathbb R}^3)}\right) \ .
\end{align*}
We may thus conclude that
there exists $T'_2\in (0,T_0)$ such that for any $T\in (0, T'_2)$ we have
\begin{equation}
\label{gt2}
\left\|\ \left(\Pi_0 \left(\ \frac{a^T N}{|a^T N|^2}\cdot (\  S_{\eta(t,\cdot)}(w,r)\cdot a^T\ N\ )\ a^T N\ \right)\right)_t\ \right\|_{H^{-\frac{1}{2}}({\Gamma_0};{\mathbb R}^3)}
 \le C\ T^{\frac{1}{8}}\ \|(w,r)\|_{X_T}\ .
\end{equation}
From (\ref{gt1}) and (\ref{gt2}), we have that for all $0<T<T_{\bar{g}_t} 
=\min(T'_1,T'_2)$,
\begin{equation}
\label{gt}
\|{\bar{g}}_t\|_{L^2(0,T;H^{-\frac{1}{2}}({\Gamma_0};{\mathbb R}^3))}\le 
C\  T^{\frac{1}{16}}\ \|(w,r)\|_{X_T}\ .
\end{equation}

{\bf 7. Estimates of $g_1$, $g_2$,$\partial_t g_1$, and $\partial_t g_2$:}\hfill\break
Since $\bar{g}_2=\bar{g}\cdot N$ and $\bar{g}_1=\bar{g}- (\bar{g}\cdot N)\ N$, we infer from (\ref{g}), (\ref{gt}), and (\ref{dual}) that for any
$0<T<\min(T_{\bar{g}}, T_{\bar{g}_t})$,
\begin{align}
&\|\bar{g}_1 (w,r)\|_{L^2(0,T;H^{\frac{3}{2}}({\Gamma_0};{\mathbb R}^3))}+
\|{\bar{g}_1 (w,r)}_t\|_{L^2(0,T;H^{-\frac{1}{2}}({\Gamma_0};{\mathbb R}^3))}\nonumber\\
&+ \|\bar{g}_2 (w,r)\|_{L^2(0,T;H^{\frac{3}{2}}({\Gamma_0};{\mathbb R}))}+
\|{\bar{g}_2 (w,r)}_t\|_{L^2(0,T;H^{-\frac{1}{2}}({\Gamma_0};{\mathbb R}))}\nonumber\\
&\qquad\qquad\qquad\qquad\qquad\qquad\qquad\qquad\qquad\qquad\qquad\le C\  T^{\frac{1}{16}}\ \|(w,r)\|_{X_T} \ .
\label{g1g2}
\end{align}

{\bf 8. Estimate of $(\overline{B})_t$:}\hfill\break
As we defined in (\ref{laplacian}), 
\begin{equation*}
 \triangle_{g}w(t,x)  = g^{\alpha \beta}
\left(\frac{\partial^2 w} 
{\partial y^\alpha \partial y^\beta}(t,x)  - 
\Gamma^\gamma_{\alpha\beta}\frac{\partial w }{\partial y^\gamma}(t,x)\right)\,.
\end{equation*}
Consequently,
\begin{align}
\label{bt0}
\overline{B}_t=&(g^{\alpha \beta}(t,x)-g^{\alpha \beta}(0,x)) 
\frac{\partial^2 w} 
{\partial y^\alpha \partial y^\beta}(t,x)\nonumber\\
&  - (g^{\alpha \beta}(t,x) \Gamma^\gamma_{\alpha\beta}(t,x)-
g^{\alpha \beta}(0,x) \Gamma^\gamma_{\alpha\beta}(t,x)) 
\frac{\partial w }{\partial y^\gamma}(t,x)\nonumber\\
&+ \partial_t(g^{\alpha \beta})(t,x) \int_0^t \frac{\partial^2 w} 
{\partial y^\alpha \partial y^\beta}(t',x) dt'\nonumber\\
&  - \partial_t(g^{\alpha \beta}(t,x) \Gamma^\gamma_{\alpha\beta})(t,x) 
\int_0^t\frac{\partial w}{\partial y^\gamma}(t',x)dt'
 \ .
\end{align}
It follows that
\begin{align*}
&\|(g^{\alpha \beta}(t,\cdot)-g^{\alpha \beta}(0,\cdot)) 
\frac{\partial^2 w} 
{\partial y^\alpha \partial y^\beta}(t,\cdot) \|_{ H^{\frac{1}{2}}
(\Gamma_0;{\mathbb R})^3}\nonumber\\ 
&\qquad\le C  \|\ |\nabla\eta|^2 |\nabla\nabla\nabla w |\ 
\|_{L^2(\Omega_0;{\mathbb R}^3)}+
C  \|\ |\nabla\nabla\eta||\nabla\eta| |\nabla\nabla w|\ 
\|_{L^2(\Omega_0;{\mathbb R}^3)}\ ,
\end{align*}
and from (\ref{etah3}), and the Cauchy-Schwarz inequality,
\begin{align*}
&\|(g^{\alpha \beta}(t,\cdot)-g^{\alpha \beta}(0,\cdot)) 
\frac{\partial^2 w} 
{\partial y^\alpha \partial y^\beta}(t,\cdot)) \|_{ H^{\frac{1}{2}}
(\Gamma_0;{\mathbb R}^3)}\nonumber\\ 
&\qquad\le C \ \sqrt{T}\ \|w\|_{H^3(\Omega_0;{\mathbb R}^3)} +
 C\ \|\ |\nabla\nabla\eta||\nabla\eta|\ \|_{L^4(\Omega_0;{\mathbb R})}\ \|\ |\nabla\nabla w |\ \|_{L^4(\Omega_0;{\mathbb R})}\ .
\end{align*}
From the Sobolev embedding theorem and the Gagliardo-Nirenberg inequality, 
we then have that
\begin{align*}
&\|(g^{\alpha \beta}(t,\cdot)-g^{\alpha \beta}(0,\cdot)) 
\frac{\partial^2 w} 
{\partial y^\alpha \partial y^\beta}(t,\cdot) \|_{ H^{\frac{1}{2}}
(\Gamma_0;{\mathbb R}^3)}\nonumber\\ 
&\qquad\le C  \sqrt{T}\ \|w\|_{H^3(\Omega_0;{\mathbb R}^3)} +
 C \| \eta \|^2_{H^3(\Omega_0;{\mathbb R}^3)}\|  w\|^{\frac{1}{8}}_{H^1(\Omega_0;{\mathbb R}^3)} \| w \|^{\frac{7}{8}}_{H^3(\Omega_0;{\mathbb R}^3)}\ .
\end{align*}
Since $w(0)=0$, we use (\ref{etah3}) to deduce that
\begin{align*}
&\|(g^{\alpha \beta}(t,\cdot)-g^{\alpha \beta}(0,\cdot)) 
\frac{\partial^2 w} 
{\partial y^\alpha \partial y^\beta}(t,\cdot) \|_{ H^{\frac{1}{2}}
(\Gamma_0;{\mathbb R})}\nonumber\\ 
&\le C\  \sqrt{T}\ \|w\|_{H^3(\Omega_0;{\mathbb R}^3)} +
 C \ T^{\frac{1}{16}}\| w_t\|^{\frac{1}{8}}_{L^2(0,T; H^1(\Omega_0;{\mathbb R}^3))} \| w\|^{\frac{7}{8}}_{H^3(\Omega_0;{\mathbb R}^3)}\ ,
\end{align*}
and from Young's inequality,
\begin{align}
\label{bt1}
\|(g^{\alpha \beta}&-g^{\alpha \beta}(0,\cdot)) 
\frac{\partial^2 w} 
{\partial y^\alpha \partial y^\beta}\|_{L^2(0,T; H^{\frac{1}{2}}
(\Gamma^i_0;{\mathbb R}))}\nonumber\\  
&\le   C\ \sqrt{T} \ \|w\|_{L^2(0,T;H^3(\Omega_0;{\mathbb R}^3))} \nonumber\\ &+
 C\ T^{\frac{1}{16}}\left(\| w_t\ |\|_{L^2(0,T; H^1(\Omega_0;{\mathbb R}^3))}+ \| w\|_{L^2(0,T;H^3(\Omega_0;{\mathbb R}^3))}\right) \ .
\end{align}

Similarly,
\begin{align*}
&\|(g^{\alpha \beta}(t,\cdot) \Gamma^\gamma_{\alpha\beta}(t,\cdot)
-g^{\alpha \beta}(0,\cdot) \Gamma^\gamma_{\alpha\beta}(0,\cdot))
 \frac{\partial w} 
{\partial y^\gamma}(t,\cdot) \|_{ H^{\frac{1}{2}}
(\Gamma_0;{\mathbb R}^3)}\nonumber\\ 
&\qquad\qquad\le C  \|\ (g^{\alpha \beta}(t,\cdot) \Gamma^\gamma_{\alpha\beta}(t,\cdot)-g^{\alpha \beta} (0,\cdot) \Gamma^\gamma_{\alpha\beta}(0,\cdot))\ |\nabla\nabla w |\ \|_{L^2(\Omega_0;{\mathbb R})}\nonumber\\
&\qquad\qquad +
C  \|\ |\nabla(g^{\alpha \beta}(t,\cdot) \Gamma^\gamma_{\alpha\beta}(t,\cdot)-g^{\alpha \beta}(0,\cdot) \Gamma^\gamma_{\alpha\beta}(0,\cdot))|\ |\nabla w |\ \|_{L^2(\Omega_0;{\mathbb R})}\ ,
\end{align*}
and thus,
\begin{align*}
&\|(g^{\alpha \beta}(t,\cdot) \Gamma^\gamma_{\alpha\beta}(t,\cdot)
-g^{\alpha \beta}(0,\cdot) \Gamma^\gamma_{\alpha\beta}(0,\cdot)) 
\frac{\partial w} 
{\partial y^\gamma}(t,\cdot) \|_{ H^{\frac{1}{2}}
(\Gamma_0;{\mathbb R}^3)}\nonumber\\ 
&\le C\  \|\ g^{\alpha \beta}(t,\cdot) \Gamma^\gamma_{\alpha\beta}(t,\cdot)-g^{\alpha \beta} (0,\cdot) \Gamma^\gamma_{\alpha\beta}(0,\cdot)\|_{L^\infty (0,T;L^4(\Omega_0;{\mathbb R}))} \| w\|_{W^{2,4} (\Omega_0;{\mathbb R}^3)}\nonumber\\
&+
C\  \|\ \nabla(g^{ \alpha \beta }(t,\cdot) 
\Gamma^\gamma_{\alpha\beta}(t,\cdot)
-g^{ \alpha \beta }(0,\cdot) 
\Gamma^\gamma_{\alpha\beta}(0,\cdot))\|
_{L^\infty (0,T;L^2(\Omega_0;{\mathbb R}))} \|w \|
_{H^3(\Omega_0;{\mathbb R}^3)}\ .
\end{align*}
From (\ref{etah3}), we see that
\begin{align*}
\|(g^{\alpha \beta}(t,\cdot) \Gamma^\gamma_{\alpha\beta}(t,\cdot)-
g^{\alpha \beta}(0,\cdot) \Gamma^\gamma_{\alpha\beta}(0,\cdot)) 
\frac{\partial w} 
{\partial y^\gamma}(t,\cdot)& \|_{ H^{\frac{1}{2}}
(\Gamma^i_0;{\mathbb R}^3)}\nonumber\\
&\le C \sqrt{T} \ \|w\|_{H^3(\Omega_0;{\mathbb R}^3)}\ .
\end{align*}
Consequently,
\begin{align}
\label{bt2}
\|(g^{\alpha \beta}(t,\cdot) \Gamma^\gamma_{\alpha\beta}(t,\cdot)
-g^{\alpha \beta}(0,\cdot) \Gamma^\gamma_{\alpha\beta}(0,\cdot))& 
\frac{\partial w} 
{\partial y^\gamma}(t,\cdot) \|_{L^2(0,T; H^{\frac{1}{2}}
(\Gamma_0;{\mathbb R}^3)}\nonumber\\
\le &C\ \sqrt{T}\  \|w\|_{L^2(0,T;H^3(\Omega_0;{\mathbb R}^3))}\ .
\end{align}

Next, we have that
\begin{align*}
\|\partial_tg^{\alpha \beta} (t,\cdot) \int_0^t \frac{\partial^2 w} 
{\partial y^\alpha \partial y^\beta}(t',\cdot)& dt'\|_{ H^{\frac{1}{2}}
(\Gamma_0;{\mathbb R}^3)}\nonumber\\
&\le C  \left\|\ |\nabla\nabla v| |\nabla\eta| \left|\int_0^t \nabla\nabla w (t',\cdot) dt'\right|\ \right\|_{L^2(\Omega_0;{\mathbb R}^3)}\nonumber\\
&+ C \left\|\ |\nabla v| |\nabla\nabla\eta| \left|\int_0^t \nabla\nabla w (t',\cdot) dt'\right|\ \right\|_{L^2(\Omega_0;{\mathbb R}^3)}\nonumber\\
&+  C \left\|\ |\nabla v| |\nabla\eta| \left|\int_0^t \nabla\nabla\nabla w (t',\cdot) dt'\right|\ \right\|_{L^2(\Omega_0;{\mathbb R}^3)}\ .
\end{align*}

From (\ref{etah3}), and the Cauchy-Schwarz inequality, we then find that
\begin{align*}
\|\partial_t g^{\alpha \beta}(t,\cdot)& \int_0^t \frac{\partial^2 w} 
{\partial y^\alpha \partial y^\beta}(t',\cdot) dt' \|_{ H^{\frac{1}{2}}
(\Gamma_0;{\mathbb R}^3)}\nonumber\\
&\le C\ \| \nabla\nabla v\|_{L^4(\Omega_0;{\mathbb R}^{27})} 
\ \left\|\int_0^t \nabla\nabla w (t',\cdot) dt'\ \right\|_{L^4(\Omega_0;{\mathbb R}^{27})}\nonumber\\
&+C\ \| \nabla v\|_{L^\infty(\Omega_0;{\mathbb R}^9)} 
\ \left\|\int_0^t \nabla\nabla w (t',\cdot) dt'\ \right\|_{L^4(\Omega_0;{\mathbb R}^{27})}\nonumber\\
&+C\ \| \nabla v\|_{L^\infty(\Omega_0;{\mathbb R}^{9})} 
\ \left\|\int_0^t \nabla\nabla\nabla w (t',\cdot) dt'\ \right\|_{L^2(\Omega_0;{\mathbb R}^{27})}\ .
\end{align*}
By the Sobolev embedding theorem and Jensen's inequality, we deduce from the 
previous inequality that
\begin{align*}
\|\partial_tg^{\alpha \beta} (t,\cdot) &\int_0^t\frac{\partial^2 w} 
{\partial y^\alpha \partial y^\beta}(t',\cdot) dt'\|_{ H^{\frac{1}{2}}
(\Gamma_0;{\mathbb R})}\nonumber\\
&\le C\  
\| v\|_{H^3(\Omega_0;{\mathbb R}^{3})}
\ \int_0^t \| \nabla\nabla w (t',\cdot)\ \|_{L^4(\Omega_0;{\mathbb R}^{27})} dt'\nonumber\\
&+C\  
\| v\|_{H^3(\Omega_0;{\mathbb R}^{3})}\ \int_0^t \| w (t',\cdot) \ \|_{H^3(\Omega_0;{\mathbb R}^{3})}dt'\ ,
\end{align*}
and from another use of the Sobolev embedding theorem together with 
the Cauchy-Schwarz inequality,
\begin{align*}
\|\partial_t g^{\alpha \beta}(t,\cdot)& \int_0^t\frac{\partial^2 w} 
{\partial y^\alpha \partial y^\beta}(t',\cdot) dt'\|_{ H^{\frac{1}{2}}
(\Gamma_0;{\mathbb R}^3)}\nonumber\\
&\le C \  
\| v\|_{H^3(\Omega_0;{\mathbb R}^{3})}\ \sqrt{T}\ \| w  \ \|_{L^2(0,T;H^3(\Omega_0;{\mathbb R}^{3}))}\ .
\end{align*}
Consequently,
\begin{align}
\|\partial_t g^{\alpha \beta} (t,\cdot)& \frac{\partial^2 w} 
{\partial y^\alpha \partial y^\beta}(t,\cdot)
\|_{L^2(0,T; H^{\frac{1}{2}}(\Gamma_0;{\mathbb R}^3))}\nonumber\\
&\le C \  
\| v\|_{L^2(0,T;H^3(\Omega_0;{\mathbb R}^{3}))}\ 
\sqrt{T}\ \| w  \ \|_{L^2(0,T;H^3(\Omega_0;{\mathbb R}^{3}))}\nonumber \\
&\le C
\ \sqrt{T}\ \| w  \ \|_{L^2(0,T;H^3(\Omega_0;{\mathbb R}^{3}))}\ ,
\label{bt3}
\end{align}
where we have used that
$\| v\|_{L^2(0,T;H^3(\Omega_0;{\mathbb R}^{3}))}\le M(T)\le M(T_0)$
(since $(v,q)\in C_T$)  for the last inequality.  Recall that we allow
$C$ to depend on $T_0$.

Now, to estimate the last term of the right-hand-side of (\ref{bt0}), 
we note that
\begin{align*}
&\left\|\partial_t(g^{\alpha \beta}(t,\cdot)
\Gamma^\gamma_{\alpha\beta}(t,\cdot))  
\int_0^t \frac{\partial w} 
{\partial y^\gamma}(t',\cdot)dt' \right\|_{ H^{\frac{1}{2}}
(\Gamma_0;{\mathbb R}^3)}\nonumber\\
&\qquad\qquad\qquad\le C  \left\|\ |\nabla \eta|^3 |\nabla\nabla\nabla v| \left|\int_0^t \nabla w (t',\cdot) dt'\right|\ \right\|_{L^2(\Omega_0;{\mathbb R})}\nonumber\\
& \qquad\qquad\qquad+ C \left\|\ |\nabla \eta|^2 |\nabla\nabla\eta| |\nabla\nabla v|\left|\int_0^t \nabla w (t',\cdot) dt'\right|\ \right\|_{L^2(\Omega_0;{\mathbb R})}\nonumber\\
& \qquad\qquad\qquad+  C \left\|\ |\nabla \eta|^2 |\nabla\nabla\nabla\eta| |\nabla v|\left|\int_0^t \nabla w (t',\cdot) dt'\right|\ \right\|_{L^2(\Omega_0;{\mathbb R})} \nonumber\\
& \qquad\qquad\qquad+  C \left\|\ |\nabla \eta|^3 |\nabla\nabla v|\left|\int_0^t \nabla\nabla w (t',\cdot) dt'\right|\ \right\|_{L^2(\Omega_0;{\mathbb R})} \nonumber\\
& \qquad\qquad\qquad+  C \left\|\ |\nabla \eta|^2 |\nabla\nabla\eta| |\nabla v|\left|\int_0^t \nabla\nabla w (t',\cdot) dt'\right|\ \right\|_{L^2(\Omega_0;{\mathbb R})} \nonumber\\
& \qquad\qquad\qquad+  C \left\|\ |\nabla \eta| |\nabla\nabla\eta|^2 |\nabla v|\left|\int_0^t \nabla w (t',\cdot) dt'\right|\ \right\|_{L^2(\Omega_0;{\mathbb R})} \nonumber\ .
\end{align*}

From (\ref{etah3}), the Sobolev embedding theorem, and the Cauchy-Schwarz 
inequality, we then deduce that
\begin{align*}
& \left\|\partial_t(g^{\alpha \beta}(t,\cdot)
\Gamma^\gamma_{\alpha\beta}(t,x))
\int_0^t \frac{\partial w} 
{\partial y^\gamma}(t'\cdot)dt' \right\|_{ H^{\frac{1}{2}}
(\Gamma_0;{\mathbb R}^3)} \nonumber\\
&\qquad\qquad\le C\  \|v\|_{H^3(\Omega_0;{\mathbb R}^3)}\ \left\|\int_0^t \nabla w (t',\cdot) dt'\ \right\|_{L^\infty(\Omega_0;{\mathbb R}^9)}\nonumber\\
& \qquad\qquad + C\  \|\nabla\nabla v\|_{L^4(\Omega_0;{\mathbb R}^{27})}\ \left\|\int_0^t \nabla w (t',\cdot) dt'\ \right\|_{L^\infty(\Omega_0;{\mathbb R}^9)}\nonumber\\ 
& \qquad\qquad + C\  \|\nabla v\|_{L^\infty(\Omega_0;{\mathbb R}^{9})}\ \left\|\int_0^t \nabla w (t',\cdot) dt'\ \right\|_{L^\infty(\Omega_0;{\mathbb R}^9)}\nonumber\\ 
&  \qquad\qquad+ C\   \|\nabla\nabla v\|_{L^4(\Omega_0;{\mathbb R}^{27})}\ \left\|\int_0^t \nabla\nabla w (t',\cdot) dt'\ \right\|_{L^4(\Omega_0;{\mathbb R}^{27})}\nonumber\\ 
& \qquad\qquad+ C\  \|\nabla v\|_{L^\infty(\Omega_0;{\mathbb R}^{9})}\ \left\|\int_0^t \nabla\nabla w (t',\cdot) dt'\ \right\|_{L^4(\Omega_0;{\mathbb R}^{27})}\nonumber\\ 
& \qquad\qquad+ C\  \|\nabla v\|_{L^\infty(\Omega_0;{\mathbb R}^{9})}\ \left\|\int_0^t \nabla w (t',\cdot) dt'\ \right\|_{L^\infty(\Omega_0;{\mathbb R}^{9})}\ .
\nonumber
\end{align*}
By the Sobolev embedding theorem and Jensen's inequality, the previous 
inequality yields
\begin{align*}
& \left\|\partial_t(g^{\alpha \beta}(t,\cdot)
\Gamma^\gamma_{\alpha\beta}(t,\cdot))
 \int_0^t \frac{\partial w} 
{\partial y^\gamma}(t',\cdot)dt' \right\|_{ H^{\frac{1}{2}}
(\Gamma_0;{\mathbb R}^3)} \nonumber\\
&\qquad\qquad\le C\   \|v\|_{H^3(\Omega_0;{\mathbb R}^{3})}\ \int_0^t \| w (t',\cdot)\|_{H^3(\Omega_0;{\mathbb R}^3)} dt'\nonumber \ .
\nonumber
\end{align*}
Consequently, by the Cauchy-Schwarz inequality,
\begin{align*}
& \left\|\partial_t(g^{\alpha \beta}(t,\cdot)
\Gamma^\gamma_{\alpha\beta}(t,\cdot)) 
\int_0^t \frac{\partial w} 
{\partial y^\gamma}(t',\cdot)dt' \right\|_{H^{\frac{1}{2}}
(\Gamma^i_0;{\mathbb R}^3)} \nonumber\\
&\qquad\qquad\qquad\qquad\qquad \le  C\ \sqrt{T}\   \|v\|_{H^3(\Omega_0;{\mathbb R}^{3})}\ \| w \|_{L^2(0,T;H^3(\Omega_0;{\mathbb R}^3))} \nonumber \ ,
\nonumber
\end{align*}
which leads us to
\begin{align}
& \left\|\partial_t(g^{\alpha \beta}\Gamma^\gamma_{\alpha\beta})
 \int_0^\cdot \frac{\partial w} 
{\partial y^\gamma}(t',\cdot)dt' \right\|_{L^2(0,T;H^{\frac{1}{2}}
(\Gamma_0;{\mathbb R}^3))} \nonumber\\
&\qquad\qquad\qquad\qquad\qquad\le  C\ \sqrt{T}\   \|v\|_{L^2(0,T;H^3(\Omega_0;{\mathbb R}^{3}))}\ \| w \|_{L^2(0,T;H^3(\Omega_0;{\mathbb R}^3))}\nonumber\\
&\qquad\qquad\qquad\qquad\qquad\qquad\le C\ \sqrt{T}\   \| w \|_{L^2(0,T;H^3(\Omega_0;{\mathbb R}^3))}\ . 
\label{bt4}
\end{align}
where we have used that $(v,q)\in C_T$ for the last inequality.

Finally, combining (\ref{bt0}), (\ref{bt1}), (\ref{bt2}), (\ref{bt3}),
and (\ref{bt4}), there exists $T_{\overline{B}_t}\in (0,T_0)$ such that for 
any $T\in (0, T_{\overline{B}_t})$ we have the estimate
\begin{equation}
\label{bt}
\|\overline{B}_t (w)\|_{L^2(0,T;H^{\frac{1}{2}} (\Gamma_0;{\mathbb R}))}\le C\ T^{\frac{1}{16}}\ \|(w,r)\|_{X_T}\ .
\end{equation}

{\bf 9. Estimate of $\overline{B}$:}\hfill\break
Since $\overline{B}(0)=0$, we then deduce from (\ref{bt0}) 
that for any $T\in (0, T_{\overline{B}_t})$,
\begin{align}
\label{b}
\|\overline{B}(w)\|_{L^2(0,T;H^{\frac{1}{2}} (\Gamma_0;{\mathbb R}))}\le C\ T^{\frac{9}{16}}\ \|(w,r)\|_{X_T}\ .
\end{align}

{\bf 10. End of proof of Lemma \ref{estimates}.}\hfill\break
Using the results of the first nine steps, we have that 
$${T}_1=\min\left(T_{\bar{a}_t}, T_{\bar{f}_t}, T_{\bar{g}_t}, T_{\overline{B}_t}\right) $$
satisfies the statement of Lemma \ref{estimates}.

\end{proof}

We next have the following weak continuity result.

\begin{lemma}
\label{weakcontinuity}
For $0<T\le \overline{T}_0$ the mapping $\Theta_T$ associating $(\tilde v,\tilde q)$ to $(v,q)\in C_T$ is weakly continuous from
$C_T$ into $C_T$.
\end{lemma}
\begin{proof}
Let 
Let $(v_p,q_p)_{p\in\mathbb N}$ be a given sequence of elements of $C_T$
weakly convergent (in $X_T$) towards a given element $(v,q)\in C_T$ (
$C_T$ is sequentially weakly closed as a closed convex set) and let  $(v_{\sigma(p)},q_{\sigma(p)})_{p\in\mathbb N}$ be
any subsequence of this sequence.

Since $V^3(T)$ is compactly embedded into
$L^2((0,T);W^{2,5}(\Omega;{\mathbb R}^3))$ (see \cite{Evans1998}), 
we deduce the following strong convergence results in 
$L^2((0,T);L^{2}(\Omega;{\mathbb R}^3))$  as $p$ goes to $\infty$:
\begin{subequations}
\label{scv}
\begin{align}
(a^j_l)_p (a^k_l)_p &\rightarrow a^j_l a^k_l\,, 
\label{scv.a}\\
[(a^j_l)_p (a^k_l)_p]_{,j} &\rightarrow (a^j_l a^k_l)_{,j}\,,
\label{scv.b}\\
(a^k_i)_p &\rightarrow a^k_i \ .
\end{align}
\end{subequations}

Since $(q_{\sigma(p)})_{p\in\mathbb N}$ is bounded in $V^2 (T)$,
let  $(q_{\sigma'(p)})_{p\in\mathbb N}$ be a subsequence of $(q_{\sigma(p)})_{p\in\mathbb N}$, weakly convergent in $V^2 (T)$ to a limit
$\tilde q\in V^2 (T)$. 
Since $(v_{\sigma'(p)})_{p\in\mathbb N}$ is a bounded sequence in $V^3 (T)$, let
$(v_{\sigma''(p)})_{p\in\mathbb N}$ be a subsequence weakly convergent in $V^3 (T)$ to an element
$\tilde v \in V^3 (T)$.

From the strong convergence results (\ref{scv}), we obtain
the weak convergence results in
$L^2((0,T);L^{2}(\Omega;{\mathbb R}^3))$: 
\begin{align}
\partial_t{\tilde{v}}_{\sigma''(p)}^i - \nu (a^j_l ({\sigma''(p)}) a^k_l ({\sigma''(p)}) {\tilde{v}}_{\sigma''(p)}^i,_k),_j& + a^k_i ({\sigma''(p)}) 
{\tilde{q}}_{\sigma''(p)},_k  \nonumber \\
& \rightharpoonup 
 {\tilde{v}}^i_t - 
\nu (a^j_l a^k_l {\tilde{v}}^i,_k),_j + a^k_i \tilde{q},_k \,, \nonumber \\
 a^k_i ({\sigma''(p)}) {\tilde{v_{\sigma''(p)}}}^i,_k 
& \rightharpoonup a^k_i {\tilde{v}}^i,_k \,, \nonumber
\end{align}
and thus we have that
\begin{subequations}
\label{e1bis}
\begin{gather}
{\tilde{v}}^i_t - \nu (a^j_l a^k_l {\tilde{v}}^i,_k),_j + a^k_i \tilde{q},_k = F^i\,, \\
 a^k_i {\tilde{v}}^i,_k =0 \ . 
\end{gather}
\end{subequations}
in $L^2((0,T);L^{2}(\Omega;{\mathbb R}^3))$.

Moreover, the continuity of the embedding of $V^3 (T)$ into 
$C([0,T];H^{2}(\Omega;{\mathbb R}^3))$ ensures us that in 
$H^{2}(\Omega;{\mathbb R}^3)\}$
\begin{equation}
\label{e2bis}
\tilde v (0)=u_0\ .
\end{equation}

By the trace theorem,  
\begin{align*}
(v_{\sigma''(p)})_{p\in\mathbb N} 
&\text{ is bounded in } L^2((0,T);H^{2.5}(\Gamma_0;{\mathbb R}^3))\,, \\
(\partial_t v_{\sigma''(p)})_{p\in\mathbb N}
&\text{ is bounded in } L^2((0,T);H^{0.5}(\Gamma_0;{\mathbb R}^3)). 
\end{align*}
Then, we infer that  $(v_{\sigma''(p)})_{p\in\mathbb N}$ strongly converges 
in $L^2((0,T);H^{2}(\Gamma_0;{\mathbb R}^3))$ to $\tilde v$; furthermore, 
the compactness
of the embedding of $V^2(T)$ into $L^2((0,T);H^1(\Omega_0;{\mathbb R}^3))$ 
provides the strong convergence of $(q_{\sigma''(p)})_{p\in\mathbb N}$ to 
$\tilde q$ in $L^2((0,T);H^{1}(\Omega_0;{\mathbb R}^3))$.
We can then easily deduce that in $L^2((0,T);L^2(\Gamma_0;{\mathbb R}^3))$, 
we have as $p\rightarrow \infty$,
\begin{align*}
\Pi_0\Pi_{\eta_{\sigma''(p)}}\left( D_{\eta_{\sigma''(p)}} (\tilde{v})\cdot a^T ({\sigma''(p)})\ N\right)
&\rightarrow\Pi_0\Pi_{\eta}\left( D_{\eta} 
(\tilde{v})\cdot a^T \ N\right)\,,\\
N\cdot S_{\eta_{\sigma''(p)}}({\tilde{v}}_{\sigma''(p)},{\tilde{q}}_{\sigma''(p)})\cdot a^T({\sigma''(p)})\ N &-
\sigma\ N\cdot \triangle_{g_{\sigma''(p)}(t)}(\int_0^t {\tilde{v}}_{\sigma''(p)})\\
&\rightarrow 
N\cdot (\ S_{\eta}(\tilde{v},\tilde{q})\cdot a^T\ N\ ) -
\sigma\ N\cdot \triangle_{g(t)}(\int_0^t {\tilde{v}})\,,
\end{align*}
and thus in $L^2((0,T);L^{2}(\Gamma_0;{\mathbb R}^3))$,
\begin{subequations}
\label{e3bis}
\begin{align}
 \Pi_0\Pi_{\eta} D_{\eta} (\tilde{v})\cdot a^T \ N&=0 \,, \\
N\cdot (\ S_{\eta}(\tilde{v},\tilde{q})\cdot a^T\ N\ ) -
\sigma\ N\cdot \triangle_{g(t)}(\int_0^t {\tilde{v}})\ 
&=\sigma\ N\cdot \triangle_{g(t)}(\text{Id}) \,.
\end{align}
\end{subequations}
Concerning the forcings, we obviously have the strong convergence in the space 
$L^2(0,T;L^{2}(\Omega_0;{\mathbb R}^3))$ of $f\circ \eta_{\sigma''(p)}$ to
$f\circ \eta$. Thus, with (\ref{e1bis}), (\ref{e2bis}) and (\ref{e3bis}), we then deduce that
\begin{subequations}
 \begin{align}
{\tilde{v}}^i_t - \nu (a^j_l a^k_l {\tilde{v}}^i,_k),_j + a^k_i \tilde{q},_k 
&= F^i 
\ \ \text{in} \ \ (0,T)\times \Omega_0 \,, 
        \nonumber\\
   a^k_i {\tilde{v}}^i,_k &= 0     \ \ \text{in} \ \ (0,T)\times \Omega_0 \,, 
         \nonumber\\
\Pi_0\Pi_{\eta}\left( D_{\eta} (\tilde{v})\cdot a^T\ N\right) &=0
\ \ \text{on} \ \ (0,T)\times \Gamma_0 \,, 
        \nonumber\\
N\cdot (\ S_{\eta}(\tilde{v},\tilde{q})\cdot a\ N\ ) -
\sigma\ N\cdot \triangle_{g(t)}(\int_0^t {\tilde{v}})&= 
\sigma\ N\cdot \triangle_{g(t)}(\text{Id})
 \ \ \text{on} \ \ \Gamma_0\times \{ t=0\} \,, 
         \nonumber\\
   \tilde{v} &= u_0    
 \ \ \text{on} \ \ \Omega_0\times \{ t=0\} \,, 
         \nonumber
\end{align}
\end{subequations}

This precisely shows that 
$$(\tilde v,\tilde q)= \Theta_T (v,q).$$
Hence, we deduce that the whole sequence 
${(\Theta_T (v_n,q_n))}_{n\in\mathbb N}$ weakly converges
in $C_T$ towards $ \Theta_T (v,q) $, which proves the lemma. 
\end{proof} 

\section{Proof of the main theorem}

\subsection{Existence}
Let $T\in (0,\overline{T}_0)$. 
The mapping $\Theta$ being weakly continuous from the closed bounded convex
set $C_T$ into itself from Lemmas \ref{stability} and \ref{weakcontinuity}, we
infer from the Tychonoff fixed point theorem ({\it see} for instance 
\cite{Deimling}) that it admits (at least) one fixed point 
$(v,q)=\theta(v,q)$ in $C_T$. Hence, we have
\begin{subequations}
  \label{fp}
\begin{align}
{{v}}^i_t - \nu (a^j_l a^k_l {{v}}^i,_k),_j + a^k_i {q},_k &= F^i 
\ \ \text{in} \ \ (0,T)\times \Omega_0 \,, 
         \label{fp.a}\\
   a^k_i {v}^i,_k &= 0     \ \ \text{in} \ \ (0,T)\times \Omega_0 \,, 
         \label{nlfp.b}\\
\Pi_0\Pi_{\eta}\left( D_{\eta} ({v})\cdot a^T\ N\right) &=0
\ \ \text{on} \ \ (0,T)\times \Gamma_0 \,, 
         \label{fp.c}\\
N\cdot \left(S_{\eta}({v},{q})\cdot a^T\ N\right) -
\sigma\ N\cdot \triangle_{g(t)}(\int_0^t {{v}})&= 
\sigma\ N\cdot \triangle_{g(t)}(\text{Id})
 \ \ \text{on} \ \  (0,T)\times \Gamma_0\,, 
         \label{fp.d}\\
  {v} &= u_0    
 \ \ \text{on} \ \ \Omega_0\times \{ t=0\} \,, 
         \label{fp.e}
\end{align}
\end{subequations}
Condition (\ref{fp.d}) also reads
\begin{equation}
\label{fpbis.d}
N\cdot \left(S_{\eta}({v},{q})\cdot a^T\ N -
\sigma\ \triangle_{g(t)}(\eta)\right)= 
0 \ \ \text{on} \ \ (0,T)\times \Gamma_0 \ . 
\end{equation}
On the other hand, (\ref{fp.c}) is equivalent to
\begin{equation*}
\Pi_{\eta}\left( D_{\eta} ({v})\cdot a^T\ N \right)= \lambda N 
 \ \ \text{on} \ \ (0,T)\times \Gamma_0 \ , 
\end{equation*}
where $\lambda\in L^2((0,T);H^{0.5}(\Gamma_0;{\mathbb R}^3))$.
From (\ref{control2}), $N$ and $a^T N$ are not orthogonal. The previous identity
can then hold only if $\lambda=0$. Thus,
\begin{equation}
\label{fpbis.c}
\Pi_{\eta}\left( D_{\eta} ({v})\cdot a^T\ N\right) = 0 
 \ \ \text{on} \ \ (0,T)\times \Gamma_0 \ . 
\end{equation}
Since the vectors $q\ a^T\  N$ and $\triangle_{g(t)}(\eta)$ are also
colinear to $a^T\ N$, the direction normal to $\eta (\Gamma_0)$, 
(\ref{fpbis.c}) implies that
\begin{equation}
\label{fpter.c}
\Pi_{\eta}\left( S_{\eta} ({v},q)\cdot a^T\ N - \sigma\ \triangle_{g(t)}(\eta)\right)= 0  \ \ \text{on} \ \ (0,T)\times \Gamma_0 \ . 
\end{equation}
Adding together (\ref{fpbis.d}) and (\ref{fpter.c}), we have that 
\begin{equation}
\label{fpbis}
S_{\eta} ({v},q)\cdot a^T\ N - \sigma\ \triangle_{g(t)}(\eta)= 0  \ \ \text{on} \ \ (0,T)\times \Gamma_0 \ . 
\end{equation}
This finally shows that
\begin{align*}
{{v}}^i_t - \nu (a^j_l a^k_l {{v}}^i,_k),_j + a^k_i {q},_k &= F^i 
\ \ \text{in} \ \ (0,T)\times \Omega_0 \,, 
         \\
   a^k_i {v}^i,_k &= 0     \ \ \text{in} \ \ (0,T)\times \Omega_0 \,, 
         \\
S_{\eta}({v},{q})\cdot a^T\ N&= 
\sigma\  \triangle_{g(t)}(\eta)
 \ \ \text{on} \ \ (0,T)\times \Gamma_0 \,, 
         \\
  {v} &= u_0    
 \ \ \text{on} \ \ \Omega_0\times \{ t=0\} \,, 
\end{align*}
so that $v$ is a solution of (\ref{nsl}).

\subsection{Uniqueness}
Now, to prove the uniqueness of a solution to (\ref{nsl}) in $X_T$, 
we use the Lipschitz condition (\ref{Lip}).
Let $(\tilde v,\tilde q)$ be another solution of (\ref{nsl}). Then, 
\begin{subequations}
\label{unique}
\begin{align}
({{v-\tilde{v}}})^i_t - \nu (a^j_l a^k_l ({{v}}^i,_k-{\tilde{v}}^i,_k )),_j
 + a^k_i ({q},_k-\tilde{q},_k ) 
&= \delta f\ \ \text{in} \ \ (0,T)\times \Omega_0 \,, \\
 a^k_i {(v-\tilde{v})}^i,_k &= \delta a  \ \ \text{in}
 \ \ (0,T)\times \Omega_0 \,, \\
\Pi_0\Pi_{\eta} \left(D_{\eta} ({v-\tilde{v}})\cdot a^T\ N\right) 
&=\delta g_1 \ \ \text{on} \ \ (0,T)\times \Gamma_0 \,, \\
N\cdot (\ S_{\eta}({v-\tilde{v}},q-\tilde{q})\cdot a^T\ N\ ) 
-\sigma\ N\cdot \triangle_{g(t)}&(\int_0^t {{v-\tilde{v}}}) \\
& = \delta g_2 +\sigma \delta B \ \ \text{on} \ \ (0,T)\times \Gamma_0 \,, \\
  {v}-\tilde{v} &= 0    
 \ \ \text{on} \ \ \{0\}\times \Omega_0 \,, 
\end{align}
\end{subequations}
with
\begin{align}
\delta f&=- \nu ((a^j_l a^k_l -\tilde{a}^j_l \tilde{a}^k_l)\ {\tilde{v}}^i,_k )),_j + f\circ\eta-f\circ\tilde{\eta}
\ \text{in} \ \ (0,T)\times \Omega_0 \,, \label{df} \\
\delta a&= (\tilde{a}^k_i-a^k_i) {\tilde{v}}^i,_k 
\ \text{in} \ \ (0,T)\times \Omega_0 \,,\nonumber\\
\delta g_1&= \Pi_0\left(\ \Pi_{\tilde\eta}(\ D_{\tilde\eta}(\tilde{v})\cdot \tilde{a}^T\ N\ )-\Pi_{\eta}(\  D_{\eta}(\tilde{v})\cdot a^T\ N\ )\ \right)
\ \text{on} \ \ (0,T)\times \Gamma_0 \,,\nonumber\\
\delta g_2&=N\cdot \left(S_{\tilde\eta}(\tilde{v},\tilde{q})\cdot \tilde{a}^T\ N -S_{\eta}(\tilde{v},\tilde{q})\cdot a^T\ N\right)\text{on} \ \ (0,T)\times \Gamma_0 \,, \nonumber\\
\delta B&= -
\sigma\ N\cdot (\triangle_{\tilde{g}(t)}-\triangle_{g(t)})(\int_0^t {\tilde{v}})\ \text{on} \ \ (0,T)\times \Gamma_0 \ .
\label{u6}
\end{align}
We note that these forcing terms are similar to those appearing in (\ref{ilfp}),
with $(v_n,q_n)$ replaced by $(\tilde v,\tilde q)$, and Id replaced by 
$\tilde a$, and with the two addition terms on the right-hand-side of 
(\ref{df}):
\begin{equation}\label{2terms}
 f \circ \eta  - f \circ  \tilde \eta  \,.
 \end{equation}	
 
Since the initial boundary forcings satisfy $\delta g_1 (0,\cdot)=0$, $\delta g_1 (0,\cdot)=0$ and $\delta B (0,\cdot)=0$, and the 
initial data satisfies the compatibility condition 
$$\Pi_0 \left(\operatorname{Def} (0)\cdot N\right)=0\ ,$$
associated to the condition $\delta g_1(0,\cdot)=0$, we  can apply the
same type of iteration procedure that we used to establish the solvability
of the linear problem in the proof of Lemma \ref{stability} to solve
(\ref{unique}).  We find that for any $T'\in (0,T)$,
\begin{align}
&\| (v-\tilde{v},q-\tilde{q})\|_{X_{T'}} \nonumber\\
&\le   C \ \left(
\| \delta f \|_ {L^2(0,T';H^1(\Omega_0;{\mathbb R}^3))}+ 
\|\delta f_t \|_{L^2(0,T';H^1(\Omega_0;{\mathbb R}^3)')} \right.\nonumber \\
& \qquad\qquad \qquad 
+ \|\delta a \|_ {L^2(0,T';H^2(\Omega_0;{\mathbb R}^9))} 
+ \|\delta a_t \|_ {L^2(0,T';L^2(\Omega_0;{\mathbb R}^9))} \nonumber \\
& \qquad\qquad \qquad 
+ \|\delta g_1\|_ {L^2(0,T';H^\frac{3}{2}(\Gamma_0;{\mathbb R}^3))}
+ \|{\delta {g_1}}_t \|_ {L^2(0,T';H^{-\frac{1}{2}}(\Gamma_0;{\mathbb R}^3))}\nonumber \\
& \qquad\qquad \qquad 
+ \|\delta g_2\|_ {L^2(0,T';H^\frac{3}{2}(\Gamma_0;{\mathbb R}))}
+ \|{\delta {g_2}}_t \|_ {L^2(0,T';H^{-\frac{1}{2}}(\Gamma_0;{\mathbb R}))}\nonumber \\
&\qquad\qquad \qquad\left. + \|\delta B \|_ {L^2(0,T';H^\frac{1}{2}(\Gamma_0;{\mathbb R}))}
+ \|\delta B_t \|_ {L^2(0,T';H^{\frac{1}{2}}(\Gamma_0;{\mathbb R}))}\right) \ .
\label{u1}
\end{align}

Now, by using (\ref{deltaah2}) together with the remark made after (\ref{u6}),
the same method as we used for proving Lemma \ref{estimates} can be applied 
for estimating the right-hand-side of (\ref{u1}) (with (\ref{deltaah2}) 
playing the role of (\ref{ah2}) in the present case); 
we find that for any $T'\in (0,T)$,
\begin{align}
\label{u2}
& \| \delta f \|_ {L^2(0,T';H^1(\Omega_0;{\mathbb R}^3))}+ 
\|\delta f_t \|_{L^2(0,T';H^1(\Omega_0;{\mathbb R}^3)')} \nonumber \\
& + \|\delta a \|_ {L^2(0,T';H^2(\Omega_0;{\mathbb R}^9))} 
+ \|\delta a_t \|_ {L^2(0,T;L^2(\Omega_0;{\mathbb R}^9))} \nonumber \\
&  + \|\delta g_1\|_ {L^2(0,T';H^\frac{3}{2}(\Gamma_0;{\mathbb R}^3))}
+ \|{\delta g_1}_t \|_ {L^2(0,T';H^{-\frac{1}{2}}(\Gamma_0;{\mathbb R}^3))}\nonumber \\
&  + \|\delta g_2\|_ {L^2(0,T';H^\frac{3}{2}(\Gamma_0;{\mathbb R}))}
+ \|{\delta g_2}_t \|_ {L^2(0,T';H^{-\frac{1}{2}}(\Gamma_0;{\mathbb R}))}\nonumber \\
& + \|\delta B \|_ {L^2(0,T';H^\frac{1}{2}(\Gamma_0;{\mathbb R}))}
+ \|\delta B_t \|_ {L^2(0,T';H^{\frac{1}{2}}(\Gamma_0;{\mathbb R}))}\nonumber \\
&\qquad\qquad\qquad\le C T'^{\frac{1}{16}}\ \|(\tilde {v},\tilde {q})\|_{X_{T'}}\ \|(v-\tilde {v},q-\tilde {q})\|_{X_{T'}}\ . 
\end{align}
Note that the additional terms (\ref{2terms}) arising in (\ref{df}) 
do not cause any difficulties because of the estimate (\ref{fcirceta}).

Thus, from (\ref{u1}) and (\ref{u2}), we have that
\begin{equation}
\label{u3}
\|(v-\tilde v,q-\tilde q)\|_{X_{T'}}\le C \ \|(\tilde {v},\tilde {q})\|_{X_{T'}}\ T'^{\frac{1}{16}}\ \|(v-\tilde v,q-\tilde q)\|_{X_{T'}}\ .
\end{equation}
Now, let $T'\in (0,T)$ be such that $C\ \|(\tilde {v},\tilde {q})\|_{X_{T}}\ T'^{\frac{1}{16}} <1$. Then, from (\ref{u3}), we have $\|(v-\tilde v,q-\tilde q)\|_{X_{T'}}=0$ and thus $(v-\tilde v,q-\tilde q)=0$ on $[0,T']$. 
Let us define $$T_s=\sup\{t\in (0,T)|\ \forall t'\le t,\ (v-\tilde v,q-\tilde q)(t',\cdot)=0\}\ .$$
We have $0<T'\le T_s\le T$. let us assume that $0<T_s< T$. Since $V^3(T)\times V^2(T)\subset
C([0,T]; H^2(\Omega_0;{\mathbb R}^3))\times C([0,T]; H^1(\Omega_0;{\mathbb R}^3))$, we have by continuity $$(v-\tilde v,q-\tilde q)(T_s,\cdot)=(0,0)\ .$$ 
By an integration in time from 0 to $T_s$, one also immediately gets $\eta(T_s,\cdot)=\tilde{\eta} (T_s,\cdot)$ and as a consequence $\delta g_1(T_s,\cdot)=0$,
$\delta g_2(T_s,\cdot)=0$, $\delta B(T_s,\cdot)=0$. Thus, we see that
\begin{subequations}
\label{unique2}
\begin{align}
({{v-\tilde{v}}})^i_t - \nu (a^j_l a^k_l ({{v}}^i,_k-{\tilde{v}}^i,_k )),_j
 + a^k_i ({q},_k-\tilde{q},_k ) 
&= \delta f\ \ \text{in} \ \ (T_s,T)\times \Omega_0 \,, \\
a^k_i {(v-\tilde{v})}^i,_k
& = \delta a     \ \ \text{in} \ \ (T_s,T)\times \Omega_0 \,, \\
\Pi_0\Pi_{\eta} \left( D_{\eta} ({v-\tilde{v}})\cdot a^T\ N \right)
& =\delta g_1
\ \ \text{on} \ \ (T_s,T)\times \Gamma_0 \,, \\
N\cdot(\ S_{\eta}({v-\tilde{v}},q-\tilde{q})\cdot a^T\ N\ )
-\sigma\ N\cdot \triangle_{g(t)}(\int_{T_s}^t {{v-\tilde{v}}})
&= \delta g_2+ \delta B \ \ \text{on} \ \ (T_s,T)\times \Gamma_0 \,, \\
 {v}-\tilde{v} &= 0    
 \ \ \text{on} \ \ \Omega_0\times \{ t=T_s\} \,. 
\end{align}
\end{subequations}

Since the new 
initial data $(v-\tilde{v})({T_s},\cdot)=(0,0)$ satisfies the compatibility
condition 
$$\Pi_0 \left(\operatorname{Def} ((v-\tilde{v})({T_s},\cdot))\cdot a^T({T_s},\cdot)\ N\right)
=0\ ,$$ 
associated to the conditions $\delta g_1(T_s,\cdot)=0$, $\delta g_2(T_s,\cdot)=0$ and $\delta B(T_s,\cdot)=0$, we can then use the
same estimate as (\ref{u1}) for the system (\ref{unique2}); thus, 
for any $T''\in (0,T)$,
\begin{align}
\label{u4}
&\| (v-\tilde{v},q-\tilde{q})(T_s+\cdot,\cdot)\|_{X_{T''}} \nonumber\\
&\le   C \ \left(
\| \delta f \|_ {L^2(T_s,T'';H^1(\Omega_0;{\mathbb R}^3))}+ 
\|\delta f_t \|_{L^2(T_s,T'';H^1(\Omega_0;{\mathbb R}^3)')} \right.\nonumber \\
& \qquad\qquad \qquad 
+ \|\delta a \|_ {L^2(T_s,T'';H^2(\Omega_0;{\mathbb R}^9))} 
+ \|\delta a_t \|_ {L^2(T_s,T'';L^2(\Omega_0;{\mathbb R}^9))} \nonumber \\
& \qquad\qquad \qquad 
+ \|\delta g_1\|_ {L^2(T_s,T'';H^\frac{3}{2}(\Gamma_0;{\mathbb R}^3))}
+ \|{\delta g_1}_t \|_ {L^2(T_s,T'';H^{-\frac{1}{2}}(\Gamma_0;{\mathbb R}^3))}\nonumber \\
&\qquad\qquad \qquad + \|\delta g_2\|_ {L^2(T_s,T'';H^\frac{3}{2}(\Gamma_0;{\mathbb R}))}
+ \|{\delta g_2}_t \|_ {L^2(T_s,T'';H^{-\frac{1}{2}}(\Gamma_0;{\mathbb R}))}\nonumber \\
& \qquad\qquad \qquad\left. 
+ \|\delta B \|_ {L^2(T_s,T'';H^\frac{1}{2}(\Gamma_0;{\mathbb R}))}
+ \|\delta B_t \|_ {L^2(T_s,T'';H^{\frac{1}{2}}(\Gamma_0;{\mathbb R}))}\right)
\ . 
\end{align}
Since the same methods that we applied for proving Lemma \ref{estimates} 
can also be used for estimating the right-hand-side of (\ref{u4}), 
we find that for any $T''\in (T_s,T)$,
\begin{align}
\label{u5}
& \| \delta f \|_ {L^2(T_s,T'';H^1(\Omega_0;{\mathbb R}^3))}+ 
\|\delta f_t \|_{L^2(T_s,T'';H^1(\Omega_0;{\mathbb R}^3)')} \nonumber \\
& + \|\delta a \|_ {L^2(T_s,T'';H^2(\Omega_0;{\mathbb R}^9))} 
+ \|\delta a_t\|_ {L^2(0,T;L^2(\Omega_0;{\mathbb R}^9))} \nonumber \\
&  + \|\delta g_1\|_ {L^2(T_s,T'';H^\frac{3}{2}(\Gamma_0;{\mathbb R}^3))}
+ \|{\delta g_1}_t \|_ {L^2(T_s,T'';H^{-\frac{1}{2}}(\Gamma_0;{\mathbb R}^3))}\nonumber \\
&  + \|\delta g_2\|_ {L^2(T_s,T'';H^\frac{3}{2}(\Gamma_0;{\mathbb R}^3))}
+ \|{\delta g_2}_t \|_ {L^2(T_s,T'';H^{-\frac{1}{2}}(\Gamma_0;{\mathbb R}^3))}\nonumber \\
& + \|\delta B \|_ {L^2(T_s,T'';H^\frac{1}{2}(\Gamma_0;{\mathbb R}))}
+ \|\delta B_t \|_ {L^2(T_s,T'';H^{\frac{1}{2}}(\Gamma_0;{\mathbb R}))}\nonumber \\
&\qquad\qquad\qquad\le C (T''-T_s)^{\frac{1}{16}}\ \|(\tilde {v},\tilde {q})(T_s+\cdot,\cdot)\|_{X_{T''}}\ \|(v-\tilde {v},q-\tilde {q})(T_s+\cdot,\cdot)\|_{X_{T''}}\ . 
\end{align}
Consequently, from (\ref{u4}) and (\ref{u5}), we have that
\begin{align*}
&\|(v-\tilde v,q-\tilde q)(T_s+\cdot,\cdot)\|_{X_{T''}}\nonumber\\
&\qquad\qquad\qquad \le C \ \|(\tilde {v},\tilde {q})(T_s+\cdot,\cdot)\|_{X_{T''}}\ (T''-T_s)^{\frac{1}{16}}\ \|(v-\tilde v,q-\tilde q)(T_s+\cdot,\cdot)\|_{X_{T''}}\ .
\end{align*}
Now, we see from this inequality that for  $T''\in (T_s,T)$ sufficiently 
close to $T_s$,  $(v-\tilde v,q-\tilde q)(t',\cdot)=0$ in $[T_s,T'']$. Thus,
$$T''\le T_s\ ,$$ 
so that we arrive at a contradiction, showing that $T_s=T$, 
{\it i.e.} the uniqueness of the solution to (\ref{nsl}) in $X_T$.

\section*{Acknowledgments}
We are grateful to Ching-Hsiao Cheng for carefully reading the manuscript and
correcting a number of typographical errors.
SS and DC were partially supported by NSF ATM-98-73133.
SS was partially supported by  NSF DMS-0105004 and the Alfred P.  
Sloan Foundation Research Fellowship.

\end{document}